\documentclass[a4paper,11pt,twoside]{scrartcl}
\usepackage{german,amssymb,amsmath,epsfig,graphics,graphicx,mathrsfs,xcolor,amsthm,cite} 
\usepackage[latin1]{inputenc}
\usepackage{fancyhdr}
\usepackage{makeidx}  

\newtheorem{theorem}{Theorem}[section]

\newtheorem{definition}[theorem]{Definition}
\newtheorem{lemma}[theorem]{Lemma}

\newtheorem{beispiel}[theorem]{Beispiel}
\newtheorem{satz}[theorem]{Satz}
\newtheorem{bemerkung}[theorem]{Bemerkung}

\DeclareOldFontCommand{\it}{\normalfont\itshape}{\mathit}
\DeclareOldFontCommand{\rm}{\normalfont\rmfamily}{\mathrm}
\DeclareOldFontCommand{\bf}{\normalfont\bfseries}{\mathbf}

\numberwithin{equation}{section}
\newcommand{\R}{{\mathbb R}}
\newcommand{\N}{{\mathbb N}}

\newcommand{\Zt}{{\mathscr{Z}^k([t_0,t_1])}}

\begin{document}
\thispagestyle{empty}
\vspace*{10mm}
{\bf{\LARGE Nico Tauchnitz}} \\[2cm]
{\bf{\Huge Steuerungsprobleme mit \\[5mm] freiem rechten Endpunkt}} \\[2cm]
{\bf{\Large Die Herleitung des Pontrjaginschen Maximumprinzips \\[5mm]
            mit der Methode der einfachen Nadelvariation}} \\[20mm]
  
\fbox{
\begin{minipage}{0.40\textwidth}Einfache Nadelvariation der \\[1mm] Optimalen Steuerung:
\begin{eqnarray*}
u_{\lambda}(t) &=& 
  \left\{ \begin{array}{ll}
          u_*(t) & \mbox{ f"ur } t \not\in [\tau-\lambda,\tau) \\
          v      & \mbox{ f"ur } t     \in [\tau-\lambda,\tau) 
          \end{array} \right. \qquad \\
\dot{x}_\lambda(t) &=& \varphi\big(t,x_\lambda(t),u_\lambda(t)\big)
\end{eqnarray*}
Marginale Zustandsvariation: \\[2mm]
$\displaystyle y(t)=\lim_{\lambda \to 0^+}\frac{x_{\lambda}(t) - x_*(t)}{\lambda}$
\end{minipage}
\begin{minipage}{0.48\textwidth}
\centering
\includegraphics[width=6cm]{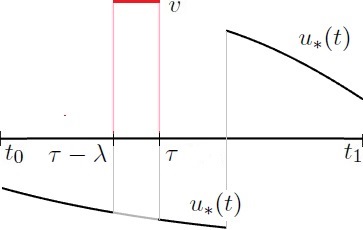}
\end{minipage}}

\cleardoublepage


\newpage
\setcounter{page}{1}
\pagenumbering{roman}

\lhead[\thepage \hspace*{1mm} Inhaltsverzeichnis]{}
\rhead[]{Inhaltsverzeichnis \hspace*{1mm} \thepage}
\tableofcontents
\cleardoublepage


\newpage
\setcounter{page}{1}
\pagenumbering{arabic}
\addcontentsline{toc}{section}{Einleitung}
\lhead[\thepage \hspace*{1mm} Einleitung]{}
\rhead[]{Einleitung \hspace*{1mm} \thepage}
\section*{Einleitung}
Diese Ausarbeitung behandelt notwendige Optimalitätsbedingungen für ein starkes lokales Optimum in
gewissen Klassen von Steuerungsproblemen mit freiem rechten Endpunkt.
Die Aufgaben mit freiem rechten Endpunkt setzen sich aus einem Optimierungskriterium in Form eines Zielfunktionals,
aus einer dynamischen Nebenbedingung in Form eines Systems von Differentialgleichungen mit Anfangsbedingungen
und aus Steuerbeschränkungen zusammen.
Dieser Standardtyp eines Steuerungsproblems lässt sich konstruktiv mit Hilfe einer einfachen Nadelvariation behandeln.
Im Zuge dieser Ausarbeitung wird dieser Zugang auf Steuerungsprobleme mit individuellen Charakteristiken erweitert. \\[2mm]
Der Begriff eines starken lokalen Optimum ist der Klassischen Variationsrechnung,
aus der sich die Theorie der Optimalen Steuerungen entwickelte, entlehnt. 
Wenn wir von einem Optimum sprechen,
so bezeichnet dies eine Lösung der Aufgabe,
die sich über einer Menge von Konkurrenten als die Bestmögliche behauptet.
Es ist somit die Untersuchungen auf einer möglichst großen Konkurrentenmenge anzustreben.
In diesem Zusammenhang kann man einer Optimalstelle eine Qualtitäts- bzw. Gütestufe zuweisen,
nämlich in Bezug auf den Umfang der Konkurrenten. \\
In der Klassischen Variationsrechnung haben sich die Begriffe eines schwachen bzw. starken lokalen Optimum etabliert.
Während bei einer schwachen lokalen Minimalstelle (``Stelle'' bezeichnet hier eine gewisse Funktion) 
sowohl die konkurriende Funktion und auch deren Ableitung sich wenig von einem Optimum und der Ableitung des Optimum zu unterscheiden haben,
so schränkt man bei einem starken lokalen Optimum lediglich die konkurrierende Funktion, aber nicht deren Ableitung ein.
Das bedeutet, die Konkurrentenmenge bei der Untersuchung auf starke lokale Optimalität ist umfassender als diejenige bei der
Behandlung von schwachen lokalen Optimalstellen.
Folglich besitzt ein starkes Optimum eine höhere Güte als ein schwaches
und die Herleitung von Optimalitätsbedingungen ist für ein starkes Optimum ist anzustreben. \\[2mm]
Im Rahmen der Optimalen Steuerung sucht man nicht nur eine Lösungsfunktion, 
sondern ein Funktionenpaar (den sogenannten Steuerungsprozess) aus Zustand(sfunktion) und Steuerung(sfunktion).
Von zentraler mathematischer und auch anwendungsbezogener Bedeutung ist die freie Wahl an Steuerungen,
denn in der Praxis ist man oft mit der Situation konfrontiert, 
dass man sprichwörtlich ``den Schalter umlegen'' oder ``das Steuer herumreißen'' muss.
Ein solcher abrupter Wechsel wird gerade durch Abläufe mit sprunghaftem Änderungsverhalten widergegeben
und ist in den Anwendungen von besonderem Interesse.
Um dies in den Rahmen der Optimalen Steuerungen einzubinden, ergibt sich die Anforderung,
dass man konkurrierenden Steuerungsprozessen zwar Einschränkungen an den Zustand, aber nicht an die Steuerung auferlegt.
Die so gebildete Konkurrentenmenge führt (in Analogie zur Klassischen Variationsrechnung)
zur starken lokalen Optimalität im Rahmen der Theorie Optimaler Steuerungen.
\newpage
Die notwendigen Optimalitätsbedingungen für ein starkes lokales Optimum in der Steuerungstheorie bilden einen Verbund,
dessen Zusammenstellung man als Pontrjaginsches Maximumprinzip bezeichnet.
Im Rahmen dieser Ausarbeitung wollen wir bei freiem rechten Endpunkt das Pontrjaginsche Maximumprinzip mit Hilfe
der elementaren Beweismethode der einfachen Nadelvariation für
\begin{itemize}
\item[$\cdot$] eine Standardaufgabe mit Wiedergewinnungswert,
\item[$\cdot$] optimale Multiprozesse,
\item[$\cdot$] einer Aufgabe mit unendlichem Zeithorizont,
\item[$\cdot$] zeitverzögerte Systeme,
\item[$\cdot$] die Steuerung Volterrascher Integralgleichungen
\end{itemize}
vorstellen.
Diese Methode haben wir für die Standardaufgabe bei Ioffe \& Tichomirov \cite{Ioffe}
entnommen und für die weiteren Klassen ausgebaut. \\[2mm]
Wir haben uns entschieden,
die Methode der einfachen Nadelvariation in dieser Ausarbeitung ausführlich darzustellen,
denn der Nachweis der notwendigen Optimalitätsbedingungen für Aufgaben der Optimalen Steuerung mit freiem rechten Endpunkt mit Hilfe dieser Methode
fällt ``vergleichsweise einfach'' aus: \\
Die ``Einfachheit'' besteht in einem konstruktiven Nachweis des Pontrjaginschen Maximumprinzips.
Auf Basis der einfachen Nadelvariation lassen sich notwendige Optimalitätsbedingungen direkt ermitteln.
Die Herausforderungen des Beweises liegen in den mathematischen Instrumenten,
nämlich in Existenz-, Eindeutigkeits- und Abhängigkeitssätzen für Differential- und Integralgleichungen mit stückweise stetigen rechten Seiten.
Aufgrund des speziellen Rahmens werden die benötigten Resultate im Anhang ausführlich erarbeitet.
Dazu sind aber gewisse Grundbausteine der Funktionalanalysis,
z.\,B. der Begriff des Banachraumes, Fixpunktsätze und der Satz über implizite Funktionen, heranzuziehen.
Somit sind die nötigen Instrumente eher von ``abstrakter'' als von ``einfacher'' Natur. \\
Im ``Vergleich'' zu Steuerungsproblemen mit allgemeineren Nebenbedingungen ist bei der elementaren Beweismethode der Nachweis der
Gültigkeit des Lagrangeschen Prinzips in unendlich-dimensionalen Funktionenräumen nicht nötig.
Dadurch wird das Beweisverfahren deutlich abgekürzt,
da die Maß- und Integrationstheorie oder die Grundprinzipien der Funktionalanalysis und der Satz von Ljusternik oder
die Rieszschen Darstellungssätze oder die benötigten Elemente der Konvexen Analysis umgangen werden. \\[2mm]
Anhand einer Vielzahl an Beispielen demonstrieren wir die erzielten Resultate.
Zu unseren wichtigsten Quellen zählen dabei die Lehrb"ucher von Feichtinger \& Hartl \cite{Feichtinger},
Kamien \& Schwartz \cite{Kamien} und Seierstad \& Syds\ae ter \cite{Seierstad}.
Diese Lehrbücher sind nicht auf die Aufgaben mit freiem rechten Endpunkt beschränkt und
zum weiteren Studium der Theorie Optimaler Steuerungen wärmstens empfohlen.
Außerdem möchten wir den interessierten Leser auf die Arbeiten Pesch \& Plail \cite{PeschPlail} und Plail \cite{Plail} zur historischen Entwicklung
der Optimalen Steuerung zu einem eigenst"andigen mathematischen Fachgebiet hinweisen.

\newpage
\addcontentsline{toc}{section}{Wichtige Bezeichungen}
\lhead[\thepage \hspace*{1mm} Wichtige Bezeichungen]{}
\rhead[]{Wichtige Bezeichungen \hspace*{1mm} \thepage}
\section*{Wichtige Bezeichungen}
Die Methode der einfachen Nadelvariation erlaubt es uns die Aufgaben in dem Rahmen der Räume von
stückweise stetigen und stückweise stetig differenzierbaren Funktionen zu betrachten.
Diese Funktionenklassen kennzeichnen wir mit \\[2mm]
\begin{tabular}{lll}
$PC([a,b],\R)$ &--& Raum der stückweise stetigen Funktionen, \\[1mm]
$PC_1([a,b],\R)$ &--& Raum der stückweise stetig differenzierbaren Funktionen. \\[2mm]
\end{tabular}

Für das Skalarprodukt der Funktionswerte von $f,g: [a,b] \to \R^n$ schreiben wir $\langle f(t),g(t) \rangle$. \\[2mm]
Der gewählte Rahmen dieser Ausarbeitung gestattet die auftretenden Differentationen und Integrationen im klassischen Sinn aufzufassen.
Weiterhin bezeichnen $\dot{x}(t)$ die Zeitableitung und $\displaystyle \int_a^b \,dt$ das (für $b=\infty$ uneigentliche) Riemann-Integral. \\[2mm]
Der Aufbau der Steuerungsprobleme, die wir behandeln werden,
stimmen grundsätzlich im Auftreten von Zielfunktional, Dynamik und Steuerungsbeschränkungen überein.
Deswegen verwenden wir meistens die Bezeichnungen \\[2mm]
\begin{tabular}{lll}
$J$ &--& für das Zielfunktional, \\[1mm]
$f$ &--& für den Integranden im Zielfunktional, \\[1mm]
$S$ &--& für den Wiedergewinnungswert, \\[1mm]
$\varphi$ &--& für die rechte Seite des dynamischen Systems, \\[1mm]
$U$ &--& für den Steuerbereich, \\[1mm]
$x(\cdot)$ &--& für den Zustand, \\[1mm]
$u(\cdot)$ &--& für die Steuerung. \\[2mm]
\end{tabular}

Im Rahmen unserer Untersuchungen treten einseitige Grenzwert wiederholt auf.
Für den links- bzw. rechtseitigen Grenzwert einer Abbildung $\psi$ an der Stelle $x_0$ schreiben wir
$$\lim_{x \to x_0^-} \psi(x) \quad\mbox{bzw.}\quad \lim_{x \to x_0^+} \psi(x).$$  
Die verschiedenen Aufgabentypen besitzen ihre eigenen Charakteristiken und erfordern teils individuelle Voraussetzungen.
Deswegen werden die Pontrjagin-Funktion $H$ und die Hamilton-Funktion $\mathscr{H}$, sowie die Annahmen an die Aufgabenklasse,
wie z.\,B. die Menge der zul"assigen Steuerungsprozesse $\mathscr{A}_{\rm adm}$,
durch einen zur Aufgabe gehörenden Index versehen: \\[2mm]
\begin{tabular}{lll}
$\mathcal{S}\,$tandardaufgabe &--& $H^{\mathcal{S}}$, $\mathscr{H}^{\mathcal{S}}$,
                                 $\mathscr{A}^{\mathcal{S}}_{\rm adm}$, $\mathscr{A}^{\mathcal{S}}_{\rm Lip}$, \\[2mm]
$\mathcal{M}\,$ultiprozesse &--& $H^{\mathcal{M}}$, $\mathscr{H}^{\mathcal{M}}$,
                               $\mathscr{A}^{\mathcal{M}}_{\rm adm}$, $\mathscr{A}^{\mathcal{M}}_{\rm Lip}$, \\[2mm]
$\mathcal{U}\,$nendlicher Zeithorizont &--& $H^{\mathcal{U}}$, $\mathscr{H}^{\mathcal{U}}$,
                                          $\mathscr{A}^{\mathcal{U}}_{\rm adm}$, $\mathscr{A}^{\mathcal{U}}_{\rm Lip}$,
                                          $\mathscr{A}^{\mathcal{U}}_{\lim}$, \\[2mm]
$\mathcal{Z}\,$eitverzögerte Systeme &--& $H^{\mathcal{Z}}$, $\mathscr{H}^{\mathcal{Z}}$,
                                        $\mathscr{A}^{\mathcal{Z}}_{\rm adm}$, $\mathscr{A}^{\mathcal{Z}}_{\rm Lip}$, \\[2mm]
$\mathcal{I}\,$ntegralgleichungen &--& $H^{\mathcal{I}}$, $\mathscr{H}^{\mathcal{I}}$,
                                     $\mathscr{A}^{\mathcal{I}}_{\rm adm}$, $\mathscr{A}^{\mathcal{I}}_{\rm Lip}$.
\end{tabular}

\newpage
\addcontentsline{toc}{section}{Stückweise Stetigkeit und Differenzierbarkeit}
\lhead[\thepage \hspace*{1mm} Stückweise Stetigkeit und Differenzierbarkeit]{}
\rhead[]{Stückweise Stetigkeit und Differenzierbarkeit \hspace*{1mm} \thepage}
\section*{Stückweise Stetigkeit und Differenzierbarkeit}
Die vorliegende Untersuchung von Steuerungsproblemen mit freiem rechten Endpunkt betten wir in den
Rahmen stückweise stetiger und stückweise stetig differenzierbarer Funktionen ein.
\index{Raum, beschränkter Funktionen!ststetig@-- stückweise stetiger Funktionen}
\index{Raum, beschränkter Funktionen!ststetigdiff@-- stückweise stetig differenzierbarer Funktionen}
Der Begriff ``stückweise'' bezieht sich dabei auf das Vorhandensein von höchstens endlich vielen Unstetigkeiten
in Form von Sprungstellen.
Wir räumen diesen -- für diese Ausarbeitung zentralen -- Begriffen eigenständige Definitionen ein:

\begin{definition}[Stückweise Stetigkeit] \label{DefinitionSTSTFunktion}
Die Funktion $f(\cdot): [a,b] \to \R$ heißt stückweise stetig, \index{Funktion, stückweise stetige}
wenn sie in endlich vielen Stellen $a<s_1<s_2<...<s_N<b$ Sprünge besitzt,
d.\,h. in diesen Stellen existieren beide einseitigen Grenzwerte von $f$ im eigentlichen Sinn.
In den Stellen $s_1,...,s_N$ wählen die Werte der Funktion $f$ so, dass $f$ rechtsseitig stetig ist.
\end{definition}

\begin{definition}[Stückweise stetige Differenzierbarkeit] \label{DefinitionSTSTDFunktion}
Die Funktion $f(\cdot): [a,b] \to \R$ heißt stückweise stetig differenzierbar,
\index{Funktion, stückweise stetige!ststetigdiff@--, stückweise stetig differenzierbare}
wenn sie auf $[a,b]$ stetig und in den endlich vielen Teilintervallen $(a,s_1),(s_1,s_2),...,(s_N,b)$ stetig differenzierbar ist,
sowie ihre Ableitung eine stückweise stetige und in den Stellen $s_1,...,s_N$ rechtsseitig stetige Funktion über $[a,b]$ ist.
\end{definition}

Unter den wichtigen Bezeichnungen haben wir die Funktionenräume \\[2mm]
\begin{tabular}{lll}
$PC([a,b],\R^n)$ &--& Raum der stückweise stetigen Funktionen, \\[1mm]
$PC_1([a,b],\R^n)$ &--& Raum der stückweise stetig differenzierbaren Funktionen \\[2mm]
\end{tabular}

über $[a,b]$ bereits aufgeführt.
In der Aufgabe mit unendlichem Zeithorizont wird das Intervall $[a,b]$ durch $[0,\infty)$ ersetzt.
Es ergeben sich die Funktionenräume  \\[2mm]
\begin{tabular}{lll}
$PC([0,\infty),\R^n)$ &--& Raum der stückweise stetigen Funktionen, \\[1mm]
$PC_1([0,\infty),\R^n)$ &--& Raum der stückweise stetig differenzierbaren Funktionen \\[2mm]
\end{tabular}

über $[0,\infty)$.
In dem Fall des unbeschränkten Intervalls $[0,\infty)$ definieren wir:
\begin{definition}[Stückweise Stetigkeit] \label{DefinitionSTSTFunktion2}
Die Funktion $f(\cdot): [0,\infty) \to \R$ heißt stückweise stetig, 
wenn sie über $[0,\infty)$ beschränkt und über jedem endlichen Intervall $[0,T]$ stückweise stetig ist.
\end{definition}

\begin{definition}[Stückweise stetige Differenzierbarkeit] \label{DefinitionSTSTDFunktion2}
Die Funktion $f(\cdot): [0,\infty) \to \R$ heißt stückweise stetig differenzierbar,
wenn sie über $[0,\infty)$ beschränkt und über jedem endlichen Intervall $[0,T]$ stückweise stetig differenzierbar ist.
\end{definition}
\cleardoublepage

\lhead[ ]{}
\rhead[]{Standardaufgabe \hspace*{1mm} \thepage}
\section{Eine Standardaufgabe mit Wiedergewinnungswert}       
       gIn diesem Abschnitt behandeln wir eine ``einfache'' Aufgabe der Optimalen Steuerung.  \label{AbschnittPMPeinfach}
Darin besteht das Optimierungskriterium aus einem Zielfunktional in Integralform und einem zusätzlichen Terminalfunktional im Endpunkt des Zustandes.
Das Terminalfunktional \index{Terminalfunktional}
bezeichnet in ökonomischen Anwendungen häufig einen Schrotterlös (``scrap value'')
oder den Wert von wiederverwendbaren Materialien oder Teilen (``salvage value'') einer nicht weiter nutzbaren Maschine
am Ende des Planungszeitraumes.
Weiterhin liegt bezüglich dem Zustand ``lediglich'' eine Nebenbedingung in Form einer Differentialgleichung vorliegt.
Diese nennen wir die Dynamik.
Ferner weist die Aufgabenstellung Steuerbeschränkungen auf.
Da keine weiteren Einschränkungen, insbesondere im Endzeitpunkt an den Zustand, vorliegen,
ist die Methode der einfachen Nadelvariation zur Herleitung von notwendigen Optimalitätsbedingungen für ein starkes lokales Minimum dieser Aufgabe gut geeignet.
Die vorliegende Herangehensweise ist Ioffe \& Tichomirov \cite{Ioffe} entnommen.
Aber im Gegensatz zu \cite{Ioffe} beschränken wir uns auf den sogenannten ``normalen'' Fall notwendiger Bedingungen,
welcher in Aufgaben mit freiem Endpunkt vorliegt.
       \lhead[\thepage \hspace*{1mm} Pontrjaginsches Maximumprinzip]{ }
       \subsection{Die Aufgabenstellung und das Pontrjaginsche Maximumprinzip} 
Über dem gegebenen Intervall $[t_0,t_1]$ betrachten wir zu $x_0 \in \R^n$ die Aufgabe 
\begin{eqnarray}
&&\label{PMPeinfach1} J\big(x(\cdot),u(\cdot)\big) = \int_{t_0}^{t_1} f\big(t,x(t),u(t)\big) \, dt +S\big(x(t_1)\big) \to \inf, \\
&&\label{PMPeinfach2} \dot{x}(t) = \varphi\big(t,x(t),u(t)\big), \quad x(t_0)=x_0,\\
&&\label{PMPeinfach3} u(t) \in U \subseteq \R^m, \quad U \not= \emptyset.
\end{eqnarray}
Die Aufgabe (\ref{PMPeinfach1})--(\ref{PMPeinfach3}) untersuchen wir bez"uglich dem Zustand $x(\cdot) \in PC_1([t_0,t_1],\R^n)$
und der Steuerung $u(\cdot) \in PC([t_0,t_1],U)$.
Das Paar $\big(x(\cdot),u(\cdot)\big)$ nennen wir einen Steuerungsprozess. \\[2mm]
Mit $\mathscr{A}^{\,\mathcal{S}}_{\rm Lip}$ bezeichnen wir die Menge aller Paare $\big(x(\cdot),u(\cdot)\big)$,
für die es ein $\gamma>0$ derart gibt,
dass die Abbildungen $f(t,x,u)$, $\varphi(t,x,u)$ auf der Menge aller $(t,x,u) \in \R \times \R^n \times \R^m$ mit
$$t \in [t_0,t_1], \qquad \|x-x(t)\| < \gamma, \qquad u \in \R^m$$
stetig in der Gesamtheit aller Variablen und stetig differenzierbar bezüglich $x$ sind. \\[2mm]
Das Paar $\big(x(\cdot),u(\cdot)\big) \in PC_1([t_0,t_1],\R^n) \times PC([t_0,t_1],U)$
hei"st ein zul"assiger Steuerungsprozess in der Aufgabe (\ref{PMPeinfach1})--(\ref{PMPeinfach3}),
falls $\big(x(\cdot),u(\cdot)\big)$ der Dynamik (\ref{PMPeinfach2}) zu $x(t_0)=x_0$ gen"ugt.
Mit $\mathscr{A}^{\mathcal{S}}_{\rm adm}$ bezeichnen wir die Menge der zul"assigen Steuerungsprozesse.

\begin{bemerkung}
{\rm Die Reduktion auf einen Umgebungsstreifen um eine Trajektrie $x(\cdot)$ ist angebracht,
wenn die Abbildungen nicht willkürliche Werte für den Zustand annehmen dürfen. 
Dies ist im Beispiel \ref{BeispielKonInv} erforderlich,
in dem die Abbildung $x \to x^\alpha$ mit $\alpha \in (0,1)$ nur für $x \geq 0$ wohldefiniert und die Ableitung $x \to \alpha x^{\alpha-1}$
nur für $x \geq \delta >0$ mit einem gegebenen $\delta>0$ beschränkt sind. \hfill $\square$}
\end{bemerkung}

Ein zul"assiger Steuerungsprozess $\big(x_*(\cdot),u_*(\cdot)\big)$ ist eine
starke lokale Minimalstelle\index{Minimum, starkes lokales!elementar@-- Standardaufgabe}
der Aufgabe (\ref{PMPeinfach1})--(\ref{PMPeinfach3}),
falls eine Zahl $\varepsilon > 0$ derart existiert, dass die Ungleichung 
$$J\big(x(\cdot),u(\cdot)\big) \geq J\big(x_*(\cdot),u_*(\cdot)\big)$$
f"ur alle $\big(x(\cdot),u(\cdot)\big) \in \mathscr{A}^{\mathcal{S}}_{\rm adm}$ mit $\|x(\cdot)-x_*(\cdot)\|_\infty < \varepsilon$ gilt.\\[2mm]
Es bezeichnet $H: \R \times \R^n \times \R^m \times \R^n \to \R$ die Pontrjagin-Funktion
\begin{equation} \label{PontrjaginFunktionGA}
H^{\mathcal{S}}(t,x,u,p) = \langle p, \varphi(t,x,u) \rangle - f(t,x,u).
\end{equation}

\begin{theorem}[Pontrjaginsches Maximumprinzip] \label{SatzPMPeinfach}
\index{Pontrjaginsches Maximumprinzip!Standard@-- Standardaufgabe} 
Sei $\big(x_*(\cdot),u_*(\cdot)\big) \in \mathscr{A}^{\mathcal{S}}_{\rm adm} \cap \mathscr{A}^{\mathcal{S}}_{\rm Lip}$. 
Ist $\big(x_*(\cdot),u_*(\cdot)\big)$ ein starkes lokales Minimum der Aufgabe (\ref{PMPeinfach1})--(\ref{PMPeinfach3}),
dann existiert eine Vektorfunktion $p(\cdot) \in PC_1([t_0,t_1],\R^n)$ derart, dass
\begin{enumerate}
\item[(a)] die adjungierte Gleichung
           \index{adjungierte Gleichung!Standard@-- Standardaufgabe}
           \begin{equation}\label{PMPeinfach4} 
           \dot{p}(t) = -H_x^{\mathcal{S}}\big(t,x_*(t),u_*(t),p(t)\big),
           \end{equation} 
\item[(b)] in $t=t_1$ die Transversalitätsbedingung
           \index{Transversalitätsbedingungen!Standard@-- Standardaufgabe}
           \begin{equation}\label{PMPeinfach5} 
           p(t_1)=-S'\big(x_*(t_1)\big)
           \end{equation} 
\item[(c)] und in fast allen Punkten $t \in [t_0,t_1]$ die Maximumbedingung
           \index{Maximumbedingung!Standard@-- Standardaufgabe}
           \begin{equation}\label{PMPeinfach6} 
           H^{\mathcal{S}}\big(t,x_*(t),u_*(t),p(t)\big) = \max_{u \in U}H^{\mathcal{S}}\big(t,x_*(t),u,p(t)\big).
           \end{equation}
\end{enumerate}
erfüllt sind.
\end{theorem}

\begin{bemerkung}
{\rm Im dem Fall, dass in (\ref{PMPeinfach1}) das Terminalfunktional $S$ nicht vorkommt,
besitzt die Transversalitätsbedingung in (\ref{PMPeinfach4}) die Form $p(t_1)=0$. \hfill $\square$}
\end{bemerkung}

Die nachstehenden Investitionsmodelle sind Seierstad \& Syds\ae ter \cite{Seierstad} entnommen:
\begin{beispiel} \label{BeispielLinInv}
{\rm Es bezeichne $K(t)$ das Kapital und
$u(t)$ die Investitions- bzw. $\big(1-u(t)\big)$ die Konsumptionsrate zum Zeitpunkt $t$.
Wir betrachten das lineare Investitionsmodell
\begin{eqnarray}
&&\label{Beispiel1einfach1} J\big(K(\cdot),u(\cdot)\big) = \int_0^T \big( 1-u(t)\big) \cdot K(t) \, dt \to \sup, \\
&&\label{Beispiel1einfach2} \dot{K}(t) = u(t) \cdot K(t), \quad K(0)=K_0 >0,\\
&&\label{Beispiel1einfach3} u(t) \in [0,1], \quad T > 1 \mbox{ fest}.
\end{eqnarray}
Um zu einem Minimierungsproblem "uberzugehen, 
multiplizieren wir den Integranden in (\ref{Beispiel1einfach1}) mit $-1$.
Wir erhalten damit in der Aufgabe (\ref{Beispiel1einfach1})--(\ref{Beispiel1einfach3}) die Abbildungen
$$f(t,K,u) = -( 1-u) \cdot K, \quad \varphi(t,K,u) = u \cdot K.$$
Dann gilt für die Pontrjaginsche Funktion (\ref{PontrjaginFunktionGA}) f"ur $u = u_*(t)$ nach (\ref{PMPeinfach6})
\begin{equation} \label{Beispiel1einfach4}
H^{\mathcal{S}}\big(t,K_*(t),u_*(t),p(t)\big) = \max_{u \in [0,1]} \big[ \big( p(t) -1 \big) u + 1 \big] \cdot K_*(t)
\end{equation}
und $p(\cdot)$ ist die L"osung der adjungierten Gleichung zur Transversalit"atsbedingung in (\ref{PMPeinfach4}):
\begin{equation} \label{Beispiel1einfach5}
\dot{p}(t) = - u_*(t) \cdot p(t) - \big( 1 - u_*(t)\big) = -1 - u_*(t) \big( p(t) - 1 \big), \quad p(T)=0.
\end{equation}
Aus der Maximumbedingung (\ref{Beispiel1einfach4}) k"onnen wir unmittelbar
$$u_*(t) = 0 \mbox{ f"ur } p(t) < 1 \quad\mbox{und}\quad u_*(t) = 1 \mbox{ f"ur } p(t) > 1$$
entnehmen. 
Für die adjungierte Gleichung (\ref{Beispiel1einfach5}) folgt damit sofort
$$\dot{p}(t) = \left\{ \begin{array}{ll}
               -1 ,   & p(t) < 1 \\
               -p(t), & p(t) > 1 \end{array} \right\} = - \max \{1, p(t) \} \leq -1 <0 \mbox{ für alle } t \in (0,T).$$
Die Stelle $\tau_1$ mit $p(\tau_1)=1$ ist wegen $\dot{p}(t)<-1$ und $p(T)=0$ eindeutig bestimmt. \\
Aus formalen Gründen betrachten wir zur Bestimmung von $\tau_1$ die Funktion $g$ mit
$$\dot{g}(t) = -1 \mbox{ und } g(T) = 0 \quad\Rightarrow\quad g(t)=T-t \mbox{ für } t \leq T.$$
Also gibt es genau eine L"osung $\tau_1 \in (0,T)$ der Gleichung $g(t)=1$ und es gilt 
\begin{equation} \label{Beispiel1einfach6} \tau_1 = T-1. \end{equation}
Da die adjungierte Funktion $p(\cdot)$ streng monoton fallend ist, folgt f"ur diese
\begin{equation} \label{Beispiel1einfach7}
\dot{p}(t) = \left\{ \begin{array}{ll}
               -p(t), & t \in (0,T-1), \\
               -1,    & t \in (T-1,T), \end{array} \right.
\qquad
p(t) = \left\{ \begin{array}{ll}
               e^{T-(1+t)}, & t \in [0,T-1), \\
               T-t,         & t \in [T-1,T]. \end{array} \right.
\end{equation}
Mit der Steuerung
\begin{equation} \label{Beispiel1einfach8}
u_*(t) = \left\{ \begin{array}{ll}
               1 & \mbox{ f"ur } t \in [0,T-1), \\
               0 & \mbox{ f"ur } t \in [T-1,T], \end{array} \right.
\end{equation}
erhalten wir aus (\ref{Beispiel1einfach2}) f"ur das Kapital
\begin{equation} \label{Beispiel1einfach9}
\dot{K}_*(t) = \left\{ \begin{array}{ll}
                K_*(t), & t \in (0,T-1), \\
                0,     & t \in (T-1,T), \end{array} \right.
\qquad
K_*(t) = \left\{ \begin{array}{ll}
               K_0 \cdot e^t     & t \in [0,T-1), \\
               K_0 \cdot e^{T-1} & t \in [T-1,T], \end{array} \right.
\end{equation}
und f"ur den Wert des Zielfunktionals (\ref{Beispiel1einfach1})
\begin{equation} \label{Beispiel1einfach10}
J\big(K_*(\cdot),u_*(\cdot)\big) = \int_0^T \big( 1-u_*(t)\big) \cdot K_*(t) \, dt = \int_{T-1}^{T} K_0 \cdot e^{T-1} \, dt = K_0 \cdot e^{T-1}.
\end{equation}
In diesem Beispiel ist die optimale Steuerung $u_*(\cdot)$ unstetig und nimmt ausschlie"slich Werte auf dem Rand des kompakten
Steuerbereiches $U=[0,1]$ an. \hfill $\square$}
\end{beispiel}

\begin{beispiel} \label{BeispielKonInv}
{\rm Im Gegensatz zum vorherigen linearen Modell nehmen wir an,
dass wir mit einer einprozentigen Erh"ohung des Kapitals $K$ eine konstante prozentuale Erh"ohung des Nutzens $U$ um $\alpha \in (0,1)$ erreichen:
$$\frac{\Delta U}{U} = \alpha \frac{\Delta K}{K} \mbox{ bzw. } \frac{dU}{U} = \alpha \frac{dK}{K}
  \quad \Longrightarrow \quad U(K) = c \cdot K^\alpha.$$
Damit erhalten wir folgende Aufgabe mit der Cobb-Douglas-Funktion $U(K)=K^\alpha$:
\begin{eqnarray}
&&\label{Beispiel2einfach1} J\big(K(\cdot),u(\cdot)\big) = \int_0^T \big( 1-u(t)\big) \cdot K^\alpha(t) \, dt \to \sup, \\
&&\label{Beispiel2einfach2} \dot{K}(t) = u(t) \cdot K^\alpha(t), \quad K(0)=K_0 >0,\\
&&\label{Beispiel2einfach3} u(t) \in [0,1], \quad \alpha \in (0,1) \mbox{ konstant}, \quad T \mbox{ fest mit } \alpha T -K_0^{1-\alpha} > 0.
\end{eqnarray}
Wir multiplizieren den Integranden in (\ref{Beispiel2einfach1}) mit $-1$ und erhalten
$$f(t,K,u) = -( 1-u) \cdot K^\alpha, \quad \varphi(t,K,u) = u \cdot K^\alpha.$$
Dann gilt für die Pontrjaginsche Funktion (\ref{PontrjaginFunktionGA}) nach (\ref{PMPeinfach6}) f"ur $u = u_*(t)$
\begin{equation}\label{Beispiel2einfach5}
H^{\mathcal{S}}\big(t,K_*(t),u_*(t),p(t)\big) =  \max_{u \in [0,1]} \big[ \big( p(t) -1 \big) u + 1 \big] \cdot K_*^\alpha(t)
\end{equation}
und $p(\cdot)$ ist die L"osung der adjungierten Gleichung zur Transversalit"atsbedingung (\ref{PMPeinfach4}):
\begin{equation}\label{Beispiel2einfach6}
\dot{p}(t) = - \alpha K_*^{\alpha-1}(t) \big[ 1 + u_*(t) \big( p(t) - 1 \big) \big], \quad p(T)=0.
\end{equation}
Aus der Maximumbedingung (\ref{Beispiel2einfach5}) ergibt sich unmittelbar
$$u_*(t) = 0 \mbox{ f"ur } p(t) < 1 \quad\mbox{und}\quad u_*(t) = 1 \mbox{ f"ur } p(t) > 1.$$
Für die adjungierte Gleichung (\ref{Beispiel2einfach6}) folgt weiterhin
$$\dot{p}(t) = \left\{ \begin{array}{ll}
               -\alpha \, K_*^{\alpha-1}(t),            & p(t) < 1 \\
               -\alpha \, K_*^{\alpha-1}(t) \cdot p(t), & p(t) > 1 \end{array} \right\}
             = -\alpha \, K_*^{\alpha-1}(t) \cdot \max \{ 1,p(t) \}<0 \mbox{ über } (0,T).$$
Die Adjungierte $p(\cdot)$ ist stetig und streng monoton fallend über $[0,T]$ mit $p(T)=0$. \\[2mm]
Wir zeigen nun die Existenz einer Zahl $\tau_2 \in (0,T)$ mit $p(\tau_2)=1$.
Es sei $\tau_2 \in [0,T]$ mit
$$\tau_2 = \inf \{ t \in [0,T] \,|\, p(t)<1 \}.$$
F"ur $t > \tau_2$ gilt $u_*(t)=0$ und es ist $K_*(t)$ auf $[\tau_2,T]$ konstant.
Es bezeichne $K=K_*(\tau_2)$ diesen konstanten Wert.
Anstelle der Adjungieren $p(\cdot)$ betrachten wir mit Hilfe der Konstanten $K$ zun"achst die Funktion $g$ mit
$$\dot{g}(t) = -\alpha \, K^{\alpha-1}, \quad g(T) = 0, \quad \alpha \in (0,1), \quad K>0 \mbox{ konstant}.$$
Es ergibt sich $g(t) = \alpha K^{\alpha-1} (T-t)$.
Nun bestimmen wir $\tau_2$.
F"ur $t < \tau_2$ ist $p(t) >1$ und damit $u_*(t) = 1$.
In der Dynamik (\ref{Beispiel2einfach2}) sei $u(t) \equiv 1$.
Es ergibt sich die Funktion $h(\cdot)$ mit $\dot{h}(t) = h^\alpha(t)$, $\alpha \in (0,1)$ und $h(0) = K_0$.
Wir erhalten nach Trennung der Ver"anderlichen und anschlie"sender Integration
$$\int_{K_0}^h x^{-\alpha}\, dx = \int_0^t 1 \, ds \quad\Rightarrow\quad 
  \frac{1}{1-\alpha}\big[h^{1-\alpha}- K_0^{1-\alpha}\big] = t\quad\Rightarrow\quad 
  h(t) = \big[ (1-\alpha) t + K_0^{1-\alpha} \big]^{\frac{1}{1-\alpha}}.$$
Aus $K=h(\tau_2)$ und $g(\tau_2)=1$ ergibt sich f"ur $\tau_2$ die Gleichung
$$1 = \underbrace{\alpha K^{\alpha-1} (T-\tau_2)}_{=g(\tau_2)}
    = \alpha \cdot \underbrace{\big[ (1-\alpha) \tau_2 + K_0^{1-\alpha} \big]^{\frac{\alpha-1}{1-\alpha}}}_{=h^{\alpha-1}(\tau_2)} \cdot (T-\tau_2)
    = \frac{\alpha (T-\tau_2)}{(1-\alpha) \tau_2 + K_0^{1-\alpha}}.$$
Aus dieser Gleichung folgt zusammen mit der Annahme $\alpha T -x_0^{1-\alpha} > 0$ in (\ref{Beispiel2einfach3})
\begin{equation} \label{Beispiel2einfach7} \tau_2 = \alpha T - K_0^{1-\alpha}>0. \end{equation}
F"ur die Analyse der Aufgabe liefert dies die Steuerung
\begin{equation} \label{Beispiel2einfach9}
u_*(t) = \left\{ \begin{array}{ll}
               1 & \mbox{ f"ur } t \in [0,\tau_2), \\
               0 & \mbox{ f"ur } t \in [\tau_2,T]. \end{array} \right.
\end{equation}
Mit $u_*(\cdot)$ erhalten wir in (\ref{Beispiel2einfach2}) f"ur die "Anderung des Kapitalstocks
$$\dot{K}_*(t) = \left\{ \begin{array}{ll}
                K_*^{\alpha}(t) & \mbox{ f"ur } t \in (0,\tau_2), \\
                0      & \mbox{ f"ur } t \in (\tau_2,T). \end{array} \right.$$
Wenden wir über $[0,\tau_2]$ wieder die Trennung der Ver"anderlichen wie für die Funktion $h$ an,
so ergibt sich f"ur die zeitliche Entwicklung des optimalen Kapitalstocks
\begin{equation} \label{Beispiel2einfach10}
K_*(t) = \left\{ \begin{array}{ll}
               \big[ (1-\alpha)t +  K_0^{1-\alpha} \big]^\frac{1}{1-\alpha} & \mbox{ f"ur } t \in [0,\tau_2), \\
               \big[ \alpha(1-\alpha)T + \alpha K_0^{1-\alpha} \big]^\frac{1}{1-\alpha} & \mbox{ f"ur } t \in [\tau_2,T]. \end{array} \right.
\end{equation}
Ferner liefert dies in der adjungierten Gleichung
$$\dot{p}(t) = \left\{ \begin{array}{ll}
               -\alpha \, K_*^{\alpha-1}(t) \cdot p(t) & \mbox{ f"ur } t \in (0,\tau_2), \\
               -\alpha \, K_*^{\alpha-1}(\tau_2)       & \mbox{ f"ur } t \in (\tau_2,T). \end{array} \right.$$
Zur Berechnung von $p(\cdot)$ auf $[0,\tau_2]$ verwenden wir $p(\tau_2)=1$ und erhalten wieder mit
einer Trennung der Ver"anderlichen und anschlie"sender Integration
$$\ln\left(\frac{p(t)}{p(\tau_2)}\right) = -\alpha \int_{\tau_2}^t \frac{1}{(1-\alpha)s +K_0^{1-\alpha}} \, ds
  = -\frac{\alpha}{1-\alpha} \ln\left( \frac{(1-\alpha)t +K_0^{1-\alpha}}{(1-\alpha)\tau_2 +K_0^{1-\alpha}}\right).$$
Daraus folgt
$$p(t) = \left[ \frac{(1-\alpha)\tau_2 +K_0^{1-\alpha}}{(1-\alpha)t + K_0^{1-\alpha}} \right]^\frac{\alpha}{1-\alpha}
       = \left[ \frac{\alpha(1-\alpha)T + \alpha K_0^{1-\alpha}}{(1-\alpha)t + K_0^{1-\alpha}} \right]^\frac{\alpha}{1-\alpha}$$
auf $[0,\tau_2]$. Auf $[0,T]$ gilt damit f"ur die adjungierte Funktion
\begin{equation} \label{Beispiel2einfach8}
p(t) = \left\{ \begin{array}{ll}  
         \displaystyle \left[ \frac{\alpha(1-\alpha)T + \alpha K_0^{1-\alpha}}{(1-\alpha)t + K_0^{1-\alpha}} \right]^\frac{\alpha}{1-\alpha}
                                                        & \mbox{ f"ur } t \in [0,\tau_2), \\[4mm]
         \displaystyle \frac{T-t}{(1-\alpha)T + K_0^{1-\alpha}} & \mbox{ f"ur } t \in [\tau_2,T]. \end{array} \right.
\end{equation}
Für den Wert des Zielfunktionals (\ref{Beispiel2einfach1}) ergibt sich in diesem Beispiel
\begin{eqnarray}
J\big(K_*(\cdot),u_*(\cdot)\big) &=& \int_0^T \big( 1-u_*(t)\big) \cdot K^{\alpha}_*(t) \, dt
                                     = \int_{\tau_2}^{T} \big[ \alpha(1-\alpha)T + \alpha K_0^{1-\alpha} \big]^\frac{\alpha}{1-\alpha} \, dt \nonumber \\
                                 &=& (T- \alpha T +K_0^{1-\alpha}) \cdot \big[ \alpha(1-\alpha)T + \alpha K_0^{1-\alpha} \big]^\frac{\alpha}{1-\alpha}
                                     \nonumber \\
\label{Beispiel2einfach11}       &=& \alpha^{\frac{\alpha}{1-\alpha}} \cdot \big[ (1-\alpha)T + K_0^{1-\alpha} \big]^\frac{1}{1-\alpha}.
\end{eqnarray}
Die Auswertung des Pontrjaginschen Maximumprinzip ist damit abgeschlossen. \hfill $\square$}
\end{beispiel}

\begin{beispiel} \label{Beispiel2Sektoren}
{\rm Wir betrachten nach Seierstad \& Syds\ae ter \cite{Seierstad} ein Zwei-Sektoren-Modell,
das aus der Produktion von Investitons- und Konsumg"uter besteht.
Wir bezeichnen mit $x(t)$ bzw. $y(t)$ die Rate der G"uterproduktion im Investitions- bzw. Konsumsektor
und es beschreibe $u(t) \in [0,1]$ die Aufteilung der Investitionsg"uter auf beide Sektoren zur Produktion.
Es entsteht damit die Aufgabe
\begin{equation}
\left.\begin{array}{l}
\hspace*{-5mm}
\displaystyle J\big(x(\cdot),y(\cdot),u(\cdot)\big) =\int_0^T y(t) \, dt \to \sup,  \\[3mm]
\hspace*{-5mm}   
\displaystyle \dot{x}(t)=u(t)x(t), \quad x(0)=x_0>0,  \\[2mm]
\hspace*{-5mm}   
\displaystyle \dot{y}(t)=\big(1-u(t)\big)x(t), \quad y(0)=y_0>0,  \\[2mm]
\hspace*{-5mm}   
u \in [0,1], \qquad T>2.
\end{array} \right\}
\end{equation}
Wir gehen zu einem Minimierungsproblem über und wenden Theorem \ref{SatzPMPeinfach} an. \\
Die Pontrjagin-Funktion besitzt die Gestalt
$$H^{\mathcal{S}}(t,x,y,u,p,q) = uxp+(1-u)xq+y.$$
F"ur die Adjungierten $p(\cdot)$ und $q(\cdot)$ ergeben sich die Gleichungen
$$\dot{p}(t)=-u_*(t)p(t)-\big(1-u_*(t)\big)q(t),\quad p(T)=0,\qquad \dot{q}(t)=-1, \quad q(T)=0.$$
Es folgt unmittelbar $q(t)=T-t$.
Die Maximumbedingung ist äquivalent zu
$$\max_{u \in [0,1]} u \cdot \big[p(t)-q(t)\big] \cdot x(t)$$
und wir erhalten für die optimale Steuerung $u_*(\cdot)$ 
$$u_*(t) = 1 \mbox{ f"ur } p(t) > q(t) \quad\mbox{und}\quad u_*(t) = 0 \mbox{ f"ur } p(t) < q(t).$$
Wegen $p(T)=0$ und $q(t)=T-t$ ist $\dot{p}(t)>-1$ über einem gewissen Intervall $(\tau,T)$.
Demzufolge gilt $p(t)<q(t)$ und $\dot{p}(t)=-(T-t)$ über $(\tau,T)$. \\[2mm]
Zur Bestimmung des Zeitpunktes $t$ mit $p(t)=q(t)$ nutzen wir wieder eine Nebenrechnung.
Dazu betrachten wir die Funktion $g(t)$ mit $\dot{g}(t)=-(T-t)$ und $g(T)=0$.
Es ergeben sich
$$g(t)=\frac{1}{2}(T-t)^2 \quad\mbox{ und }\quad g(t)=q(t) \;\Leftrightarrow\; t =T-2.$$
Die Adjungierten $p(\cdot)$, $q(\cdot)$ sind streng monoton fallen und es gilt $p(T-2)=q(T-2)=2$.
Es sind damit $\dot{p}(t)<-2$ für $t \in (0,T-2)$ und $p(t)> q(t)$ über $[0,T-2)$.
Wir erhalten zusammenfassend für die adjungierten Funktionen 
$$p(t)= \left\{\begin{array}{ll} 2e^{T-(t+2)}, & t \in [0,T-2), \\[1mm] \frac{1}{2}(T-t)^2, & t \in [T-2,T], \end{array}\right.
  \qquad q(t)= T-t,$$
für die optimale Steuerung 
$$u_*(t) = 1 \mbox{ f"ur } t \in [0,T-2) \quad\mbox{und}\quad u_*(t) = 0 \mbox{ f"ur } t \in [T-2,T],$$
für die zugeh"origen Zustandstrajektorien
\begin{eqnarray*}
x_*(t) &=& \left\{\begin{array}{ll} x_0 e^t, & t \in [0,T-2), \\[1mm] x_0e^{T-2}, & t \in [T-2,T], \end{array}\right. \\
y_*(t) &=& \left\{\begin{array}{ll} y_0, & t \in [0,T-2), \\[1mm] y_0+(t-T+2)x_0e^{T-2}, & t \in [T-2,T], \end{array}\right.
\end{eqnarray*}
und für den Wert des Zielfunktionals
$$J\big(x_*(\cdot),y_*(\cdot),u_*(\cdot)\big) = \int_0^T y_*(t) \, dt= y_0 T+ \frac{1}{2}(T-t-2)^2x_0e^{T-2}.$$
Unsere Untersuchung der Aufgabe ist abgeschlossen. \hfill $\square$}
\end{beispiel}
       \newpage
       \lhead[\thepage \hspace*{1mm} Beweis des Maximumprinzips]{ }
       \subsection{Der Beweis des Maximumprinzips} \label{AbschnittPMPBeweiseinfach}
Der nachstehende Beweis ist \cite{Ioffe} entnommen.
Wir haben lediglich die Transversalitätsbedingung (\ref{PMPeinfach4}) bezüglich dem Terminalfunktional $S$ hinzugefügt. \\
Da s"amtliche Abbildungen stetig und stetig differenzierbar bezüglich $x$ sind, und die Funktion $u_*(\cdot)$ dem Raum $PC([t_0,t_1],U)$ angeh"ort,
sind die Abbildungen $t \to f\big(t,x_*(t),u_*(t)\big)$ und $t \to \varphi\big(t,x_*(t),u_*(t)\big)$ über $[t_0,t_1]$ lediglich stückweise stetig und nicht stetig.
Deswegen verweisen wir auf die Ergebnisse im Anhang \ref{AnhangDGL} über Differentialgleichungen mit stückweise stetigen rechten Seiten. \\
Ferner erf"ullt die adjungierte Gleichung (\ref{PMPeinfach4}) die Voraussetzungen von Lemma \ref{LemmaDGL3} und Lemma \ref{LemmaDGL5}.
Daher gibt es eine eindeutige L"osung $p(\cdot)$ der Gleichung (\ref{PMPeinfach4}) zur Randbedingung $p(t_1)=-S'\big(x_*(t_1)\big)$,
die dem Raum $PC_1([t_0,t_1],\R^n)$ angehört. 
Die behauptete Existenz einer Lösung $p(\cdot) \in PC_1([t_0,t_1],\R^n)$ der adjungierten Gleichung (\ref{PMPeinfach4})
zur Transversalitätsbedingung (\ref{PMPeinfach5}) ist damit bereits gezeigt. \\[2mm]
Es sei $\tau \in (t_0,t_1)$ ein Stetigkeitspunkt der Steuerung $u_*(\cdot)$.
Dann ist $u_*(\cdot)$ auch in einer gewissen hinreichend kleinen Umgebung von $\tau$ stetig und wir w"ahlen ein festes
$\lambda$ positiv und hinreichend klein, so dass sich $\tau-\lambda$ in dieser Umgebung befindet. \\
Weiter sei nun $v$ ein beliebiger Punkt aus $U$.
Wir setzen\index{Nadelvariation, einfache}\\
\begin{minipage}{0.59\textwidth}
$$u(t;v,\tau,\lambda) = u_{\lambda}(t) = 
  \left\{ \begin{array}{ll}
          u_*(t) & \mbox{ f"ur } t \not\in [\tau-\lambda,\tau), \\
          v      & \mbox{ f"ur } t     \in [\tau-\lambda,\tau), 
          \end{array} \right.$$
und es bezeichne $x_\lambda(\cdot)$, $x_\lambda(t)=x(t;v,\tau,\lambda)$, die eindeutige L"osung der Gleichung
$$\dot{x}(t) = \varphi\big(t,x(t),u_\lambda(t)\big), \qquad x(t_0)=x_0.$$
\end{minipage}
\begin{minipage}{0.4\textwidth}
\centering
\includegraphics[width=4.5cm]{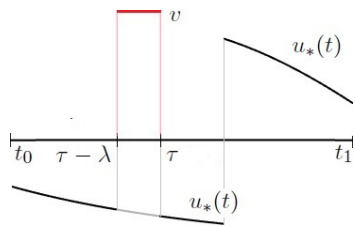}
\captionof{figure}[Einfache Nadelvariation der Steuerungstheorie]{Variation $u(t;v,\tau,\lambda)$.}
\end{minipage} \\[3mm]
Dann ist $x_{\lambda}(t) = x_*(t)$ f"ur $t_0 \leq t \leq  \tau - \lambda$.
F"ur $t \geq \tau$ betrachten wir den Grenzwert
$$y(t)=\lim_{\lambda \to 0^+}\frac{x_{\lambda}(t) - x_*(t)}{\lambda},$$
und werden zeigen, dass dieser existiert, die Funktion $y(\cdot)$ der Integralgleichung
\begin{equation} \label{BeweisPMPeinfach1}
y(t)=y(\tau) + \int_{\tau}^{t} \varphi_x\big(s,x_*(s),u_*(s)\big)\,y(s) \, ds
\end{equation}
zur Anfangsbedingung
\begin{equation} \label{BeweisPMPeinfach2}
y(\tau)=\varphi\big(\tau,x_*(\tau),v\big) - \varphi\big(\tau,x_*(\tau),u_*(\tau)\big)
\end{equation} 
gen"ugt und die Beziehung
\begin{equation} \label{BeweisPMPeinfach3}
\langle p(\tau) , y(\tau) \rangle = -\big\langle S'\big(x_*(t_1)\big), y(t_1) \big\rangle
                                    - \int_{\tau}^{t_1} \big\langle f_x\big(t,x_*(t),u_*(t)\big) , y(t) \big\rangle \, dt
\end{equation}
erf"ullt ist.

\newpage
Herleitung von (\ref{BeweisPMPeinfach2}):
Nach Wahl des Punktes $\tau$ ist $\dot{x}_*(\cdot)$ in einer Umgebung dieses Punktes
stetig und f"ur hinreichend kleine positive $\lambda$ gelten
\begin{eqnarray*}
    x_*(\tau)
&=& x_*(\tau - \lambda) + \lambda 
    \varphi\big(\tau - \lambda,x_*(\tau - \lambda),u_*(\tau - \lambda)\big) + o(\lambda), \\
    x_{\lambda}(\tau)
&=& x_*(\tau - \lambda) + \lambda
    \varphi\big(\tau - \lambda,x_*(\tau - \lambda),v\big) + o(\lambda). \hspace{13,8mm}
\end{eqnarray*}
Hieraus folgt
$$\frac{x_{\lambda}(\tau) - x_*(\tau)}{\lambda}
= \varphi\big(\tau - \lambda,x_*(\tau - \lambda),v\big) 
- \varphi\big(\tau - \lambda,x_*(\tau - \lambda),u_*(\tau - \lambda)\big)
+ \frac{o(\lambda)}{\lambda}.$$
D.\,h., dass der Grenzwert
$$y(\tau)=\lim_{\lambda \to 0^+}\frac{x_{\lambda}(\tau) - x_*(\tau)}{\lambda}$$
existiert und gleich (\ref{BeweisPMPeinfach2}) ist. \\[2mm]
Herleitung von (\ref{BeweisPMPeinfach1}):
Auf dem Intervall $[\tau,t_1]$ gen"ugen sowohl $x_*(\cdot)$ als auch $x_{\lambda}(\cdot)$ der Gleichung
$$\dot{x}(t) = \varphi\big(t,x(t),u_*(t)\big).$$
Aus den Sätzen \ref{SatzEEglobal} und \ref{SatzDGLDifferenzierbarkeit} "uber die Stetigkeit und Differenzierbarkeit
der L"osung eines Differentialgleichungssystems in Abh"angigkeit von den Anfangsdaten folgt,
dass f"ur hinreichend kleine positive $\lambda$ die Vektorfunktionen $x_{\lambda}(\cdot)$ auf $[\tau,t_1]$ definiert sind,
dass sie f"ur $\lambda \to 0^+$ gleichm"a"sig gegen $x_*(\cdot)$ konvergieren und dass der Grenzwert
$$y(t) = \lim_{\lambda \to 0^+}\frac{x_{\lambda}(t) - x_*(t)}{\lambda}$$
f"ur jedes $t \in [\tau,t_1]$ existiert. 
Weiterhin liefert Satz \ref{SatzDGLDifferenzierbarkeit}, dass
$y(\cdot)$ für alle $t \in [\tau,t_1]$ der Gleichung
$$y(t)=y(\tau) + \int_{\tau}^{t} \varphi_x\big(s,x_*(s),u_*(s)\big)\,y(s) \, ds$$
genügt, d.\,h. $y(\cdot)$ die Gleichung (\ref{BeweisPMPeinfach1}) zum Anfangswert (\ref{BeweisPMPeinfach2}) erfüllt. \\[2mm]
Herleitung von (\ref{BeweisPMPeinfach3}):
Mit (\ref{PMPeinfach4}) und (\ref{BeweisPMPeinfach1}) erhalten wir f"ur $t \geq \tau$: 
$$\frac{d}{dt} \langle p(t),y(t) \rangle 
  = \langle \dot{p}(t),y(t) \rangle + \langle p(t),\dot{y}(t) \rangle
  = \big\langle f_x\big(t,x_*(t),u_*(t)\big) , y(t) \big\rangle.$$
Daher gilt für $t \in [\tau,t_1]$ die Gleichung
$$\langle p(t) ,y(t) \rangle - \langle p(\tau) ,y(\tau) \rangle = \int_{\tau}^{t} \big\langle f_x\big(s,x_*(s),u_*(s)\big) , y(s) \big\rangle \, ds.$$
Für $t=t_1$ folgt (\ref{BeweisPMPeinfach3}) mit der Transversalitätsbedingung $p(t_1)=-S'\big(x_*(t_1)$:
$$\langle p(\tau) ,y(\tau) \rangle= \langle -S'\big(x_*(t_1) ,y(t_1) \rangle - \int_{\tau}^{t_1} \big\langle f_x\big(s,x_*(s),u_*(s)\big) , y(s) \big\rangle \, ds.$$
Beweisschluss:
Da $\big(x_*(\cdot),u_*(\cdot)\big)$ ein starkes lokales Minimum ist, ist f"ur alle hinreichend kleine und positive $\lambda$
$$\frac{J\big(x_\lambda(\cdot),u_\lambda(\cdot)\big) - J\big(x_*(\cdot),u_*(\cdot)\big)}{\lambda} \geq 0.$$
In diesem Ausdruck ergeben sich für $\lambda \to 0^+$
\begin{eqnarray*}
&& \lim_{\lambda \to 0^+}\frac{1}{\lambda} \int_{\tau - \lambda}^{\tau} 
    \Big[ f\big(t,x_\lambda(t),u_\lambda(t)\big) - f\big(t,x_*(t),u_*(t)\big) \Big] \, dt \\
&& + \lim_{\lambda \to 0^+} \frac{1}{\lambda} \int_{\tau}^{t_1} 
    \Big[ f\big(t,x_\lambda(t),u_\lambda(t)\big) - f\big(t,x_*(t),u_*(t)\big) \Big] \, dt \\
&=& f\big(\tau,x_*(\tau),v\big) - f\big(\tau,x_*(\tau),u_*(\tau)\big)
    + \int_{\tau}^{t_1} \big\langle f_x\big(t,x_*(t),u_*(t)\big) , y(t) \big\rangle \, dt
\end{eqnarray*}
und für das Terminalfunktional
$$\lim_{\lambda \to 0^+}\frac{S\big(x_\lambda(t_1)\big)-S\big(x_*(t_1)\big)}{\lambda}=\big\langle S'\big(x_*(t_1)\big), y(t_1) \big\rangle.$$
Für das Zielfunktional liefert der Grenzübergang $\lambda \to 0^+$ dies den Ausdruck.
\begin{eqnarray*}
0 &\leq& \lim_{\lambda \to 0^+} \frac{J\big(x_\lambda(\cdot),u_\lambda(\cdot)\big)- J\big(x_*(\cdot),u_*(\cdot)\big)}{\lambda} \\
  &=   & f\big(\tau,x_*(\tau),v\big) - f\big(\tau,x_*(\tau),u_*(\tau)\big) \\
  &    & + \int_{\tau}^{t_1} \big\langle f_x\big(t,x_*(t),u_*(t)\big),y(t) \big\rangle \, dt
            + \big\langle S'\big(x_*(t_1)\big), y(t_1) \big\rangle.
\end{eqnarray*}
Nach (\ref{BeweisPMPeinfach3}) besteht in dieser Ungleichung der Zusammenhang
$$-\langle p(\tau) , y(\tau) \rangle =
  \int_{\tau}^{t_1} \big\langle f_x\big(t,x_*(t),u_*(t)\big),y(t) \big\rangle \, dt + \big\langle S'\big(x_*(t_1)\big), y(t_1) \big\rangle.$$
Bei Anwendung der Gleichungen (\ref{BeweisPMPeinfach2}) und (\ref{BeweisPMPeinfach3}) ergibt sich weiter
\begin{eqnarray*}
\lefteqn{\big\langle p(\tau) , \varphi\big(\tau,x_*(\tau),u_*(\tau)\big) \big\rangle - f\big(\tau,x_*(\tau),u_*(\tau)\big)} \\
&\geq& \big\langle p(\tau) , \varphi\big(\tau,x_*(\tau),v\big) \big\rangle - f\big(\tau,x_*(\tau),v\big).
\end{eqnarray*}
Damit folgt aus der Definition der Pontrjagin-Funktion:
\begin{equation*}
H^{\mathcal{S}}\big(\tau,x_*(\tau),u_*(\tau),p(\tau)\big) \geq H^{\mathcal{S}}\big(\tau,x_*(\tau),v,p(\tau)\big).
\end{equation*}
Nun ist $\tau$ ein beliebiger Stetigkeitspunkt von $u_*(\cdot)$ und $v$ ein beliebiger Punkt der Menge $U$. 
Demzufolge ist die Beziehung (\ref{PMPeinfach6}) in allen Stetigkeitspunkten von $u_*(\cdot)$ wahr und damit ist das
Maximumprinzip bewiesen. \hfill $\blacksquare$
       \newpage
       \lhead[\thepage \hspace*{1mm} Ökonomische Deutung]{ }
       \subsection{\"Okonomische Deutung des Maximumprinzips} \label{AbschnittDeutung}
Die "okonomische Interpretation des Pontrjaginschen Maximumprinzips
\index{Pontrjaginsches Maximumprinzip!oekonomische@-- ökonomische Interpretation}
wurde von Dorfman \cite{Dorfman} angeregt.
Wir geben eine Variante an,
die sich auf den Beweis von Theorem \ref{SatzPMPeinfach} im letzten Abschnitt bezieht. \\[2mm]
Der Kern des Beweises im letzten Abschnitt war,
die Auswirkungen der Nadelvariation
$$u_{\lambda}(t) = 
  \left\{ \begin{array}{ll}
          u_*(t) & \mbox{ f"ur } t \not\in [\tau-\lambda,\tau), \\
          v      & \mbox{ f"ur } t     \in [\tau-\lambda,\tau), 
          \end{array} \right.$$
auf das Steuerungsproblem (\ref{PMPeinfach1})--(\ref{PMPeinfach3}) zu untersuchen.
Die erzeugte marginale "Anderung
$$y(t)=\lim_{\lambda \to 0^+}\frac{x_{\lambda}(t) - x_*(t)}{\lambda}$$
des optimalen Zustandes $x_*(\cdot)$ bewirkt nach (\ref{BeweisPMPeinfach3}) im Zielfunktional
\begin{equation} \label{Interpretation1}
\langle p(\tau) , y(\tau) \rangle = -\big\langle S'\big(x_*(t_1)\big), y(t_1) \big\rangle
                      - \int_{\tau}^{t_1} \big\langle f_x\big(t,x_*(t),u_*(t)\big) , y(t) \big\rangle \, dt,
\end{equation}
wobei $p(\cdot)$ die L"osung der adjungierten Gleichung (\ref{PMPeinfach4}) ist.
Au"serdem gilt f"ur den marginalen Profit die Beziehung
\begin{eqnarray}
\lefteqn{\lim_{\lambda \to 0^+} \frac{J\big(x_\lambda(\cdot),u_\lambda(\cdot)\big)- J\big(x_*(\cdot),u_*(\cdot)\big)}{\lambda}}
         \nonumber \\
  &=   & f\big(\tau,x_*(\tau),v\big) - f\big(\tau,x_*(\tau),u_*(\tau)\big) \nonumber \\
  &    & + \int_{\tau}^{t_1} \big\langle f_x\big(t,x_*(t),u_*(t)\big),y(t) \big\rangle \, dt + \big\langle S'\big(x_*(t_1)\big), y(t_1) \big\rangle \nonumber \\
  &=   & \label{Interpretation2} H^{\mathcal{S}}\big(\tau,x_*(\tau),u_*(\tau),p(\tau)\big) - H^{\mathcal{S}}\big(\tau,x_*(\tau),v,p(\tau)\big).
\end{eqnarray}
Die Entscheidung zum Zeitpunkt $\tau$ mit dem Parameter $v$ statt $u_*(\tau)$ zu steuern hat einen direkten und indirekten Effekt
(Feichtinger \& Hartl \cite{Feichtinger}, S.\,29):
\begin{enumerate}
\item[--] Der unmittelbare Effekt besteht darin, dass die Profitrate $f\big(\tau,x_*(\tau),v\big)$ erzielt wird.
\item[--] Die indirekte Wirkung manifestiert sich in der "Anderung der Kapitalbest"ande um $y(\tau)$.
          Der Kapitalstock wird durch die in $\tau$ getroffene Entscheidung transformiert,
          was bemessen mit der Adjungierten $p(\cdot)$ in (\ref{Interpretation1}) ausgedr"uckt wird.
\end{enumerate}
Die Gleichung (\ref{Interpretation1}) beschreibt damit die Opportunit"atskosten,
die sich durch die Entscheidung $v$ zum Zeitpunkt $\tau$ ergeben.
Wegen der Darstellung der Opportunit"atskosten in der Form $\langle p(\tau) , y(\tau) \rangle$ in (\ref{Interpretation1}),
wird die Adjungierte $p(\cdot)$ in der "Okonomie h"aufig als Schattenpreis bezeichnet. \\
Ferner bemisst die Funktion $v \to H^{\mathcal{S}}\big(\tau,x_*(\tau),v,p(\tau)\big)$ nach Gleichung (\ref{Interpretation2}) die Profitrate aus
direktem und indirektem Gewinn,
die sich aus der "Anderung der Kapitalbest"ande ergibt.
Die Maximumbedingung sagt daher aus,
dass die Instrumente zu jedem Zeitpunkt so eingesetzt werden sollen,
dass die totale Profitrate maximal wird.
       \newpage
       \lhead[\thepage \hspace*{1mm} Arrow-Bedingungen]{ }
       \subsection{Hinreichende Bedingungen nach Arrow} \label{AbschnittHBPMP}
Im Gegensatz zu notwendigen Bedingungen,
die Kandidaten für Optimalstellen liefern,
geben hinreichende Bedingungen Gewissheit über das Vorliegen eines Optimums.
Deswegen erweitern wir die notwendigen Optimalitätsbedingungen in Form des Pontrjaginschen Maximumprinzips
um die hinreichenden Bedingungen nach Arrow. \\
Unser\index{hinreichende Bedingungen nach Arrow!Standard@-- Standardaufgabe}
Vorgehen zur Herleitung der hinreichenden Bedingungen basiert auf der Darstellung in 
Seierstad \& Syds\ae ter \cite{Seierstad},
die wir in der Beweisführung um Argumente in Aseev \& Kryazhimskii \cite{AseKry} erg"anzt haben. \\[2mm]
Für $t \in [t_0,t_1]$ führen wir die Menge $V_\gamma(t)=\{ x \in \R^n \,|\, \|x-x_*(t)\| < \gamma\}$ ein.
Au"serdem bezeichnet $\mathscr{H}^{\mathcal{S}}$ die Hamilton-Funktion
\begin{equation} \label{PMPHamilton}
\mathscr{H}^{\mathcal{S}}(t,x,p) = \sup_{u \in U} H^{\mathcal{S}}(t,x,u,p).
\end{equation}

\begin{theorem} \label{SatzHBPMP}
In der Aufgabe (\ref{PMPeinfach1})--(\ref{PMPeinfach3}) sei
$\big(x_*(\cdot),u_*(\cdot)\big) \in \mathscr{A}^{\mathcal{S}}_{\rm adm} \cap \mathscr{A}^{\mathcal{S}}_{\rm Lip}$
und es sei $p(\cdot) \in PC_1([t_0,t_1],\R^n)$. Ferner gelte:
\begin{enumerate}
\item[(a)] Das Tripel $\big(x_*(\cdot),u_*(\cdot),p(\cdot)\big)$
           erf"ullt (\ref{PMPeinfach4})--(\ref{PMPeinfach6}) in Theorem \ref{SatzPMPeinfach}.        
\item[(b)] F"ur jedes $t \in [t_0,t_1]$ ist die Funktion $\mathscr{H}^{\mathcal{S}}\big(t,x,p(t)\big)$ konkav in $x$ auf $V_\gamma(t)$.
\item[(c)] Die Abbildung $S:\R^n\to \R$ im Zielfunktional (\ref{PMPeinfach1}) ist konvex in $x$ auf $V_\gamma(t_1)$.
\end{enumerate}
Dann ist $\big(x_*(\cdot),u_*(\cdot)\big)$ ein starkes lokales Minimum der Aufgabe (\ref{PMPeinfach1})--(\ref{PMPeinfach3}).
\end{theorem}

{\bf Beweis} Es sei $t \in [t_0,t_1]$ gegeben.
Da die Abbildung $x \to -\mathscr{H}^{\mathcal{S}}\big(t,x,p(t)\big)$ auf $V_\gamma(t)$ konvex ist, ist die Menge
$Z=\big\{ (\alpha,x) \in \R \times \R^n \,\big|\, x \in V_\gamma(t), \alpha \geq -\mathscr{H}^{\mathcal{S}}\big(t,x,p(t)\big) \big\}$
konvex und besitzt ein nichtleeres Inneres.
Wir setzen $\alpha_*= -\mathscr{H}^{\mathcal{S}}\big(t,x_*(t),p(t)\big)$.
Dann ist $\big(\alpha_*,x_*(t)\big)$ mit ein Randpunkt der Menge $Z$.
Daher existiert nach dem Trennungssatz \ref{Trennungssatz} ein nichttrivialer Vektor $\big(a_0(t),a(t)\big) \in \R \times \R^n$ mit
\begin{equation} \label{BeweisHBPMP1}
a_0(t) \alpha + \langle a(t),x\rangle \geq a_0(t) \alpha_* + \langle a(t),x_*(t)\rangle \quad \mbox{ f"ur alle } (\alpha,x) \in Z.
\end{equation}
Es ist $x_*(t)$ ein innerer Punkt der Menge $V_\gamma(t)$.
Weiterhin folgt aus den elementaren Eigenschaften konvexer Funktionen,
dass $x \to -\mathscr{H}^{\mathcal{S}}\big(t,x,p(t)\big)$ in $V_\gamma(t)$ stetig ist,
da sie auf $V_\gamma(t)$ konvex und nach oben durch $x \to -H^{\mathcal{S}}\big(t,x,u_*(t),p(t)\big)$ beschr"ankt ist. \\
Deswegen existiert ein $\delta>0$ mit $x_*(t)+\xi \in V_\gamma(t)$ und $\big(\alpha_*+1,x_*(t)+\xi\big) \in Z$
f"ur alle $\|\xi\| \leq \delta$.
Aus (\ref{BeweisHBPMP1}) folgt daher $a_0(t) +\langle a(t),\xi\rangle \geq 0$ f"ur alle $\|\xi\| \leq \delta$.
Dies zeigt $a_0(t) >0$ und wir k"onnen ohne Einschr"ankung $a_0(t)=1$ annehmen.
Wiederum (\ref{BeweisHBPMP1}) liefert damit
\begin{equation} \label{BeweisHBPMP2}
\langle a(t),x-x_*(t)\rangle \geq  \mathscr{H}^{\mathcal{S}}\big(t,x,p(t)\big) - \mathscr{H}^{\mathcal{S}}\big(t,x_*(t),p(t)\big)
  \quad \mbox{ f"ur alle } x \in V_\gamma(t).
\end{equation}
Es sei nun $t \in [t_0,t_1]$ so gew"ahlt,
dass die Maximumbedingung (\ref{PMPeinfach6}) zu diesem Zeitpunkt erf"ullt ist.
Dann folgt aus (\ref{BeweisHBPMP2}), dass
\begin{eqnarray*}
       \langle a(t),x-x_*(t)\rangle
&\geq& \mathscr{H}^{\mathcal{S}}\big(t,x,p(t)\big) - \mathscr{H}^{\mathcal{S}}\big(t,x_*(t),p(t)\big) \\
&=&    \sup_{u \in U} H^{\mathcal{S}}\big(t,x,u,p(t)\big) - \mathscr{H}^{\mathcal{S}}\big(t,x_*(t),p(t)\big) \\
&\geq& \big\langle p(t), \varphi\big(t,x,u_*(t)\big) \big\rangle - f\big(t,x,u_*(t)\big) \\
&    & -\big[ \big\langle p(t), \varphi\big(t,x_*(t),u_*(t)\big) \big\rangle - f\big(t,x_*(t),u_*(t)\big)\big]
\end{eqnarray*}
f"ur alle $x \in V_\gamma(t)$ gilt.
Wir setzen
\begin{eqnarray*}
\Phi(x) &=& \big\langle p(t), \varphi\big(t,x,u_*(t)\big)-\varphi\big(t,x_*(t),u_*(t)\big) \big\rangle \\
        & & - \big[f\big(t,x,u_*(t)\big)-f\big(t,x_*(t),u_*(t)\big)\big] - \langle a(t),x-x_*(t)\rangle.
\end{eqnarray*}
Die Funktion $\Phi(x)$ ist stetig differenzierbar auf $V_\gamma(t)$.
Ferner gelten $\Phi(x) \leq 0$ f"ur alle $x \in V_\gamma(t)$ und $\Phi(x_*(t))=0$.
Damit nimmt die Funktion $\Phi$ in dem inneren Punkt $x_*(t)$ der Menge $V_\gamma(t)$ ihr globales Maximum an.
Also gilt $0=\Phi'(x_*(t))$, d.\,h.
\begin{equation} \label{BeweisHBPMP3}
-a(t)= -\varphi_x^T\big(t,x_*(t),u_*(t)\big) p(t) + f_x\big(t,x_*(t),u_*(t)\big).
\end{equation}
Die Gleichung (\ref{BeweisHBPMP3}) wurde unter der Annahme erzielt,
dass die Maximumbedingung (\ref{PMPeinfach6}) in dem Zeitpunkt $t \in [t_0,t_1]$ erf"ullt ist.
Da (\ref{PMPeinfach6}) f"ur fast alle $t \in [t_0,t_1]$ gilt,
stimmt $-a(t)$ mit der Ableitung $\dot{p}(t)$ "uberein.
Also gilt auf $V_\gamma(t)$ die Ungleichung
\begin{equation} \label{BeweisHBPMP4}
\langle \dot{p}(t),x-x_*(t)\rangle \leq \mathscr{H}^{\mathcal{S}}\big(t,x_*(t),p(t)\big)- \mathscr{H}^{\mathcal{S}}\big(t,x,p(t)\big)
\end{equation}
f"ur fast alle $t \in [t_0,t_1]$.
Die Abbildung $S:\R^n \to \R$ ist nach Voraussetzung (c) von Theorem \ref{SatzHBPMP} über $V_\gamma(t_1)$.
Deswegen gilt für alle $x \in V_\gamma(t_1)$ die Ungleichung
\begin{equation} \label{BeweisHBPMP5}
S(x)-S\big(x_*(t_1)\big) \geq \big\langle S'\big(x_*(t_1)\big), x-x_*(t_1) \big\rangle.
\end{equation}
Es sei $\big(x(\cdot),u(\cdot)\big) \in \mathscr{A}^{\mathcal{S}}_{\rm adm}$ mit $\|x(\cdot)-x_*(\cdot)\|_\infty < \gamma$.
Dann erhalten wir
\begin{eqnarray*}
    \lefteqn{J\big(x(\cdot),u(\cdot)\big)-J\big(x_*(\cdot),u_*(\cdot)\big)} \\
&=& \int_{t_0}^{t_1} \big[f\big(t,x(t),u(t)\big)-f\big(t,x_*(t),u_*(t)\big)\big] \, dt +S\big(x(t_1)\big)-S\big(x_*(t_1)\big) \\
&\geq& \int_{t_0}^{t_1} \big[\mathscr{H}^{\mathcal{S}}\big(t,x_*(t),p(t)\big)-\mathscr{H}^{\mathcal{S}}\big(t,x(t),p(t)\big)\big] \, dt 
    + \int_{t_0}^{t_1} \langle p(t), \dot{x}(t)-\dot{x}_*(t) \rangle dt \\
& & + S\big(x(t_1)\big)-S\big(x_*(t_1)\big) \\
&\geq& \int_{t_0}^{t_1} \big[\langle \dot{p}(t),x(t)-x_*(t)\rangle + \langle p(t), \dot{x}(t)-\dot{x}_*(t) \rangle \big] \, dt
    + \big\langle S'\big(x_*(t_1)\big), x(t_1)-x_*(t_1) \big\rangle \\
&=& \langle p(t_1)+S'\big(x_*(t_1)\big),x(t_1)-x_*(t_1)\rangle-\langle p(t_0),x(t_0)-x_*(t_0)\rangle.
\end{eqnarray*}
Wegen (\ref{PMPeinfach5}) ist $p(t_1)=-S'\big(x_*(t_1)\big)$ erfüllt.
Ferner betrachten wir Aufgaben mit festem Anfangswert,
d.\,h. $x(t_0)=x_*(t_0)=x_0$.
Also gilt $J\big(x(\cdot),u(\cdot)\big)-J\big(x_*(\cdot),u_*(\cdot)\big) \geq 0$ für alle
$\big(x(\cdot),u(\cdot)\big) \in \mathscr{A}^{\mathcal{S}}_{\rm adm}$ mit $\|x(\cdot)-x_*(\cdot)\|_\infty < \gamma$. \hfill $\blacksquare$

\begin{beispiel}
{\rm Im Beispiel \ref{BeispielLinInv} des linearen Investitionsmodells lautet die Pontrjagin-Funktion
$$H^{\mathcal{S}}(t,K,u,p) = puK+(1-u)K= \big((p-1)u+1\big)K.$$
Ferner ergab sich nach (\ref{Beispiel1einfach7}) die adjungierte Funktion
$$p(t) = \left\{ \begin{array}{ll}
               e^{T-(1+t)}, & t \in [0,T-1), \\
               T-t,         & t \in [T-1,T]. \end{array} \right.$$
Dies führt auf die Hamilton-Funktion
\begin{eqnarray*}
\mathscr{H}^{\mathcal{S}}\big(t,K,p(t)\big)
&=& \sup_{u \in [0,1]} \Big(\big(p(t)-1\big)u+1\Big)K \\
&=& \left\{\begin{array}{rl} e^{T-(1+t)}K, & t \in [0,T-1), \\ K,& t \in [T-1,T], \end{array}\right.
\end{eqnarray*}
die für alle $t \in [0,T]$ konkav in der Variablen $K$ ist.
Damit liefert der Steuerungsprozess $\big(K_*(\cdot),u_*(\cdot)\big)$,
der in (\ref{Beispiel1einfach8}) und (\ref{Beispiel1einfach9}) angegeben ist,
ein starkes lokales Maximum. \hfill $\square$}
\end{beispiel}

\begin{beispiel}
{\rm In dem Beispiel \ref{Beispiel2Sektoren} eines Zwei-Sektoren-Modells besitzt die Pontrjagin-Funktion die Gestalt
$$H(t,x,y,u,p,q) = uxp+(1-u)xq+y= \big(u(p-q)+1\big)x+y.$$
Ferner lieferte die Analyse der Aufgabe die Adjungierten
$$p(t)= \left\{\begin{array}{ll} 2e^{T-(t+2)}, & t \in [0,T-2), \\[1mm] \frac{1}{2}(T-t)^2, & t \in [T-2,T], \end{array}\right.
  \qquad\qquad q(t)= T-t.$$
Daraus resultiert
$$p(t) \left\{ \begin{array}{ll} >q(t) & \mbox{für } t \in [0,T-2), \\ \leq q(t) & \mbox{für } t \in [T-2,T].
  \end{array}\right.$$
Damit erhalten wir die Hamilton-Funktion
\begin{eqnarray*}
    \mathscr{H}^{\mathcal{S}}\big(t,x,y,p(t),q(t)\big)
&=& \sup_{u \in [0,1]} \Big\{ \Big(u\big(p(t)-q(t)\big)+1\Big)x+y\Big\} \\
&=& \left\{\begin{array}{ll} \big(p(t)-q(t)+1\big)x+y, & t \in [0,T-2), \\ x+y, & t \in [T-2,T]. \end{array}\right.
\end{eqnarray*}
F"ur jedes $t \in [0,T]$ ist $\mathscr{H}^{\mathcal{S}}$ linear in $(x,y)$, also eine konkave Funktion in $(x,y)$.
Damit ist der im Beispiel  \ref{Beispiel2Sektoren} ermittelte Kandidat $\big(x_*(\cdot),y_*(\cdot),u_*(\cdot)\big)$,
der alle Bedingungen des Maximumprinzips erf"ullt,
ein starkes lokales Maximum. \hfill $\square$}
\end{beispiel}
       \newpage
       \lhead[\thepage \hspace*{1mm} Instandhaltungsmanagement]{ }
       \subsection{Instandhaltungsmanagement}
Vorbeugende Instandhaltung wirkt sich positiv auf Lebensdauer und Funktionalität von Produktionsanlagen aus,
woraus sich au"serdem Einsparungen beim Ersatz von Anlagen ergeben.
Die Aufgabe besteht hier in der Ermittlung eines optimierten Zusammenspiels von Instandhaltungsintensität im Laufe der
Produktionsdauer und dem Verkaufswert bzw. dem Schrotterlös der Anlage am Ende der Produktionsdauer. \\
Das nachstehende Modell wurde durch die Betrachtungen in Feichtinger \cite{Feichtinger} angeregt. \\[2mm]
Es sei $x(t)$ der Zustand einer Maschine,
die im Produktionseinsatz den Erlös $E \cdot x(t)$ pro Zeiteinheit einbringt.
Der Einsatz der Maschine führt zum Verschleiß mit der Rate $\delta \cdot x(t)$,
welchem durch instandhaltende Maßnahmen entgegengewirkt werden kann.
Die Kosten für die Instandhaltung wird durch $u(t)$ angegeben und die Effizienz der Maßnahmen durch $C\big(u(t)\big)$.
Zusammenfassend bezeichnen:
\begin{center}\begin{tabular}{rcl}
$x(t)$ &--& den Zustand der Maschine zur Zeit $t$, \\
$E \cdot x(t)$ &--& den Erlös der Maschine pro Zeiteinheit aus dem Betrieb der Maschine, \\
$\delta \cdot x(t)$ &--& den Verschleiß der Maschine bei Produktion, \\
$W \cdot x(T)$ &--& den Wiederverkaufswert am Ende der Planungsperiode, \\
$u(t)$ &--& die Instandshaltungskosten, \\
$C(u)$ &--& die Effizienz der Instandhaltung zur Rate $u$.
\end{tabular}\end{center}
Der zu erwartende Nettoerlös aus dem Betrieb der Maschine, den anfallenden Instandhaltungskosten und dem Barwert des Verkaufes
berechnet sich gemäß
$$\int_0^T e^{-\varrho t}\big[E x(t)-u(t)\big] \, dt + e^{-\varrho T} \cdot W x(T).$$
Die instandhaltenden Maßnahmen $C(u)$ seien mit wachsendem Kostenaufwand $u$ weniger effizient.
Daher besitze die Funktion $C$ die Eigenschaften: 
$$C(0)=0, \quad C(\infty)=\infty, \quad C'(u) >0, \quad C'(0)=\infty, \quad C'(\infty)=0, \quad C''(u)<0, \quad u \geq 0.$$
Diese Eigenschaften spiegelt zum Beispiel die Funktion $C(u)=u^\sigma$ mit $\sigma \in (0,1)$ wider.
Über dem Planungszeitraum werden außerdem die instandhaltenden Maßnahmen durch den betrieblichen Verschleiß und durch wiederkehrende Reparaturen weniger effektiv.
Wir machen für diese Beobachtung den Ansatz $e^{-\alpha t} \cdot C(u)$. \\[2mm]  
Zusammenfassend ergibt sich für unser Instandhaltungsmodell die Aufgabe
\begin{equation} \label{Instand} 
\left.\begin{array}{l}
\hspace*{-5mm}
\displaystyle J\big(x(\cdot),u(\cdot)\big) = \int_0^T e^{-\varrho t}\big[E x(t)-u(t)\big] \, dt + e^{-\varrho T} \cdot W x(T) \to \sup,  \\[4mm]
\hspace*{-5mm}   
\displaystyle \dot{x}(t)=e^{-\alpha t} \cdot C\big(u(t)\big)- \delta x(t), \qquad x(0)=x_0>0,  \\[4mm]
\hspace*{-5mm}   
u \geq 0, \qquad \varrho,\,\delta,\, \alpha,\, E,\, W >0, \qquad\displaystyle W < \frac{E}{\varrho + \delta}.
\end{array} \right\}
\end{equation}
Wir wenden Theorem \ref{SatzPMPeinfach} an:
Die Pontrjagin-Funktion $H^{\mathcal{S}}$ der Aufgabe (\ref{Instand}) hat die Form
$$H^{\mathcal{S}}(t,x,u,p)=p (e^{-\alpha t}C(u)-\delta x) + e^{-\varrho t}(E x -u).$$
Sie ist linear in $x$, womit die Hamilton-Funktion $\mathscr{H}^{\mathcal{S}}(t,x,p)$ konkav ist.
Daher sind die Bedingungen des Maximumprinzips hinreichend für ein starkes lokales Maximum. 
Die Maximumbedingung (\ref{PMPeinfach6}) lautet
$$H^{\mathcal{S}}\big(t,x_*(t),u_*(t),p(t)\big)=\max_{u \geq 0} \Big[ p(t) \big(e^{-\alpha t}C(u)-\delta x_*(t)\big) + e^{-\varrho t}\big(E x_*(t) -u\big)\Big].$$
Aufgrund der Eigenschaften der Funktion $C(u)$ können nur positive Instandhaltungskosten $u_*(t)$ optimal sein.
Deswegen führt die Maximumbedingung nach Ausschluss der Randlösung $u=0$ auf die Gleichung
$$H^{\mathcal{S}}_u\big(t,x_*(t),u_*(t),p(t)\big)=p(t)e^{-\alpha t}C'\big(u_*(t)\big) - e^{-\varrho t} =0
  \Leftrightarrow p(t)= e^{-(\varrho - \alpha) t} \frac{1}{C'\big(u_*(t)\big)}.$$
Wir beachten den Vorzeichenwechsel beim Übergang zu einem Minimierungsproblem in der Transversalitätsbedingung.
Dann lautet die adjungierte Gleichung (\ref{PMPeinfach4})
$$\dot{p}(t)=\delta p(t)-E \cdot e^{-\varrho t}, \quad p(T)=W \cdot  e^{-\varrho T}.$$
Für die adjungierte Funktion $p(\cdot)$ erhalten wir die Abbildung
$$p(t)=W e^{-(\varrho + \delta)T} e^{\delta t} + \frac{E}{\varrho + \delta}\big[e^{-\varrho t} - e^{-(\varrho + \delta)T} e^{\delta t}\big].$$
Den Ausdruck in der letzten Klammer formen wir um und erhalten
$$e^{-\varrho t} - e^{-(\varrho + \delta)T} e^{\delta t} = e^{-\varrho t} - e^{-\varrho T} \cdot e^{-\delta (T-t)} > e^{-\varrho t} - e^{-\varrho T} >0$$
für $t \in [0,T)$.
Daher nimmt die Adjungierte $p(t)$ nur positive Werte über $[0,T]$ an und
die optimale Strategie $u_*(\cdot)$ ist durch $p(t)= e^{-(\varrho - \alpha) t} / C'\big(u_*(t)\big)$ sinnvoll festgelegt. \\[2mm]
In $t=0$ ergibt sich für die Adjungierte $p(\cdot)$ der Wert
$$p(0)= W e^{-(\varrho + \delta)T} + \frac{E}{\varrho + \delta}\big[1 - e^{-(\varrho + \delta)T}\big] \approx \frac{E}{\varrho + \delta}.$$
Im Sinne des Schattenpreises bemisst $p(\cdot)$ zu Beginn den gesamten Erlös unter Beachtung von Diskontierung und Verschleiß,
$$p(0) \approx \int_0^T E e^{-(\varrho + \delta)t} \, dt \approx \int_0^\infty E e^{-(\varrho + \delta)t} \, dt = \frac{E}{\varrho + \delta},$$
und nimmt am Ende des Planungszeitraumes den Barwert $p(T)=W e^{-\varrho T}$ des Verkaufspreises an. \hfill $\square$
       \newpage
       \lhead[\thepage \hspace*{1mm} Optimale Werbestrategien]{ }
       \subsection{Optimale Werbestrategien} \label{AbschnittWerbungDGL}
Es bezeichne $x(\cdot)$ den Bekanntheitsgrad eines Produktes, der auf $x(t) \in [0,1]$ skaliert sei.
Den Umsatz, den das Produkt zum Zeitpunkt $t$ erbringt, sei $\pi x(t)$.
Der Einsatz von Werbung wird durch die Steuerung $u(\cdot)$ widergegeben.
Es ergibt sich in Anlehnung an Feichtinger \& Hartl \cite{Feichtinger} (vgl. Abschnitt \ref{AbschnittWerbungIGL}) das Modell 
\begin{eqnarray}
&& \label{WerbungDGL1} J\big(x(\cdot),u(\cdot)\big) = \int_0^T e^{-\varrho t} \big[\pi x(t)-u(t)\big] \, dt
                      + e^{-\varrho T} S\big(x(T)\big) \to \sup, \\
&& \label{WerbungDGL2} \dot{x}(t) = -(\alpha+\delta) x(t) + \delta x_0 + \mu\big(u(t)\big), \quad x(0)=x_0, \\
&& \label{WerbungDGL3} x_0 \in [0,1],\quad u(t) \geq 0, \quad \varrho >0, \quad \alpha, \delta > 0.
\end{eqnarray}
Dabei sei $S$ stetig differenzierbar und konkav mit $S(0)=0$ und $S'(x)>0$ für $x\geq 0$. \\[2mm]
Die Auswirkungen der Werbestrategie auf den Bekanntheitsgrad wird durch (\ref{WerbungDGL2}) festgelegt.
Dabei sei $\mu$ über $[0,\infty)$ stetig, differenzierbar und streng monoton wachsend mit
$$\mu(0)=0, \quad \lim_{u \to \infty} \mu(u)=\gamma \leq \alpha, \quad \lim_{u \to 0+} \mu'(u)=\infty, \quad \lim_{u \to \infty} \mu'(u)=0.$$
Für $0< \sigma < 1$ genügen die Funktionen
$\displaystyle \mu_\sigma(u)=\gamma\frac{u^\sigma}{1+u^\sigma}$
diesen Bedingungen. \\[2mm]
Im Fall $\delta =0$ bekommt die Differentialgleichung (\ref{WerbungDGL2}) die Gestalt
$$\dot{x}(t)= \mu\big(u(t)\big)-\alpha x(t), \quad x(0)=x_0,$$
die die Form eines Ebbinghausenschen Vergessensmodell trägt.
Hermann Ebbinghausen (1850--1909) gilt als Begründer der experimentellen Gedächtnisforschung. \\[2mm]
Der Fall $u(t) \equiv 0$ überführt (\ref{WerbungDGL2}) in die Gestalt $\dot{x}(t)=-(\alpha+\delta)x(t)+ \delta x_0$
zum Anfangswert $x(0)=x_0$ und besitzt die Lösung
\begin{equation} \label{WerbungDGL4}
x(t)=\frac{\alpha}{\alpha+\delta} x_0e^{-(\alpha+\delta)t}+\frac{\delta}{\alpha+\delta} x_0.
\end{equation}
Die Popularität $x(t)$ ergibt sich anhand Gleichung (\ref{WerbungDGL4}) durch das Zusammenspiel der Parameter
$\alpha$ und $\delta$.
Dabei spiegelt der Parameter $\alpha$ einem der Natur der Sache innewohnenden Attraktivitätsverlust für das Produkt wider,
während $\delta$ eine Faszination,
d.\,h. einen eigendynamischen Erhalt der Popularität für das Produkt innerhalb der Interessengruppe, beschreibt.
Speziell ergeben sich  $x(t) \equiv x_0$ für $\alpha=0$ und es liegt kein Attraktivitätsverlust für das Produkt vor.
Im Fall $\delta =0$ stellt sich der Popularitätsverlust $x(t)=x_0e^{-\alpha t}$ zum Parameter $\alpha >0$ ein.
Ferner liegt für $\delta =\infty$ eine unbegrenzte Faszination für das Produkt vor und es ist
$x(t) \equiv x_0$ unabhängig vom Paramter $\alpha \geq 0$.
\newpage
In der Aufgabe (\ref{WerbungDGL1})--(\ref{WerbungDGL3}) ohne Wiedergewinnungswert $S$ liefert das Pontrjaginsche Maximumprinzip in Form von
Theorem \ref{SatzPMPeinfach} die adjungierte Gleichung
\begin{equation} \label{WerbungDGL5}
\dot{p}(t)=(\alpha+\delta)p(t)-\pi e^{-\varrho t}, \quad p(T)=0,
\end{equation}
welche die Lösung
\begin{equation} \label{WerbungDGL6}
p(t)=-\frac{\pi e^{-\varrho T}}{\alpha + \delta + \varrho} \cdot e^{(\alpha+\delta)(t-T)}
       +\frac{\pi e^{-\varrho t}}{\alpha + \delta + \varrho} 
\end{equation}
besitzt.
Die Maximumbedingung führt nun weiter zur der Beziehung
\begin{equation} \label{WerbungDGL7}
\max_{u \geq 0} \big[p(t)\cdot \mu(u)- e^{-\varrho t} \cdot u\big],
\end{equation}
aus der die Gültigkeit der Gleichung $p(t)\cdot \mu'\big(u_*(t)\big)-e^{-\varrho t} =0$ bzw.
\begin{equation} \label{WerbungDGL8}
\mu'\big(u_*(t)\big) = \frac{1}{p(t) \cdot e^{\varrho t}} =\frac{\alpha + \delta + \varrho}{\pi ( 1-e^{(\alpha + \delta + \varrho)(t-T)})}=\psi(t)
\end{equation}
resultiert.
Die Funktion $\psi(t)$ ist über $[0,T)$ positiv und streng monton wachsend.
Im Grenzwert $t \to T^-$ ergibt sich $\psi(t) \to \infty$.
Aus den Eigenschaften von $\mu$ folgt,
dass $u_*(\cdot)$ über $[0,T]$ stetig, streng monoton fallend mit $u_*(t) \in [0,1)$ und $u_*(T)=0$ ist. \\[2mm]
In dem Modell (\ref{WerbungDGL1})--(\ref{WerbungDGL3}) ist es nun vernünftig anzunehmen,
dass zum Abschluss der Werbekampagne für das Produkt ein weiterer Absatz mit Gesamtumsatz $S\big(x(T)\big)$ in
Abhängigkeit vom Bekanntheitsgrad $x(T)$ erwartet werden darf.
Das Pontrjaginsche Maximumprinzip für die Aufgabe führt wieder
auf die adjungierte Gleichung (\ref{WerbungDGL5}), jedoch zur Transversalitätsbedingung $p(T)=e^{-\varrho T} S'\big(x_*(T)\big)$:
\begin{equation} \label{WerbungDGL9}
\dot{p}(t)=(\alpha+\delta)p(t)-\pi e^{-\varrho t}, \quad p(T)=e^{-\varrho T} S'\big(x_*(T)\big).
\end{equation}
Da eine Maximierungsaufgabe vorliegt, tritt in der Transversalitätsbedingung ein Vorzeichenwechsel ein. 
Es ergibt sich daraus die Adjungierte
\begin{equation} \label{WerbungDGL10}
p(t)=\bigg(e^{-\varrho T}S'\big(x_*(T)\big)-\frac{\pi e^{-\varrho T}}{\alpha + \delta + \varrho}\bigg)
       \cdot e^{(\alpha+\delta)(t-T)} +\frac{\pi e^{-\varrho t}}{\alpha + \delta + \varrho}.
\end{equation}
Die Maximumbedingung führt erneut auf die Beziehung $\mu'\big(u_*(t)\big) = 1/ [p(t) \cdot e^{\varrho t}]$, d.\,h.
$$\mu'\big(u_*(t)\big)
  =\frac{1}{ \displaystyle \frac{\pi}{\alpha + \delta + \varrho}( 1-e^{(\alpha + \delta + \varrho)(t-T)})
             + S'\big(x_*(T)\big)\cdot e^{(\alpha + \delta + \varrho)(t-T)}}.$$
Die Werbestrategie $u_*(\cdot)$ ist erneut stetig, aber am Ende der Planungsperiode gilt
$$\mu'\big(u_*(T)\big) =\frac{1}{S'\big(x_*(T)\big)}\in (0,\infty), \quad\mbox{ d.\,h. } u_*(T)>0.$$
Da die Variable $x$ linear in den Integranden und in die Dynamik einfließt,
und ferner die Abbildung $S$ konkav in $x$ ist,
sind die Optimalitätsbedingungen des Maximumprinzips nach Theorem \ref{SatzHBPMP} hinreichend für ein starkes lokales Maximum. \hfill $\square$ 
       \newpage
       \lhead[\thepage \hspace*{1mm} Kapitalismusspiel]{ }
       \subsection{Kapitalismusspiel} \label{AbschnittKapitalismusspiel}
In vielen Problemstellungen sind verschiedene Entscheidungsträger involviert,
deren Interessen nicht im Einklang stehen müssen.
Eine derartige Situation bezeichnet man als mathematisches Spiel. 
Unterliegen die Zielkriterien der Konfliktgruppen dynamischen Nebenbedingungen,
so spricht man von einem Differentialspiel. \\[2mm]
Es liege ein dynamisches Spiel vor,
in dem die Steuerungen $u_1(\cdot)$ und $u_2(\cdot)$ die Einflussnahmen durch den ersten bzw. zweiten Spieler widergeben.
Durch die Einführung verschiedener Zielfunktionale in der Standardaufgabe,
welche die Zielstellungen der einzelnen Spieler darstellen,
gelangen zur dynamischen Spielsituation:
\begin{eqnarray*}
&& J_1\big(x(\cdot),u_1(\cdot),u_2(\cdot)\big) = \int_{t_0}^{t_1} f_1\big(t,x(t),u_1(t),u_2(t)\big) \, dt \to \inf, \\
&& J_2\big(x(\cdot),u_1(\cdot),u_2(\cdot)\big) = \int_{t_0}^{t_1} f_2\big(t,x(t),u_1(t),u_2(t)\big) \, dt \to \inf, \\
&& \dot{x}(t) = \varphi\big(t,x(t),u_1(t),u_2(t)\big), \quad x(t_0)=x_0, \quad u_1(t) \in U_1, \; u_2(t) \in U_2.
\end{eqnarray*}
Ein L"osungskonzept zur Behandlung eines solchen Spiels ist das
Nash-Gleichgewicht\index{Nash-Gleichgewicht} f"ur nichtkooperative Spiele.
In unserem Differentialspiel ist ein Nash-Gleichgewicht ein Tripel
$\big(x^*(\cdot),u_1^*(\cdot),u_2^*(\cdot)\big)$ mit
\begin{eqnarray*}
    J_1\big(x^*(\cdot),u_1^*(\cdot),u_2^*(\cdot)\big)
&=& \min_{(x(\cdot),u_1(\cdot))} J_1\big(x(\cdot),u_1(\cdot),u_2^*(\cdot)\big), \\
    J_2\big(x^*(\cdot),u_1^*(\cdot),u_2^*(\cdot)\big)
&=& \min_{(x(\cdot),u_2(\cdot))} J_2\big(x(\cdot),u_1^*(t),u_2(\cdot)\big).
\end{eqnarray*}
Befindet sich das Spiel also im Nash-Gleichgewicht $\big(x^*(\cdot),u_1^*(\cdot),u_2^*(\cdot)\big)$,
so hat kein Spieler einen Anreiz seine Strategie zu ändern,
denn das bestmögliche Ergebnis zur Strategie des Konkurrenten wurde gefunden. \\
Ausf"uhrlich geschrieben ist das Nash-Gleichgewicht ein System von zwei gekoppelten Steuerungsproblemen.
Für den ersten Spieler ergibt sich das Steuerungsproblem
\begin{eqnarray*}
&& J_1\big(x(\cdot),u_1(\cdot),u^*_2(\cdot)\big) = \int_{t_0}^{t_1} f_1\big(t,x(t),u_1(t),u^*_2(t)\big) \, dt \to \inf, \\
&& \dot{x}(t) = \varphi\big(t,x(t),u_1(t),u^*_2(t)\big), \quad x(t_0)=x_0, \quad u_1(t) \in U_1
\end{eqnarray*}
bei gegebener Strategie $u_2^*(\cdot)$ des zweiten Spielers; entsprechend für den zweiten Spieler:
\begin{eqnarray*}
&& J_2\big(x(\cdot),u_1^*(\cdot),u_2(\cdot)\big) = \int_{t_0}^{t_1} f_2\big(t,x(t),u^*_1(t),u_2(t)\big) \, dt \to \inf, \\
&& \dot{x}(t) = \varphi\big(t,x(t),u^*_1(t),u_2(t)\big), \quad x(t_0)=x_0, \quad u_2(t) \in U_2.
\end{eqnarray*}
\newpage
Formal sind die Abbildungen $f_1\big(t,x,u_1,u^*_2(t)\big)$, $f_2\big(t,x,u^*_1(t),u_2\big)$ und
$\varphi\big(t,x,u_1,u^*_2(t)\big)$, $\varphi\big(t,x,u^*_1(t),u_2\big)$ nicht stetig bezüglich der Variablen $t$
und verletzen die Anforderungen an die Standardaufgabe (\ref{PMPeinfach1})--(\ref{PMPeinfach3}).
Sind die Steuerungen $u_1^*(\cdot)$ und $u_2^*(\cdot)$ zumindest stückweise stetig,
so besitzen diese Abbildungen für die gegebene Strategie des Gegenspielers höchstens endlich viele Unstetigkeiten in der Variablen $t$.
Bei der Setzung der einfachen Nadelvariation im Beweis des Maximumprinzips sind für die Wahl von $\tau \in (t_0,t_1)$ diese Stellen auszuschließen.
Auf diese Weise bleiben das Beweisschema und die Optimalitätsbedingungen in Form des Pontrjaginschen Maximumprinzips gültig. \\[2mm]
Wir betrachten im Folgenden ein Differentialspiel zwischen Arbeitern und Unternehmern.
Beide Parteien sind bestrebt ihr Konsumverlangen zu befriedigen.
Zudem zeichnen sich die Unternehmer für Investitionen verantwortlich.
F"ur die Arbeiter stellt sich dabei das Problem,
inwieweit sie den produzierten Ertrag konsumieren oder den Unternehmern "uberlassen sollen,
damit eine k"unftige hohe G"uterproduktion erm"oglicht wird, von der auch die Arbeiter wieder profitieren.
Das Dilemma, das sich dem Arbeiter stellt, ist,
dass sie keine Garantien "uber ausreichende Neuinvestitionen seitens der Unternehmer haben.
Die Unternehmer stehen ihrerseits vor der Frage,
wie sie mit dem verbliebenen Teil des Ertrages, der nicht dem Arbeiter zugesprochen wird, umgehen:
Sollen sie diesen investieren oder konsumieren? \\
Die Spielsituation ensteht dadurch,
dass der Nutzen f"ur den Arbeiter und den Unternehmern durch den aufgeteilten Ertrag miteinander gekoppelt ist.
D.\,h. der Gewinn eines jeden Spielers h"angt von der Entscheidung des anderen Spielers ab. \\[2mm]
In diesem Modell werden die Geld- und G"uterwerte im Kapitalstock $K(\cdot)$ zusammengefasst.
Weiter bezeichne $u(t)$ denjenigen relativen Anteil an der Produktion, der dem Arbeiter zum Zeitpunkt $t$ zugesprochen wird.
Von dem verbliebenen Teil kann der Unternehmer mit der Rate $v(t)$ zur Zeit $t$ investieren.
Mit den Zielfunktionalen $W$ bzw. $C$ f"ur den Arbeiter bzw. Unternehmer entsteht damit die folgende Aufgabe
\begin{eqnarray}
&&\label{KapSpiel1} J_W\big(K(\cdot),u(\cdot),v(\cdot)\big) = \int_0^T u(t)K(t) \, dt \to \sup,\\
&&\label{KapSpiel2} J_C\big(K(\cdot),u(\cdot),v(\cdot)\big) = \int_0^T \big(1-v(t)\big)\big(1-u(t)\big)K(t) \, dt \to \sup, \\
&&\label{KapSpiel3} \dot{K}(t) = v(t)\big(1-u(t)\big)K(t), \qquad K(0)=K_0>0, \quad K(T) \mbox{ frei}, \\
&&\label{KapSpiel4} \big(u(t),v(t)\big) \in [a,b] \times [0,1], \quad 0<a<b<1, \quad b > \frac{1}{2},\quad T>\frac{1}{1-b}.\hspace*{5mm}
\end{eqnarray}
Das Modell geht auf Lancaster \cite{Lancaster} zur"uck.
Wir folgen der Darstellung und untersuchen das Nash-Gleichgewicht und die Kollusionsl"osung,
die wir abschlie"send miteinander vergleichen.
Wir verweisen au"serdem auf die umfassenden Ausf"uhrungen in
Feichtinger \& Hartl \cite{Feichtinger} und Seierstad \& Syds\ae ter \cite{Seierstad}.

\newpage
Wir platzieren die zwei Spieler Arbeiter und Unternehmer in einer Spielsituation und benutzen als L"osungskonzept das
Nash-Gleichgewicht\index{Nash-Gleichgewicht} f"ur nichtkooperative Spiele.
Im vorliegenden Differentialspiel ist ein Nash-Gleichgewicht ein Tripel
$$\big(K_*(\cdot),u_*(\cdot),v_*(\cdot)\big) \in PC_1([0,T],\R) \times PC([0,T],[a,b]) \times PC([0,T],[0,1]),$$
f"ur das die folgenden Gleichungen erf"ullt sind:
\begin{eqnarray*}
J_W\big(K_*(t),u_*(t),v_*(t)\big) &=& \max_{(K(\cdot),u(\cdot))} J_W\big(K(\cdot),u(\cdot),v_*(t)\big), \\
J_C\big(K_*(t),u_*(t),v_*(t)\big) &=& \max_{(K(\cdot),v(\cdot))} J_C\big(K(\cdot),u_*(t),v(\cdot)\big).
\end{eqnarray*}
Ausf"uhrlich geschrieben ist das Nash-Gleichgewicht das folgende System von zwei gekoppelten Steuerungsproblemen
f"ur den Arbeiter,
\begin{equation} \label{KapSpiel5} \left. \begin{array}{l}
J_W\big(K(\cdot),u(\cdot),v_*(\cdot)\big) = \displaystyle\int_0^T u(t)K(t) \, dt  \to \sup, \\[1mm]
\dot{K} =   v_*(t)\big(1-u(t)\big)K(t), \quad K(0)=K_0>0, \quad u \in [a,b], \\[1mm]
\end{array} \right\}
\end{equation}
und für den Unternehmer,
\begin{equation} \label{KapSpiel6} \left. \begin{array}{l}
J_C\big(K(\cdot),u_*(\cdot),v(\cdot)\big)
 = \displaystyle\int_0^T \big(1-v(t)\big)\big(1-u_*(t)\big)K(t) \, dt \to \sup, \\[1mm]
\dot{K} = v(t)\big(1-u_*(t)\big)K(t), \quad K(0)=K_0>0, \quad v \in [0,1].
\end{array} \right\}
\end{equation} 
Auf beide Steuerungsprobleme wenn wir die Optimalitätsbedingungen im Pontrjaginsche Maximumprinzip
unter Beachtung der Strategie des Gegenspielers an: \\
In der Aufgabe der Arbeiter (\ref{KapSpiel5}) erhalten wir die Maximumbedingung
$$H_W^{\mathcal{S}}\big(t,K_*(t),u_*(t),p(t)\big) = \max_{u \in [a,b]} \Big\{ \big(1-p(t)v_*(t)\big)uK_*(t) \Big\} + p(t)v_*(t)K_*(t).$$
Da jede zul"assige Trajektorie $K(\cdot)$ stets positiv ist, gilt
\begin{equation} \label{KapSpiel7} \left. 
   \begin{array}{ll} 
      u_*(t)=a,         & \mbox{ wenn } p(t)v_*(t)>1, \\
      u_*(t)=b,         & \mbox{ wenn } p(t)v_*(t)<1, \\
      u_*(t) \in [a,b], & \mbox{ wenn } p(t)v_*(t)=1.    
    \end{array} \right\}
\end{equation}
Weiter gen"ugt $p(\cdot)$ der adjungierten Gleichung (\ref{PMPeinfach4}) und Transversalitätsbedingung (\ref{PMPeinfach5}):
\begin{equation} \label{KapSpiel8}
\dot{p}(t) = -v_*(t)\big(1-u_*(t)\big)p(t) - u_*(t), \qquad p(T)=0.   
\end{equation}
In Gleichung (\ref{KapSpiel8}) erkennt man unmittelbar $\dot{p}(t)\leq  - u_*(t)\leq -a$ f"ur $t \in (0,T)$.
Daher gilt $p(t)>0$ f"ur $t \in [0,T)$. \\[2mm]
Entsprechend ergibt sich in der Aufgabe des Unternehmers (\ref{KapSpiel6}) die Maximumbedingung
$$H_C^{\mathcal{S}}\big(t,K_*(t),v_*(t),q(t)\big)=
  \max_{v \in [0,1]} \Big\{ \big(q(t)-1\big)v\big(1-u_*(t)\big)K_*(t) \Big\} + \big(1-u_*(t)\big)K_*(t)$$
und es gilt
\begin{equation} \label{KapSpiel9} \left. 
   \begin{array}{ll} 
      v_*(t)=0         & \mbox{ wenn } q(t)<1, \\
      v_*(t)=1         & \mbox{ wenn } q(t)>1, \\
      v_*(t) \in [0,1] & \mbox{ wenn } q(t)=1.    
    \end{array} \right\}
\end{equation}
Ferner gen"ugt $q(\cdot)$ der adjungierten Gleichung zur Transversalit"atsbedingung:
\begin{equation} \label{KapSpiel10}
\dot{q}(t) = -v_*(t)\big(1-u_*(t)\big)q(t) - \big(1-v_*(t)\big)\big(1-u_*(t)\big),\qquad q(T)=0.
\end{equation}
Zusammen erhalten wir aus (\ref{KapSpiel9}) und (\ref{KapSpiel10}) 
$$\dot{q}(t) = -\big(1-u_*(t)\big) \cdot \max\{1,\,q(t)\}<0.$$
Aufgrund der strengen Monotonie existiert eine eindeutig bestimmte L"osung von
$$\tau= \min \{ t \in [0,T] \,|\, q(t) \leq 1 \}.$$
Wegen (\ref{KapSpiel9}) gilt $v_*(t) = 0$ f"ur alle $t \in [\tau,T]$ und aus (\ref{KapSpiel7}) folgt nun $u_*(t) = b$ auf $[\tau,T]$.
Auf $[\tau,T]$ erhalten wir daraus
$$q(t)=(1-b)(T-t), \quad q(\tau)=1, \quad p(t)=b(T-t),\quad p(\tau)=\frac{b}{1-b}.$$
Weiterhin folgt nach (\ref{KapSpiel4}) und aus der Bedingung $q(\tau)=1$
\begin{equation} \label{KapSpiel11}
\tau=T-\frac{1}{1-b}>0.
\end{equation}
Aufgrund $b \geq 1/2$ in (\ref{KapSpiel4}) ist $p(\tau)\geq 1$.
Da $p(\cdot)$ streng monoton fallend ist,
erhalten wir $u_*(t)=a$ und $v_*(t)=1$ f"ur alle $t \in [0,\tau)$. \\[2mm]
Zusammenfassend lauten im Nash-Gleichgewicht die jeweiligen Strategien
$$u_*(t)=\left\{\begin{array}{l} a,\; t \in [0,\tau), \\ b,\; t \in [\tau,T], \end{array}\right. \quad
  v_*(t)=\left\{\begin{array}{l} 1,\; t \in [0,\tau), \\ 0,\; t \in [\tau,T], \end{array}\right.$$
zum Kapitalstock 
$$K_*(t)=\left\{\begin{array}{l} K_0e^{(1-a)t},\; t \in [0,\tau), \\ K_0e^{(1-a)\tau},\; t \in [\tau,T]. \end{array}\right.$$
F"ur die Zielfunktionale ergeben sich
\begin{equation} \label{KapSpiel12} \left. 
\begin{array}{l}
J_W\big(K_*(t),u_*(t),v_*(t)\big)
   = \displaystyle\left(\frac{a}{1-a}+\frac{b}{1-b}\right)K_0 e^{(1-a)\tau} - \frac{a}{1-a}K_0, \\[5mm]
J_C\big(K_*(t),u_*(t),v_*(t)\big) = \displaystyle K_0 e^{(1-a)\tau}. 
\end{array} \right\}
\end{equation}
Die Untersuchung des Nash-Gleichgewichtes ist damit abgeschlossen.

\newpage
Im Gegensatz zur Spielsituation besprechen sich die Arbeiter und Unternehmer,
wie sie gemeinsam am besten agieren k"onnen.
F"ur das gemeinsame Zielfunktional ergibt sich
$$\int_0^T u(t)K(t) \, dt + \int_0^T \big(1-v(t)\big)\big(1-u(t)\big)K(t) \, dt =
     \int_0^T \big[1- v(t)\big(1-u(t)\big)\big]K(t) \, dt.$$
F"uhren wir die neue Steuerungsvariable $w(\cdot)$ durch
$$w(t)=v(t)\big(1-u(t)\big), \qquad w(t) \in [0,1-a],$$
ein, so erhalten wir die Aufgabe 
\begin{equation}
\left.\begin{array}{l}
\hspace*{-5mm}
\displaystyle J\big(K(\cdot),w(\cdot)\big) = \int_0^T \big(1-w(t)\big)K(t) \, dt \to \sup,  \\[3mm]
\hspace*{-5mm}   
\displaystyle \dot{K}(t) = w(t)K(t), \qquad K(0)=K_0>0,  \\[2mm]
\hspace*{-5mm}   
w(t) \in [0,1-a], \quad 0<a<b<1, \quad b > \frac{1}{2},\quad T>\frac{1}{1-b}.
\end{array} \right\}
\end{equation}
Die kooperative Lösung ist (sie ergibt sich analog zu Beispiel \ref{BeispielLinInv})
$$w_*(t)=\left\{\begin{array}{ll} 1-a,& t \in [0,T-1), \\ b,& 0 \in [T-1,T], \end{array}\right. \qquad
  K_*(t)=\left\{\begin{array}{ll} K_0e^{(1-a)t},& t \in [0,T-1), \\ K_0e^{(1-a)(T-1)},& t \in [T-1,T], \end{array}\right.$$
mit folgendem Optimalwert f"ur das Zielfunktional:
\begin{equation} \label{KapSpiel13}
J\big(K_*(t),w_*(t)\big) = \frac{1}{1-a}K_0 e^{(1-a)(T-1)} - \frac{a}{1-a}K_0.
\end{equation}
Der Gesamtnutzen im Nash-Gleichgewicht nach (\ref{KapSpiel12}) ist gleich
\begin{eqnarray}
&& \hspace*{-2.5cm} J_W\big(K_*(t),u_*(t),v_*(t)\big) + J_C\big(K_*(t),u_*(t),v_*(t)\big) \nonumber \\
&&\label{KapSpiel14} = \left(\frac{1}{1-a}+\frac{b}{1-b}\right)K_0 e^{(1-a)\tau} - \frac{a}{1-a}K_0.
\end{eqnarray}
Zum Vergleich der Werte (\ref{KapSpiel13}) und (\ref{KapSpiel14}) verwenden wir die Beziehung
$$T-1=\tau+\frac{b}{1-b}.$$
Damit erhalten wir
\begin{eqnarray*}
J\big(K_*(t),w_*(t)\big) &>& J_W\big(K_*(t),u_*(t),v_*(t)\big) + J_C\big(K_*(t),u_*(t),v_*(t)\big) \\
\Leftrightarrow \qquad \frac{1}{1-a}e^{(1-a)\frac{b}{1-b}} &>& \frac{1}{1-a}+\frac{b}{1-b}.
\end{eqnarray*}
Benutzen wir die elementare Ungleichung $e^x>1+x$ f"ur $x>0$, so folgt abschlie"send
$$\frac{1}{1-a}e^{(1-a)\frac{b}{1-b}}>\frac{1}{1-a}\left(1+(1-a)\frac{b}{1-b}\right)= \frac{1}{1-a}+\frac{b}{1-b}.$$
D.\,h. der Gesamtnutzen ist im sozialen Optimum gr"o"ser als im Wettbewerb.
An dieser Stelle bleibt aber die spannende Frage nach einer ``gerechten'' Aufteilung des Nutzens auf Unternehmer und Arbeiter offen.  \hfill $\square$
\cleardoublepage

\lhead[ ]{}
\rhead[]{Multiprozesse \hspace*{1mm} \thepage}
\section{Optimale Multiprozesse}    
       Der Begriff des Multiprozesses\index{Multiprozess} l"asst sich am Beispiel des Autofahrens illustrieren:
Neben den kontinuierlichen Steuerungen ``beschleunigen'' und ``bremsen''
bildet die Wahl des konkreten Ganges einen wesentlichen Beitrag zur Minimierung des Treibstoffverbrauchs oder
zum Erreichen des Zielortes in k"urzester Zeit.
Denn die Auswahl des jeweiligen Ganges beeinflusst ma"sgeblich das dynamische Beschleunigungsverhalten, den momentanen Benzinverbrauch
und die Geschwindigkeit.
Das Fahrverhalten und der Treibstoffverbrauch wird in jedem einzelnen Gang durch eine eigene Dynamik und Verbrauchsfunktion beschrieben.
Dementsprechend setzt sich dieses Optimierungsproblem aus verschiedenen,
dem jeweilig ausgew"ahlten Gang zugeordneten Steuerungsproblemen zusammen.
Die Schaltfolge zwischen den einzelnen G"angen wird durch eine Wechselstrategie\index{Wechselstrategie}
\index{Multiprozess!Wechselstrategie@-- Wechselstrategie} beschrieben.
Das Beispiel des Autofahrens verdeutlicht dabei den speziellen Charakter der Wechselstrategie:
Im Vergleich zu den Pedalen, die stufenlos gesteuert werden k"onnen,
ist die Auswahl des Ganges eine rein diskrete Gr"o"se. \\[2mm]
Ein Multiprozess besteht aus einer gewissen endlichen Anzahl von einzelnen Steuerungssystemen,
welche sich jeweils aus individuellen Dynamiken, Zielfunktionalen und Steuerungsbereichen zusammensetzen.
Neben der Suche nach der optimalen Steuerung f"ur das jeweils gew"ahlte Steuerungssystem
liegt das Hauptaugenmerk bei der Untersuchung von Multiprozessen auf der Bestimmung der optimalen Wechselstrategie.
Dabei wird eine Wechselstrategie durch die Anzahl der Wechsel zwischen den einzelnen Steuerungssystemen,
durch die konkrete Auswahl des jeweiligen Steuerungssystems und durch diejenigen Zeitpunkte,
zu denen diese Wechsel stattfinden, beschrieben.
Wir betrachten ausschlie"slich Multiprozesse mit stetigen Zust"anden.
Au"serdem d"urfen keine Wechselkosten anfallen. \\[2mm]
Im Vergleich zur Standardaufgabe liegt die besondere Charakteristik eines Multiprozesses in der Wechselstrategie.
Die wesentliche Herausforderung in der folgenden Untersuchung von Multiprozessen besteht somit darin,
die Wechselstrategien so in eine Form zu gießen,
dass die Aufgabe eines Multiprozesses die Gestalt eines Steuerungsproblems erhält.
       \lhead[\thepage \hspace*{1mm} Vorbereitungen]{ }
       \subsection{Vorbereitende Betrachtungen}
Es bezeichne $k \in \N$ die Anzahl an Steuerungssystemen,
aus denen sich der Multiprozess zusammen setzt.
Für ein $i \in \{1,...,k\}$ bezeichnen $f_i(t,x,u_i)$ den Integranden im Zielfunktional,
$\varphi_i(t,x,u_i)$ die rechte Seite der Dynamik und $U_i \subseteq \R^{m_i}$ den Steuerbereich des $i-$ten Steuerungssystems.
Dann ist es naheliegend einen Multiprozess wie folgt aufzustellen:
Wir zerlegen das Zeitintervall $[t_0,t_1]$ mit Hilfe der (gegebenen) Wechselzeitpunkte $t_0=s_0<s_1<...<s_N=t_1$
in $N$ Teilintervalle $[s_{j-1},s_j]$ mit $j =1,...,N$.
Über einem Zeitabschnitt $[s_{j-1},s_j]$ ist ein gewisses Steuerungsproblem ausgewählt, nämlich $i_j \in \{1,...,k\}$.
Auf diese Weise ergeben sich über einem Teilabschnitt die Elemente
$$\int_{s_{j-1}}^{s_j} f_{i_j}\big(t,x(t),u(t)\big) \, dt, \quad
  \dot{x}(t) = \varphi_{i_j}\big(t,x(t),u(t)\big) \mbox{ f"ur } t \in (s_{j-1},s_j),\quad u(t) \in U_{i_j} \subseteq \R^{m_{i_j}}.$$
Über dem gesamten Zeitintervall $[t_0,t_1]$ erhält der Multiprozess dadurch die Gestalt
\begin{equation} \label{MultiprozessSt}
\left.\begin{array}{l}
\hspace*{-5mm}
\displaystyle J\big(x(\cdot),u(\cdot)\big) = \sum_{j=1}^N \int_{s_{j-1}}^{s_j} f_{i_j}\big(t,x(t),u(t)\big) \, dt \to \inf,  \\[3mm]
\hspace*{-5mm}   
\displaystyle \dot{x}(t) = \varphi_{i_j}\big(t,x(t),u(t)\big) \quad \mbox{ f"ur } t \in (s_{j-1},t_j), \quad j=1,...,N,\quad x(t_0)=x_0,  \\[2mm]
\hspace*{-5mm} 
u(t) \in U_{i_j} \subseteq \R^{m_{i_j}},\quad {i_j} \in \{1,...,k\}, \quad j=1,...,N.
\end{array} \right\}
\end{equation}
Die Darstellung (\ref{MultiprozessSt}) eines Multiprozesses ist für die Behandlung als Aufgabe der Optimalen Steuerung weniger geeignet,
denn die Folge der Indizes $\{i_j\}$,
welche Anzahl der Wechsel und die Abfolge der aktiven Steuerungssysteme beschreibt,
ist bereits determiniert.
Somit kann die wesentliche Charakteristik eines Multiprozesses, nämlich dessen optimale Wechselstrategie,
nicht in die Optimierung integriert werden.
Die Beschreibung einer Wechselstrategie als ``echte'' Steuerungsvariable wird in dieser Form nicht erreicht. \\[2mm]
Unsere Zielstellung ist also eine Formulierung eines Multiprozesses zu finden,
so dass die Wechselstrategien die Form einer Steuervariablen erhalten und damit alle möglichen Abfolgen von aktiven
Steuerunssystemen mit beliebiger Anzahl an Wechseln und zu beliebigen Wechselzeitpunkten zur Konkurrenz zugelassen sind.
Dazu erarbeiten wir einen alternativen Ansatz, welcher auf den Zerlegungen des Zeitintervall $[t_0,t_1]$ beruht.

\begin{definition}
Ein stückweise zusammengesetztes Intervall ist die endliche Vereinigung von
rechtsseitig halboffenen und abgeschlossenen Intervallen $[a,b)$ bzw. $[a,b]$ mit $a<b$.
\end{definition}

\begin{definition}[$k$-fache Zerlegung] \label{DefinitionZerlegung}\index{Multiprozess!kfacheZerlegung@-- $k$-fache Zerlegung}\index{kfacheZerlegung@$k$-fache Zerlegung}
Eine $k$-fache Zerlegung des Intervalls $I \subseteq \R$ ist ein endliches System
$\mathscr{A} = \{\mathscr{A}_1,...,\mathscr{A}_k\}$ von stückweise zusammengesetzten Intervallen  mit
$$\bigcup_{1\leq s\leq k} \mathscr{A}_s = I, \qquad \mathscr{A}_{s} \cap \mathscr{A}_{s'} = \emptyset \quad\mbox{f"ur } s \not=s'.$$
Es bezeichnet $\mathscr{Z}^k(I) = \{ \mathscr{A} \} = \big\{\{ \mathscr{A}_s\}_{1\leq s\leq k}\big\}$ die Menge der
$k$-fachen Zerlegungen von $I$.
\end{definition}

In unserem nächsten Schritt verknüpfen wir die $k$-fachen Zerlegungen $\mathscr{Z}^k(I)$ mit den Steuerungssystemen des Multiprozesses.
Dazu identifizieren wir ein Element $\mathscr{A} \in \mathscr{Z}^k(I)$ mit der Vektorfunktion
der charakteristischen Funktionen der Mengen ${\mathscr A}_s$:
$$\chi_{\mathscr{A}}(t) = \big(\chi_{\mathscr{A}_1}(t),...,\chi_{\mathscr{A}_k}(t) \big), \quad t \in I.$$
Da die Mengen ${\mathscr A}_s$ stückweise zusammengesetzte Intervalle sind,
ist jede Funktion $\chi_{\mathscr{A}}(\cdot)$ stückweise stetig und geh"ort zu der Menge
$$\bigg\{ y(\cdot) \in PC(I,\R^k)\,\bigg|\, y(t)=\big(y_1(t),...,y_k(t)\big), y_s(t) \in \{0,1\},
                             \sum_{s=1}^k y_s(t)=1\bigg\}.$$
Weiterhin betrachten wir die $k$ Funktionen $h_s(t,x,u_s): I \times \R^n \times \R^{m_s} \to \R^m$ und fassen diese zur
Vektorfunktion $h=(h_1,...,h_k)$ zusammen.
Außerdem wurden darin die Variablen $u=(u_1,...,u_k) \in \R^{m_1+...+m_k}$ vereint.
Die Verkn"upfung der Funktion $h$ mit einer $k$-fachen Zerlegung setzen wir wie folgt fest:
$$\mathscr{A} \circ h(t,x,u) = \chi_{\mathscr{A}}(t) \circ h(t,x,u)
                             = \sum_{s=1}^k \chi_{\mathscr{A}_s}(t) \cdot h_s(t,x,u_s).$$
Sind die Funktionen $h_{i}$ nach $x$ differenzierbar, dann setzen wir f"ur die Ableitung:
$$\mathscr{A} \circ h_x(t,x,u) = \chi_{\mathscr{A}}(t) \circ h_x(t,x,u) = \sum_{s=1}^k \chi_{\mathscr{A}_s}(t) \cdot h_{s,x}(t,x,u_s).$$
Die Verknüpfung $\mathscr{A} \circ h$ hat zur Folge,
dass die $k$-fache Zerlegung $\mathscr{A}$ zu jedem Zeitpunkt $t \in I$ einerseits in eindeutiger Weise eine der
Funktionen $h_s(t,x,u_s)$ und au"serdem die entsprechende Komponente $u_s$ des Vektors $u=(u_1,...,u_k)$ auswählt.
Auf diese Weise erhält $\mathscr{A}$ das Wesen einer Auswahlvariable.
Die Herausforderung, einen Mulitiprozess unabhängig von einer Schaltfolge zu formulieren,
wird damit durch $k$-fache Zerlegung gemeistert.
Damit sind wir bereit die Multiprozesse in eine geeignete Form zu gießen.
       \lhead[\thepage \hspace*{1mm} Pontrjaginsches Maximumprinzip]{ }
       \subsection{Die Aufgabenstellung und das Pontrjaginsche Maximumprinzip} \label{AbschnittPMPeinfachMP}
Es sei $k \in \N$ die Anzahl der verschiedenen gegebenen Steuerungssysteme
mit ihren Integranden $f_s(t,x,u_s): \R \times \R^n \times \R^{m_s} \to \R$,
rechten Seiten $\varphi_s(t,x,u_s): \R \times \R^n \times \R^{m_s} \to \R^n$
und Steuerbereichen $U_s \subseteq \R^{m_s}$.
Wir fassen die Steuerungssysteme durch die Setzungen
$$f=(f_1,...,f_k),\qquad \varphi=(\varphi_1,...,\varphi_k), \qquad U=U_1 \times ... \times U_k$$
zusammen.
Dann besitzt mit Hilfe der $k$-fachen Zerlegungen die Aufgabe eines optimalen Multiprozesses die Form eines Steuerungsproblems:
\begin{eqnarray}
&&\label{PMPeinfach1MP} J\big(x(\cdot),u(\cdot),\mathscr{A}\big)
               = \int_{t_0}^{t_1} \chi_{\mathscr{A}}(t) \circ f\big(t,x(t),u(t)\big) dt \to \inf, \\
&&\label{PMPeinfach2MP} \dot{x}(t) = \chi_{\mathscr{A}}(t) \circ \varphi\big(t,x(t),u(t) \big), \quad x(t_0)=x_0, \\
&&\label{PMPeinfach3MP} u(t) \in U= U_1 \times ... \times U_k, \quad U_s \not= \emptyset, \quad \mathscr{A} \in \mathscr{Z}^k([t_0,t_1]).
\end{eqnarray}
Die Aufgabe (\ref{PMPeinfach1MP})--(\ref{PMPeinfach3MP}) untersuchen wir bez"uglich der Tripel
$$\big(x(\cdot),u(\cdot),\mathscr{A}\big) \in PC_1([t_0,t_1],\R^n) \times PC([t_0,t_1],U) \times \mathscr{Z}^k([t_0,t_1]).$$
Mit $\mathscr{A}^{\,\mathcal{M}}_{\rm Lip}$ bezeichnen wir die Menge aller Tripel $\big(x(\cdot),u(\cdot),\mathscr{A}\big)$,
für die es ein $\gamma>0$ derart gibt,
dass die Abbildungen $f_s(t,x,u_s)$, $\varphi_s(t,x,u_s)$ für $s=1,...,k$ auf der Menge aller Punkte $(t,x,u) \in \R \times \R^n \times \R^m$ mit 
$t \in [t_0,t_1]$, $\|x-x(t)\| < \gamma$, $u \in \R^m$
stetig in der Gesamtheit aller Variablen und stetig differenzierbar bezüglich $x$ sind. \\[2mm]
Ein Tripel $\big(x(\cdot),u(\cdot),\mathscr{A}\big) \in PC_1([t_0,t_1],\R^n) \times PC([t_0,t_1],U) \times \mathscr{Z}^k([t_0,t_1])$
ist ein zul"assiger Steuerungsprozess der Aufgabe (\ref{PMPeinfach1MP})--(\ref{PMPeinfach3MP}),
falls $\big(x(\cdot),u(\cdot),\mathscr{A}\big)$ dem System (\ref{PMPeinfach2MP}) zu $x(t_0)=x_0$ gen"ugt.
Mit $\mathscr{A}^{\mathcal{M}}_{\rm adm}$ bezeichnen wir die Menge der zul"assigen Steuerungsprozesse. \\[2mm]
Ein zul"assiger Steuerungsprozess $\big(x_*(\cdot),u_*(\cdot),\mathscr{A}_*\big)$ ist eine
starke lokale Minimalstelle\index{Minimum, starkes lokales!elementar@-- Multiprozesse}
der Aufgabe (\ref{PMPeinfach1MP})--(\ref{PMPeinfach3MP}),
falls eine Zahl $\varepsilon > 0$ derart existiert, dass die Ungleichung 
$$J\big(x(\cdot),u(\cdot),\mathscr{A}\big) \geq J\big(x_*(\cdot),u_*(\cdot),\mathscr{A}_*\big)$$
f"ur alle $\big(x(\cdot),u(\cdot),\mathscr{A}\big) \in \mathscr{A}^{\mathcal{M}}_{\rm adm}$ mit $\|x(\cdot)-x_*(\cdot)\|_\infty < \varepsilon$ gilt. \\[2mm]
Es bezeichnen $H_s$ die Pontrjagin-Funktionen der einzelnen Steuerungssysteme:
$$H_s(t,x,u_s,p)= \langle p, \varphi_s(t,x,u_s) \rangle - f(t,x,u_s), \quad s=1,...,k.$$
Die Pontrjagin-Funktionen fügen wir zu der Vektorfunktion $H=(H_1,...,H_k)$ zusammen:
$$H(t,x,u,p)=\big(H_1(t,x,u_1,p),...,H_k(t,x,u_k,p)\big).$$
Damit definieren wir die Pontrjagin-Funktion $H^{\mathcal{M}}=\mathscr{A} \circ H$ auf folgende Weise:
$$H^{\mathcal{M}}(t,x,u,p,\mathscr{A})=\chi_{\mathscr{A}}(t) \circ H(t,x,u,p) = \sum_{s=1}^k \chi_{\mathscr{A}_s}(t) \cdot H_s(t,x,u_s,p).$$

\begin{theorem}[Pontrjaginsches Maximumprinzip] \label{SatzPMPeinfachMP}
\index{Pontrjaginsches Maximumprinzip!Multi@-- Multiprozesse} 
Sei $\big(x_*(\cdot),u_*(\cdot),\mathscr{A}_*\big) \in \mathscr{A}^{\mathcal{M}}_{\rm adm} \cap \mathscr{A}^{\mathcal{M}}_{\rm Lip}$. 
Ist $\big(x_*(\cdot),u_*(\cdot),\mathscr{A}_*\big)$ ein starkes lokales Minimum der Aufgabe (\ref{PMPeinfach1MP})--(\ref{PMPeinfach3MP}),
dann existiert eine Vektorfunktion $p(\cdot) \in PC_1([t_0,t_1],\R^n)$ derart, 
dass
\begin{enumerate}
\item[(a)] die adjungierte Gleichung
           \index{adjungierte Gleichung!Multi@-- Multiprozesse}
           \begin{equation}\label{PMPeinfach4MP} 
           \dot{p}(t) =- H^{\mathcal{M}}_x\big(t,x_*(t),u_*(t),p(t),\mathscr{A}_*\big),
           \end{equation} 
\item[(b)] in $t=t_1$ die Transversalitätsbedingung
           \index{Transversalitätsbedingungen!Multi@-- Multiprozesse}
           \begin{equation}\label{PMPeinfach5MP} 
           p(t_1)=0
           \end{equation} 
\item[(c)] und in fast allen Punkten $t \in [t_0,t_1]$ die Maximumbedingung
           \index{Maximumbedingung!Multi@-- Multiprozesse}
           \begin{equation}\label{PMPeinfach6MP} 
           H^{\mathcal{M}}\big(t,x_*(t),u_*(t),p(t),\mathscr{A}_*\big)
           = \max_{s \in\{1,...,k\}}\Big\{ \max_{u_s \in U_s} H_s\big(t,x_*(t),u_s,p(t)\big)\Big\}
           \end{equation}
\end{enumerate}
erfüllt sind.
\end{theorem}

\begin{beispiel} \label{AbschnittBeispielhybrid}
{\rm Wir untersuchen ein Investitionsmodell,
dass sich aus dem linearen und dem konkaven Beispielen \ref{BeispielLinInv} und \ref{BeispielKonInv} zusammensetzt:
\begin{equation}
\left.\begin{array}{l}
\hspace*{-5mm}
\displaystyle J\big(x(\cdot),u(\cdot),\mathscr{A} \big)  = \int_0^T \chi_{\mathscr{A}}(t) \circ f\big(t,x(t),u(t)\big) dt \to \sup,  \\[3mm]
\hspace*{-5mm}   
\displaystyle \dot{x}(t) = \chi_{\mathscr{A}}(t) \circ \varphi\big(t,x(t),u(t)\big), \quad x(0)=x_0>0,  \\[2mm]
\hspace*{-5mm}   
u(t) \in [0,1] \times [0,1], \quad \mathscr{A}=\{ \mathscr{A}_1,\mathscr{A}_2 \} \in \mathscr{Z}^2([0,T]).
\end{array} \right\}
\end{equation}
In den einzelnen Steuerungssystemen sind
$$\begin{array}{ll}
     f_1(t,x,u_1) = (1-u_1)x,        & \varphi_1(t,x,u_1) = u_1 x, \\[1mm]
     f_2(t,x,u_2) = (1-u_2)x^\alpha, & \varphi_2(t,x,u_2) = u_2 x^\alpha,
  \end{array} \quad U_1=U_2=[0,1]$$
und es gelten für die Modellparameter
$$\alpha \in (0,1) \mbox{ konstant}, \qquad T \mbox{ fest mit } T > \max \left\{1, \frac{x_0^{1-\alpha}}{\alpha} \right\}.$$
Da $U_1=U_2$ ist, unterscheiden wir nicht zwischen den Steuervariablen $u_1=u_2$ und schreiben $u$.
Wir wenden Theorem \ref{SatzPMPeinfachMP} an.
Es gilt nach (\ref{PMPeinfach6MP}):
\begin{eqnarray*}
\lefteqn{H^{\mathcal{M}}\big(t,x_*(t),u_*(t),p(t),\mathscr{A}_*\big)
         =\max_{s \in \{1,2\}} H_s\big(t,x_*(t),u_*(t),p(t)\big)} \\
&=& \max_{s \in \{1,2\}} \Big\{ \max_{u \in [0,1]} \big[ \big( p(t) -1 \big) u + 1 \big] \cdot x_*(t),
                \max_{u \in [0,1]} \big[ \big( p(t) -1 \big) u + 1 \big] \cdot x_*^\alpha(t) \Big\}.
\end{eqnarray*}
Dies k"onnen wir weiterhin in die Form
$$\max_{u \in [0,1]} \big[ \big( p(t) -1 \big) u \big] \cdot \max_{s \in \{1,2\}}\{ x_*(t) , x_*^\alpha(t) \}$$
bringen.
Daraus erhalten wir f"ur die optimale Investitionsrate und Wechselstrategie
$$u_*(t) = \left\{ \begin{array}{ll} 1, & p(t) < 1, \\ 0, & p(t) > 1, \end{array} \right. \qquad
  \chi_{\mathscr{A}_*}(t) = \left\{ \begin{array}{ll} ( 1,0 ), & x_*(t) > 1, \\ ( 0,1 ), & 0 < x_*(t) < 1. \end{array} \right.$$
Die Funktion $p(\cdot)$ ist die L"osung der adjungierten Gleichung (\ref{PMPeinfach4MP}):
$$\dot{p}(t) = \left\{ \begin{array}{ll} - \big[ \big( p(t) -1 \big) u_*(t) + 1 \big],    & t \in \mathscr{A}_1, \\[1mm]
          - \big[ \big( p(t) -1 \big) u_*(t) + 1 \big] \cdot \alpha x_*^{\alpha-1}(t), & t \in \mathscr{A}_2, \end{array} \right.
   \quad p(T)=0.$$
Betrachten wir die einzelnen Steuerungssysteme,
so sind nach den Beispielen \ref{BeispielLinInv} und \ref{BeispielKonInv} durch $\tau_1=T-1$ bzw. $\tau_2=\alpha T -x_0^{1-\alpha}$
die Zeitpunkte f"ur den optimalen Wechsel von vollst"andiger Investition in komplette Kosumption gegeben.
Au"serdem ist f"ur die L"osung $x(\cdot)$ der Differentiagleichung $\dot{x}(t)=x^\alpha(t)$ mit $x_0 \in (0,1)$ der
Zeitpunkt $\sigma$, in dem $x(\sigma)=1$ gilt, durch $\sigma=\frac{1-x_0^{1-\alpha}}{1-\alpha}$ bestimmt.
Im Weiteren seien also
$$\tau_1=T-1,  \qquad \tau_2=\alpha T -x_0^{1-\alpha}, \qquad \sigma=\frac{1-x_0^{1-\alpha}}{1-\alpha}.$$
Wir diskutieren die m"oglichen Szenarien:
\begin{enumerate}
\item[(a)] Sei $x_0 \geq 1$: In diesem Fall lautet der Kandidat
           \begin{eqnarray*}
           && x_*(t) = \left\{ \begin{array}{ll} x_0 \cdot e^t, & t \in [0,\tau_1), \\
                                                 x_0 \cdot e^{\tau_1}, & t \in [\tau_1,T], \end{array} \right.
              \qquad u_*(t) = \left\{ \begin{array}{ll} 1, & t \in [0,\tau_1), \\
                                                         0, & t \in [\tau_1,T], \end{array} \right. \\
           && \chi_{\mathscr{A}_*}(t) = (1,0), \qquad J\big(x_*(\cdot),u_*(\cdot),\mathscr{A}_*\big) = x_0 \cdot e^{\tau_1}.          
           \end{eqnarray*}
\item[(b)] Seien $x_0<1$ und $T < \frac{1 - \alpha x_0^{1-\alpha}}{\alpha (1-\alpha)}$:
           Dann ist $\tau_2<\sigma$ und der Kandidat lautet
           \begin{eqnarray*}
           && y_*(t) = \left\{ \begin{array}{ll}
                       \big[ (1-\alpha)t +  x_0^{1-\alpha} \big]^\frac{1}{1-\alpha}, & t \in [0,\tau_2), \\
                       \big[ \alpha(1-\alpha)T + \alpha x_0^{1-\alpha} \big]^\frac{1}{1-\alpha}, & t \in [\tau_2,T],
                       \end{array} \right.
              \quad v_*(t) = \left\{ \begin{array}{ll} 1, & t \in [0,\tau_2), \\
                                                         0, & t \in [\tau_2,T], \end{array} \right. \\
           && \chi_{\mathscr{B}_*}(t) = (0,1), \qquad 
              J\big(y_*(\cdot),v_*(\cdot),\mathscr{B}_*\big) = 
              \alpha^\frac{\alpha}{1-\alpha} \cdot \big[ (1-\alpha)T + x_0^{1-\alpha} \big]^\frac{1}{1-\alpha}.          
           \end{eqnarray*}
\item[(c)] Seien $x_0<1$ und $T > \frac{1 - \alpha x_0^{1-\alpha}}{\alpha (1-\alpha)}$:
           Dann ist $\sigma < \tau_1$ und wir erhalten den Kandidaten
           \begin{eqnarray*}
           && z_*(t) = \left\{ \begin{array}{ll}
                       \big[ (1-\alpha)t +  x_0^{1-\alpha} \big]^\frac{1}{1-\alpha}, & t \in [0,\sigma), \\
                       e^{t-\sigma}, & t \in [\sigma,\tau_1), \\
                       e^{\tau_1-\sigma}, & t \in [\tau_1,T], \end{array} \right.                      
              \quad w_*(t)= \left\{ \begin{array}{ll} 1, & t \in [0,\tau_1), \\
                                                      0, & t \in [\tau_1,T], \end{array} \right. \\
           && \chi_{\mathscr{C}_*}(t) = \left\{ \begin{array}{ll} (0,1), & t \in [0,\sigma), \\
                                                         (1,0), & t \in [\sigma,T], \end{array} \right.
              \qquad J\big(z_*(\cdot),w_*(\cdot),\mathscr{C}_*\big) = e^{\tau_1-\sigma}.        
           \end{eqnarray*}
\item[(d)] Seien $x_0<1$ und $T = \frac{1 - \alpha x_0^{1-\alpha}}{\alpha (1-\alpha)}$:
           Wegen $\sigma=\tau_2$ sind alle Bedingungen von Theorem \ref{PMPeinfach6MP} f"ur die Multiprozesse
           $\big(y_*(\cdot),v_*(\cdot),\mathscr{B}_*\big)$ und $\big(z_*(\cdot),w_*(\cdot),\mathscr{C}_*\big)$ erf"ullt. \\
           Sei $n_0 \in \N$ mit $n_0 >\alpha/(1-\alpha)$.
           F"ur $n \geq n_0$ betrachten wir die Steuerungen
           $$v_n(t)= v_*(t) + \chi_{[\sigma,\sigma+\frac{1}{n})}(t) \cdot \big(w_*(t)-v_*(t)\big),
             \quad \chi_{\mathscr{B}_n}(t)= \chi_{\mathscr{C}_*}(t),$$
           die wie im Fall (c) einen Wechsel des Steuerungssystems zum Zeitpunkt $t=\sigma=\tau_2$ und eine verl"angerte Investitionsphase
           vorgeben.
           F"ur den zugeh"origen Kapitalbestand $y_n(\cdot)$ gilt $\|y_n(\cdot)-y_*(\cdot)\|_\infty=e^{\frac{1}{n}}-1$
           und f"ur den Wert des Zielfunktionals
           \begin{eqnarray*}
           \lefteqn{J\big(y_n(\cdot),v_n(\cdot),\mathscr{B}_n\big)-J\big(y_*(\cdot),v_*(\cdot),\mathscr{B}_*\big)
              =e^{\frac{1}{n}}\bigg(T-\sigma-\frac{1}{n}\bigg)- (T-\sigma)} \\
           && >\bigg(1+\frac{1}{n}\bigg)\bigg(\frac{1}{\alpha}-\frac{1}{n}\bigg)-\frac{1}{\alpha}
              =\frac{1}{n}\bigg(\frac{1}{\alpha}-1-\frac{1}{n}\bigg)>0 \quad\mbox{ f"ur alle } n \geq n_0.
           \end{eqnarray*}
           Daher stellt $\big(y_*(\cdot),v_*(\cdot),\mathscr{B}_*\big)$ in diesem Fall kein starkes lokales Maximum dar.
           Der Multiprozess $\big(z_*(\cdot),w_*(\cdot),\mathscr{C}_*\big)$ ist der einzige der Kandidat.
\end{enumerate}
Die Diskussion der notwendigen Optimalitätsbedingungen ist damit abgeschlossen. \hfill $\square$}
\end{beispiel}
   
       \newpage
       \lhead[\thepage \hspace*{1mm} Beweis des Maximumprinzips]{ }   
       \subsection{Der Beweis des Maximumprinzips}
Wir passen den Beweis für die elementare Standardaufgabe in Abschnitt \ref{AbschnittPMPBeweiseinfach}
an die zusätzliche Steuerung $\mathscr{A} \in \mathscr{Z}^k([t_0,t_1])$ an:
Es sei $\tau \in (t_0,t_1)$ ein Stetigkeitspunkt der Steuerungen $u_*(\cdot)$ und $\chi_{\mathcal{A}_*}(\cdot)$.
Weiterhin seien $v=(v_1,...,v_k) \in U=U_1 \times ... \times U_k$, $s \in \{1,...,k\}$ fest gewählt und
es bezeichne $\mathscr{A}(s) \in \mathscr{Z}^k([t_0,t_1])$ diejenige Zerlegung, für die $\mathscr{A}_s(s)=[t_0,t_1]$ ist.
Damit definieren wir die einfachen Nadelvariationen $u_\lambda(\cdot)$, $\chi_{\mathcal{A}_\lambda}(\cdot)$ durch \index{Nadelvariation, einfache}
$$u_{\lambda}(t) = 
  \left\{ \begin{array}{ll}
          u_*(t) & \mbox{ f"ur } t \not\in [\tau-\lambda,\tau), \\
          v      & \mbox{ f"ur } t     \in [\tau-\lambda,\tau),
          \end{array} \right. \qquad
   \chi_{\mathcal{A}_\lambda}(t) =
   \left\{ \begin{array}{ll}
            \chi_{\mathcal{A}_*}(t) & \mbox{ f"ur } t \not\in [\tau-\lambda,\tau), \\
            \chi_{\mathcal{A}(s)}(t)     & \mbox{ f"ur } t     \in [\tau-\lambda,\tau).
          \end{array} \right.$$ 
Dabei wählt $\chi_{\mathcal{A}_\lambda}(\cdot)$ über $[\tau-\lambda,\tau)$ von den gegebenen Steuerungssystemen
dasjenige zum Index $s \in \{1,...,k\}$ aus.
Weiter sei $x_\lambda(\cdot)$ die eindeutige L"osung der Gleichung
$$\dot{x}_\lambda(t) = \chi_{\mathcal{A}_\lambda}(t) \circ \varphi\big(t,x(t),u_\lambda(t)\big), \quad x(t_0)=x_0.$$
F"ur $t \geq \tau$ untersuchen wir den Grenzwert
$\displaystyle y(t)=\lim_{\lambda \to 0^+}\frac{x_{\lambda}(t) - x_*(t)}{\lambda}$.
Auf die gleiche Weise wie im Abschnitt \ref{AbschnittPMPBeweiseinfach} zeigt sich,
dass dieser Grenzwert für alle $t \in [\tau,t_1]$ existiert,
in $t =\tau$ die Gleichung
$$y(\tau) = \chi_{\mathcal{A}_\lambda}(\tau) \circ \varphi\big(\tau,x_*(\tau),u_\lambda(\tau)\big) 
            - \chi_{\mathcal{A}_*}(\tau) \circ \varphi\big(\tau,x_*(\tau),u_*(\tau)\big)$$
gilt und ferner $y(t)$ über $[\tau,t_1]$ der Integralgleichung
$$y(t) = y(\tau) + \int_{\tau}^{t} \big[\chi_{\mathcal{A}_*}(s) \circ \varphi_x\big(s,x_*(s),u_*(s)\big) \big] \, y(s) \, ds$$
genügt.
Weiterhin ergibt sich die Beziehung
$$\langle p(\tau) , y(\tau) \rangle
 = - \int_{\tau}^{t_1} \big\langle \chi_{\mathcal{A}_*}(t) \circ f_x\big(t,x_*(t),u_*(t)\big), y(t) \big\rangle dt.$$
Abschließend zeigt sich
\begin{eqnarray*}
0 &\leq& \lim_{\lambda \to 0^+} \frac{J\big(x_\lambda(\cdot),u_\lambda(\cdot),\mathcal{A}_\lambda\big)- J\big(x_*(\cdot),u_*(\cdot),\mathcal{A}_*\big)}{\lambda} \\
  &=& \chi_{\mathcal{A}_\lambda}(\tau) \circ f\big(\tau,x_*(\tau),u_\lambda(\tau)\big) - \chi_{\mathcal{A}_*}(\tau) \circ f\big(\tau,x_*(\tau),u_*(\tau)\big) \\
  & & + \int_{\tau}^{t_1} \big\langle \chi_{\mathcal{A}_*}(t) \circ f_x\big(t,x_*(t),u_*(t)\big) , y(t) \big\rangle \, dt \\
  &=& \chi_{\mathcal{A}_\lambda}(\tau) \circ f\big(\tau,x_*(\tau),u_\lambda(\tau)\big) - \chi_{\mathcal{A}_*}(\tau) \circ f\big(\tau,x_*(\tau),u_*(\tau)\big) \\
  & & - \big\langle p(\tau), \chi_{\mathcal{A}_\lambda}(\tau) \circ \varphi\big(\tau,x_*(\tau),u_\lambda(\tau)\big) 
            - \chi_{\mathcal{A}_*}(\tau) \circ \varphi\big(\tau,x_*(\tau),u_*(\tau)\big) \big\rangle,
\end{eqnarray*}
welche zur Gültigkeit der Maximumbedingung (\ref{PMPeinfach6MP}) führt.  \hfill $\blacksquare$
       \newpage
       \lhead[\thepage \hspace*{1mm} Arrow-Bedingungen]{ }
       \subsection{Hinreichende Bedingungen nach Arrow} \index{hinreichende Bedingungen nach Arrow!Multi@-- Multiprozesse}
Die Aufgabe (\ref{PMPeinfach1MP})--(\ref{PMPeinfach3MP}) eines Multiprozesses enth"alt den Steuerbereich $\Zt$,
der die Wechselstrategien charakterisiert und keine konvexe Menge darstellt.
Es wird sich im Beispiel \ref{AbschnittBeispielhybrid1} zeigen,
dass die hinreichenden Arrow-Bedingungen im Fall optimaler Multiprozesse recht eingeschränkt anwendbar sind. \\[2mm]
Für $t \in [t_0,t_1]$ sei wieder $V_\gamma(t) = \{ x \in \R^n \,|\, \|x-x_*(t)\| < \gamma\}$.
Au"serdem bezeichnet $\mathscr{H}^{\mathcal{M}}$ die Hamilton-Funktion
\begin{equation} \label{PMPHamiltonMP}
\mathscr{H}^{\mathcal{M}}(t,x,p) = \max_{s \in\{1,...,k\}}\big\{ \max_{u_s \in U_s} H_s(t,x,u_s,p)\big\}.
\end{equation}

\begin{theorem} \label{SatzHBMP}
In der Aufgabe (\ref{PMPeinfach1MP})--(\ref{PMPeinfach3MP}) sei
$\big(x_*(\cdot),u_*(\cdot),\mathscr{A}_*\big) \in \mathscr{A}^{\mathcal{M}}_{\rm adm} \cap \mathscr{A}^{\mathcal{M}}_{\rm Lip}$
und es sei $p(\cdot) \in PC_1([t_0,t_1],\R^n)$. Ferner gelte:
\begin{enumerate}
\item[(a)] Das Quadrupel $\big(x_*(\cdot),u_*(\cdot),\mathscr{A}_*,p(\cdot)\big)$
           erf"ullt (\ref{PMPeinfach4MP})--(\ref{PMPeinfach6MP}) in Theorem \ref{SatzPMPeinfachMP}.        
\item[(b)] F"ur jedes $t \in [t_0,t_1]$ ist die Funktion $\mathscr{H}^{\mathcal{M}}\big(t,x,p(t)\big)$ konkav in $x$ auf $V_\gamma(t)$.
\end{enumerate}
Dann ist $\big(x_*(\cdot),u_*(\cdot),\mathscr{A}_*\big)$ ein starkes lokales Minimum der Aufgabe (\ref{PMPeinfach1MP})--(\ref{PMPeinfach3MP}).
\end{theorem}

{\bf Beweis} Auf die gleiche Weise wie die Ungleichungen (\ref{BeweisHBPMP2}), (\ref{BeweisHBPMP4}) im Beweis von Theorem \ref{SatzHBPMP}
in Abschnitt \ref{AbschnittHBPMP} ergeben sich:
Für alle $x \in V\gamma(t)$ und alle $t \in [t_0,t_1]$, in denen die Maximumbedingung (\ref{PMPeinfach6MP}) erfüllt ist, gelten
\begin{eqnarray*}
       \langle a(t),x-x_*(t)\rangle
&\geq& -\mathscr{H}^{\mathcal{M}}\big(t,x_*(t),p(t)\big) + \mathscr{H}^{\mathcal{M}}\big(t,x,p(t)\big), \\
       \langle \dot{p}(t),x-x_*(t)\rangle 
&\leq& \mathscr{H}^{\mathcal{M}}\big(t,x_*(t),p(t)\big)- \mathscr{H}^{\mathcal{M}}\big(t,x,p(t)\big).
\end{eqnarray*}
Es sei $\big(x(\cdot),u(\cdot),\mathscr{A}\big) \in \mathscr{A}^{\mathcal{M}}_{\rm adm}$ mit $\|x(\cdot)-x_*(\cdot)\|_\infty < \gamma$.
Dann erhalten wir
\begin{eqnarray*}
    \lefteqn{J\big(x(\cdot),u(\cdot),\mathscr{A}\big)-J\big(x_*(\cdot),u_*(\cdot),\mathscr{A}_*\big)} \\
&=& \int_{t_0}^{t_1} \big[\chi_{\mathscr{A}}(t) \circ f\big(t,x(t),u(t)\big)-\chi_{\mathscr{A}_*}(t) \circ f\big(t,x_*(t),u_*(t)\big)\big] \, dt \\
&\geq& \int_{t_0}^{t_1} \big[\mathscr{H}^{\mathcal{M}}\big(t,x_*(t),p(t)\big)-\mathscr{H}^{\mathcal{M}}\big(t,x(t),p(t)\big)\big] \, dt 
    + \int_{t_0}^{t_1} \langle p(t), \dot{x}(t)-\dot{x}_*(t) \rangle dt \\
&\geq& \int_{t_0}^{t_1} \big[\langle \dot{p}(t),x(t)-x_*(t)\rangle + \langle p(t), \dot{x}(t)-\dot{x}_*(t) \rangle \big] \, dt \\
&=& \langle p(t_1),x(t_1)-x_*(t_1)\rangle-\langle p(t_0),x(t_0)-x_*(t_0)\rangle.
\end{eqnarray*}
Mit $p(t_1)=0$ und $x(t_0)=x_*(t_0)=x_0$ ergibt sich
$J\big(x(\cdot),u(\cdot),\mathscr{A}\big)-J\big(x_*(\cdot),u_*(\cdot),\mathscr{A}_*\big) \geq 0$ für alle
$\big(x(\cdot),u(\cdot),\mathscr{A}\big) \in \mathscr{A}^{\mathcal{M}}_{\rm adm}$ mit $\|x(\cdot)-x_*(\cdot)\|_\infty < \gamma$. \hfill $\blacksquare$ \\[2mm]
Die Aufgabe eines Multiprozesses enth"alt neben der "ublichen Steuerung noch mittels der Wechselstrategie die M"oglichkeit auf die Auswahl eines
bestimmten Steuerungssystems.
Diese Vermischung von kontinuierlichen und diskreten Steuervariablen erschwert den Nachweis der Konkavit"at der Hamilton-Funktion.

\begin{beispiel} \label{AbschnittBeispielhybrid1}
{\rm Das Investitionsmodell im Beispiel \ref{AbschnittBeispielhybrid} lauten die Pontrjagin-Funktionen $H_1,H_2$
zu den beiden Steurungssystemen
$$H_1(t,x,u_1,p)=[(p -1) u_1 + 1 ] \cdot x, \quad H_2(t,x,u_2,p)=[(p -1) u_2 + 1 ] \cdot x^\alpha.$$
Wegen $u_1,u_2 \in [0,1]$ beachten wir nur eine Steuervariable $u$.
Dies führt zu der Hamilton-Funktion
$$\mathscr{H}^{\mathcal{M}}(t,x,p) = \max_{\substack{u_s \in [0,1] \\ s \in \{1,2\}}} H_s(t,x,u_s,p)= \max\{p,1\} \cdot \max\{ x , x^\alpha \},$$
welche die Funktion
$$f(x)= \max\{p,1\} \cdot \left\{\begin{array}{ll} x^\alpha, & x \in (0,1), \\ x, & x\geq 1, \end{array} \right.$$
enthält, die in $x=1$ nicht konkav ist. \\[2mm]
Damit ist die Hamilton-Funktion $\mathscr{H}^{\mathcal{M}}\big(t,x,p(t)\big)$ bez"uglich der Trajektorie $x_*(\cdot)$
genau dann f"ur jedes $t \in [t_0,t_1]$ konkav in der Variablen $x$ auf $V_\gamma(t)$,
wenn $x_*(t) \not=1$ f"ur alle $t \in [t_0,t_1]$ gilt.
F"ur die einzelnen F"alle im Beispiel \ref{AbschnittBeispielhybrid} ergeben sich damit:
\begin{enumerate}
\item[(a)] Theorem \ref{SatzHBMP} ist nur dann anwendbar, wenn $x_0>1$ gilt.
           In diesem Fall ist der Kandidat $\big(x_*(\cdot),u_*(\cdot),\mathscr{A}_*\big)$ optimal.
\item[(b)] Da in diesem Fall $x_0\leq y_*(t)<1$ auf $[0,T]$ gilt, ist Theorem \ref{SatzHBMP} anwendbar und
           der Kandidat $\big(y_*(\cdot),v_*(\cdot),\mathscr{B}_*\big)$ optimal.
\item[(c)] Wegen $z_*(\sigma)=1$ ist Theorem \ref{SatzHBMP} nicht anwendbar.
\item[(d)] Wegen $y_*(\sigma)=1$ und $z_*(\sigma)=1$ ist Theorem \ref{SatzHBMP} auf keinen der ermittelten Kandidaten
           $\big(y_*(\cdot),v_*(\cdot),\mathscr{B}_*\big)$, $\big(z_*(\cdot),w_*(\cdot),\mathscr{C}_*\big)$ anwendbar.
           Die lokale Optimalität des Kandidaten $\big(y_*(\cdot),v_*(\cdot),\mathscr{B}_*\big)$ wurde bereits
           ausgeschlossen. \hfill $\square$
\end{enumerate}}
\end{beispiel}
\cleardoublepage

\lhead[ ]{}
\rhead[]{Unendlicher Zeithorizont \hspace*{1mm} \thepage}
\section{Aufgaben mit unendlichem Zeithorizont} 
       Steuerungsprobleme mit einem unendlichen Zeithorizont finden z.\,B. bei dynamischen System,
welche nur asymptotische Stabilität aufweisen,
oder in der Ökonomischen Wachstumstheorie Anwendung.
Die Einbindung dieser Klasse in die Wirtschaftstheorie geht auf Frank P. Ramsey (1903--1930) zurück,
der in der Arbeit \cite{Ramsey} aus dem Jahr 1928 der Frage einer optimalen Sparquote nachgeht,
die einer "Okonomie langfristiges und wohlfahrtsoptimiertes Wachstum garantiert.
Das besondere an der Modellierung des Variationsproblems war die Einf"uhrung des unendlichen Zeithorizontes.
Die Philosophie dahinter ist die Vorstellung,
dass kein nat"urliches Ende f"ur den Betrachtungszeitraum existiert.
M"ochte man s"amtlichen nachfolgenden Generationen Beachtung schenken, dann ist die Idealisierung in Form des
zeitlich unbeschr"ankten Rahmens die einzig m"ogliche Konsequenz. \\[2mm]
In der Literatur werden die Steuerungsprobleme mit unendlichem Zeithorizont häufig als eine Folgerung
der Standardaufgabe dargestellt.
Dabei sind die Hürden, die das unbeschränkte Zeitintervall mit sich führt,
von grundlegender Natur und mit den bisher zur Verfügung gestellten Methoden nicht erfassbar.
Anhand der elementaren Aufgabe mit freiem Endpunkt lassen sich diese Schwierigkeiten bereits sehr gut verdeutlichen.
Letztendlich zeigt sich,
dass ohne zusätzliche Argumente mit der einfachen Nadelvariation kein vollständiger Beweis eines Pontrjaginschen Maximumprinzips
für die Aufgabe mit unendlichem Zeithorizont erbracht werden kann. 
       \lhead[\thepage \hspace*{1mm} Vorbereitungen]{ }     
       \subsection{Pathologien und Variierbarkeit}
Steuerungsprobleme mit unendlichem Zeithorizont werden häufig als Aufgabe über einem sehr langen Planungszeitraum aufgefasst.
Auf diese Weise wird die Aufgabenstellung mit unbeschränktem Zeitintervall auf eine Aufgabe mit endlichem Zeithorizont reduziert,
für die die bekannten Ergebnisse angewendet werden können.
Diese Methode der Approximation mit endlichem Horizont ist in der Literatur weit verbreitet. \\
Allerdings kann dieser Ansatz zu pathologischen Situationen führen.
Einige der entarteten Fälle sind bei Aseev \& Kryazhimskii \cite{AseKry} dokumentiert.
Die wesentliche Ursache für das Auftreten von Pathologien besteht bei der Methode der Approximation mit endlichem Horizont in der Erwartung,
dass im Grenzübergang von endlichen zum unendlichen Zeithorizont die Struktur des optimalen Steuerungsprozesses und dessen Optimalwert,
aber auch die Optimalitätsbedingungen des Pontrjaginschen Maximumprinzips ein stetiges Verhalten aufweisen. 
Im Allgemeinen darf auf diese Hoffnung nicht gebaut werden. \\[2mm]
Um diese Schilderungen zu untermauern, betrachten wir illustrative Beispiele.
Wir beginnen mit einer Variante von Beispiel \ref{BeispielLinInv}:

\begin{beispiel} \label{BeispielInvestmentUH}
{\rm Wir betrachten die Aufgabe
$$\int_0^T e^{-\varrho t} \big(1-u(t)\big)x(t) \, dt \to \sup, \quad
  \dot{x}(t) = u(t)x(t), \quad x(0)=1, \quad u(t) \in [0,1], \quad \varrho \in (0,1).$$
Für jedes feste $T>-\ln (1-\varrho)/\varrho$ erhalten wir mit den gleichen Argumenten wie im Beispiel \ref{BeispielLinInv} den global
optimalen Steuerungsprozess
$$x^T_*(t)= \left\{ \begin{array}{ll} e^t,& t \in [0,\tau), \\[1mm] e^\tau,& t \in [\tau,T], \end{array} \right. \quad
  u^T_*(t)= \left\{ \begin{array}{ll} 1,& t \in [0,\tau), \\[1mm] 0,& t \in [\tau,T], \end{array} \right. \quad
  \tau = T + \frac{\ln (1-\varrho)}{\varrho}.$$
Betrachten wir den Grenz"ubergang $T \to \infty$,
dann konvergiert die Familie $\big\{\big(x_*^T(\cdot),u_*^T(\cdot)\big)\big\}$, $T \in \R_+$, punktweise gegen den
Steuerungsprozess
$$x_*(t) = e^t, \qquad u_*(t)=1, \qquad t \in \R_+.$$
Dieses Paar ist das globale Minimum der Aufgabe mit unendlichem Zeithorizont. \hfill $\square$}
\end{beispiel}

Im Beispiel \ref{BeispielInvestmentUH} bleibt der Wert der Zielfunktionals für die optimalen Steuerungsprozesse $\big(x_*^T(\cdot),u_*^T(\cdot)\big)$
im Grenzwert für $T \to \infty$ nicht endlich.
Die Aufgabe ist demnach im Unendlichen unstetig.
Deswegen erweist sich die Approximation mit endlichem Zeithorizont an dieser Stelle als ungeeignet.
Demgegenüber erfüllen die folgenden Beispiele gewisse Anforderungen an Endlichkeit und Stetigkeit im Unendlichen.
Trotzdem besitzen die notwendigen Optimalitätsbedingungen nicht die erwartete Form. \\[2mm]
Im Vorwort zur Standardaufgabe (\ref{PMPeinfach1})--(\ref{PMPeinfach3}) im Abschnitt \ref{AbschnittPMPeinfach} haben wir angekündigt,
uns lediglich auf den ``normalen'' Fall notwendiger Optimilitätsbedingungen zu beschränken.
In Verallgemeinerung zum normalen Fall\index{Pontrjaginsches Maximumprinzip!normal@-- normaler Fall} besitzt die Pontrjagin-Funktion die Gestalt
$$H(t,x,u,p,\lambda_0) = \langle p, \varphi(t,x,u) \rangle - \lambda_0 f(t,x,u)$$
mit dem zusätzlichen Faktor $\lambda_0 \in \R$.
In diesem Fall liefert das Pontrjaginsche Maximumprinzip für die Aufgabe ohne Terminalfunktional $S$ die Existenz nichttrivialer Multiplikatoren
$\lambda_0 \geq 0$ und $p(\cdot) \in PC_1([t_0,t_1],\R^n)$ derart,
dass die adjungierte Gleichung
$$\dot{p}(t) = -\varphi_x^T \big(t,x_*(t),u_*(t)\big) \, p(t) + \lambda_0 f_x\big(t,x_*(t),u_*(t)\big)$$
zur Transversalitätsbedingung $p(t_1)=0$ und in fast allen Punkten $t \in [t_0,t_1]$ die Maximumbedingung
$$H\big(t,x_*(t),u_*(t),p(t),\lambda_0\big) = \max_{u \in U}H\big(t,x_*(t),u,p(t),\lambda_0\big)$$
erfüllt sind.
In unseren Betrachtungen zu Aufgaben mit freiem rechten Endpunkt ergibt sich,
dass stets der normale Fall mit $\lambda_0=1$ vorliegt.
Die nachfolgenden Beispiele von Halkin \cite{Halkin} zeigen nun,
dass im Übergang zum unendlichen Zeithorizont einerseits der normale Fall mit $\lambda_0=1$ annulliert werden kann
und andererseits auf Basis von $p(t_1)=0$ die erwartete Transversalitätsbedingung $p(t) \to 0$ für $t \to \infty$ nicht eintreten muss.

\begin{beispiel} {\rm Wir diskutieren nach Halkin \cite{Halkin} die Aufgabe
$$\int_0^\infty u(t)\big(1-x(t)\big)\, dt \to \sup, \quad \dot{x}(t)=u(t)\big(1-x(t)\big), \quad x(0)=0,\quad u \in [0,\infty).$$
Die Dynamik lässt sich mittels der Trennung der Veränderlichen behandeln.
Auf diese Weise ergibt sich für den Zustand des Steuerungsprozesses $\big(x(\cdot),u(\cdot)\big)$ die Darstellung
$$x(t)= 1-\exp\bigg\{ -\int_0^t u(s)\,ds \bigg\}.$$
Da $u(t) \geq 0$ ausfällt, existiert für jede uneigentlich integrierbare Funktion $u(\cdot)$ mit endlichem oder unendlichem Integralwert
der Grenzwert der Zustandes $x(t)$ für $t \to \infty$ und besitzt einen Wert aus $[0,1]$. 
Ferner ergibt sich im Zielfunktional
$$\int_0^\infty u(t)\big(1-x(t)\big)\, dt = \int_0^\infty \dot{x}(t)\, dt = \lim_{t \to \infty} x(t) \in [0,1].$$
Daher stellt jeder zulässige Steuerungsprozess $\big(x(\cdot),u(\cdot)\big)$ mit $x(t) \to 1$ für $t \to \infty$ ein globales Maximum dar. \\[2mm]
Wir diskutieren die notwendigen Optimalitätsbedingungen für den optimalen Steuerungsprozess
$\big(x_*(\cdot),u_*(\cdot)\big)$ mit $u_*(t)=1$ und $x_*(t)=1-e^{-t}$ über $[0,\infty)$:
In dieser Aufgabe besitzt die Pontrjagin-Funktion die Form $H(t,x,u,p)=u(1+p)(1-x)$.
Die Anwendung der Maximumbedingung (\ref{PMPeinfach6}) in Theorem \ref{SatzPMPeinfach} führt auf die Beziehung
$$H\big(t,x_*(t),u_*(t),p(t)\big) = \max_{u \in [0,\infty)} u \big(1+p(t)\big)\big(1-x_*(t)\big).$$
Diese Bedingung kann für den inneren Wert $u_*(t)=1 \in [0,\infty)$ nur dann erfüllt sein, 
wenn $p(t) = -1$ für alle $t \in [0,\infty)$ gilt.
Damit ergibt sich $p(t) \not\to 0$ für $t \to \infty$ und die Gültigkeit einer gleichartigen Transversalitätsbedingung
zu (\ref{PMPeinfach4}) stellt sich nicht ein. \hfill $\square$}
\end{beispiel}

\begin{beispiel} {\rm In der folgenden Variante eines Beispieles nach Halkin \cite{Halkin},
\begin{eqnarray*}
&& J\big(x(\cdot),u(\cdot)\big) = \int_0^\infty e^{-\varrho t}\big(u(t)-x(t)\big) \, dt \to \sup, \\
&& \dot{x}(t) = x(t)+u(t), \quad x(0)=0, \quad  u(t) \in [0,1], \quad \varrho \in (0,1),
\end{eqnarray*}
erhalten wir für einen zulässigen Steuerungsprozess $\big(x(\cdot),u(\cdot)\big)$ für den Zustand
$$x(t)=e^t \cdot \int_0^t u(s)e^{-s} \, ds, \quad t \geq 0.$$
Damit liefert der Steuerungsprozess $\big(x_*(t),u_*(t)\big) \equiv (0,0)$ das globale Maximum der Aufgabe.
Denn f"ur jeden anderen zulässigen Steuerungsprozess mit einer stückweise stetigen Steuerung,
die nicht identisch verschwindet, besitzt das Zielfunktional den Wert $-\infty$. 

\newpage
Wir werten die Bedingungen des Pontrjaginschen Maximumprinzips im allgemeinen Fall aus:
Die Pontrjagin-Funktion lautet $H(t,x,u,p,\lambda_0)=p[x+u]+\lambda_0 e^{-\varrho t}[u-x]$.
Für die Lösung $p(\cdot)$ der adjungierten Gleichung ergibt sich
$$\dot{p}(t)=-p(t)+\lambda_0 e^{-\varrho t} \quad\Rightarrow\quad
  p(t)=\bigg(p(0)-\frac{\lambda_0}{1-\varrho}\bigg) e^{-t}+\frac{\lambda_0}{1-\varrho} e^{-\varrho t},\; p(0) \in \R,$$
und es gilt $p(t) \to 0$ für $t \to \infty$.
Die Maximumbedingung lautet
$$H\big(t,x_*(t),u_*(t),p(t),\lambda_0\big) =
  \max_{u \in [0,1]} \Big[p(t)\big(x_*(t)+u\big)+\lambda_0 e^{-\varrho t}\big(u-x_*(t)\big)\Big]$$
und ist äquivalent zu der Maximierungsaufgabe
$$\max_{u \in [0,1]} u\big[p(t)+\lambda_0 e^{-\varrho t}\big]
  =\max_{u \in [0,1]} u \cdot \bigg[ \bigg(p(0)-\frac{\lambda_0}{1-\varrho}\bigg) e^{-t}+\frac{2-\varrho}{1-\varrho} \lambda_0e^{-\varrho t}\bigg].$$
Wäre nun $\lambda_0 >0$,
so fällt wegen $\varrho \in (0,1)$ der Ausdruck in der eckigen Klammer für alle hinreichend große $t$ positiv aus,
und $u_*(t) \equiv 0$ kann nicht der Maximumbedingung für fast alle $t \geq 0$ genügen.
Deswegen muss der anormale Fall der notwendigen Optimalitätsbedingungen mit $\lambda_0=0$ eintreten. \hfill $\square$}
\end{beispiel}

Unsere Untersuchungen der elementaren Aufgaben auf starke lokale Optimalstellen basieren auf der einfachen Nadelvariation\index{Nadelvariation, einfache}
$$u(t;v,\tau,\lambda) = u_{\lambda}(t) = 
  \left\{ \begin{array}{ll}
          u_*(t) & \mbox{ f"ur } t \not\in [\tau-\lambda,\tau), \\
          v      & \mbox{ f"ur } t     \in [\tau-\lambda,\tau), 
          \end{array} \right. \qquad v \in U.$$
Im Vergleich zum endlichen Zeitintervall ist die Variation mittels $u_\lambda(\cdot)$ nicht in jedem Fall möglich.
Zur Veranschaulichung betrachten wir die Dynamik
$$\dot{x}(t)=x(t)+x^2(t)+u(t), \qquad x(0)=0, \qquad u(t)\in U = [0,\infty).$$
Die Trajektorie $x_0(t) \equiv 0$ zur Steuerung $u_0(t)\equiv 0$ ist nicht variierbar.
Denn jede Steuerung $u(t;v,\tau,\lambda)$ mit $\lambda >0$ und $v>0$ liefert eine Trajektorie $x_\lambda(\cdot)$,
die als Lösung der Differentialgleichung nicht über der gesamten Halbachse $[0,\infty)$ existiert. \\
Im Beweis des Pontrjaginschen Maximumprinzips für die elementare Aufgabe sind die Abhängigkeitssätze im Anhang \ref{AnhangDGL} zentraler Bestandteil.
Auf das unbeschränkte Zeitintervall sind diese Resultate nicht unmittelbar übertragbar. \\[2mm]
In der nachfolgenden Untersuchung legen wir uns auf beschränkte Steuerungsprozesse fest,
für die das Zielfunktional endlich ausfällt und beschränken uns außerdem auf Variationen $x_\lambda(\cdot)$,
die gleichmäßig gegen $x_*(\cdot)$ über $[0,\infty)$ für $\lambda \to 0^+$ konvergieren.
Die Existenz derartiger Variationen bleibt an dieser Stelle offen,
da wir auf adäquate Abhängigkeitssätze nicht zugreifen können.  
       \lhead[\thepage \hspace*{1mm} Pontrjaginsches Maximumprinzip]{ }    
       \subsection{Die Aufgabenstellung und das Pontrjaginsche Maximumprinzip}
Die Aufgabe mit unendlichem Zeithorizont besitzt die Gestalt
\begin{eqnarray}
&& \label{UHA1} J\big(x(\cdot),u(\cdot)\big) = \int_0^\infty e^{-\varrho t} f\big(t,x(t),u(t)\big) \, dt \to \inf, \\
&& \label{UHA2} \dot{x}(t) = \varphi\big(t,x(t),u(t)\big), \quad x(0)=x_0, \\
&& \label{UHA3} u(t) \in U \subseteq \R^m, \quad U \not= \emptyset, \quad \varrho > 0.
\end{eqnarray}
Die Aufgabe (\ref{UHA1})--(\ref{UHA3}) mit unendlichem Zeithorizont untersuchen wir bez"uglich der Steuerungsprozesse
$\big(x(\cdot),u(\cdot)\big) \in PC_1([0,\infty),\R^n) \times PC([0,\infty),U)$.
Dabei bezeichnen wir mit $PC([0,\infty),\R)$ bzw. $PC_1([0,\infty),\R)$ die Räume derjenigen Funktionen,
die über $[0,\infty)$ beschränkt und die über jedem endlichen Intervall $[0,T]$ stückweise stetig bzw. stückweise stetig differenzierbar sind. \\[2mm]
Mit $\mathscr{A}^{\,\mathcal{U}}_{\rm Lip}$ bezeichnen wir die Menge aller Paare $\big(x(\cdot),u(\cdot)\big)$,
für die es ein $\gamma>0$ derart gibt,
dass die Abbildungen $f(t,x,u)$, $\varphi(t,x,u)$ für jede kompakte Teilmenge $U_1 \subset \R^m$
auf der Menge aller Punkte $(t,x,u) \in \R \times \R^n \times \R^m$ mit
$0 \leq t < \infty$, $\|x-x(t)\| < \gamma$ und $u \in U_1$
beschränkt, stetig in der Gesamtheit aller Variablen und stetig differenzierbar bezüglich $x$
mit beschränkten Ableitungen $f_x(t,x,u)$, $\varphi_x(t,x,u)$ sind. \\[2mm]
Das Paar $\big(x(\cdot),u(\cdot)\big) \in PC_1([0,\infty),\R^n) \times PC([0,\infty),U)$
hei"st ein zul"assiger Steuerungsprozess in der Aufgabe (\ref{UHA1})--(\ref{UHA3}),
falls $\big(x(\cdot),u(\cdot)\big)$ über jedem endlichen Intervall der Dynamik (\ref{UHA2}) zu $x(t_0)=x_0$ gen"ugt und
die Steuerbeschränkungen (\ref{UHA3}) erfüllt.
Mit $\mathscr{A}^{\mathcal{U}}_{\rm adm}$ bezeichnen wir die Menge der zul"assigen Steuerungsprozesse. \\[2mm]
Ein zul"assiger Steuerungsprozess $\big(x_*(\cdot),u_*(\cdot)\big)$ ist eine
starke lokale Minimalstelle\index{Minimum, starkes lokales!elementar@-- unendlicher Zeithorizont}
der Aufgabe (\ref{UHA1})--(\ref{UHA3}),
falls eine Zahl $\varepsilon > 0$ derart existiert, dass die Ungleichung 
$$J\big(x(\cdot),u(\cdot)\big) \geq J\big(x_*(\cdot),u_*(\cdot)\big)$$
f"ur alle $\big(x(\cdot),u(\cdot)\big) \in \mathscr{A}^{\mathcal{U}}_{\rm adm}$ mit $\|x(\cdot)-x_*(\cdot)\|_\infty < \varepsilon$ gilt. \\[2mm]
Eine wesentliche Herausforderung ist der Nachweis einer adjungierten Funktion,
die sowohl die adjungierte Gleichung als auch eine Transversalitätsbedingung im Unendlichen erfüllt.
Die Existenz der Adjungierten wird im Folgenden mit Hilfe der Bedingung
\begin{equation} \label{PMPBedingung}
\int_0^\infty \big\|\varphi_x\big(t,x(t),u(t)\big)\big\| \, dt < \infty
\end{equation}
nach Lemma \ref{LemmaDGL4} und \ref{LemmaDGL6} gesichert.
Allerdings ist diese Bedingung nicht unkritisch,
da sie zum Beispiel lineare Dynamiken mit konstanten Koeffizienten ausschließt.
Außerdem wird durch (\ref{PMPBedingung}) die Problematik der fehlenden Variierbarkeit einer Trajektorie nicht behoben. \\
Ungeachtet dessen bezeichnet $\mathscr{A}^{\mathcal{U}}_{\lim}$ die Menge aller Paare $\big(x(\cdot),u(\cdot)\big)$,
welche der Bedingung (\ref{PMPBedingung}) genügen. \\ \\[2mm] 
F"ur die Aufgabe (\ref{UHA1})--(\ref{UHA3}) bezeichnet
$H^{\mathcal{U}}: \R \times \R^n \times \R^m \times \R^n \to \R$ die Pontrjagin-Funktion
$$H^{\mathcal{U}}(t,x,u,p) = \langle p, \varphi(t,x,u) \rangle - e^{-\varrho t} f(t,x,u).$$

\begin{theorem}[Pontrjaginsches Maximumprinzip] \label{SatzPMPUHA}
\index{Pontrjaginsches Maximumprinzip!unendlich@-- unendlicher Zeithorizont} 
Sei $\big(x_*(\cdot),u_*(\cdot)\big) \in \mathscr{A}^{\mathcal{U}}_{\rm adm} \cap \mathscr{A}^{\mathcal{U}}_{\rm Lip} \cap \mathscr{A}^{\mathcal{U}}_{\lim}$.
Ist $\big(x_*(\cdot),u_*(\cdot)\big)$ ein starkes lokales Minimum der Aufgabe (\ref{UHA1})--(\ref{UHA3}),
dann existiert eine Vektorfunktion $p(\cdot) \in PC_1([0,\infty),\R^n)$ derart, dass
\begin{enumerate}
\item[(a)] die adjungierte Gleichung
           \index{adjungierte Gleichung!Standard@-- unendlicher Zeithorizont}
           \begin{equation}\label{PMPUHA1} 
           \dot{p}(t) = -H_x^{\mathcal{U}}\big(t,x_*(t),u_*(t),p(t)\big),
           \end{equation} 
\item[(b)] für $t \to \infty$ die Transversalitätsbedingung
           \index{Transversalitätsbedingungen!Standard@-- unendlicher Zeithorizont}
           \begin{equation}\label{PMPUHA2} 
           \lim_{t \to \infty} p(t)=0
           \end{equation} 
\item[(c)] und in fast allen Punkten $t \in [0,\infty)$ die Maximumbedingung           
           \index{Maximumbedingung!Standard@-- unendlicher Zeithorizont}
           \begin{equation}\label{PMPUHA3} 
           H^{\mathcal{U}}\big(t,x_*(t),u_*(t),p(t)\big) = \max_{u \in U} H^{\mathcal{U}}\big(t,x_*(t),u,p(t)\big)
           \end{equation}
\end{enumerate}
erfüllt sind.
\end{theorem}

\begin{bemerkung} {\rm
Im Folgenden sind die Wohldefiniertheit von $x_\lambda(\cdot)$ auf $[0,\infty)$ und die gleichmäßige Konvergenz von $x_\lambda(\cdot)$ gegen $x_*(\cdot)$
nicht gesichert,
da uns adäquate Abhängigkeitssätze über $[0,\infty)$ nicht zur Verfügung stehen.
Aus diesem Grund ist der nachstehende ``Beweis'' bezüglich diesem Argument unvollständig. \hfill $\square$}
\end{bemerkung}

{\bf Beweis} Da der Steuerungsprozess $\big(x_*(\cdot),u_*(\cdot)\big)$ der Menge
$\mathscr{A}^{\,\mathcal{U}}_{\rm adm} \cap \mathscr{A}^{\,\mathcal{U}}_{\rm Lip} \cap \mathscr{A}^{\,\mathcal{U}}_{\lim}$ angehört,
ist er beschränkt.
Ferner sind $t \to e^{-\varrho t} f_x\big(t,x_*(t),u_*(t)\big)$ und
$t \to \varphi_x\big(t,x_*(t),u_*(t)\big)$ über $[0,\infty)$ integrierbar, beschränkt und stückweise stetig.
Damit sind die Voraussetzungen von Lemma \ref{LemmaDGL4} und von Lemma \ref{LemmaDGL6} erfüllt
und es existiert eine eindeutige stetige L"osung $p(\cdot)$ der adjungierten Gleichung (\ref{PMPUHA1})
zur Transversalitätsbedingung  (\ref{PMPUHA2}).
Wegen der stückweisen Stetigkeit der einfließenden Abbildungen ist $\dot{p}(\cdot)$ ebenfalls stückweise stetig
und die Adjungierte $p(\cdot)$ gehört dem Raum $PC_1([0,\infty),\R^n)$ an.  \\[2mm]
Ferner besitzt nach Lemma \ref{LemmaDGL4} die Integralgleichung
\begin{equation} \label{BeweisPMPUH1}
y(t)=y(\tau) + \int_{\tau}^{t} \varphi_x\big(s,x_*(s),u_*(s)\big)\,y(s) \, ds, \quad t \in [0,\infty),
\end{equation}
für jedes $\tau \in [0,\infty)$ und jedes $y(\tau) \in \R^n$ eine eindeutige Lösung,
da $t \to \varphi_x\big(t,x_*(t),u_*(t)\big)$ über $[0,\infty)$ integrierbar ist. \\[2mm]
Wir verwenden wieder die Nadelvariation\index{Nadelvariation, einfache} des Beweises in Abschnitt \ref{AbschnittPMPBeweiseinfach}:
$$u(t;v,\tau,\lambda) = u_{\lambda}(t) = 
  \left\{ \begin{array}{ll}
          u_*(t) & \mbox{ f"ur } t \not\in [\tau-\lambda,\tau), \\
          v      & \mbox{ f"ur } t     \in [\tau-\lambda,\tau) 
          \end{array} \right.$$
und es bezeichne $x_\lambda(\cdot)$, $x_\lambda(t)=x(t;v,\tau,\lambda)$, die eindeutige L"osung der Gleichung
$$\dot{x}(t) = \varphi\big(t,x(t),u_\lambda(t)\big), \qquad x(t_0)=x_0.$$
Wir wählen $T \in (\tau,\infty)$.
Mit den gleichen Argumenten wie in Abschnitt \ref{AbschnittPMPBeweiseinfach} folgt,
dass für $t \in [\tau,T]$ der Grenzwert
$$y(t;T)=\lim_{\lambda \to 0^+}\frac{x_{\lambda}(t) - x_*(t)}{\lambda}$$
existiert und die Funktion $y(\cdot;T)$ der Integralgleichung
$$y(t;T) = y(\tau) + \int_{\tau}^{t} \varphi_x\big(s,x_*(s),u_*(s)\big)\,y(s;T) \, ds$$
zur Anfangsbedingung
$$y(\tau) = \varphi\big(\tau,x_*(\tau),v\big) - \varphi\big(\tau,x_*(\tau),u_*(\tau)\big)$$
gen"ugt.
Da die Integralgleichung (\ref{BeweisPMPUH1}) eine eindeutige Lösung $y(\cdot)$ über $[0,\infty)$ besitzt,
gilt $y(t;T)= y(t)$ für alle $[\tau,T]$, und dies für beliebiges $T > \tau$.
Ferner ist die Beziehung
$$\langle p(\tau) , y(\tau) \rangle-\langle p(T) , y(T) \rangle = - \int_{\tau}^T e^{-\varrho t} \big\langle f_x\big(t,x_*(t),u_*(t)\big) , y(t) \big\rangle \, dt$$
erf"ullt.
Da $y(\cdot)$ nach Lemma \ref{LemmaDGL4} einen Grenzwert im Unendlichen besitzt, ergibt sich
\begin{equation} \label{BeweisPMPUH2}
\langle p(\tau) , y(\tau) \rangle = - \int_{\tau}^\infty e^{-\varrho t} \big\langle f_x\big(t,x_*(t),u_*(t)\big) , y(t) \big\rangle \, dt.
\end{equation}
Es sei $x_\lambda(\cdot)$ die Variation von $x_*(\cdot)$ zur Steuerung $u_\lambda(\cdot)$,
die auf $[\tau,\infty)$ gleichmäßig konvergent gegen $x_*(\cdot)$ sei.
Es zeigt sich abschließend wie in Abschnitt \ref{AbschnittPMPBeweiseinfach}
\begin{eqnarray*}
0 &\leq& \lim_{\lambda \to 0^+} \frac{J\big(x_\lambda(\cdot),u_\lambda(\cdot)\big)- J\big(x_*(\cdot),u_*(\cdot)\big)}{\lambda} \\
  &=   & e^{-\varrho \tau} \big[f\big(\tau,x_*(\tau),v\big) - f\big(\tau,x_*(\tau),u_*(\tau)\big)\big]
         + \int_{\tau}^\infty e^{-\varrho t} \big\langle f_x\big(t,x_*(t),u_*(t)\big),y(t) \big\rangle \, dt,
\end{eqnarray*}
d.\,h. es ergibt sich zusammen mit (\ref{BeweisPMPUH2}) die Ungleichung
\begin{eqnarray*}
\lefteqn{\big\langle p(\tau) \,,\, \varphi\big(\tau,x_*(\tau),v\big) - e^{-\varrho \tau} f\big(\tau,x_*(\tau),v\big) } \\
&\leq& \big\langle p(\tau) \,,\, \varphi\big(\tau,x_*(\tau),u_*(\tau)\big) - e^{-\varrho \tau} f\big(\tau,x_*(\tau),u_*(\tau)\big).
\end{eqnarray*}
Da $\tau$ als ein Stetigkeitspunkt von $u_*(\cdot)$ und $v \in U$ willkürlich ausgewählt wurden,
folgt die Gültigkeit der Maximumbedingung (\ref{PMPUHA3}). \hfill $\blacksquare$ 
       \newpage
       \lhead[\thepage \hspace*{1mm} Fischerei-Differentialspiel]{ }   
       \subsection{Ein Fischerei-Differentialspiel}
\label{BeispielDockner} {\rm Wir betrachten nach Dockner et\,al. \cite{Dockner} das Differentialspiel
\begin{equation*}
\left.\begin{array}{l}
\hspace*{-5mm}
\displaystyle \tilde{J}_i\big(x(\cdot),u_1(\cdot),u_2(\cdot)\big) =\int_0^\infty e^{-\varrho t}\big(p x(t)-c_i\big)u_i(t) \, dt \to \sup,  \\[3mm]
\hspace*{-5mm}   
\displaystyle \dot{x}(t)=e^{-\delta t} \big[x(t) \big(\alpha-r\ln x(t) \big) -u_1(t)x(t)-u_2(t)x(t)\big], \quad x(0)=x_0>0, \\[2mm]
\hspace*{-5mm} 
u_i > 0, \quad \alpha,c_i,p,r,\varrho, \delta >0, \quad \alpha> \frac{1}{c_1+c_2}, \quad i=1,2.
\end{array} \right\}
\end{equation*}
Im Vergleich zu \cite{Dockner} haben wir den Faktor $e^{-\delta t}$ in der Dynamik hinzugefügt,
um die Stabilität des dynamischen Systems gemäß der Bedingungen (\ref{PMPBedingung}) zu sichern.
In dieser Aufgabe sei der Preis $p$ nicht konstant, sondern umgekehrt proportional zum Angebot $u_1x+u_2x$:
$$p= \frac{1}{u_1x+u_2x}, \quad px=\frac{1}{u_1+u_2}.$$
Dieser Ansatz spiegelt die "Okonomie einer ``Eskimo''-Gesellschaft wider,
in welcher der Fischbestand die wichtigste Nahrungsgrundlage darstellt und kein echtes Ersatzprodukt existiert.
Unter diesen Umst"anden f"uhrt eine prozentuale Preissteigerung zu einem Umsatzr"uckgang in gleicher Relation. \\[2mm]
Wenden wir die Transformation $z=\ln x$ an, so gilt
$$\dot{z}(t)=\frac{d}{dt} \ln\big(x(t)\big)=\frac{\dot{x}(t)}{x(t)}
            = e^{-\delta t} \big[\alpha-r\ln x(t)  -u_1(t)-u_2(t)\big] = e^{-\delta t} \big[\alpha-rz(t)  -u_1(t)-u_2(t)\big]$$
und wir erhalten das Spielproblem
\begin{equation*}
\left.\begin{array}{l}
\hspace*{-5mm}
\displaystyle J_i\big(z(\cdot),u_1(\cdot),u_2(\cdot)\big) =\int_0^\infty e^{-\varrho t}\bigg(\frac{1}{u_1(t)+u_2(t)}-c_i\bigg)u_i(t) \, dt \to \sup,  \\[3mm]
\hspace*{-5mm}   
\displaystyle \dot{z}(t)=e^{-\delta t} \big[-r z(t)+\alpha-u_1(t)-u_2(t)\big], \quad z(0)=\ln x_0>0, \\[2mm]
\hspace*{-5mm} 
u_i>0, \quad \alpha,c_i,p,r,\varrho >0, \quad \alpha> \frac{1}{c_1+c_2}, \quad i=1,2.
\end{array} \right\}
\end{equation*}
Ein zul"assiger Steuerungsprozess $\big(z^*(\cdot),u^*_1(\cdot),u^*_2(\cdot)\big)$ ist ein Nash-Gleichgewicht\index{Nash-Gleichgewicht} des Spiels,
falls f"ur alle anderen zul"assigen Steuerungsprozesse $\big(z(\cdot),u_1(\cdot),u^*_2(\cdot)\big)$,
$\big(z(\cdot),u^*_1(\cdot),u_2(\cdot)\big)$ die Ungleichungen
$$\left\{\begin{array}{l}
  J_1\big(z^*(\cdot),u^*_1(\cdot),u^*_2(\cdot)\big) \geq J_1\big(z(\cdot),u_1(\cdot),u^*_2(\cdot)\big) \\[2mm]
  J_2\big(z^*(\cdot),u^*_1(\cdot),u^*_2(\cdot)\big) \geq J_2\big(z(\cdot),u^*_1(\cdot),u_2(\cdot)\big)
  \end{array}\right.$$
gelten.
Halten wir die optimale Strategie des Gegenspielers fest,
dann ergeben sich f"ur $i,j \in \{1,2\}$, $i \not= j$, folgende miteinander gekoppelte Steuerungsprobleme:
\begin{equation*}
\left.\begin{array}{l}
\hspace*{-5mm}
\displaystyle J_i\big(z(\cdot),u_i(\cdot),u^*_j(\cdot)\big) =
      \int_0^\infty e^{-\varrho t}\bigg(\frac{1}{u_i(t)+u^*_j(t)}-c_i\bigg)u_i(t) \, dt \to \sup,  \\[3mm]
\hspace*{-5mm}   
\displaystyle \dot{z}(t)=e^{-\delta t} \big[ -r z(t)+\alpha-u_i(t)-u^*_j(t)\big], \quad z(0)=\ln x_0>0, \quad u_i>0.
\end{array} \right\}
\end{equation*}
Die zul"assigen Steuerungen geh"oren dem Raum $PC([0,\infty),U)$ an.
Deshalb sind im Ansatz des Nash-Gleichgewichtes die Abbildungen
$$f_i(t,z,u_i)=e^{-\varrho t}\bigg(\frac{1}{u_i+u^*_j(t)}-c_i\bigg)u_i, \quad \varphi_i(t,z,u_i)=e^{-\delta t}[-r z+\alpha-u_i-u^*_j(t)]$$
nicht stetig bez"uglich der Variable $t$.
An dieser Stelle verweisen auf die Bemerkungen in der Analyse des Kapitalismusspiels in Abschnitt \ref{AbschnittKapitalismusspiel},
dass unter diesen Rahmenbedingungen Theorem \ref{SatzPMPUHA} seine G"ultigkeit beh"alt. \\[2mm]
Die Pontrjagin-Funktionen lauten
\begin{eqnarray*}
H^{\mathcal{U}}_1(t,z,u_1,p_1) &=& p_1 e^{-\delta t} \big[-r z+\alpha-u_1-u^*_2(t)\big]+ e^{-\varrho t}\bigg(\frac{1}{u_1+u^*_2(t)}-c_1\bigg)u_1, \\
H^{\mathcal{U}}_2(t,z,u_2,p_2) &=& p_2 e^{-\delta t} \big[-r z+\alpha-u^*_1(t)-u_2\big]+ e^{-\varrho t}\bigg(\frac{1}{u^*_1(t)+u_2}-c_2\bigg)u_2.
\end{eqnarray*}  
Mit Lemma \ref{LemmaDGL6} ergibt sich die eindeutige Lösung der adjungierten Gleichung (\ref{PMPUHA1}):
$$\dot{p}_i(t)=r e^{-\delta t} p_i(t), \quad \lim_{t \to \infty} p_i(t)=0,\qquad\Rightarrow\qquad p_i(t) \equiv 0.$$
Die Maximumbedingung (\ref{PMPUHA2}) liefert das Gleichungssystem
$$\frac{\partial}{\partial u_i} H^{\mathcal{U}}_i\big(t,z_*(t),u_i^*(t),p_i(t)\big)=0 \quad\Rightarrow\quad
  \frac{1}{u_1^*(t)+u_2^*(t)}- \frac{u_i^*(t)}{\big(u_1^*(t)+u_2^*(t)\big)^2}=c_i, \; i=1,2.$$
Aus der Summe beider Gleichungen ergibt sich
$$\frac{2}{u_1^*(t)+u_2^*(t)}- \frac{u_1^*(t)+u_2^*(t)}{\big(u_1^*(t)+u_2^*(t)\big)^2}=\frac{1}{u_1^*(t)+u_2^*(t)}=c_1+c_2$$
und wir erhalten die optimalen Steuerungen
$$u_1^*(t)\equiv \frac{c_2}{(c_1+c_2)^2}, \qquad u_2^*(t)\equiv \frac{c_1}{(c_1+c_2)^2}.$$
Die optimale Trajektorie besitzt die Gestalt
$$z_*(t)=\bigg(z_0-\frac{c_0}{r}\bigg)\exp\bigg\{\frac{r}{\delta}( e^{-\delta t}-1)\bigg\} + \frac{c_0}{r}, \qquad c_0=  \alpha-\frac{1}{c_1+c_2}.$$
Für diese gilt im Unendlichen
$$z_\infty=\lim_{t \to \infty}z_*(t)=\bigg(z_0-\frac{c_0}{r}\bigg)e^{-\frac{r}{\delta}} + \frac{c_0}{r}
                           = z_0 e^{-\frac{r}{\delta}}+\bigg(\alpha-\frac{1}{c_1+c_2}\bigg)(1-e^{-\frac{r}{\delta}})>0.$$
Die Funktion $z_*(\cdot)$ ist streng monoton und nimmt nur Werte zwischen $z_0$ und $z_\infty$ an.
Somit ist $x_*(t)=\exp\big(z_*(t)\big)$ "uber $[0,\infty)$ wohldefiniert, besitzt eine untere positive Schranke und
$\big(x_*(\cdot),u^*_1(\cdot),u^*_2(\cdot)\big)$ liefert einen Kandidaten f"ur das urspr"ungliche Differentialspiel. \hfill $\square$
       \newpage
       \lhead[\thepage \hspace*{1mm} Arrow-Bedingungen]{ }
       \subsection{Hinreichende Bedingungen nach Arrow}
Wir leiten nun für die Aufgabe (\ref{UHA1})--(\ref{UHA3})
die hinreichenden Bedingungen nach Arrow ab. \index{hinreichende Bedingungen nach Arrow!unendlich@-- unendlicher Zeithorizont}
Das Vorgehen erschließt sich unmittelbar aus Abschnitt \ref{AbschnittHBPMP}. 
Für die Standardaufgabe (\ref{PMPeinfach1})--(\ref{PMPeinfach3}) liefern die hinreichenden Bedingungen nach Arrow
Gewissheit über die Optimalität eines Kandidaten.
Für die Anwendung des Maximumprinzips \ref{SatzHBPMPUH} müssen die Einschr"ankungen (\ref{PMPBedingung}) erfüllt sein.
In der Herleitung der Arrow-Bedingungen können wir auf diese Einschränkungen verzichten,
da sie wesentlich auf der Konkavit"at der Hamilton-Funktion $\mathscr{H}^{\,\mathcal{U}}$ basieren.
Dadurch erreichen wir eine deutlich umfassendere Anwendbarkeit der Bedingungen des (unvollständig bewiesenen) Maximumprinzips. \\[2mm]
Für $t \in [0,\infty)$ sei wieder $V_\gamma(t)=\{ x \in \R^n \,|\, \|x-x_*(t)\| < \gamma\}$.
Au"serdem bezeichnet $\mathscr{H}^{\mathcal{U}}$ die Hamilton-Funktion
$$\mathscr{H}^{\mathcal{U}}(t,x,p) = \sup_{u \in U} H^{\mathcal{U}}(t,x,u,p).$$

\begin{theorem} \label{SatzHBPMPUH}
In der Aufgabe (\ref{UHA1})--(\ref{UHA3}) sei
$\big(x_*(\cdot),u_*(\cdot)\big) \in \mathscr{A}^{\mathcal{U}}_{\rm adm} \cap \mathscr{A}^{\mathcal{U}}_{\rm Lip}$
und es sei $p(\cdot) \in PC_1([0,\infty),\R^n)$. Ferner gelte:
\begin{enumerate}
\item[(a)] Das Tripel $\big(x_*(\cdot),u_*(\cdot),p(\cdot)\big)$
           erf"ullt (\ref{PMPUHA1})--(\ref{PMPUHA3}) in Theorem \ref{SatzPMPUHA}.        
\item[(b)] F"ur jedes $t \in [0,\infty)$ ist die Funktion $\mathscr{H}^{\mathcal{U}}\big(t,x,p(t)\big)$ konkav in $x$ auf $V_\gamma(t)$.
\end{enumerate}
Dann ist $\big(x_*(\cdot),u_*(\cdot)\big)$ ein starkes lokales Minimum der Aufgabe (\ref{UHA1})--(\ref{UHA3}).
\end{theorem}

{\bf Beweis} Wir wählen $T \in (0,\infty)$.
Auf die gleiche Weise wie im Beweis von Theorem \ref{SatzHBPMP} in Abschnitt \ref{AbschnittHBPMP}
ergibt sich die Beziehung
\begin{eqnarray*}
\lefteqn{\int_0^T e^{-\varrho t}\big[f\big(t,x(t),u(t)\big)-f\big(t,x_*(t),u_*(t)\big)\big] \, dt} \\
&\geq & \int_0^T\big[\mathscr{H}^{\mathcal{U}}\big(t,x_*(t),p(t)\big)-\mathscr{H}^{\mathcal{U}}\big(t,x(t),p(t)\big)\big] \, dt 
        + \int_0^T \langle p(t), \dot{x}(t)-\dot{x}_*(t) \rangle \, dt \\
&\geq& \langle p(T),x(T)-x_*(T)\rangle-\langle p(0),x(0)-x_*(0)\rangle.
\end{eqnarray*}
Im Anfangszeitpunkt $t=0$ gelten $x(0)=x_*(0)=x_0$.
Für $T \to \infty$ sind die Trajektorien $x(\cdot),x_*(\cdot) \in PC_1([0,\infty),\R^n)$ beschränkt und es gilt $p(T) \to 0$
nach der Transversalitätsbedingung (\ref{PMPUHA2}).
Daher folgt die Beziehung
\begin{eqnarray*}
    J\big(x(\cdot),u(\cdot)\big)-J\big(x_*(\cdot),u_*(\cdot)\big)
&=& \lim_{T \to \infty} \int_0^T e^{-\varrho t}\big[f\big(t,x(t),u(t)\big)-f\big(t,x_*(t),u_*(t)\big)\big] \, dt \\
&\geq& \lim_{T \to \infty} \langle p(T),x(T)-x_*(T)\rangle-\langle p(0),x(0)-x_*(0)\rangle=0
\end{eqnarray*}
f"ur alle zul"assigen $\big(x(\cdot),u(\cdot)\big)$ mit $\|x(\cdot)-x_*(\cdot)\|_\infty < \gamma$. \hfill $\blacksquare$

\begin{beispiel} 
{\rm In der Aufgabe eines linear-quadratischen Reglers, 
\begin{equation} \label{BeispielRegler1}
\left.\begin{array}{l}
\displaystyle J\big(x(\cdot),u(\cdot)\big) = \frac{1}{2}\int_0^\infty e^{-t} \cdot \Big(\big( x(t)-1\big)^2+u^2(t)\big)\Big) \, dt \to \inf, \\[2mm]
\dot{x}(t) = x(t)+u(t), \qquad x(0)=0, \qquad u(t) \in \R,
\end{array}\right\}
\end{equation}
ist die restrikitive Bedingung (\ref{PMPBedingung}) nicht erf"ullt und es ist das
Pontrjaginsche Maximumprinzip \ref{SatzPMPUHA} nicht anwendbar.
Andererseits lautet die Pontrjagin-Funktion
$$H^{\mathcal{U}}(t,x,u,p) = p (x+u) - \frac{1}{2}e^{-t} \big((x-1)^2+u^2\big),$$
welche für jedes $u \in \R$ und $p \in \R$ konvex in $x$ ist.
Damit ist die Hamilton-Funktion
$$\mathscr{H}^{\mathcal{U}}(t,x,p) = \sup_{u \in U} \big\{p (x+u) - \frac{1}{2}e^{-t} \big((x-1)^2+u^2\big)\big\}$$
konvex, denn das Supremum wird durch $u_{\max}=pe^t$ angenommen.
Um die Lösung der Aufgabe (\ref{BeispielRegler1}) zu bestimmen,
müssen wir das Tripel $\big(x_*(\cdot),u_*(\cdot),p(\cdot)\big)$ ermitteln,
das die Bedingungen (\ref{PMPUHA1})--(\ref{PMPUHA3}) in Theorem \ref{SatzPMPUHA} erfüllt. \\[2mm]
Da der Steuerbereich offen ist und die Pontrjagin-Funktion streng konvex in $u$ ist,
darf man in der Maximumbedingung zur Ableitung nach $u$ übergehen.
Damit liefern die Dynamik, die adjungierte Gleichung und die Maximumbedingung das Gleichungssystem
$$\dot{x}(t) = x(t)+u(t), \quad \dot{p}(t)=-p(t)+e^{-t}\big(x(t)-1\big), \quad p(t)=e^{-t}u(t).$$
Das Ersetzen von $u(t)$ durch $e^tp(t)$ und Ableiten der Dynamik führt auf
$$\ddot{x}(t)=\dot{x}(t)+e^{t}\big(p(t)+\dot{p}(t)\big)=\dot{x}(t) +x(t)-1.$$
Die charakteristische Gleichung der Differentialgleichung für $x(\cdot)$ besitzt die Eigenwerte
$$\lambda_1=\frac{1 - \sqrt{5}}{2}, \qquad  \lambda_2=\frac{1 + \sqrt{5}}{2}.$$
Damit ergeben sich die Funktionen
$$x(t) = c_1e^{\lambda_1 t} + c_2e^{\lambda_2 t} + 1, \quad
  p(t) = (\lambda_1-1) c_1e^{(\lambda_1-1) t} + (\lambda_2-1) c_2e^{(\lambda_2-1) t}-e^{-t}.$$
Aufgrund der Transversalitätsbedingung $p(t) \to 0$ für $t \to \infty$ gilt $c_2=0$.
Damit führt die Anfangsbedingung $x(0)=0$ zu $c_1=-1$.
Wir erhalten das zulässige Tripel
$$x_*(t)=1-e^{\lambda_1 t}, \quad u_*(t)=(1-\lambda_1) e^{\lambda_1 t}-1 , \quad p(t)=(1-\lambda_1) e^{(\lambda_1-1) t}-e^{-t},$$
welches die Bedingungen (\ref{PMPUHA1})--(\ref{PMPUHA3}) in Theorem \ref{SatzPMPUHA} erfüllt.
Damit ist nach Theorem \ref{SatzHBPMPUH} der Steuerungsprozess
$\big(x_*(\cdot),u_*(\cdot)\big) \in PC_1([0,\infty),\R) \times PC([0,\infty),\R)$
ein starkes lokales Minimum der Aufgabe (\ref{BeispielRegler1}). \hfill $\square$}
\end{beispiel}

\begin{beispiel} 
{\rm Wir modifizieren das vorhergehende Beispiel und betrachten wir den linear-quadratischen Regler 
\begin{equation} \label{BeispielRegler2}
\left.\begin{array}{l}
\displaystyle J\big(x(\cdot),u(\cdot)\big) = \frac{1}{2}\int_0^\infty e^{-t} \cdot \Big(\big( x(t)-\cos t\big)^2+u^2(t)\big)\Big) \, dt \to \inf, \\[2mm]
\dot{x}(t) = x(t)+u(t), \qquad x(0)=0, \qquad u(t) \in \R.
\end{array}\right\}
\end{equation}
Die Pontrjagin-Funktion der Aufgabe \ref{BeispielRegler2} lautet
$$H^{\mathcal{U}}(t,x,u,p) = p (x+u) - \frac{1}{2}e^{-t} \big((x-\cos t)^2+u^2\big).$$
Damit ist wie im vorgehenden Beispiel die Hamilton-Funktion
$$\mathscr{H}^{\mathcal{U}}(t,x,p) = \sup_{u \in U} \big\{p (x+u) - \frac{1}{2}e^{-t} \big((x-\cos t)^2+u^2\big)\big\}$$
konvex und das Supremum wird durch $u_{\max}=pe^t$ angenommen.
Die Dynamik, die adjungierte Gleichung und die Maximumbedingung führen auf das Gleichungssystem
$$\dot{x}(t) = x(t)+u(t), \quad \dot{p}(t)=-p(t)+e^{-t}\big(x(t)-\cos t\big), \quad p(t)=e^{-t}u(t)$$
und liefern die Differentialgleichung
$$\ddot{x}(t)=\dot{x}(t)+e^{t}\big(p(t)+\dot{p}(t)\big)=\dot{x}(t) +x(t)- \cos t,$$
welche die Eigenwerte
$$\lambda_1=\frac{1 - \sqrt{5}}{2}, \qquad  \lambda_2=\frac{1 + \sqrt{5}}{2}$$
besitzt.
Damit ergeben sich die Funktionen
\begin{eqnarray*}
x(t) &=& c_1e^{\lambda_1 t} + c_2e^{\lambda_2 t} + \frac{1}{5}\sin t + \frac{2}{5}\cos t, \\
p(t) &=& (\lambda_1-1) c_1e^{(\lambda_1-1) t} + (\lambda_2-1) c_2e^{(\lambda_2-1) t}- \bigg( \frac{3}{5}\sin t + \frac{1}{5}\cos t\bigg) e^{-t}.
\end{eqnarray*}
Mit Hilfe der Transversalitätsbedingung $p(t) \to 0$ für $t \to \infty$ und der Anfangsbedingung $x(0)=0$
erhalten wir den Steuerungsprozess
$$x_*(t)=-\frac{2}{5} e^{\lambda_1 t} + \frac{1}{5}\sin t + \frac{2}{5}\cos t, \quad
  u_*(t)=\frac{2}{5}(1-\lambda_1) e^{\lambda_1 t}- \frac{3}{5}\sin t - \frac{1}{5}\cos t$$
und die Adjungierte
$$p(t)=\frac{2}{5}(1-\lambda_1) e^{(\lambda_1-1) t}- \bigg( \frac{3}{5}\sin t + \frac{1}{5}\cos t\bigg) e^{-t},$$
welche die Bedingungen (\ref{PMPUHA1})--(\ref{PMPUHA3}) in Theorem \ref{SatzPMPUHA} erfüllen. \\[2mm]
Der Steuerungsprozess $\big(x_*(\cdot),u_*(\cdot)\big)$ ist über $[0,\infty)$ stetig und beschränkt und gehört damit
der Menge $PC_1([0,\infty),\R) \times PC([0,\infty),\R)$ an -- obwohl $x_*(\cdot)$ und $u_*(\cdot)$ im Unendlichen divergieren!
Dennoch ist nach Theorem \ref{SatzHBPMPUH} der Steuerungsprozess
$\big(x_*(\cdot),u_*(\cdot)\big)$
ein starkes lokales Minimum der Aufgabe (\ref{BeispielRegler2}). \hfill $\square$}
\end{beispiel}
\cleardoublepage

\lhead[ ]{}
\rhead[]{Zeitverzögerte Systeme \hspace*{1mm} \thepage}
\section{Zeitverz\"ogerte dynamische Systeme}
       Dynamische Systeme mit Zeitverzögerungen sind wesentliche Elemente in der Modellierung und Auswertung 
lebensnaher Phänomene.
Anwendungsfelder sind in der Biologie oder der Biomedizin zu finden. 
So zählt zu den gängigen Beispielen einer retardierten Arznei der Typus,
bei dem der Wirkstoff zeitlich verzögert freigesetzt wird. \\[2mm]
Das verzögerte System entsteht dadurch,
dass die Zustandsänderung $\dot{x}(t)$ nicht ausschließlich aus dem vorliegenden Zustand $x(t)$ und
den simultanen Ma"snahmen $u(t)$ resultiert,
sondern gleichzeitig mit einer zeitlichen Verzögerungen $\Delta >0$ auf den vorherigen Status $x(t-\Delta)$ und auf die getroffenen
Maßnahmen $u(t-\Delta)$ reagiert.
Dadurch erhält das dynamische System die Form
$$\dot{x}(t) = \varphi\big(t,x(t),x(t-\Delta),u(t),u(t-\Delta)\big),$$
welche zentraler Gegenstand der Untersuchungen dieses Abschnitts ist. \\[2mm]
Die Auswertung von Steuerungsproblemen mit einer Zeitverzögerung ist schwierig und erfolgt in der Literatur meistens numerisch.
Als eine Anwendung auf eine reale Problemstellung verweisen wir die Untersuchung einer Chemoimmuntherapie,
einer Kombination von Chemo- und Immuntherapie,
bei Rihan\,et.\,al. \cite{Rihan} oder bei Göllmann \& Maurer in \cite{GoMa}. \\[2mm]
Bei einer Chemotherapie verwendet man verschiedene Medikamente (Zytostatika),
um Krebszellen abzutöten oder das Wachstum der Krebszellen zu verlangsamen.
Bei der Behandlung bösartiger Tumorerkrankungen nutzen die meisten Zytostatika die schnelle Teilungsfähigkeit der Tumorzellen,
da diese empfindlicher als gesunde Zellen auf Störungen der Zellteilung reagieren.
Auf gesunde Zellen mit ähnlich guter Teilungsfähigkeit üben sie allerdings eine ähnliche Wirkung aus,
wodurch sich Nebenwirkungen wie Haarausfall oder Durchfall einstellen können.
Immuntherapien beinhalten Behandlungen,
die das Immunsystem im Kampf gegen Krebs stimulieren oder stärken.
Effektorzellen wie B-Lymphozyten (kurz B-Zellen) gehören zu den Leukozyten (weiße Blutkörperchen),
die Plasmazellen  bilden, die wiederum Antikörper bilden.
       \lhead[\thepage \hspace*{1mm} Pontrjaginsches Maximumprinzip]{ }         
       \subsection{Die Aufgabenstellung und das Pontrjaginsche Maximumprinzip} \label{AbschnittPMPeinfachRet}
Die Aufgabe eines verzögerten Steuerungsproblems besitzt die Form
\begin{eqnarray}
&& \label{Ret1} \hspace*{-5mm} J\big(x(\cdot),u(\cdot)\big) = \int_{t_0}^{t_1} f\big(t,x(t),x(t-\Delta),u(t),u(t-\Delta)\big) \, dt +S\big(x(t_1)\big)\to \inf, \\
&& \label{Ret2} \hspace*{-5mm} \dot{x}(t) = \varphi\big(t,x(t),x(t-\Delta),u(t),u(t-\Delta)\big),\quad x(t)=x_0 \mbox{ für } t \in [t_0-\Delta,t_0], \\
&& \label{Ret3} \hspace*{-5mm} u(t) \in U \subseteq \R^m, \quad u(t)=u_0 \in U \mbox{ für } t \in [t_0-\Delta,t_0),\quad U\not= \emptyset.
\end{eqnarray} 
Im Vergleich zur Aufgabe (\ref{PMPeinfach1})--(\ref{PMPeinfach3}) ergänzen wir in den Abbildungen $f$ und $\varphi$
die Variablen $y$ und $v$ bez"uglich des verzögerten Zustandes bzw. der verzögerten Steuerung:

\begin{eqnarray*}
&& f=f(t,x,y,u,v):\R \times \R^n \times \R^n \times \R^m \times \R^m \to \R, \\
&& \varphi=\varphi(t,x,y,u,v):\R \times \R^n \times \R^n \times \R^m \times \R^m \to \R^n.
\end{eqnarray*}
Bei der Behandlung von Steuerungsproblemen mit einer Zeitverzögerung stellt sich die Frage,
wie das dynamische System (\ref{Ret2})
zu behandeln ist.
Für den Teilabschnitt $[t_0,t_0+\Delta]$ reduziert es sich auf die Gleichung
$$\dot{x}(t) = \varphi\big(t,x(t),x_0,u(t),u_0\big), \quad x(t_0)=x_0,$$
welche sich für eine stückweise stetige rechte Seite $\varphi$ mit den Werkzeugen im Anhang \ref{AnhangDGL} über Differentialgleichungen handhaben lässt.
Dieser Gedankengang kann anschließend schrittweise über den Teilabschnitten $[t_0+\Delta,t_0+2\Delta]$, $[t_0+2\Delta,t_0+3\Delta]$, ... wiederholt werden.
Daher kann die Dynamik mit Zeitverzögerung abschnittsweise wie eine Differentialgleichung mit stückweise stetigen rechten Seiten
aufgefasst werden. 
In ähnlicher Weise verfährt man mit der adjungierten Gleichung in der umgekehrten Zeitrichtung. \\[2mm]
Die Aufgabe (\ref{Ret1})--(\ref{Ret3}) betrachten wir bez"uglich der Paare
$$\big(x(\cdot),u(\cdot)\big) \in PC_1([t_0,t_1],\R^n) \times PC([t_0,t_1],U).$$
Mit $\mathscr{A}^{\,\mathcal{Z}}_{\rm Lip}$ bezeichnen wir die Menge aller Paare $\big(x(\cdot),u(\cdot)\big)$,
für die es ein $\gamma>0$ derart gibt,
dass die Abbildungen $f(t,x,y,u,v)$, $\varphi(t,x,y,u,v)$ auf der Menge aller Punkte
$(t,x,y,u,v) \in \R \times \R^n \times \R^n \times \R^m \times \R^m$ mit
$$t \in [t_0,t_1], \qquad \|x-x(t)\| < \gamma, \qquad \|y-x(t-\delta)\| < \gamma, \qquad u,v \in \R^m$$
stetig in der Gesamtheit aller Variablen und stetig differenzierbar bezüglich $x$ und $y$ sind. \\[2mm]
Das Paar $\big(x(\cdot),u(\cdot)\big) \in PC_1([t_0,t_1],\R^n) \times PC([t_0,t_1],U)$
ist zul"assig in der Aufgabe (\ref{Ret1})--(\ref{Ret3}),
falls $\big(x(\cdot),u(\cdot)\big)$ dem System (\ref{Ret2}) gen"ugt und die Steuerbeschränkungen (\ref{Ret3}) erfüllt.
Mit $\mathscr{A}^{\mathcal{Z}}_{\rm adm}$ bezeichnen wir die Menge der zul"assigen Steuerungsprozesse. \\[2mm]
Ein zul"assiger Steuerungsprozess $\big(x_*(\cdot),u_*(\cdot)\big)$ ist eine
starke lokale Minimalstelle\index{Minimum, starkes lokales!ZDelay@-- zeitverzögerte Systeme}
der Aufgabe (\ref{Ret1})--(\ref{Ret3}),
falls eine Zahl $\varepsilon > 0$ derart existiert, dass die Ungleichung 
$$J\big(x(\cdot),u(\cdot)\big) \geq J\big(x_*(\cdot),u_*(\cdot)\big)$$
f"ur alle $\big(x(\cdot),u(\cdot)\big) \in \mathscr{A}^{\mathcal{Z}}_{\rm adm}$
mit $\|x(\cdot)-x_*(\cdot)\|_\infty < \varepsilon$ gilt. \\[2mm]
Es bezeichnet $H^{\mathcal{Z}}: \R \times \R^n \times \R^n \times \R^m \times \R^m \times \R^n \to \R$ die Pontrjagin-Funktion
$$H^{\mathcal{Z}}(t,x,y,u,v,p) = \langle p,\varphi(t,x,y,u,v)\rangle - f(t,x,y,u,v).$$
Weiterhin führen wir die folgenden abkürzenden Bezeichnungen ein:
\begin{eqnarray*}
H^{\mathcal{Z}}[t] &=& H^{\mathcal{Z}}\big(t,x_*(t),x_*(t-\Delta), u_*(t),u_*(t-\Delta),p(t)\big), \\
H^{\mathcal{Z}}[t,u,v] &=& H^{\mathcal{Z}}\big(t,x_*(t),x_*(t-\Delta), u,v,p(t)\big).
\end{eqnarray*}

Damit formulieren wir das Pontrjaginsche Maximumprinzip der Aufgabe (\ref{Ret1})--(\ref{Ret3}): 
\begin{theorem}[Pontrjaginsches Maximumprinzip] \label{SatzRetPMP}
\index{Pontrjaginsches Maximumprinzip!ZDelay@-- zeitverzögerte Systeme} 
Es sei $\big(x_*(\cdot),u_*(\cdot)\big) \in \mathscr{A}^{\mathcal{Z}}_{\rm adm} \cap \mathscr{A}^{\mathcal{Z}}_{\rm Lip}$.
Ist $\big(x_*(\cdot),u_*(\cdot)\big)$ ein starkes lokales Minimum der Aufgabe (\ref{Ret1})--(\ref{Ret3}),
dann existiert eine Vektorfunktion $p(\cdot) \in PC_1([t_0,t_1],\R^n)$ derart, dass
\begin{enumerate}
\item[(a)] fast überall in $[t_0,t_1]$ die adjungierte Gleichung
           \index{adjungierte Gleichung!ZDelay@-- zeitverzögerte Systeme}
           \begin{equation}\label{SatzRetPMP1}
           \dot{p}(t) = - H^{\mathcal{Z}}_x[t] - \chi_{[t_0,t_1-\Delta]}(t) \cdot H^{\mathcal{Z}}_y[t+\Delta],
           \end{equation}
\item[(b)] in $t =t_1$ die Transversalitätsbedingung
           \index{Transversalitätsbedingungen!ZDelay@-- zeitverzögerte Systeme}
           \begin{equation}\label{SatzRetPMP2} 
           p(t_1)=-S'\big(x_*(t_1)\big)
           \end{equation}
\item[(c)] und in fast allen Punkten $t \in [t_0,t_1]$ die Maximumbedingung
           \index{Maximumbedingung!ZDelay@-- zeitverzögerte Systeme}
           \begin{eqnarray}
           && \hspace*{-10mm}
           H^{\mathcal{Z}}[t,u_*(t),u_*(t-\Delta)]+\chi_{[t_0,t_1-\Delta]}(t) \cdot H^{\mathcal{Z}}[t+\Delta,u_*(t+\Delta),u_*(t)] \nonumber \\
           \label{SatzRetPMP3}
           && \hspace*{-10mm}
              = \max_{u \in U} \Big\{H^{\mathcal{Z}}[t,u,u_*(t-\Delta)]+\chi_{[t_0,t_1-\Delta]}(t) \cdot H^{\mathcal{Z}}[t+\Delta,u_*(t+\Delta),u] \Big\}
           \end{eqnarray}
\end{enumerate}
erfüllt sind.
\end{theorem}

In Theorem \ref{SatzRetPMP} sind die Bedingungen kompakt formuliert.
Die adjungierte Gleichung (\ref{SatzRetPMP1}) besitzt ausführlich aufgeschrieben die Gestalt
\begin{eqnarray*}
\dot{p}(t) &=& - \varphi_x^T\big(t,x_*(t),x_*(t-\Delta), u_*(t),u_*(t-\Delta)\big)p(t) \\
           & & \hspace*{2cm} +f_x\big(t,x_*(t),x_*(t-\Delta), u_*(t),u_*(t-\Delta)\big) \\
           & & + \chi_{[t_0,t_1-\Delta]}(t) \cdot
                 \Big[ - \varphi_y^T\big(t+\Delta,x_*(t+\Delta),x_*(t), u_*(t+\Delta),u_*(t)\big)p(t+\Delta) \\
           & &\hspace*{4.5cm}  +f_y\big(t+\Delta,x_*(t+\Delta),x_*(t), u_*(t+\Delta),u_*(t)\big)\Big].
\end{eqnarray*}
Ähnlich wie das System (\ref{Ret2}) besitzt die adjungierte Gleichung (\ref{SatzRetPMP1}) über $[t_1-\Delta,t_1]$ die
``gewöhnliche'' Form $\dot{p}(t) = - H^{\mathcal{Z}}_x[t]$ zur Transversalitätsbedingung (\ref{SatzRetPMP2}) 
und lässt sich dann abschnittsweise über $[t_1-2\Delta,t_1-\Delta]$, $[t_1-3\Delta,t_1-2\Delta]$, ... behandeln.
Ferner hat die Maximumbedingung (\ref{SatzRetPMP3}) die Gestalt
\begin{eqnarray*}
& & \big\langle \varphi\big(t,x_*(t),x_*(t-\Delta), u_*(t),u_*(t-\Delta)\big),p(t) \big\rangle \\
& & \hspace*{2cm} f\big(t,x_*(t),x_*(t-\Delta), u_*(t),u_*(t-\Delta)\big) \\
& & + \chi_{[t_0,t_1-\Delta]}(t) \cdot
      \Big[ \big\langle \varphi\big(t+\Delta,x_*(t+\Delta),x_*(t), u_*(t+\Delta),u_*(t)\big) , p(t+\Delta)\big\rangle \\
& &\hspace*{4.5cm}  -f\big(t+\Delta,x_*(t+\Delta),x_*(t), u_*(t+\Delta),u_*(t)\big)\Big] \\
&=& \max_{u \in U}\Big\{ \big\langle \varphi\big(t,x_*(t),x_*(t-\Delta), u,u_*(t-\Delta)\big),p(t) \big\rangle \\
& & \hspace*{3cm} f\big(t,x_*(t),x_*(t-\Delta), u,u_*(t-\Delta)\big) \\
& & \hspace*{1cm}+ \chi_{[t_0,t_1-\Delta]}(t) \cdot
      \Big[ \big\langle \varphi\big(t+\Delta,x_*(t+\Delta),x_*(t), u_*(t+\Delta),u\big) , p(t+\Delta)\big\rangle \\
& &\hspace*{5.5cm}  -f\big(t+\Delta,x_*(t+\Delta),x_*(t), u_*(t+\Delta),u\big)\Big]\Big\}.
\end{eqnarray*}

\begin{beispiel} \label{BeispielInvDelay}
{\rm In Anlehnung an Beispiel \ref{BeispielLinInv} betrachten wir eine Aufgabe,
in welcher sich die Investitionen mit der Verzögerung von einer Zeiteinheit $\Delta = 1$ auf den Kapitalbestand auswirken.
Es ergibt sich das Investitionsmodell
\begin{eqnarray}
&&\label{BeispielRet1} J\big(K(\cdot),u(\cdot)\big) = \int_0^T \big( 1-u(t)\big) \cdot K(t) \, dt \to \sup, \\
&&\label{BeispielRet2} \dot{K}(t) = u(t-1) \cdot K(t-1), \quad K(0)=K_0 >0,\\
&&\label{BeispielRet3} u(t) \in [0,1] \mbox{ für } t \in [0,T], \quad u(t)=0 \mbox{ für } t \in [-1,0), \quad T > 2 \mbox{ fest}.
\end{eqnarray}
In der Aufgabe (\ref{BeispielRet1})--(\ref{BeispielRet3}) verwenden wir für $K(t),K(t-1),u(t),u(t-1)$ die Variablenbezeichnungen $K,L,u,v$.
Damit lauten die Abbildungen der Aufgabe
$$f(t,K,L,u,v) = -( 1-u) \cdot K, \quad \varphi(t,K,L,u,v) = v \cdot L$$
und die Pontrjaginsche Funktion
$$H^{\mathcal{Z}}(t,K,L,u,v,p) = p vL +(1-u)K.$$
Die Maximumbedingung (\ref{SatzRetPMP3}) besitzt die Form
\begin{eqnarray*}
\lefteqn{\max_{u \in U} \Big\{H^{\mathcal{Z}}[t,u,u_*(t-1)]+\chi_{[0,T-1]}(t) \cdot H^{\mathcal{Z}}[t+1,u_*(t+1),u] \Big\}} \\
&=& \max_{u \in [0,1]}\Big\{p(t)u_*(t-1)K_*(t-1) +(1-u)K_*(t) \\
& & \hspace*{1.5cm} + \chi_{[0,T-1]}(t)\cdot\big[p(t+1)u K_*(t) +\big(1-u_*(t+1)\big)K_*(t+1)\Big\}
\end{eqnarray*}
und ist äquivalent zu folgenden Maximierungen über Teilabschnitte:
\begin{equation} \label{BeispielRet4}
\max_{u \in [0,1]} \left\{\begin{array}{ll}
u \big(p(t+1)-1\big) K_*(t), & t \in [0,T-1], \\[1mm] 
-u K_*(t),& t \in [T-1,T].
\end{array} \right.
\end{equation}
Die adjungierte Gleichung (\ref{SatzRetPMP1}) erhält die Gestalt
$$\dot{p}(t) = - H^{\mathcal{Z}}_K[t] - \chi_{[0,T-1]}(t) \cdot H^{\mathcal{Z}}_L[t+1] = -\big(1-u_*(t)\big) - \chi_{[0,T-1]}(t)\cdot p(t+1)u_*(t)$$
und lautet abschnittsweise
\begin{equation} \label{BeispielRet5}
\dot{p}(t)=\left\{\begin{array}{ll} 
\displaystyle -\big(1-u_*(t)\big) - p(t+1)u_*(t), & t \in (0,T-1), \\[1mm]
\displaystyle -\big(1-u_*(t)\big),& t \in (T-1,T),
\end{array} \right.
\end{equation}
zur Transversalitätsbedingung $p(T)=0$. \\[2mm]
Da $\dot{K}(t) \geq 0$ und $\dot{p}(t)<0$ in $(0,T)$ gelten, führen (\ref{BeispielRet4}) und (\ref{BeispielRet5}) zu:
\begin{enumerate}
\item[$\cdot$] Über $[T-1,T]$ lauten (\ref{BeispielRet4}) und (\ref{BeispielRet5})
               $$\max_{u \in [0,1]} -u K_*(t), \quad \dot{p}(t)=-\big(1-u_*(t)\big), \quad p(T)=0.$$
               Daraus ist unmittelbar $u_*(t)=0$ und $p(t)=T-t$ über $[T-1,T]$ einzusehen.
\item[$\cdot$] Damit folgt für (\ref{BeispielRet4}) über $[T-2,T-1)$:
               $\displaystyle \max_{u \in [0,1]} u \big(p(t+1)-1\big) K_*(t)$.
               Da $p(t)<1$ für $t \in (T-1,T]$ ausfällt, ist $u_*(t)=0$ für $t \in [T-2,T-1)$.
               In (\ref{BeispielRet5}) ergibt sich ferner $\dot{p}(t)=-1$ und es gilt weiterhin $p(t)=T-t$ für $t \in [T-2,T-1)$.
\item[$\cdot$] Da stets $\dot{p}<0$ ausfällt, ist $p(t+1)>2$ über $[0,T-2)$.
               Dies hat anhand der Maximumbedingung (\ref{BeispielRet4}) die Investitionsrate
               $u_*(t)=1$ für $t \in [0,T-2)$ zur Folge.
               Weiterhin ergibt sich $\dot{p}(t)=-p(t+1)$ in (\ref{BeispielRet5}).              
\end{enumerate}
In der Aufgabe (\ref{BeispielRet1})--(\ref{BeispielRet3}) mit verzögerten Investitionen konnten wir die Steuerung
$$u_*(t)=\left\{\begin{array}{rl} 0, & t \in [-1,0), \\ 1, & t \in [0,T-2), \\ 0, & t \in [T-2,T], \end{array} \right.$$
ermitteln.
Mit dieser ergeben sich für den Kaptialbestand und für die Adjungierte die dynamischen Entwicklungen
$$\dot{K}_*(t)=\left\{\begin{array}{rl} 0, & t \in (0,1), \\ K_*(t-1), & t \in (1,T-1), \\ 0, & t \in (T-1,T), \end{array} \right. \qquad
  \dot{p}(t)=\left\{\begin{array}{rl} -p(t+1), & t \in (0,T-2), \\ -1 , & t \in (T-2,T). \end{array} \right.$$
Zusammen mit $K(0)=K_0$ und $p(T)=0$ führt dies abschnittsweise auf die Darstellungen
\begin{eqnarray*}
K_*(t) &=& \left\{\begin{array}{rl} K_0, & t \in [0,1), \\ K_0+K_0(t-1), & t \in [1,2), \\ 
                                    2K_0 + K_0(t-2)+ \frac{1}{2}K_0\big((t-1)^2-1\big), & t \in [2,3), \\ 
                                    \vdots \\ K_*(T-1), & t \in [T-1,T], \end{array} \right. \\[2mm]
p(t) &=& \left\{\begin{array}{rl} \vdots \\ 2+\frac{1}{2}\big[\big(T-(t+1)\big)^2-1\big], & t \in [T-3,T-2), \\ T-t , & t \in [T-2,T]. \end{array} \right.
\end{eqnarray*}
Aufgrund der abschnittsweise rekursiven Bildung des optimalen Kapitalbestandes $K_*(\cdot)$ verzichten wir an dieser Stelle auf die Berechnung
des Optimalwertes. \hfill $\square$}
\end{beispiel}     
       \lhead[\thepage \hspace*{1mm} Beweis des Maximumprinzips]{ }  
       \subsection{Der Beweis des Maximumprinzips}
Im Zuge des Beweises gelangen wir zu Ausdrücken,
welche Verzögerungen in verschiedene Zeitrichtungen aufweisen.
Mit der Substitution $t = s+ \Delta$ gilt $f(t)g(t-\Delta)=f(s+\Delta)g(s)$ und es ergibt sich
\begin{equation}\label{Zeitshift}
\int_a^b \chi_{[c+\Delta,d]}(t) \cdot f(t)g(t-\Delta) \, dt = \int_a^b \chi_{[c,d-\Delta]}(t) \cdot f(t+\Delta)g(t) \, dt.
\end{equation}
Beachten wir die unsere Bemerkungen zur Lösung von Differentialgleichungen mit zeitlichen Verzögerungen,
so gibt es nach Lemma \ref{LemmaDGL3} und Lemma \ref{LemmaDGL5} eine eindeutige L"osung $p(\cdot) \in PC_1([t_0,t_1],\R^n)$
der Gleichung (\ref{SatzRetPMP1}) zur Randbedingung (\ref{SatzRetPMP2}). \\[2mm]
Für $\lambda < \Delta$ betrachten wir wieder die einfache Nadelvariation\index{Nadelvariation, einfache}
$$u(t;v,\tau,\lambda) = u_{\lambda}(t) = 
  \left\{ \begin{array}{ll}
          u_*(t) & \mbox{ f"ur } t \not\in [\tau-\lambda,\tau), \\
          v      & \mbox{ f"ur } t     \in [\tau-\lambda,\tau), 
          \end{array} \right.$$
und untersuchen f"ur $t \geq \tau$ den Grenzwert (vgl. Abschnitt \ref{AbschnittPMPBeweiseinfach})
\begin{equation} \label{BeweisRet1}
y(t)=\lim_{\lambda \to 0^+}\frac{x_{\lambda}(t) - x_*(t)}{\lambda}.
\end{equation}
Im weiteren Vorgehen müssen wir beachten,
dass in der Aufgabe (\ref{Ret1})--(\ref{Ret3}) die Nadelvariation $u_\lambda(\cdot)$ über den beiden Intervall
$[\tau-\lambda,\tau)$ und $[\tau+\Delta-\lambda,\tau+\Delta)$ in die Abbildungen $f$ und $\varphi$ einfließt.
Für den Quotienten auf der rechten Seite in (\ref{BeweisRet1}) ergeben sich für $t \geq \tau$ abschnittsweise
über dem Intervall $[\tau-\lambda,t_1]$ die folgenden Ausdrücke:
\begin{eqnarray*}
&& \frac{1}{\lambda} \int_{\tau - \lambda}^{\tau} \chi_{[t_0,t]}(s) \cdot
   \Big[ \varphi\big(s,x_\lambda(s),x_*(s-\Delta),v,u_*(s-\Delta)\big) \nonumber\\
&& \hspace*{4cm} - \varphi\big(s,x_*(s),x_*(s-\Delta),u_*(s),u_*(s-\Delta)\big) \Big] \, ds \nonumber\\
&& + \frac{1}{\lambda} \int_{\tau}^{\tau+\Delta - \lambda} \chi_{[t_0,t]}(s) \cdot
   \Big[ \varphi\big(s,x_\lambda(s),x_*(s-\Delta),u_*(s),u_*(s-\Delta)\big) \nonumber\\
&& \hspace*{4cm} - \varphi\big(s,x_*(s),x_*(s-\Delta),u_*(s),u_*(s-\Delta)\big) \Big] \, ds \\
&& + \frac{1}{\lambda} \int_{\tau+\Delta - \lambda}^{\tau+\Delta} \chi_{[t_0,t]}(s) \cdot
   \Big[ \varphi\big(s,x_\lambda(s),x_*(s-\Delta),u_*(s),v\big) \nonumber\\
&& \hspace*{4cm} - \varphi\big(s,x_*(s),x_*(s-\Delta),u_*(s),u_*(s-\Delta)\big) \Big] \, ds \\
&& + \frac{1}{\lambda} \int_{\tau+\Delta}^{t_1} \chi_{[t_0,t]}(s) \cdot
   \Big[ \varphi\big(s,x_\lambda(s),x_\lambda(s-\Delta),u_*(s),u_*(s-\Delta)\big) \nonumber\\
&& \hspace*{4cm} - \varphi\big(s,x_*(s),x_*(s-\Delta),u_*(s),u_*(s-\Delta)\big) \Big] \, ds. \nonumber
\end{eqnarray*}
Wir erweitern diese Summe um die Terme
$$\pm \frac{1}{\lambda} \int_{\tau+\Delta - \lambda}^{\tau+\Delta} \chi_{[t_0,t]}(s) \cdot
  \varphi\big(s,x_\lambda(s),x_*(s-\Delta),u_*(s),u_*(s-\Delta)\big)\, ds,$$
formen damit den zweiten und dritten Summanden um und bringen die gesamte Summe in die Gestalt
\begin{eqnarray*}
&& \frac{1}{\lambda} \int_{\tau - \lambda}^{\tau} \chi_{[t_0,t]}(s) \cdot
   \Big[ \varphi\big(s,x_\lambda(s),x_*(s-\Delta),v,u_*(s-\Delta)\big) \nonumber\\
&& \hspace*{4cm} - \varphi\big(s,x_*(s),x_*(s-\Delta),u_*(s),u_*(s-\Delta)\big) \Big] \, ds \nonumber\\
&& + \frac{1}{\lambda} \int_{\tau}^{\tau+\Delta} \chi_{[t_0,t]}(s) \cdot
   \Big[ \varphi\big(s,x_\lambda(s),x_*(s-\Delta),u_*(s),u_*(s-\Delta)\big) \nonumber\\
&& \hspace*{4cm} - \varphi\big(s,x_*(s),x_*(s-\Delta),u_*(s),u_*(s-\Delta)\big) \Big] \, ds \\
&& + \frac{1}{\lambda} \int_{\tau+\Delta - \lambda}^{\tau+\Delta} \chi_{[t_0,t]}(s) \cdot
   \Big[ \varphi\big(s,x_\lambda(s),x_*(s-\Delta),u_*(s),v\big) \nonumber\\
&& \hspace*{4cm} - \varphi\big(s,x_\lambda(s),x_*(s-\Delta),u_*(s),u_*(s-\Delta)\big) \Big] \, ds \nonumber\\
&& + \frac{1}{\lambda} \int_{\tau+\Delta}^{t_1} \chi_{[t_0,t]}(s) \cdot
   \Big[ \varphi\big(s,x_\lambda(s),x_\lambda(s-\Delta),u_*(s),u_*(s-\Delta)\big) \nonumber\\
&& \hspace*{4cm} - \varphi\big(s,x_*(s),x_*(s-\Delta),u_*(s),u_*(s-\Delta)\big) \Big] \, ds.
\end{eqnarray*}
Nun ergänzen wir in dieser Summe die Terme
$$\pm \frac{1}{\lambda} \int_{\tau+\Delta}^{t_1} \chi_{[t_0,t]}(s) \cdot
  \varphi\big(s,x_\lambda(s),x_*(s-\Delta),u_*(s),u_*(s-\Delta)\big)\, ds,$$
formen damit den zweiten und vierten Summanden um und bringen die Summe in diejenige Form,
für die wir den Grenzübergang $\lambda \to 0^+$ bilden wollen:
\begin{eqnarray}
&& \frac{1}{\lambda} \int_{\tau - \lambda}^{\tau}\chi_{[t_0,t]}(s) \cdot
   \Big[ \varphi\big(s,x_\lambda(s),x_*(s-\Delta),v,u_*(s-\Delta)\big) \nonumber\\
&& \hspace*{4cm} - \varphi\big(s,x_*(s),x_*(s-\Delta),u_*(s),u_*(s-\Delta)\big) \Big] \, ds  \nonumber\\
&& + \frac{1}{\lambda} \int_{\tau}^{t_1} \chi_{[t_0,t]}(s) \cdot
   \Big[ \varphi\big(s,x_\lambda(s),x_*(s-\Delta),u_*(s),u_*(s-\Delta)\big) \nonumber\\
&& \hspace*{4cm} - \varphi\big(s,x_*(s),x_*(s-\Delta),u_*(s),u_*(s-\Delta)\big) \Big] \, ds  \nonumber\\
&& + \frac{1}{\lambda} \int_{\tau+\Delta - \lambda}^{\tau+\Delta} \chi_{[t_0,t]}(s) \cdot
   \Big[ \varphi\big(s,x_\lambda(s),x_*(s-\Delta),u_*(s),v\big) \nonumber\\
&& \hspace*{4cm} - \varphi\big(s,x_\lambda(s),x_*(s-\Delta),u_*(s),u_*(s-\Delta)\big) \Big] \, ds  \nonumber\\
&& + \frac{1}{\lambda} \int_{\tau+\Delta}^{t_1} \chi_{[t_0,t]}(s) \cdot
   \Big[ \varphi\big(s,x_\lambda(s),x_\lambda(s-\Delta),u_*(s),u_*(s-\Delta)\big) \nonumber\\
&& \label{BeweisRet2} \hspace*{4cm} - \varphi\big(s,x_\lambda(s),x_*(s-\Delta),u_*(s),u_*(s-\Delta)\big) \Big] \, ds.
\end{eqnarray}
Um uns von den klobigen Ausdrücken der Art (\ref{BeweisRet2}) zu lösen,
verwenden wir im Weiteren die Kurzschreibweisen
\begin{eqnarray*}
\varphi[t] &=& \varphi\big(t,x_*(t),x_*(t-\Delta), u_*(t),u_*(t-\Delta)\big), \\
\varphi[t,u,v] &=& \varphi\big(t,x_*(t),x_*(t-\Delta), u,v\big), \\
f[t] &=& f\big(t,x_*(t),x_*(t-\Delta), u_*(t),u_*(t-\Delta)\big), \\
f[t,u,v] &=& f\big(t,x_*(t),x_*(t-\Delta), u,v\big).
\end{eqnarray*}
Mit diesen Schreibweisen ergibt der Grenzübergang $\lambda \to 0^+$ in (\ref{BeweisRet2}):
\begin{eqnarray*}
y(t)&=& \underbrace{\varphi[\tau,v,u_*(\tau-\Delta)] - \varphi[\tau,u_*(\tau),u_*(\tau-\Delta)]}_{=y(\tau)} 
        + \int_{\tau}^{t_1} \chi_{[t_0,t]}(s) \cdot \varphi_x[s] y(s) \,ds \\
    & & + \chi_{[t_0,t]}(\tau+\Delta) \cdot
        \Big[ \underbrace{\varphi[\tau+\Delta,u_*(\tau+\Delta),v] - \varphi[\tau+\Delta,u_*(\tau+\Delta),u_*(\tau)]}_{=\eta(\tau+\Delta)} \Big] \\
    & & + \int_{\tau+\Delta}^{t_1}\chi_{[t_0,t]}(s) \cdot \varphi_y[s] y(s-\Delta) \,ds.
\end{eqnarray*}
Das erste Integral beschränken wir auf das Intervall $[\tau,t]$ und entfernen dafür die charakteristische Funktion im Integranden.
Im zweiten Integral stimmt der Geltungsbereich $[\tau+\Delta,t_1] \cap [t_0,t]$ mit $[\tau,t] \cap [\tau+\Delta,t_1]$ überein.
Ferner dürfen wir $\chi_{[t_0,t- \Delta]}(\tau)$ anstelle von $\chi_{[t_0,t]}(\tau+\Delta)$ schreiben.
Es ergibt sich
\begin{eqnarray}
y(t) &=& y(\tau) +  \chi_{[t_0,t-\Delta]}(\tau) \cdot \eta(\tau+\Delta) \nonumber\\
     & & \label{BeweisRet21} + \int_{\tau}^t \varphi_x[s] y(s) \,ds + \int_{\tau}^t \chi_{[\tau+\Delta,t_1]}(s) \cdot \varphi_y[s] y(s-\Delta) \,ds. 
\end{eqnarray}
Es kann $y(t)$ in $t=\tau+\Delta \in (t_0,t_1)$ unstetig sein und dann gilt
$$\lim_{\varepsilon \to 0^+} [y(\tau+\Delta+\varepsilon) -y(\tau+\Delta-\varepsilon)] =\eta(\tau+\Delta).$$
Wir werten nun die Ableitung
$\displaystyle \frac{d}{dt} \langle p(t), y(t) \rangle = \langle \dot{p}(t), y(t) \rangle + \langle p(t), \dot{y}(t) \rangle$
aus.
Es gelten für $y(\cdot)$ und nach (\ref{SatzRetPMP1}) für die Adjungierte $p(\cdot)$ in $[\tau,t_1]$:
\begin{eqnarray*}
\dot{y}(t) &=& \varphi_x[t] y(t) + \chi_{[\tau+\Delta,t_1]}(t) \cdot \varphi_y[t] y(t-\Delta), \\
\dot{p}(t) &=& - \varphi^T_x[t] p(t)+f_x[t] + \chi_{[\tau,t_1-\Delta]}(t) \cdot \big( - \varphi^T_y[t+\Delta] p(t+\Delta)+f_y[t+\Delta] \big).
\end{eqnarray*}
Für die $\varphi-$Terme mit Zeitverschiebung $\pm\Delta$ ergibt sich mit Hilfe von (\ref{Zeitshift})
$$\int_\tau^{t_1} \chi_{[\tau+\Delta,t_1]}(t) \cdot \big\langle p(t) , \varphi_y[t] y(t-\Delta) \big\rangle \, dt
   =\int_\tau^{t_1} \chi_{[\tau,t_1-\Delta]}(t) \cdot \big\langle \varphi^T_y[t+\Delta] p(t+\Delta) , y(t) \big\rangle \, dt.$$
Da außerdem $\langle - \varphi^T_x[t] p(t), y(t) \rangle = - \langle p(t), \varphi_x[t] y(t) \rangle$ gilt, erhalten wir
\begin{eqnarray*}
    \langle p(t_1), y(t_1) \rangle
&=& \langle p(\tau), y(\tau) \rangle + \chi_{[t_0,t_1-\Delta]}(\tau) \cdot \langle p(\tau+\Delta), \eta(\tau+\Delta) \rangle
     + \int_\tau^{t_1} \frac{d}{dt} \langle p(t), y(t) \rangle \, dt \\
&=& \langle p(\tau), y(\tau) \rangle + \chi_{[t_0,t_1-\Delta]}(\tau) \cdot \langle p(\tau+\Delta), \eta(\tau+\Delta) \rangle \\
& & \hspace*{1cm} + \int_{\tau}^{t_1} \langle f_x[t] , y(t) \rangle \, dt 
    + \int_{\tau}^{t_1} \chi_{[\tau,t_1-\Delta]}(t) \cdot \langle f_y[t+\Delta], y(t) \rangle \, dt.
\end{eqnarray*}
Im zweiten Integral ebenso die Umformung (\ref{Zeitshift}) angewendet liefert dann
\begin{eqnarray}
    \hspace*{-5mm} \langle p(t_1), y(t_1) \rangle
&=& \langle p(\tau), y(\tau) \rangle + \chi_{[t_0,t_1-\Delta]}(\tau) \cdot \langle p(\tau+\Delta), \eta(\tau+\Delta) \rangle \nonumber \\
\hspace*{-5mm} & & \label{BeweisRet3} + \int_{\tau}^{t_1} \langle f_x[t] , y(t) \rangle \, dt 
    + \int_{\tau}^{t_1} \chi_{[\tau+\Delta,t_1]}(t) \cdot \langle f_y[t], y(t-\Delta) \rangle \, dt, \\
\hspace*{-5mm} y(\tau) &=& \label{BeweisRet4} \varphi[\tau,v,u_*(\tau-\Delta)] - \varphi[\tau,u_*(\tau),u_*(\tau-\Delta)], \\
\hspace*{-5mm} \eta(\tau+\Delta) &=& \label{BeweisRet5} \varphi[\tau+\Delta,u_*(\tau+\Delta),v] - \varphi[\tau+\Delta,u_*(\tau+\Delta),u_*(\tau)].
\end{eqnarray}
Die Auswertung der dynamischen Wirkung der Nadelvariation ist (endlich) abgeschlossen und wir kommen zum Abschluss des Beweises:
Da $\big(x_*(\cdot),u_*(\cdot)\big)$ ein starkes lokales Minimum ist, gilt für alle hinreichend kleine $\lambda >0$ die Ungleichung
$$0 \leq \frac{J\big(x_\lambda(\cdot),u_\lambda(\cdot)\big)- J\big(x_*(\cdot),u_*(\cdot)\big)}{\lambda}.$$
Den Quotienten auf der rechten Seite der Ungleichung formen wir auf die gleiche Weise wie zur Herleitung von (\ref{BeweisRet2}) bzw. (\ref{BeweisRet21}) um
und erhalten im Grenzübergang $\lambda \to 0^+$:
\begin{eqnarray}
\hspace*{-10mm} 0 &\leq& f[\tau,v,u_*(\tau-\Delta)] -f[\tau,u_*(\tau),u_*(\tau-\Delta)]  \nonumber\\
\hspace*{-10mm} & & + \chi_{[t_0,t_1-\Delta]}(\tau) \cdot \Big[ f[\tau+\Delta,u_*(\tau+\Delta),v] - f[\tau+\Delta,u_*(\tau+\Delta),u_*(\tau)] \Big]  \nonumber\\
\hspace*{-10mm} & & \label{BeweisRet6}
       \hspace*{-1cm} + \int_{\tau}^{t_1} \langle f_x[t] , y(t) \rangle \,dt
        +  \int_{\tau}^{t_1} \chi_{[\tau+\Delta,t_1]}(t) \cdot \langle f_y[t] , y(t-\Delta) \rangle \,dt
        + \big\langle S'\big(x_*(t_1)\big), y(t_1) \big\rangle. 
\end{eqnarray}
Die Bedingungen (\ref{BeweisRet3})--(\ref{BeweisRet6}) und die Transversalitätsbedingung (\ref{SatzRetPMP1}) ergeben 
\begin{eqnarray}
0 &\leq& f[\tau,v,u_*(\tau-\Delta)] -f[\tau,u_*(\tau),u_*(\tau-\Delta)]  \nonumber\\
    & & + \chi_{[t_0,t_1-\Delta]}(\tau) \cdot\Big[ f[\tau+\Delta,u_*(\tau+\Delta),v] - f[\tau+\Delta,u_*(\tau+\Delta),u_*(\tau)] \Big]  \nonumber\\
    & & \label{BeweisRet7} - \langle p(\tau), y(\tau) \rangle - \chi_{[t_0,t_1-\Delta]}(\tau)\cdot\langle p(\tau+\Delta), \eta(\tau+\Delta) \rangle. 
\end{eqnarray}
Da $\tau \in (t_0,t_1)$ als ein Stetigkeitspunkt der Steuerung $u_*(\cdot)$ und zudem $v \in U$ beliebig gewählt waren,
liefert (\ref{BeweisRet6}) die Gültigkeit der Maximumbedingung (\ref{SatzRetPMP3}).
Der Beweis von Theorem \ref{SatzRetPMP} ist damit abgeschlosssen. \hfill $\blacksquare$   
       \newpage
       \lhead[\thepage \hspace*{1mm} Arrow-Bedingungen]{ }    
       \subsection{Hinreichende Bedingungen nach Arrow}
Die Ableitung von hinreichenden Optimalitästbedingungen nach
Arrow \index{hinreichende Bedingungen nach Arrow!Standard@-- Zeitverzögerte Systeme}
gestaltet sich durch das Auftreten der zeitlichen Verschiebung um $\Delta>0$ schwieriger.
Dennoch führt uns das Vorgehen in Abschnitt \ref{AbschnittHBPMP} zum Ziel. \\[2mm]
Im Rahmen der zeitverzögerten Systeme treten die Zustände $x(t)$ und $x(t-\Delta)$ auf.
Aus diesem Grund verwenden wir für $t \in [t_0,t_1]$ neben der aus Abschnitt \ref{AbschnittHBPMP}
bekannten Menge $V_\gamma(t)=\{ x \in \R^n \,|\, \|x-x_*(t)\| < \gamma\}$
außerdem die Menge $U_\gamma(t)=V_\gamma(t) \times V_\gamma(t-\Delta)$, 
$$U_\gamma(t) = \{ (x,y) \in \R^n \times \R^n \,|\, \|x-x_*(t)\| < \gamma, \|y-x_*(t-\Delta)\| < \gamma\}.$$
Ferner bezeichnet $\mathscr{H}^{\mathcal{Z}}$ die Hamilton-Funktion
$$\mathscr{H}^{\mathcal{Z}}(t,x,y,p) = \sup_{u,v \in U} H^{\mathcal{Z}}\big(t,x,y,u,v,p).$$

\begin{theorem} 
In der Aufgabe (\ref{Ret1})--(\ref{Ret3}) sei
$\big(x_*(\cdot),u_*(\cdot)\big) \in \mathscr{A}^{\mathcal{Z}}_{\rm adm} \cap \mathscr{A}^{\mathcal{Z}}_{\rm Lip}$
und es sei $p(\cdot) \in PC_1([t_0,t_1],\R^n)$. Ferner gelte:
\begin{enumerate}
\item[(a)] Das Tripel $\big(x_*(\cdot),u_*(\cdot),p(\cdot)\big)$
           erf"ullt (\ref{SatzRetPMP1})--(\ref{SatzRetPMP3}) in Theorem \ref{SatzRetPMP}.        
\item[(b)] F"ur jedes $t \in [t_0,t_1]$ ist die Hamilton-Funktion $\mathscr{H}^{\mathcal{Z}}\big(t,x,y,p(t)\big)$ konkav auf $U_\gamma(t)$.
\item[(c)] Die Abbildung $S:\R^n\to \R$ im Zielfunktional (\ref{Ret1}) ist konvex in $x$ auf $V_\gamma(t_1)$.
\end{enumerate}
Dann ist $\big(x_*(\cdot),u_*(\cdot)\big)$ ein starkes lokales Minimum der Aufgabe (\ref{Ret1})--(\ref{Ret3}).
\end{theorem}

{\bf Beweis} Wie im Abschnitt \ref{AbschnittHBPMP} folgt,
dass $\alpha_*= -\mathscr{H}^{\mathcal{Z}}\big(t,x_*(t),x_*(t-\Delta),p(t)\big)$ ein Randpunkt der konvexen Menge 
$$Z=\big\{ (\alpha,x,y) \in \R \times \R^n \times \R^n \,\big|\, 
           (x,y) \in U_\gamma(t), \alpha \geq -\mathscr{H}^{\mathcal{Z}}\big(t,x,y,p(t)\big) \big\}$$
ist.
Daher existiert ein nichttrivialer Vektor $\big(a_0(t),a_1(t),a_2(t)\big) \in \R \times \R^n \times \R^n$ mit
$$a_0(t) \alpha + \langle a_1(t),x\rangle + \langle a_2(t),y\rangle
  \geq a_0(t) \alpha_* + \langle a_1(t),x_*(t)\rangle+ \langle a_2(t),x_*(t-\Delta)\rangle$$
f"ur alle $(\alpha,x,y) \in Z$.
Wählen wir hierin $y=x_*(t-\Delta)$,
dann ergibt sich aus
$$a_0(t) \alpha + \langle a_1(t),x\rangle \geq a_0(t) \alpha_* + \langle a_1(t),x_*(t)\rangle$$
mit der gleichen Argumentation wie in Abschnitt \ref{AbschnittHBPMP},
dass wir $a_0(t)=1$ annehmen können und erhalten
\begin{eqnarray}
\lefteqn{\langle a_1(t),x-x_*(t)\rangle + \langle a_2(t),y-x_*(t-\Delta)\rangle} \nonumber\\
&\geq& \label{BeweisHBPMPZ1} \mathscr{H}^{\mathcal{Z}}\big(t,x,y,p(t)\big) - \mathscr{H}^{\mathcal{Z}}\big(t,x_*(t),x_*(t-\Delta),p(t)\big).
\end{eqnarray}
Wählen wir in dieser Ungleichung nacheinander einerseits $x=x_*(t)$ und andererseits $y=x_*(t-\Delta)$,
so ergeben sich aus (\ref{BeweisHBPMPZ1}) für alle $x \in V_\gamma(t)$ und alle
$y \in V_\gamma(t-\Delta)$ die Relationen
\begin{eqnarray*}
       \langle a_1(t),x-x_*(t)\rangle
&\geq& \mathscr{H}^{\mathcal{Z}}\big(t,x,x_*(t-\Delta),p(t)\big) - \mathscr{H}^{\mathcal{Z}}\big(t,x_*(t),x_*(t-\Delta),p(t)\big), \\
       \langle a_2(t),y-x_*(t-\Delta)\rangle
&\geq& \mathscr{H}^{\mathcal{Z}}\big(t,x_*(t),y,p(t)\big) - \mathscr{H}^{\mathcal{Z}}\big(t,x_*(t),x_*(t-\Delta),p(t)\big).
\end{eqnarray*}
In der zweiten Ungleichung nehmen wir eine Zeitverschiebung um $\Delta$ vor und überführen den Zeitpunkt $t$ in $t + \Delta$.
Damit gehört $y$ der Menge $V_\gamma(t)$ an und wir dürfen nach der Verschiebung $y=x \in V_\gamma(t)$ wählen.
So ergibt sich für alle $x \in V_\gamma(t)$ die Ungleichung
\begin{eqnarray}
      \lefteqn{\langle a_1(t)+\chi_{[t_0,t_1-\Delta]}(t) \cdot a_2(t+\Delta),x-x_*(t)\rangle} \nonumber \\
&\geq& \mathscr{H}^{\mathcal{Z}}\big(t,x,x_*(t-\Delta),p(t)\big) - \mathscr{H}^{\mathcal{Z}}\big(t,x_*(t),x_*(t-\Delta),p(t)\big) \nonumber \\
&& + \chi_{[t_0,t_1-\Delta]}(t) \cdot \Big(\mathscr{H}^{\mathcal{Z}}\big(t+\Delta,x_*(t+\Delta),x,p(t+\Delta)\big) \nonumber \\
&& \label{BeweisHBPMPZ2} \hspace*{3,5cm} - \mathscr{H}^{\mathcal{Z}}\big(t+\Delta,x_*(t+\Delta),x_*(t),p(t+\Delta)\big)\Big).
\end{eqnarray} 
Es sei in $t \in [t_0,t_1]$ die Maximumbedingung (\ref{SatzRetPMP3}) erf"ullt.
Mit Hilfe der Pontrjagin-Funktion können wir die Ungleichung (\ref{BeweisHBPMPZ2}) weiterführen und erhalten
\begin{eqnarray}
\lefteqn{\langle a_1(t)+\chi_{[t_0,t_1-\Delta]}(t) \cdot a_2(t+\Delta),x-x_*(t)\rangle} \nonumber \\
&\geq& H^{\mathcal{Z}}\big(t,x,x_*(t-\Delta),u_*(t),u_*(t-\Delta),p(t)\big) - \mathscr{H}^{\mathcal{Z}}\big(t,x_*(t),x_*(t-\Delta),p(t)\big) \nonumber \\
&& + \chi_{[t_0,t_1-\Delta]}(t) \cdot \Big(H^{\mathcal{Z}}\big(t+\Delta,x_*(t+\Delta),x,u_*(t+\Delta),u_*(t),p(t+\Delta)\big) \nonumber \\
&& \label{BeweisHBPMPZ3} \hspace*{3,5cm} - \mathscr{H}^{\mathcal{Z}}\big(t+\Delta,x_*(t+\Delta),x_*(t),p(t+\Delta)\big)\Big)=\Psi(x)
\end{eqnarray} 
f"ur alle $x \in V_\gamma(t)$.
Mit Hilfe der rechten Seite $\Psi(x)$ von (\ref{BeweisHBPMPZ3}) bilden wir die Funktion
$$\Phi(x) = \Psi(x) - \langle a_1(t)+\chi_{[t_0,t_1-\Delta]}(t) \cdot a_2(t+\Delta),x-x_*(t)\rangle,$$
welche in dem inneren Punkt $x_*(t)$ der Menge $V_\gamma(t)$ ihr globales Maximum annimmt.
Also gilt $0=\Phi'(x_*(t))$, d.\,h.
\begin{eqnarray}
\lefteqn{-a_1(t)-\chi_{[t_0,t_1-\Delta]}(t) \cdot a_2(t+\Delta)} \nonumber \\
&=& -H_x^{\mathcal{Z}}\big(t,x_*(t),x_*(t-\Delta),u_*(t),u_*(t-\Delta),p(t)\big) \nonumber \\
& & \label{BeweisHBPMPZ4} - \chi_{[t_0,t_1-\Delta]}(t) \cdot H_y^{\mathcal{Z}}\big(t+\Delta,x_*(t+\Delta),x_*(t),u_*(t+\Delta),u_*(t),p(t+\Delta)\big).
\end{eqnarray}
Die adjungierte Gleichung (\ref{SatzRetPMP1}) zeigt nun
\begin{equation} \label{BeweisHBPMPZ5}
\dot{p}(t) = -a_1(t)-\chi_{[t_0,t_1-\Delta]}(t) \cdot a_2(t+\Delta)\quad \mbox{ für fast alle } t \in [t_0,t_1].
\end{equation}
Es sei $\big(x(\cdot),u(\cdot)\big) \in \mathscr{A}^{\mathcal{Z}}_{\rm adm}$ mit $\|x(\cdot)-x_*(\cdot)\|_\infty < \gamma$.
Wir erhalten mit (\ref{BeweisHBPMPZ1}) und (\ref{BeweisHBPMPZ5}):
\begin{eqnarray*}
\lefteqn{\int_{t_0}^{t_1} \big[H^{\mathcal{Z}}\big(t,x_*(t),x_*(t-\Delta),u_*(t),u_*(t-\Delta),p(t)\big)} \\
& & \hspace*{3cm} -H^{\mathcal{Z}}\big(t,x(t),x(t-\Delta),u(t),u(t-\Delta),p(t)\big)\big] \, dt \\
&=& \int_{t_0}^{t_1} \big[\mathscr{H}^{\mathcal{Z}}\big(t,x_*(t),x_*(t-\Delta),p(t)\big)
                           -H^{\mathcal{Z}}\big(t,x(t),x(t-\Delta),u(t),u(t-\Delta),p(t)\big)\big] \, dt \\
&\geq& \int_{t_0}^{t_1} \big[\mathscr{H}^{\mathcal{Z}}\big(t,x_*(t),x_*(t-\Delta),p(t)\big)
              -\mathscr{H}^{\mathcal{Z}}\big(t,x(t),x(t-\Delta),p(t)\big)\big] \, dt \\
&\geq& \int_{t_0}^{t_1} \big(\langle -a_1(t),x(t)-x_*(t)\rangle + \langle -a_2(t),x(t-\Delta)-x_*(t-\Delta)\rangle \big) \, dt \\
&=& \int_{t_0}^{t_1} \langle \underbrace{-a_1(t)-\chi_{[t_0,t_1-\Delta]}(t) \cdot a_2(t+\Delta)}_{=\dot{p}(t)},x(t)-x_*(t)\rangle \, dt.
\end{eqnarray*}
Abschließend können wir festhalten:
\begin{eqnarray*}
\lefteqn{J\big(x(\cdot),u(\cdot)\big)-J\big(x_*(\cdot),u_*(\cdot)\big)=S\big(x(t_1)\big)-S\big(x_*(t_1)\big)} \\
&& \hspace*{-5mm} + \int_{t_0}^{t_1} \big[ f\big(t,x(t),x(t-\Delta),u(t),u(t-\Delta)\big)  -f\big(t,x_*(t),x_*(t-\Delta),u_*(t),u_*(t-\Delta)\big)\big] \, dt \\
&=& \int_{t_0}^{t_1} \big[H^{\mathcal{Z}}\big(t,x_*(t),x_*(t-\Delta),u_*(t),u_*(t-\Delta),p(t)\big) \\
& & \hspace*{3cm} -H^{\mathcal{Z}}\big(t,x(t),x(t-\Delta),u(t),u(t-\Delta),p(t)\big)\big] \, dt  \\
& & + \int_{t_0}^{t_1} \langle p(t), \dot{x}(t)-\dot{x}_*(t) \rangle dt + S\big(x(t_1)\big)-S\big(x_*(t_1)\big) \\
&\geq& \int_{t_0}^{t_1}  \big[\langle \dot{p}(t),x(t)-x_*(t)\rangle + \langle p(t), \dot{x}(t)-\dot{x}_*(t) \rangle \big] \, dt
    + \big\langle S'\big(x_*(t_1)\big), x(t_1)-x_*(t_1) \big\rangle \\
&=& \langle p(t_1)+S'\big(x_*(t_1)\big),x(t_1)-x_*(t_1)\rangle-\langle p(t_0),x(t_0)-x_*(t_0)\rangle \geq 0
\end{eqnarray*}
für alle $\big(x(\cdot),u(\cdot)\big) \in \mathscr{A}^{\mathcal{Z}}_{\rm adm}$
mit $\|x(\cdot)-x_*(\cdot)\|_\infty < \gamma$. \hfill $\blacksquare$

\begin{beispiel}
{\rm Die Hamilton-Funktion im Beispiel \ref{BeispielInvDelay},
$$\mathscr{H}^{\mathcal{Z}}(t,K,L,u,v,p) = \sup_{u,v \in [0,1]} \{ p vL +(1-u)K\} =pL+K,$$
ist für $p,K,L \geq 0$ konkav und die notwendigen Optimalitätsbedingungen sind in diesem Beispiel gleichzeitig hinreichend. \hfill $\square$}
\end{beispiel}  
\cleardoublepage

\lhead[ ]{}
\rhead[]{Volterrasche Integralgleichungen \hspace*{1mm} \thepage}
\section{Steuerung Volterrascher Integralgleichungen}
       Integralgleichungen treten in natürlicher Weise in dynamischen Problemen auf,
in denen das System eine Form von ``Erinnerungsvermögen'' besitzt.
Da sich im Rahmen der Standardaufgabe und deren Erweiterungen der Einfluss einer Steuervariable stets
unmittelbar auf den Zustand auswirkt,
kann ein Effekt, der sich im Laufe der Zeit entwickelt, meist nicht modelliert werden.
Hier ist die Beschreibung des dynamischen Systems mit Hilfe einer Integralgleichung ein probates Mittel. \\[2mm]
Ein Effekt über Zeit kann bei der Zusammenstellung der aggregierten Produktionskapazitäten vorliegen.
Die Kapazitäten ergeben sich als die gesamten Investitionen in Produktionsanlagen in den vorhergehenden Jahrgängen.
Durch Verschlei"s, Wartung oder technologischen Fortschritt wird die Produktionsfähigkeit beeinflusst.
Zur Zeit $t$ sei die Effizienz der Anlagen des Jahrgangs $s \leq t$ durch eine Funktion $\pi(t,s)$ beschrieben.
Zum Zeitpunkt $t$ ergeben sich demnach die gesamten vorliegenden Produktionskapazitäten durch 
$\displaystyle P(t)=\int_{t_0}^t \pi(t,s) \, ds$. \\
Im Rahmen der Werbeindustrie wird durch den strategischen Auf- und Ausbau der Produktbekanntheit dem Konsument ein gewisses
Produktimage und eine Verbundenheit zu dem Produkt suggeriert.
Dabei ist das Management mit der Problemstellung konfrontiert,
dass der Bekanntheitsgrad des Produktes durch Werbema"snahmen über eine gewisse Zeit aufgebaut werden muss.
Andererseits lässt die Produktbekanntheit (in Folge dessen auch die Nachfrage) umso stärker nach,
je länger der Werbeimpuls in der Vergangenheit liegt.
Die aggregierte Auswirkung der Werbekampagne zum einem Zeitpunkt ist demnach das Resultat der gesamten Werbeanstrengungen,
die im Vorfeld unternommen wurden.  
       \lhead[\thepage \hspace*{1mm} Vorbereitungen]{ }    
       \subsection{Vorbereitende Betrachtungen}
Der wesentliche Unterschied zur Standardaufgabe besteht darin,
dass die Dynamik
$$\dot{x}(s) = \varphi\big(s,x(s),u(s)\big)$$
nicht nur die eine rechte Seite $\varphi(s,x,u)$,
sondern eine Familie $\{\varphi(t,s,x,u),\, t \in [t_0,t_1]\}$ rechter Seiten besitzt:
$$x(t) = x_0 + \int_{t_0}^t \varphi\big(t,s,x(s),u(s)\big) \, ds, \quad t \in [t_0,t_1].$$
Anstelle der einen Zeitvariablen liegen nun eine ``innere'' Zeitvariable $s$ und eine ``äußere'' Zeitvariable $t$ vor.
Die Konsequenz daraus ist,
dass die Dynamik statt über dem Intervall $\{t \in \R \,|\, t_0 \leq t \leq t_1\}$ nun
über dem Zeitbereich $\{(s,t) \in \R^2 \,|\, t_0 \leq s, t \leq t_1\}$ agiert.

\begin{bemerkung}
{\rm Damit man Ergebnisse für $t <s$ nicht von Beginn an ausschließt,
z.\,B. in Satz \ref{SatzFixpunktlokal} das Teilintervall $[t_0-\varepsilon,t_0)$,
ist es günstiger die Betrachtungen nicht auf den
Zeitbereich $\{(s,t) \in \R^2 \,|\, t_0 \leq s \leq t \leq t_1\}$ einzuschränkungen. \hfill $\square$}
\end{bemerkung}

Im Vergleichung zur Steuerung von Differentialgleichungen besitzt die Steuerung von Integralgleichungen
einen allgemeineren Charakter.
Insbesondere die vorliegende Situation eines dynamischen Systems mit einer Familie von rechten Seiten hat uns dazu bewegt,
die Lösungstheorie zu Differential- und Integralgleichungen im Anhang \ref{AnhangDGL} ausführlich darzustellen. \\[2mm]
Die Herleitung der Optimalitätsbedingungen führt uns zu Stieltjes-Integralen
$$\int_a^b f(t) \, dg(t)$$
mit dem Integranden $f(\cdot) \in PC([a,b],\R)$ und dem Integrator $g(\cdot) \in PC_1([a,b),\R)$.
Die Funktion $g(\cdot)$ ist dabei nur über dem halboffenen Intervall $[a,b)$ stetig und kann im Endpunkt unstetig sein.
Die Zwischensummen für Stieltjes-Integrale besitzen mit den Zerlegungspunkten $a=t_0<t_1<...<t_{N+1}=b$,
den Zwischenpunkten $s_i \in [t_i,t_{i+1}]$ die Gestalt
$$\sum_{i=0}^N f(s_i)[g(t_{i+1})-g(t_i)].$$
Da $g(\cdot)$ über $[a,b)$ stückweise stetig differenzierbar ist,
gilt -- bis auf endlich viele Stellen -- $g(t_{i+1})-g(t_i)=\dot{g}(\sigma_i) \cdot (t_{i+1}-t_i)$ mit den Zwischenstellen $\sigma_i \in (t_i,t_{i+1})$.
Bei Verfeinerung $\Delta=\max (t_{i+1}-t_i) \to 0$ ergibt sich
\begin{eqnarray*}
    \sum_{i=0}^N f(s_i)[g(t_{i+1})-g(t_i)]
&=& \sum_{i=0}^{N-1} f(s_i)\dot{g}(\sigma_i) \cdot (t_{i+1}-t_i) + f(s_N)[g(b)-g(t_N)] \\
&\to& \int_a^b f(t) \dot{g}(t) \, dt +f(b)[g(b)-g(b^-)].
\end{eqnarray*}
Zusammenfassend erhalten wir die Beziehung
\begin{equation} \label{StieltjesIGL}
\int_a^b f(t) \, dg(t) = \int_a^b f(t) \dot{g}(t) \, dt + f(b)[g(b)-g(b^-)].
\end{equation}
Die Gleichung (\ref{StieltjesIGL}) führt zu zwei Varianten der Optimalitätsbedingungen,
wobei die Sprungbedingung im Endzeitpunkt beachtet werden muss: \\[2mm]
Das Stieltjes-Integral auf der linken Seite der Beziehung (\ref{StieltjesIGL}) spiegelt den unverfälschten Charakter der Optimalitätsbedingungen wider.
Dabei ist die Sprungbedingung im Endzeitpunkt zentraler Bestandteil der Transversalitätsbedingungen.
Auf Grundlage dieser Darstellung geben wir eine ökonomische Interpretation des Maximumprinzips. \\[2mm]
Die rechte Seite in (\ref{StieltjesIGL}) mit einem Riemann-Integral besitzt einen zusätzlichem Randterm.
Von diesem fließt lediglich $g(b^-)$ in die Pontrjagin-Funktion und in die Transversalitätsbedingungen ein.
Da bei der nachfolgenden Festlegung der Pontrjagin-Funktion mittels des Riemann-Integrals die Sprungbedingung
im Endpunkt bereits Beachtung fand,
darf auf sie verzichtet werden.          
       \newpage
       \lhead[\thepage \hspace*{1mm} Pontrjaginsches Maximumprinzip]{ }      
       \subsection{Die Aufgabenstellung und das Pontrjaginsche Maximumprinzip} \label{AbschnittPMPeinfachIGL}
Die Steuerung einer Volterraschen Integralgleichung bezeichnet die Aufgabe 
\begin{eqnarray}
&&\label{PMPeinfachIGL1} J\big(x(\cdot),u(\cdot)\big) = \int_{t_0}^{t_1} f\big(t,x(t),u(t)\big) \, dt  +S\big(x(t_1)\big) \to \inf, \\
&&\label{PMPeinfachIGL2} x(t) = x_0 + \int_{t_0}^t \varphi\big(t,s,x(s),u(s)\big) \, ds, \quad t \in [t_0,t_1], \\
&&\label{PMPeinfachIGL3} u(t) \in U \subseteq \R^m, \quad U \not= \emptyset,
\end{eqnarray}
die wir bez"uglich $\big(x(\cdot),u(\cdot)\big) \in PC_1([t_0,t_1],\R^n) \times PC([t_0,t_1],U)$ untersuchen. \\[2mm]
Mit $\mathscr{A}^{\mathcal{I}}_{\rm Lip}$ bezeichnen wir die Menge $\big(x(\cdot),u(\cdot)\big)$,
für die es ein $\gamma>0$ derart gibt,
dass die Abbildungen $f(s,x,u)$, $\varphi(t,s,x,u)$ auf der Menge aller $(t,s,x,u) \in \R \times \R \times \R^n \times \R^m$ mit
$$t_0 \leq s, t \leq t_1, \qquad \|x-x(s)\| < \gamma, \qquad u \in \R^m$$
stetig in der Gesamtheit aller Variablen und stetig differenzierbar bezüglich $x$ sind. \\[2mm]
Das Paar $\big(x(\cdot),u(\cdot)\big) \in PC_1([t_0,t_1],\R^n) \times PC([t_0,t_1],U)$
hei"st ein zul"assiger Steuerungsprozess in der Aufgabe (\ref{PMPeinfachIGL1})--(\ref{PMPeinfachIGL3}),
falls $\big(x(\cdot),u(\cdot)\big)$ dem System (\ref{PMPeinfachIGL2}) zu $x(t_0)=x_0$ gen"ugt.
Mit $\mathscr{A}^{\mathcal{I}}_{\rm adm}$ bezeichnen wir die Menge der zul"assigen Steuerungsprozesse. \\[2mm]
Ein zul"assiger Steuerungsprozess $\big(x_*(\cdot),u_*(\cdot)\big)$ ist eine
starke lokale Minimalstelle\index{Minimum, starkes lokales!Integral@-- Integralgleichungen}
der Aufgabe (\ref{PMPeinfachIGL1})--(\ref{PMPeinfachIGL3}),
falls eine Zahl $\varepsilon > 0$ derart existiert, dass die Ungleichung 
$$J\big(x(\cdot),u(\cdot)\big) \geq J\big(x_*(\cdot),u_*(\cdot)\big)$$
f"ur alle $\big(x(\cdot),u(\cdot)\big) \in \mathscr{A}^{\mathcal{I}}_{\rm adm}$
mit $\|x(\cdot)-x_*(\cdot)\|_\infty < \varepsilon$ gilt. \\[2mm]
Ausgehend von der Gleichung (\ref{StieltjesIGL}) werden wir zwei Varianten des Pontrjaginschen Maximumprinzips formulieren.
Die erste Variante spiegelt einen ``unverfälschten'' Charakter wider.
Aufgrund (\ref{StieltjesIGL}) entsteht in $t=t_1$ eine Sprung-Transversalitätsbedingung.
Daraus resultiert als eine Folgerung die zweite Variante.
Sie ermöglicht einen Vergleich mit Optimalitätsbedingungen in Standardreferenzen.  \\[2mm]
Für die Aufgabe (\ref{PMPeinfachIGL1})--(\ref{PMPeinfachIGL3}) lautet die ``unverfälschte'' Pontrjagin-Funktion 
$$H^{\mathcal{I}}(t,x,u,p(\cdot)\big) = -\int_t^{t_1} \langle \varphi(\tau,t,x,u), dp(\tau) \rangle - f(t,x,u).$$

\begin{theorem}[Pontrjaginsches Maximumprinzip] \label{SatzPMPeinfachIGL2}
\index{Pontrjaginsches Maximumprinzip!Integral@-- Integralgleichungen} 
Sei $\big(x_*(\cdot),u_*(\cdot)\big) \in \mathscr{A}^{\mathcal{I}}_{\rm adm} \cap \mathscr{A}^{\mathcal{I}}_{\rm Lip}$. 
Ist $\big(x_*(\cdot),u_*(\cdot)\big)$ ein starkes lokales Minimum der Aufgabe (\ref{PMPeinfachIGL1})--(\ref{PMPeinfachIGL3}),
dann existiert eine Vektorfunktion $p(\cdot) \in PC_1([t_0,t_1),\R^n)$ derart, dass
\begin{enumerate}
\item[(a)] die adjungierte Gleichung
           \index{adjungierte Gleichung!Integral@-- Integralgleichungen}
           \begin{eqnarray} 
           \dot{p}(t) &=& -H_x^{\mathcal{I}}\big(t,x_*(t),u_*(t),p(\cdot)\big) \nonumber \\
                      &=& \label{SatzPMPIGL1} \int_t^{t_1} \varphi^T_x\big(\tau,t,x_*(t),u_*(t)\big)\, dp(\tau) + f_x\big(t,x_*(t),u_*(t)\big),
           \end{eqnarray}
\item[(b)] in $t=t_1$ die Sprung-Transversalitätsbedingung
           \index{Transversalitätsbedingungen!Integral@-- Integralgleichungen}
           \begin{equation} \label{SatzPMPIGL2} p(t_1^-)=-S'\big(x_*(t_1)\big), \qquad p(t_1)=0, \end{equation}
\item[(c)] in fast allen Punkten $t \in [t_0,t_1]$ die Maximumbedingung
           \index{Maximumbedingung!Integral@-- Integralgleichungen}
           \begin{equation} \label{SatzPMPIGL3}
           H^{\mathcal{I}}\big(t,x_*(t),u_*(t),p(\cdot)\big) = \max_{u \in U} H^{\mathcal{I}}\big(t,x_*(t),u,p(\cdot)\big)
           \end{equation}
\end{enumerate}
erfüllt sind.
\end{theorem}

Da die Adjungierte $p(\cdot) \in PC_1([t_0,t_1),\R^n)$ stetig und stückweise stetig differenzierbar ist, sowie
in $t=t_1$ die Sprung-Transversalitätsbedingung (\ref{SatzPMPIGL2}) erfüllt ist,
ergibt sich mit der Beziehung (\ref{StieltjesIGL}) der Zusammenhang
\begin{eqnarray*}
    \int_t^{t_1} \varphi^T_x[\tau,t] \, dp(\tau) 
&=& \int_t^{t_1} \varphi^T_x[\tau,t] \dot{p}(\tau) \, d\tau + \varphi^T_x[t_1,t]) \big(p(t_1)-p(t_1^-)\big) \\
&=& \int_t^{t_1} \varphi^T_x[\tau,t] \dot{p}(\tau) \, d\tau - \varphi^T_x[t_1,t] p(t_1^-).
\end{eqnarray*}
Setzen wir $p(\cdot)$ in $t=t_1$ stetig fort, so erhalten wir die Pontrjagin-Funktion
$$H_0^{\mathcal{I}}(t,x,u,p(\cdot)\big) = -\int_t^{t_1} \langle \varphi(\tau,t,x,u), \dot{p}(\tau) \rangle \, d\tau
                                    +\langle \varphi(t_1,t,x,u), p(t_1)\rangle - f(t,x,u).$$
Da in der Pontrjagin-Funktion $H_0^{\mathcal{I}}$ die Sprung-Bedingung in $t=t_1$ bereits zur Anwendung kam,
ergeben sich als Folgerung zu Theorem \ref{SatzPMPeinfachIGL2} die Bedingungen: 
  
\begin{theorem}[Pontrjaginsches Maximumprinzip] \label{SatzPMPeinfachIGL}
\index{Pontrjaginsches Maximumprinzip!Integral@-- Integralgleichungen} 
Sei $\big(x_*(\cdot),u_*(\cdot)\big) \in \mathscr{A}^{\mathcal{I}}_{\rm adm} \cap \mathscr{A}^{\mathcal{I}}_{\rm Lip}$. 
Ist $\big(x_*(\cdot),u_*(\cdot)\big)$ ein starkes lokales Minimum der Aufgabe (\ref{PMPeinfachIGL1})--(\ref{PMPeinfachIGL3}),
dann existiert eine Vektorfunktion $p(\cdot) \in PC_1([t_0,t_1],\R^n)$ derart, dass
\begin{enumerate}
\item[(a)] die adjungierte Gleichung
           \index{adjungierte Gleichung!Integral@-- Integralgleichungen}
           \begin{equation}\label{PMPeinfachIGL4} 
           \dot{p}(t) = -H_{0x}^{\mathcal{I}}\big(t,x_*(t),u_*(t),p(\cdot)\big),
           \end{equation}
\item[(b)] in $t=t_1$ die Transversalitätsbedingung
           \index{Transversalitätsbedingungen!Integral@-- Integralgleichungen}
           \begin{equation}\label{PMPeinfachIGL5} p(t_1)=-S'\big(x_*(t_1)\big), \end{equation}
\item[(c)] in fast allen Punkten $t \in [t_0,t_1]$ die Maximumbedingung
           \index{Maximumbedingung!Integral@-- Integralgleichungen}
           \begin{equation}\label{PMPeinfachIGL6} 
           H_0^{\mathcal{I}}\big(t,x_*(t),u_*(t),p(\cdot)\big) = \max_{u \in U} H_0^{\mathcal{I}}\big(t,x_*(t),u,p(\cdot)\big).
           \end{equation}
\end{enumerate}
erfüllt sind.
\end{theorem}    
       \newpage
       \lhead[\thepage \hspace*{1mm} Diskussion des Maximumprinzips]{ }  
       \subsection{Diskussion und Ökonomische Deutung}
Die Pontrjagin-Funktion $H_0^{\mathcal{I}}$ als Riemann-Integral mit Randterm 
stellt eine Verallgemeinerung der Pontrjagin-Funktion $H^{\mathcal{S}}$ dar,
was dem verallgemeinerten Charakter der Aufgabe (\ref{PMPeinfachIGL1})--(\ref{PMPeinfachIGL3})
zur Aufgabe (\ref{PMPeinfach1})--(\ref{PMPeinfach3}) entspricht.
Tatsächlich liefert Theorem \ref{SatzPMPeinfachIGL} eine Erweiterung von Theorem \ref{SatzPMPeinfach} in Abschnitt \ref{AbschnittPMPeinfach}: \\[2mm]
Mit der Abbildung $\varphi(\tau,t,x,u)=\varphi(t,x,u)$ ohne einer ``äußeren'' Zeitvariablen
besitzt die Aufgabe (\ref{PMPeinfachIGL1})--(\ref{PMPeinfachIGL3}) die Gestalt
der Aufgabe (\ref{PMPeinfach1})--(\ref{PMPeinfach3}).
Zusammen mit der Transversalitätsbedingung in $t=t_1$ erhält die Pontrjagin-Funktion $H_0^{\mathcal{I}}$ die Form
\begin{eqnarray*}
H_0^{\mathcal{I}}\big(t,x,u,p(\cdot)\big) 
&=& -\int_t^{t_1} \langle \varphi(t,x,u), \dot{p}(\tau) \rangle \, d\tau +\langle \varphi(t,x,u), p(t_1)\rangle - f(t,x,u) \\
&=& \Big\langle \varphi(t,x,u) , p(t_1) -\int_t^{t_1} \dot{p}(\tau) \, d\tau \Big\rangle - f(t,x,u) \\
&=& \langle p(t), \varphi(t,x,u)\rangle - f(t,x,u) = H^{\mathcal{S}}\big(t,x,u,p(\cdot)\big) 
\end{eqnarray*} 
und Theorem \ref{SatzPMPeinfachIGL} liefert die notwendigen Bedingungen (\ref{PMPeinfach4})--(\ref{PMPeinfach6}). \\[2mm]
Der Randterm in der Pontrjagin-Funktion $H_0^{\mathcal{I}}$ war die Konsequenz aus der Umformung (\ref{StieltjesIGL}) des Stieltjes-Integral.
Für die adjungierte Funktion $p(\cdot)$ steht sie im Bezug zur Transversalitätsbedingung in $t=t_1$,
über die wir nun Anmerkungen treffen: \\ 
In der Literatur, z.\,B. bei Feichtinger \& Hartl \cite{Feichtinger} oder Kamien \& Schwartz \cite{Kamien},
wird die adjungierte Gleichung (\ref{PMPeinfachIGL4}) häufig in der Form
$$\lambda(t) = \int_t^{t_1} \varphi_x^T\big(\tau,t,x_*(t),u_*(t)\big) \lambda(\tau) \, d\tau + f_x\big(t,x_*(t),u_*(t)\big)$$
angegeben.
Bis auf den Randterm korrespondiert dabei die Funktion $\lambda(\cdot)$ mit der Ableitung $\dot{p}(\cdot)$ der Adjungierten $p(\cdot)$ in Gleichung (\ref{PMPeinfachIGL4}).
Diese Fehlinterpretation von $\lambda(\cdot)$ als Adjungierte anstelle der Ableitung der Adjungierten führt zu unzutreffenden Transversalitätsbedingungen.
Bei Kamien \& Schwartz \cite{Kamien} wird keine Transversalitätsbedingung für $\lambda(\cdot)$ angegeben.
Jedoch stellen Feichtinger \& Hartl \cite{Feichtinger} bei einer Aufgabe mit dem gemischten Zielfunktional (\ref{PMPeinfachIGL1})
die Transversalitätsbedingung $\lambda(t_1)= - S_1'\big(x_*(t_1)\big)$ auf.
Da aber $\lambda(\cdot)$ mit $\dot{p}(\cdot)$ korrespondiert,
hat ein fehlender Wiedergewinnungswert $S(x)$ die Transversalitätsbedingung $\lambda(t_1)=\dot{p}(t_1)=0$ zur Folge.
Im Abschnitt \ref{AbschnittWerbungIGL} über optimale Werbestrategien erweist sich die Transversalitätsbedingung für $\lambda(\cdot)$ als fehlerhaft. \\[2mm]
Die ``unverfälschte'' Pontrjagin-Funktion $H^{\mathcal{I}}$ gibt Anlass zu der folgenden ökonomischen Deutung:
\index{Pontrjaginsches Maximumprinzip!oekonomische@-- ökonomische Interpretation}
In Verallgemeinerung zu Abschnitt \ref{AbschnittDeutung},
wo $v \to H^{\mathcal{S}}\big(\tau,x_*(\tau),v,p(\tau)\big)$ die Profitrate aus direktem $f\big(\tau,x_*(\tau),v\big)$
und indirektem Gewinn $\big\langle p(\tau) , \varphi\big(\tau,x_*(\tau),v\big) \big\rangle$ bemisst,
enthält $H^{\mathcal{I}}$ den Erwartungswert
$\displaystyle \int_\tau^{t_1} \big\langle \varphi\big(t,\tau,x_*(\tau),v\big), dp(t) \big\rangle$
zum Gewicht $p(\cdot)$ über die gesamte künftige Entwicklung des indirekten Gewinnes.    
       \newpage
       \lhead[\thepage \hspace*{1mm} Beweis des Maximumprinzips]{ }  
       \subsection{Der Beweis des Maximumprinzips} \label{AbschnittPMPBeweiseinfachIGL}
Der folgende Beweis bezieht sich auf das Pontrjaginsche Maximumprinzip \ref{SatzPMPeinfachIGL2}.
Im Weiteren schreiben wir der Kürze halber
$$f_x[t]=f_x\big(t,x_*(t),u_*(t)\big), \qquad \varphi_x[t,s]=\varphi_x\big(t,s,x_*(s),u_*(s)\big).$$
In Theorem \ref{SatzPMPeinfachIGL2} bringen wir mit Hilfe der Beziehung (\ref{StieltjesIGL})
die adjungierte Gleichung (\ref{SatzPMPIGL1}) zur Sprung-Transversalitätsbedingung (\ref{SatzPMPIGL2}),
\begin{equation}\label{BeweisPMPIGL1}
\dot{p}(t) = \int_t^{t_1} \langle \varphi_x[\tau,t], dp(\tau) \rangle + f_x[t], \quad p(t_1^-)=-S'\big(x_*(t_1)\big), \; p(t_1)=0,
\end{equation}
in die Form
\begin{eqnarray}
\dot{p}(t) &=& \int_t^{t_1} \varphi^T_x[\tau,t] \dot{p}(\tau) \, d\tau
                + \varphi^T_x[t_1,t] \big( p(t_1)-p(t_1^-)\big) + f_x[t] \nonumber \\
           &=& \label{BeweisPMPIGL2} \int_t^{t_1} \varphi^T_x[\tau,t] \dot{p}(\tau) \, d\tau
              + \varphi^T_x[t_1,t] S'\big(x_*(t_1)\big) + f_x[t].
\end{eqnarray}
Wir setzen $q(\cdot) =\dot{p}(\cdot)$.
Dann erfüllt $q(\cdot)$ die adjungierte Gleichung (\ref{PMPeinfachIGL4}) in Theorem \ref{SatzPMPeinfachIGL}:
$$q(t)= \int_t^{t_1} \varphi^T_x[\tau,t] q(\tau) \, d\tau + \varphi^T_x[t_1,t]  S'\big(x_*(t_1)\big) + f_x[t].$$
Da s"amtliche Abbildungen stetig und stetig differenzierbar bezüglich $x$ sind, und $u_*(\cdot)$ stückweise stetig ist,
existiert nach Lemma \ref{LemmaDGL2} eine eindeutige L"osung $q(\cdot) \in PC([t_0,t_1],U)$ der adjungierten Gleichung (\ref{PMPeinfachIGL4})
zur Transversalitätsbedingung (\ref{PMPeinfachIGL5}).
Daher erfüllt die Funktion $p(\cdot) \in PC_1([t_0,t_1),\R^n)$ mit $\dot{p}(\cdot) = q(\cdot)$ und $p(t_1^-)=-S'\big(x_*(t_1)\big)$, $p(t_1)=0$,
die adjungierte Gleichung (\ref{SatzPMPIGL1}) zur Sprung-Transversalitätsbedingung (\ref{SatzPMPIGL2}). \\[2mm]
Genauso wie im Abschnitt \ref{AbschnittPMPBeweiseinfach} betrachten wir die einfache Nadelvariation\index{Nadelvariation, einfache}
$$u(t;v,\tau,\lambda) = u_{\lambda}(t) = 
  \left\{ \begin{array}{ll}
          u_*(t) & \mbox{ f"ur } t \not\in [\tau-\lambda,\tau), \\
          v      & \mbox{ f"ur } t     \in [\tau-\lambda,\tau), 
          \end{array} \right.$$
und es bezeichne $x_\lambda(\cdot)$, $x_\lambda(t)=x(t;v,\tau,\lambda)$, die eindeutige L"osung der Gleichung
$$x(t) = x_0 + \int_{t_0}^t \varphi\big(t,s,x(s),u_\lambda(s)\big), \qquad t \in [t_0,t_1].$$
Im Folgenden untersuchen wir f"ur $t \geq \tau$ den Grenzwert
$$y(t)=\lim_{\lambda \to 0^+}\frac{x_{\lambda}(t) - x_*(t)}{\lambda}.$$
Nach den Sätzen \ref{SatzEElokal}--\ref{SatzDGLDifferenzierbarkeit} "uber die Stetigkeit und Differenzierbarkeit
der L"osung einer Integralgleichung in Abh"angigkeit von den Anfangsdaten ergibt sich für $y(t)$
im Grenzwert $\lambda \to 0^+$ die lineare Integralgleichung
\begin{eqnarray}
y(t) &=&  \underbrace{\varphi\big(t,\tau,x_*(\tau),v\big) - \varphi\big(t,\tau,x_*(\tau),u_*(\tau)\big)}_{=y(\tau;t)}
          + \int_{\tau}^{t} \varphi_x[t,s] y(s) \, ds \nonumber \\
     &=& \label{BeweisPMPeinfachIGL1} y(\tau;t) + \int_{\tau}^{t} \varphi_x[t,s]\,y(s) \, ds.
\end{eqnarray} 
Mit der adjungierten Gleichung (\ref{BeweisPMPIGL1}) in der Form (\ref{BeweisPMPIGL2}) erhalten wir f"ur $t \geq \tau$:
\begin{eqnarray}
  \int_{\tau}^{t_1} \langle f_x[t] , y(t) \rangle \, dt
&=& \int_{\tau}^{t_1} \langle \dot{p}(t) , y(t) \rangle \, dt 
    - \int_{\tau}^{t_1} \Big\langle \int_t^{t_1} \varphi^T_x[s,t] \dot{p}(s) \, ds, y(t) \Big\rangle \, dt\nonumber \\
&& \label{BeweisPMPeinfachIGL2} + \int_{\tau}^{t_1} \langle \varphi^T_x[t_1,t] p(t_1^-) , y(t) \rangle \, dt.
\end{eqnarray}
Im ersten Term in (\ref{BeweisPMPeinfachIGL2}) liefert (\ref{BeweisPMPeinfachIGL1}) den Zusammenhang
\begin{eqnarray*}
    \int_{\tau}^{t_1} \langle \dot{p}(t) , y(t) \rangle \, dt
&=& \int_{\tau}^{t_1} \Big\langle \dot{p}(t), y(\tau;t) + \int_{\tau}^{t} \varphi_x[t,s]y(s) \, ds \Big\rangle \, dt \\
&=& \int_{\tau}^{t_1} \langle \dot{p}(t) , y(\tau;t) \rangle \, dt + \int_{\tau}^{t_1} \bigg( \int_{\tau}^{t} \langle \dot{p}(t), \varphi_x[t,s]y(s) \rangle \, ds \bigg) \, dt.
\end{eqnarray*}
Im zweiten Term in (\ref{BeweisPMPeinfachIGL2}) ergibt sich durch Vertauschen der Integrationsreihenfolge und der Variablenbezeichnung von $t$ und $s$
\begin{eqnarray*}
    \int_{\tau}^{t_1} \Big\langle \int_t^{t_1} \varphi^T_x[s,t] \dot{p}(s) \, ds, y(t) \Big\rangle \, dt
&=& \int_{\tau}^{t_1} \bigg( \int_{\tau}^s \langle \varphi^T_x[s,t] \dot{p}(s), y(t) \rangle \, dt \bigg) \, ds \\
&=& \int_{\tau}^{t_1} \bigg( \int_{\tau}^t \langle \dot{p}(t) , \varphi_x[t,s]  y(s) \rangle \, ds \bigg) \, dt.
\end{eqnarray*}
Damit lassen sich die ersten beiden Terme zusammenfassen und führen zu
$$\int_{\tau}^{t_1} \langle \dot{p}(t) , y(t) \rangle \, dt 
    - \int_{\tau}^{t_1} \Big\langle \int_t^{t_1} \varphi^T_x[s,t] \dot{p}(s) \, ds, y(t) \Big\rangle \, dt
   = \int_{\tau}^{t_1} \langle \dot{p}(t) , y(\tau;t) \rangle \, dt.$$
In der Gleichung (\ref{BeweisPMPeinfachIGL1}) gilt für $t=t_1$
$$y(t_1)-y(\tau;t_1) = \int_{\tau}^{t_1} \varphi_x\big[t_1,t] y(t) \, dt$$
und es folgt für den dritten Term in (\ref{BeweisPMPeinfachIGL2}) die Beziehung
$$\int_{\tau}^{t_1} \langle p(t_1^-) , \varphi_x[t_1,t]y(t) \rangle \, dt
  = \Big\langle p(t_1^-), \int_{\tau}^{t_1} \varphi_x[t_1,t]y(t) \, dt \Big\rangle = \langle p(t_1^-) ,y(t_1)-y(\tau;t_1) \rangle.$$
Diese Umformungen ergeben für (\ref{BeweisPMPeinfachIGL2}) die Darstellung
$$\int_{\tau}^{t_1} \langle f_x[t] , y(t) \rangle \, dt
   = \int_{\tau}^{t_1} \langle \dot{p}(t) , y(\tau;t) \rangle \, dt + \langle p(t_1^-) ,y(t_1)-y(\tau;t_1) \rangle.$$
In dieser Gleichung bringen wir $\langle p(t_1^-) ,y(t_1) \rangle$ auf die linke Seite.
Beachten wir nun die Transversalitätsbedingung $p(t_1^-)=-S'\big(x_*(t_1)\big)$, dann folgt
\begin{equation} \label{BeweisPMPeinfachIGL3}
\int_{\tau}^{t_1} \langle f_x[t] , y(t) \rangle \, dt + \langle S'\big(x_*(t_1)\big) ,y(t_1) \rangle
   = \int_{\tau}^{t_1} \langle \dot{p}(t) , y(\tau;t) \rangle \, dt - \langle p(t_1^-) ,y(\tau;t_1) \rangle.
\end{equation}
Da $\big(x_*(\cdot),u_*(\cdot)\big)$ ein starkes lokales Minimum ist, gilt (vgl. Abschnitt \ref{AbschnittPMPBeweiseinfach})
\begin{eqnarray*}
0 &\leq& \lim_{\lambda \to 0^+} \frac{J\big(x_\lambda(\cdot),u_\lambda(\cdot)\big)- J\big(x_*(\cdot),u_*(\cdot)\big)}{\lambda} \\
  &=   & f\big(\tau,x_*(\tau),v\big) - f\big(\tau,x_*(\tau),u_*(\tau)\big)
         + \int_{\tau}^{t_1} \big\langle f_x[t],y(t) \big\rangle \, dt + \langle S'\big(x_*(t_1)\big) ,y(t_1) \rangle.
\end{eqnarray*}
Mit (\ref{BeweisPMPeinfachIGL3}) und der Setzung von $y(\tau;t)$ in (\ref{BeweisPMPeinfachIGL1}) ergibt sich hieraus die Ungleichung
\begin{eqnarray}
0 &\leq& f\big(\tau,x_*(\tau),v\big) - f\big(\tau,x_*(\tau),u_*(\tau)\big) \nonumber\\
  &    & + \int_{\tau}^{t_1} \big\langle \dot{p}(t) , \varphi\big(t,\tau,x_*(\tau),v\big) - \varphi\big(t,\tau,x_*(\tau),u_*(\tau)\big) \big\rangle \, dt \nonumber \\
  &    & \label{BeweisPMPeinfachIGL4}
         - \big\langle p(t_1^-) , \varphi\big(t_1,\tau,x_*(\tau),v\big) - \varphi\big(t_1,\tau,x_*(\tau),u_*(\tau)\big) \big\rangle.
\end{eqnarray}
Unter Verwendung von $p(t_1)=0$ und der Beziehung (\ref{StieltjesIGL}) gelten
\begin{eqnarray*}
    \int_{\tau}^{t_1} \big\langle \varphi\big(t,...,v \big),dp(t) \big\rangle
&=& \int_{\tau}^{t_1} \big\langle \varphi\big(t,...,v \big),\dot{p}(t) \big\rangle
         - \big\langle \varphi\big(t_1,...,v\big),p(t_1^-) \big\rangle \\
    \int_{\tau}^{t_1} \big\langle \varphi\big(t,...,u_*(\tau) \big),dp(t) \big\rangle
&=& \int_{\tau}^{t_1} \big\langle \varphi\big(t,...,u_*(\tau) \big),\dot{p}(t) \big\rangle
         - \big\langle \varphi\big(t_1...,u_*(\tau)\big),p(t_1^-) \big\rangle.
\end{eqnarray*}
Damit ist (\ref{BeweisPMPeinfachIGL4}) gleichbedeutend mit der Ungleichung
\begin{eqnarray*}
\lefteqn{H^{\mathcal{I}}\big(\tau,x_*(\tau),v,p(\cdot)\big)=
         -\int_\tau^{t_1} \big\langle \varphi\big(t,\tau,x_*(\tau),v\big), dp(t) \big\rangle - f\big(\tau,x_*(\tau),v\big)} \\  
&\leq& -\int_\tau^{t_1} \big\langle \varphi\big(t,\tau,x_*(\tau),u_*(\tau)\big), dp(t) \big\rangle - f\big(\tau,x_*(\tau),u_*(\tau)\big)
       = H^{\mathcal{I}}\big(\tau,x_*(\tau),u_*(\tau),p(\cdot)\big).
\end{eqnarray*}
Da $\tau$ ein beliebiger Stetigkeitspunkt von $u_*(\cdot)$ und $v$ ein beliebiger Punkt der Menge $U$ war,
zeigt diese Ungleichung die Gültigkeit der Maximumbedingung (\ref{PMPeinfachIGL6}). \hfill $\blacksquare$  
       \newpage
       \lhead[\thepage \hspace*{1mm} Arrow-Bedingungen]{ }  
       \subsection{Hinreichende Bedingungen nach Arrow} \label{AbschnittArrowIGL}
Wir gehen nun auf hinreichende Bedingungen nach Arrow\index{hinreichende Bedingungen nach Arrow!Integral@-- Integralgleichungen}
für die Aufgabe (\ref{PMPeinfachIGL1})--(\ref{PMPeinfachIGL3}) ein.
Dabei geben wir zwei Varianten an, die sich auf die Darstellungen der Optimalitätsbedingungen der Theoreme
\ref{SatzPMPeinfachIGL2} und \ref{SatzPMPeinfachIGL} beziehen. \\[2mm]
Für $t \in [t_0,t_1]$ sei wieder $V_\gamma(t)=\{ x \in \R^n \,|\, \|x-x_*(t)\| \leq \gamma\}$. \\[2mm]
Die erste Variante der hinreichenden Bedingungen bezieht sich auf Theorem \ref{SatzPMPeinfachIGL2}.
Hier lautet die zugehörige ``unverfälschte'' Pontrjagin-Funktion
$$H^{\mathcal{I}}(t,x,u,p(\cdot)\big) = -\int_t^{t_1} \langle \varphi(\tau,t,x,u), dp(\tau) \rangle - f(t,x,u)$$
und es ergibt sich die Hamilton-Funktion
$$\mathscr{H}^{\mathcal{I}}\big(t,x,p(\cdot)\big)
  = \sup_{u \in U} \bigg[-\int_t^{t_1} \langle \varphi(\tau,t,x,u), dp(\tau) \rangle - f(t,x,u)\bigg].$$

\begin{theorem} \label{SatzHBPMPIG2}
In der Aufgabe (\ref{PMPeinfachIGL1})--(\ref{PMPeinfachIGL3}) sei
$\big(x_*(\cdot),u_*(\cdot)\big) \in \mathscr{A}^{\mathcal{I}}_{\rm adm} \cap \mathscr{A}^{\mathcal{I}}_{\rm Lip}$
und es sei $p(\cdot) \in PC_1([t_0,t_1),\R^n)$. Ferner gelte:
\begin{enumerate}
\item[(a)] Das Tripel $\big(x_*(\cdot),u_*(\cdot),p(\cdot)\big)$
           erf"ullt (\ref{SatzPMPIGL1})--(\ref{SatzPMPIGL3}) in Theorem \ref{SatzPMPeinfachIGL2}.        
\item[(b)] F"ur jedes $t \in [t_0,t_1]$ ist die Funktion $\mathscr{H}^{\mathcal{I}}\big(t,x,p(t)\big)$ konkav in $x$ auf $V_\gamma(t)$.
\item[(c)] Die Abbildung $S:\R^n\to \R$ im Zielfunktional (\ref{PMPeinfachIGL1}) ist konvex in $x$ auf $V_\gamma(t_1)$.
\end{enumerate}
Dann ist $\big(x_*(\cdot),u_*(\cdot)\big)$ ein starkes lokales Minimum der Aufgabe (\ref{PMPeinfachIGL1})--(\ref{PMPeinfachIGL3}).
\end{theorem}

Die zweite Variante der hinreichenden Bedingungen beruht auf Theorem \ref{SatzPMPeinfachIGL}.
In diesem Fall hat die Pontrjagin-Funktion die Form
$$H_0^{\mathcal{I}}(t,x,u,p(\cdot)\big) = -\int_t^{t_1} \langle \varphi(\tau,t,x,u), \dot{p}(\tau) \rangle \, d\tau
                                    +\langle \varphi(t_1,t,x,u), p(t_1)\rangle - f(t,x,u).$$
und die Hamilton-Funktion $\mathscr{H}_0^{\mathcal{I}}$ erhält die Gestalt
$$\mathscr{H}_0^{\mathcal{I}}\big(t,x,p(\cdot)\big)
  = \sup_{u \in U} \bigg[-\int_t^{t_1} \langle \varphi(\tau,t,x,u), \dot{p}(\tau) \rangle \, d\tau
                                    +\langle \varphi(t_1,t,x,u), p(t_1)\rangle - f(t,x,u)\bigg].$$

\begin{theorem} \label{SatzHBPMPIG1}
In der Aufgabe (\ref{PMPeinfachIGL1})--(\ref{PMPeinfachIGL3}) sei
$\big(x_*(\cdot),u_*(\cdot)\big) \in \mathscr{A}^{\mathcal{I}}_{\rm adm} \cap \mathscr{A}^{\mathcal{I}}_{\rm Lip}$
und es sei $p(\cdot) \in PC_1([t_0,t_1],\R^n)$. Ferner gelte:
\begin{enumerate}
\item[(a)] Das Tripel $\big(x_*(\cdot),u_*(\cdot),p(\cdot)\big)$
           erf"ullt (\ref{PMPeinfachIGL4})--(\ref{PMPeinfachIGL6}) in Theorem \ref{SatzPMPeinfachIGL}.        
\item[(b)] F"ur jedes $t \in [t_0,t_1]$ ist die Funktion $\mathscr{H}_0^{\mathcal{I}}\big(t,x,p(t)\big)$ konkav in $x$ auf $V_\gamma(t)$.
\item[(c)] Die Abbildung $S:\R^n\to \R$ im Zielfunktional (\ref{PMPeinfachIGL1}) ist konvex in $x$ auf $V_\gamma(t_1)$.
\end{enumerate}
Dann ist $\big(x_*(\cdot),u_*(\cdot)\big)$ ein starkes lokales Minimum der Aufgabe (\ref{PMPeinfachIGL1})--(\ref{PMPeinfachIGL3}).
\end{theorem}

Mit den Anmerkungen in Abschnitt \ref{AbschnittPMPeinfachIGL},
die mittels der Beziehung \ref{StieltjesIGL} die Verbindungen zwischen den Pontrjagin-Funktionen
$\mathscr{H}^{\mathcal{I}}$ und $\mathscr{H}_0^{\mathcal{I}}$ herstellen,
ist es ausreichend den Nachweis der hinreichenden Bedingungen für das Theorem \ref{SatzHBPMPIG2} zu führen. \\[2mm]
{\bf Beweis}  Es sei $t \in [t_0,t_1]$ so gew"ahlt,
dass die Maximumbedingung (\ref{SatzPMPIGL3}) zu diesem Zeitpunkt erf"ullt ist.
Mit den gleichen Argumenten wie im Abschnitt \ref{AbschnittHBPMP} folgt aus der Konkavität der
Hamilton-Funktion $\mathscr{H}^{\mathcal{I}}$ die Gültigkeit der Ungleichung
\begin{eqnarray*}
       \langle a(t),x-x_*(t)\rangle
&\geq& \mathscr{H}^{\mathcal{I}}\big(t,x,q(\cdot)\big) - \mathscr{H}^{\mathcal{I}}\big(t,x_*(t),q(\cdot)\big) \\
&\geq& \int_t^{t_1} \langle \varphi\big(\tau,t,x_*(t),u_*(t)\big),dq(\tau) \rangle + f\big(t,x_*(t),u_*(t)\big) \\
&    & - \int_t^{t_1} \langle \varphi\big(\tau,t,x,u_*(t)\big),dq(\tau) \rangle - f\big(t,x,u_*(t)\big)
\end{eqnarray*}
f"ur alle $x \in V_\gamma(t)$.
Weiterhin besitzt die Funktion
\begin{eqnarray*}
\Phi(x) &=& \int_t^{t_1} \langle \varphi\big(\tau,t,x_*(t),u_*(t)\big),dq(\tau) \rangle + f\big(t,x_*(t),u_*(t)\big) \\
        & & -\int_t^{t_1} \langle \varphi\big(\tau,t,x,u_*(t)\big),dq(\tau) \rangle - f\big(t,x,u_*(t)\big) - \langle a(t),x-x_*(t)\rangle
\end{eqnarray*}
in dem inneren Punkt $x_*(t)$ ein globales Maximum über $V_\gamma(t)$ und deswegen gilt
$$0=\Phi'(x_*(t)) = - \int_t^{t_1} \langle \varphi_x\big(\tau,t,x_*(t),u_*(t)\big),dp(\tau) \rangle - f_x\big(t,x_*(t),u_*(t)\big) - a(t).$$
Nach der adjungierten Gleichung (\ref{SatzPMPIGL1}) stimmt $-a(t)$ mit $\dot{q}(t)$ f"ur fast alle $t \in [t_0,t_1]$ überein
und es gilt damit über der Menge $V_\gamma$ die Ungleichung
\begin{equation} \label{HBPMPIG1}
\langle \dot{q}(t),x-x_*(t)\rangle \leq \mathscr{H}\big(t,x_*(t),q(\cdot)\big)- \mathscr{H}\big(t,x,q(\cdot)\big).
\end{equation}
Ferner ist für einen zulässigen Steuerungsprozess $\big(x(\cdot),u(\cdot)\big)$ nach (\ref{PMPeinfachIGL2}):
\begin{equation} \label{HBPMPIG2}
x(t)-x_*(t) = x(t_0)-x_*(t_0) +\int_{t_0}^t \Big[\varphi\big(t,s,x(s),u(s)\big)-\varphi\big(t,s,x_*(s),u_*(s)\big)\Big] \, ds.
\end{equation}
Es sei $\big(x(\cdot),u(\cdot)\big) \in  \mathscr{A}^{\mathcal{I}}_{\rm adm}$ mit $\|x(\cdot)-x_*(\cdot)\|_\infty < \gamma$.
Dann gilt
\begin{eqnarray*}
\lefteqn{J\big(x(\cdot),u(\cdot)\big)- J\big(x_*(\cdot),u_*(\cdot)\big)
\geq \int_{t_0}^{t_1} \big[\mathscr{H}^{\mathcal{I}}\big(t,x_*(t),q(\cdot)\big)-\mathscr{H}^{\mathcal{I}}\big(t,x(t),q(\cdot)\big)\big] \, dt} \\
& & \hspace*{-5mm} + \int_{t_0}^{t_1} \bigg[\int_t^{t_1}
    \big\langle \varphi\big(\tau,t,x_*(t),u_*(t)\big) - \varphi\big(\tau,t,x(t),u(t)\big),dq(\tau) \big\rangle \bigg]\, dt
    + S\big(x(t_1)\big)-S\big(x_*(t_1)\big).
\end{eqnarray*}
Wir ändern im zweiten Ausdruck die Integrationsreihenfolge:
\begin{eqnarray*}
\lefteqn{\int_{t_0}^{t_1} \bigg[\int_t^{t_1}
     \big\langle \varphi\big(\tau,t,x_*(t),u_*(t)\big) - \varphi\big(\tau,t,x(t),u(t)\big),dq(\tau) \big\rangle
     \bigg]\, dt} \\
&=& \int_{t_0}^{t_1} \bigg[\int_{t_0}^t \big[ \varphi\big(t,s,x_*(s),u_*(s)\big) - \varphi\big(t,s,x(s),u(s)\big) \big] \, ds \bigg]^T \, dq(t).
\end{eqnarray*}
Ferner folgt aus der Konvexität der Abbildung $S$ für alle $x \in V_\gamma(t_1)$ die Ungleichung
$$S(x)-S\big(x_*(t_1)\big) \geq \big\langle S'\big(x_*(t_1)\big), x-x_*(t_1) \big\rangle.$$
Unter Verwendung von (\ref{HBPMPIG1}) und (\ref{HBPMPIG2}) erhateln wir
\begin{eqnarray}
    \lefteqn{J\big(x(\cdot),u(\cdot)\big)- J\big(x_*(\cdot),u_*(\cdot)\big)
             \geq \int_{t_0}^{t_1} \langle x(t)-x_*(t),\dot{q}(t)\rangle \, dt} \nonumber \\
&& \label{HBPMPIG3} \hspace*{-5mm} + \int_{t_0}^{t_1} \langle x_*(t)-x(t)-x_*(t_0)+x(t_0), dq(t) \rangle
      + \big\langle S'\big(x_*(t_1)\big), x(t_1)-x_*(t_1) \big\rangle.
\end{eqnarray}
Nach Theorem \ref{SatzPMPeinfachIGL2} ist $q(\cdot)$ stückweise stetig differenzierbar und erfüllt in $t=t_1$
die Sprung-Transversalitätsbedingungen (\ref{SatzPMPIGL2}).
Demnach ergeben sich mit Formel (\ref{StieltjesIGL}):
\begin{eqnarray*}
    \int_{t_0}^{t_1} \langle x_*(t)-x(t), dq(t) \rangle
&=& \int_{t_0}^{t_1} \langle x_*(t)-x(t),\dot{q}(t)\rangle \, dt + \langle q(t_1)-q(t_1^-),x_*(t_1)-x(t_1) \rangle \\
&=& \int_{t_0}^{t_1} \langle x_*(t)-x(t), \dot{q}(t) \rangle \, dt + \big\langle S'\big(x_*(t_1)\big),x_*(t_1)-x(t_1) \big\rangle, \\
\int_{t_0}^{t_1} \langle x(t_0) - x_*(t_0), dq(t) \rangle
&=& \langle x(t_0)-x_*(t_0), [q(t_1^-)-q(t_0)]+[q(t_1)-q(t_1^-)] \rangle \\
&=& - \langle q(t_0),x(t_0)-x_*(t_0) \rangle.
\end{eqnarray*}
Unter Beachtung dieser Gleichungen führt die Beziehung (\ref{HBPMPIG3}) zu der Ungleichung
\begin{eqnarray*}
    \lefteqn{J\big(x(\cdot),u(\cdot)\big)- J\big(x_*(\cdot),u_*(\cdot)\big)
             \geq \int_{t_0}^{t_1} \langle x(t)-x_*(t),\dot{q}(t) \rangle \, dt
                  + \int_{t_0}^{t_1} \langle x_*(t)-x(t), \dot{q}(t) \rangle \, dt} \nonumber \\
&& \hspace*{-2mm} - \langle q(t_0),x(t_0)-x_*(t_0) \rangle+ \big\langle S'\big(x_*(t_1)\big),x_*(t_1)-x(t_1) \big\rangle  
   + \big\langle S'\big(x_*(t_1)\big), x(t_1)-x_*(t_1) \big\rangle.
\end{eqnarray*}
Da ferner der Anfangswert $x(t_0)=x_0$ in (\ref{PMPeinfachIGL2}) fest vorgegeben ist, ergibt sich
$$J\big(x(\cdot),u(\cdot)\big)- J\big(x_*(\cdot),u_*(\cdot)\big) \geq - \langle q(t_0),x(t_0)-x_*(t_0) \rangle =0$$
für alle
$\big(x(\cdot),u(\cdot)\big) \in \mathscr{A}^{\mathcal{I}}_{\rm adm}$ mit $\|x(\cdot)-x_*(\cdot)\|_\infty < \gamma$.
Somit ist $\big(x_*(\cdot),u_*(\cdot)\big)$ ein starkes lokales Minimum der Aufgabe (\ref{PMPeinfachIGL1})--(\ref{PMPeinfachIGL3}). \hfill $\blacksquare$   
       \newpage
       \lhead[\thepage \hspace*{1mm} Optimale Werbestrategien]{ }  
       \subsection{Optimale Werbestrategien} \label{AbschnittWerbungIGL}
Werbekampagnen zeichnen sich dadurch aus,
dass sie den Absatz nicht unmittelbar,
sondern mit einer zeitlichen Verzögerung beeinflussen.
In Anlehnung an Feichtinger \& Hartl \cite{Feichtinger} betrachten wir das Modell (vgl. Abschnitt \ref{AbschnittWerbungDGL}):
\begin{eqnarray}
&& \label{WerbungIGL1} J\big(x(\cdot),u(\cdot)\big) = \int_0^T e^{-\varrho t} \big[\pi x(t)-u(t)\big] \, dt
                       + e^{-\varrho T} S\big(x(T)\big) \to \sup, \\
&& \label{WerbungIGL2} x(t)=x_0+\int_0^t e^{-\delta(t-s)} \big[\mu\big(u(s)\big)-\alpha x(s)\big] \, ds, \quad t \in [0,T], \\
&& \label{WerbungIGL3} x_0 \in [0,1],\quad u(t) \geq 0, \quad \varrho >0, \quad \alpha, \delta \geq 0.
\end{eqnarray}
Dabei sei $\mu$ über $[0,\infty)$ stetig, differenzierbar und streng monoton wachsend mit
$$\mu(0)=0, \quad \lim_{u \to \infty} \mu(u)=\gamma >0, \quad \lim_{u \to 0+} \mu'(u)=\infty, \quad \lim_{u \to \infty} \mu'(u)=0.$$
Für diese Aufgabe stellen wir die Optimalitätsbedingungen (\ref{PMPeinfachIGL4})--(\ref{PMPeinfachIGL6}) von Theorem \ref{SatzPMPeinfachIGL} auf.
Da die Zustandsvariable $x$ linear in die Aufgabe einfließt, sind die Bedingungen des Maximumprinzips hinreichend.
Die Pontrjagin-Funktion besitzt die Form
\begin{eqnarray}
    H^{\mathcal{I}}(t,x,u,q(\cdot)\big)
&=& -\int_t^{T} e^{-\delta(\tau-t)} [\mu(u)-\alpha x] \dot{q}(\tau) \, d\tau \nonumber \\
& & \label{WerbungIGL4} \hspace*{1cm} + e^{-\delta(T-t)} [\mu(u)-\alpha x] q(T) + e^{-\varrho t} [\pi x-u].
\end{eqnarray}
Mit (\ref{WerbungIGL4}) erhalten wir für die adjungierte Gleichung
\begin{equation}\label{WerbungIGL5}
\dot{q}(s) = - \int_t^{T} \alpha  e^{-\delta(\tau-t)} \dot{q}(\tau)  \, d\tau + \alpha e^{-\delta(T-t)} q(T)+\pi e^{-\varrho t},
\quad q(T)=e^{-\varrho T} S'\big(x_*(T)\big).
\end{equation}
Ferner führt die Maximumbedingung
$$H^{\mathcal{I}}\big(t,x_*(t),u_*(t),q(\cdot)\big) = \max_{u \in U} H^{\mathcal{I}}\big(t,x_*(t),u,q(\cdot)\big)$$
zu der Beziehung
$$\max_{u \geq 0} \bigg[-\int_t^{T} e^{-\delta(\tau-t)} \mu(u) \dot{q}(\tau) \, d\tau + e^{-\delta(T-t)} \mu(u) q(T) - e^{-\varrho t} u\bigg].$$
Damit gelangen wir zu der Gleichung
\begin{equation}\label{WerbungIGL6}
\bigg[ -\int_t^{T} e^{-\delta(\tau-t)} \dot{q}(\tau) \, d\tau + e^{-\delta(T-t)} q(T)\bigg]\cdot \mu'\big(u_*(t)\big)  - e^{-\varrho t}=0.
\end{equation}
Zur Auswertung der Optimalitätsbedingungen (\ref{WerbungIGL5}) und (\ref{WerbungIGL6}) stellen wir den Bezug zur
Aufgabe (\ref{WerbungDGL1})--(\ref{WerbungDGL3}) her:
\newpage 
Die Differentation der Integralgleichung (\ref{WerbungIGL2}) nach $t$ ergibt
\begin{eqnarray*}
\dot{x}(t) &=& -\delta \int_0^t e^{-\delta(t-s)} \big[\mu\big(u(s)\big)-\alpha x(s)\big] \, ds
               + \big[\mu\big(u(t)\big)-\alpha x(t)\big] \nonumber \\
           &=&  -\delta [x(t)-x_0]+ \big[\mu\big(u(t)\big)-\alpha x(t)\big]
               = -(\alpha+\delta) x(t) + \delta x_0 + \mu\big(u(t)\big)
\end{eqnarray*}
und die Aufgabe (\ref{WerbungIGL1})--(\ref{WerbungIGL3}) ist äquivalent zu dem Modell (\ref{WerbungDGL1})--(\ref{WerbungDGL3})
über optimale Werbestrategien in Abschnitt \ref{AbschnittWerbungDGL}.
Wir nutzen im Weiteren die Lösungstruktur der Aufgabe (\ref{WerbungDGL1})--(\ref{WerbungDGL3}),
in welcher nach (\ref{WerbungDGL8})--(\ref{WerbungDGL10}) die adjungierte Gleichung
\begin{equation} \label{WerbungIGL7}
\dot{p}(t)=(\alpha+\delta)p(t)-\pi e^{-\varrho t}, \quad p(T)=e^{-\varrho T} S'\big(x_*(T)\big),
\end{equation}
die Lösung
\begin{equation} \label{WerbungIGL8}
p(t)=\bigg(e^{-\varrho T}S'\big(x_*(T)\big)-\frac{\pi e^{-\varrho T}}{\alpha + \delta + \varrho}\bigg)
       \cdot e^{(\alpha+\delta)(t-T)} +\frac{\pi e^{-\varrho t}}{\alpha + \delta + \varrho}
\end{equation}
besitzt und aus der Maximumbedingung die Beziehung
\begin{equation} \label{WerbungIGL9}
p(t)\cdot \mu'\big(u_*(t)\big)-e^{-\varrho t} =0
\end{equation}
resultiert. \\[2mm]
Der Vergleich von (\ref{WerbungIGL6}) und (\ref{WerbungIGL9}) ergibt den Zusammenhang
$$p(t)= -\int_t^{T} e^{-\delta(\tau-t)} \dot{q}(\tau) \, d\tau + e^{-\delta(T-t)} q(T), \quad
  q(T)=p(T)=e^{-\varrho T} S'\big(x_*(T)\big).$$
Durch Differentation der ersten Gleichung nach $t$ erhalten wir weiter
$$\dot{p}(t) = -\delta \cdot \int_t^T e^{-\delta(\tau-t)} \dot{q}(\tau) \, d\tau + \dot{q}(t) +\delta e^{-\delta(T-t)} q(T) = \delta p(t) + \dot{q}(t).$$
Es gilt damit $\dot{q}(t)=\dot{p}(t)-\delta p(t)=\alpha p(t)-\pi e^{-\varrho t}$,
was der adjungierten Gleichung (\ref{WerbungIGL5}) mit Transversalitätsbedingung $q(T)=e^{-\varrho T} S'\big(x_*(T)\big)$ entspricht.
Unter Verwendung von (\ref{WerbungIGL8}) liefert die Berechnung von
$$q(t) = q(T) + \int_T^t \big[\alpha p(s)-\pi e^{-\varrho s}\big]\,ds$$
die explizite Darstellung der Adjungierten $q(\cdot)$. \\[2mm]
Abschließend entnehmen wir dem Abschnitt \ref{AbschnittWerbungDGL} über die Steuerung $u_*(\cdot)$ die Ergebnisse
$$\mu'\big(u_*(t)\big) = \frac{1}{p(t) \cdot e^{\varrho t}}, \qquad \mu'\big(u_*(T)\big) =\frac{1}{S'\big(x_*(T)\big)}.$$
Damit ist unsere Untersuchung der Aufgabe (\ref{WerbungIGL1})--(\ref{WerbungIGL3}) beendet. \hfill $\square$   
\cleardoublepage


\newpage
\addcontentsline{toc}{part}{Anhang}
\begin{appendix}
\lhead[\thepage \hspace*{1mm} Funktionalanalytische Hilfsmittel]{}
\rhead[]{Funktionalanalytische Hilfsmittel \hspace*{1mm} \thepage}
    \section{Funktionalanalytische Hilfsmittel} \label{AnhangHilfsmittel}
\subsection{Banachräume und Räume stetiger Funktionen} \label{AbschnittNormRaum}
Es sei $X$ ein reeller Vektorraum.
\begin{definition}[Normierter Raum] \index{Norm}
Eine Abbildung $\|\cdot\| : X \to [0,\infty)$ heißt eine Norm,
falls sie für alle $x,\,y \in X$ und alle $\alpha \in \R$ folgende Eigenschaften besitzt: \\[1mm]
\begin{tabular}{lll}
(1) & Definitheit: & $\|x\|=0 \; \Rightarrow \; x=0$, \\[1mm]
(2) & absolute Homogenität: & $\|\alpha x\|=|\alpha|\, \|x\|$, \\[1mm]
(3) & Dreiecksungleichung: & $\|x+y\|\leq \|x\|+\|y\|$. \\[1mm]
\end{tabular}\\
Das Paar $(X,\|\cdot\|)$ bildet einen normierten Raum.
\end{definition}

\begin{definition}[Cauchyfolge]
Die Folge $\{x_k\}$ des normierten Raumes $X$ ist eine Cauchyfolge,
falls es zu jedem $\varepsilon >0$ eine Zahl $N(\varepsilon) \in \N$ gibt mit
$\|x_m-x_n\| \leq \varepsilon$ für alle $m,\,n \geq N(\varepsilon)$.
\end{definition}

\begin{definition}[Konvergente Folge]
Die Folge $\{x_n\}$ konvergiert in $X$ gegen $x \in X$,
falls es zu jedem $\varepsilon >0$ eine Zahl $N(\varepsilon) \in \N$ gibt mit
$\|x_n-x\| \leq \varepsilon$ für alle $n \geq N(\varepsilon)$.
\end{definition}

\begin{definition}[Banachraum]
Ein normierter Raum, in dem jede Cauchyfolge konvergiert, heißt vollständig.
Ein vollständiger normierter Raum heißt Banachraum.
\end{definition}

\begin{lemma} \label{LemmaBanachraum}
Ein abgeschlossener Unterraum  eines Banachraumes ist vollständig.
\end{lemma}

{\bf Beweis} Es seien $U$ der Unterraum des Banachraumes $X$ und $\{x_n\}$ eine Cauchyfolge in $U$.
Da $X$ vollständig ist, besitzt die Folge  $\{x_n\}$ einen Grenzwert $x \in X$.
Wegen der Abgeschlossenheit von $U$ muss $x$ in $U$ liegen. \hfill $\blacksquare$ \\[2mm]
Wir geben nun Beispiele von normierten Räumen und von Banachräumen an.
Dabei verwenden wir die Schreibweise $x(\cdot)$ statt $x$, um den Charakter als Funktionen hervorzuheben.
Ferner seien die Definitionsbereiche der Funktionsklassen stets nichtleere Mengen. \\[2mm]
(1.) Der Raum $B(I,\R)$ der über dem Intervall $I \subseteq \R$ beschr"ankten Funktionen \index{Raum, beschränkter Funktionen}
    ist versehen mit der Supremumsnorm 
    $\displaystyle \|x(\cdot)\|_\infty = \sup_{t \in I} \|x(t)\|$ ein Banachraum. \\[1mm]
           Bezüglich der Supremusnorm gilt $|x(t)| \leq \|x(\cdot)\|_\infty$ für alle $t \in I$.
           Es sei $\{x_n(\cdot)\}$ eine Cauchyfolge in $B(I,\R)$.
           Dann bildet für jedes $t \in I$ die Folge $\{x_n(t)\}$ eine Cauchyfolge in $\R$ und besitzt damit einen Grenzwert $x(t) \in \R$.
           Nun ist zu zeigen, dass die resultierende Funktion $x(\cdot)$ beschränkt und Grenzwert der Cauchyfolge $\{x_n(\cdot)\}$ ist.
           Es sei $\varepsilon >0$ gegeben. Dann existiert eine Zahl $N \in \N$ mit $\|x_m(\cdot)-x_n(\cdot)\|_\infty \leq \varepsilon$ für alle $m,\,n \geq N$.
           Für $t \in I$ gibt es eine Zahl $k \geq N$ mit $|x_k(t)-x(t)| \leq \varepsilon$. Damit ergibt sich
           $|x_n(t)-x(t)| \leq |x_n(t)-x_k(t)| + |x_k(t)-x(t)| \leq 2\varepsilon$
           für alle $n \geq N$. Da $t \in I$ beliebig war, gilt
           $\|x(\cdot)\|_\infty  \leq \|x_N(\cdot)\|_\infty + \|x(\cdot)-x_N(\cdot)\|_\infty \leq \|x_N(\cdot)\|_\infty + 2\varepsilon,$
           d.\,h. $x(\cdot)$ ist beschränkt. Ferner liefert die Argumentation $\|x_n(\cdot)-x(\cdot)\|_\infty \leq 2\varepsilon$ für alle $n \geq N$,
           d.\,h. $x(\cdot)$ ist Grenzwert der Folge $\{x_n(\cdot)\}$. \hfill $\square$ \\[2mm]
(2.) Der Raum $C_b(I,\R)$ der über dem Intervall $I \subseteq \R$ beschr"ankten stetigen Funktionen ist versehen mit der Supremumsnorm
    $\|\cdot\|_\infty$ ein Banachraum. \\[1mm]
           Nach Lemma \ref{LemmaBanachraum} ist die Abgeschlossenheit von $C_b(I,\R)$ im Raum $B(I,\R)$ zu zeigen:
           der Grenzwert einer gleichmäßig konvergenten Folge stetiger Funktionen liefert eine stetige Funktion.
           Es bezeichne $x(\cdot) \in B(I,\R)$ den Grenzwert der gleichmäßig konvergenten Folge $\{x_n(\cdot)\}$ aus $C_b(I,\R)$.
           Zu $\varepsilon >0$ wählen wir ein $N \in \N$ mit $\|x_N(\cdot)-x(\cdot)\|_\infty \leq \varepsilon$.
           Sei $t_0 \in I$. Wegen der Stetigkeit von $x_N(\cdot)$ gibt es ein $\delta >0$ mit $|x_N(t)-x_N(t_0)| \leq \varepsilon$
           für alle $t \in I$ mit $|t-t_0| \leq \delta$.
           Damit folgt $|x(t)-x(t_0)| \leq |x(t)-x_N(t)| + |x_N(t)-x_N(t_0)|+ |x_N(t_0)-x(t_0)| \leq 3 \varepsilon$
           für alle $t \in I$ mit $|t-t_0| \leq \delta$, also die Stetigkeit von $x(\cdot)$ in $t_0$.\hfill $\square$ \\[2mm]
(3.) Der Raum $C([a,b],\R)$ der stetigen Funktionen\index{Raum, beschränkter Funktionen!stetiger@-- stetiger Funktionen}
    über dem kompakten Intervall $[a,b]$ ist versehen mit der Supremumsnorm $\|\cdot\|_\infty$ ein Banachraum. \\[1mm]
           Nach dem Satz von Weierstraß ist jede stetige Funktion über einer kompakten Menge beschränkt.
           Also ergibt sich die Vollständigkeit wie im vorhergehenden Beispiel.\hfill $\square$ \\[2mm]
(4.) Der Raum $C_{\lim}([0,\infty),\R)$ der stetigen und im Unendlichen konvergenten Funktionen ist versehen mit der Supremumsnorm $\|\cdot\|_\infty$
    ein Banchraum. \\[1mm]
           Eine Funktion $x(\cdot):[0,\infty) \to \R$ konvergiert im Unendlichen gegen $a \in \R$,
           wenn es zu jedem $\varepsilon >0$ eine Zahl $T>0$ gibt mit $|x(t)-a|\leq \varepsilon$ für alle $t \geq T$. \\[1mm]
           Der erweiterte reelle Halbstrahl $[0,\infty]$ ist unter der stetigen und bijektiven Abbildung $x \to x/(1+x)$
           topologisch gleichwertig zu dem kompakten Intervall $[0,1]$.
           Damit ist $C_{\lim}([0,\infty),\R)$ topologisch äquivalent zu $C([0,\infty],\R)$ bzw. $C([0,1],\R)$,
           und nach vorherigem Beispiel versehen mit der Supremumsnorm $\|\cdot\|_\infty$ ein Banachraum. \hfill $\square$ \\[2mm]
(5.) Der Raum $PC(I,\R)$ der stückweise stetigen Funktionen \index{Raum, beschränkter Funktionen!ststetig@-- stückweise stetiger Funktionen}
    über dem Intervall $I \subseteq \R$ ist versehen mit der Supremumsnorm $\|\cdot\|_\infty$
           nicht vollständig. \\[1mm]
    Nach Definition \ref{DefinitionSTSTFunktion} besitzen die Elemente des Raumes $PC(I,\R)$ höchstens endlich viele Unstetigkeitsstellen.
     Wir betrachten über $[0,2]$ die Funktionen $x_n(\cdot)$ und $x(\cdot)$ mit 
           $$x_n(t)=\left\{\begin{array}{rl} 0, & t \in [0,\frac{1}{n}], \\[1mm]
                              \displaystyle\frac{1}{m}, & t \in [\frac{1}{m},\frac{1}{m-1}], \; m=2,...,n, \\[1mm]
                                             1, & t \in [1,2], \end{array}\right. \;
             x(t)=\left\{\begin{array}{rl} 0, & t = 0, \\[1mm]
                              \displaystyle\frac{1}{m}, & t \in [\frac{1}{m},\frac{1}{m-1}], \; m=2,... \\[1mm]
                                             1, & t \in [1,2], \end{array}\right.$$
           Die Funktionen $x_n(\cdot)$ sind stückweise konstant und besitzen $n$ Sprungstellen $s_m=\frac{1}{m}$, $m=1,...,n$,
           wohingegen $x(\cdot)$ abzählbar viele Sprungstellen aufweist.
           Für $m \leq n$ gilt $\|x_m(\cdot)-x_n(\cdot)\|_\infty=\frac{1}{m+1}$, sowie $\|x_n(\cdot)-x(\cdot)\|_\infty=\frac{1}{n+1}$.
           Daher ist $\{x_n(\cdot)\}$ eine Cauchyfolge in $PC([0,2],\R)$, die gegen $x(\cdot) \not\in PC([0,2],\R)$ konvergiert. \hfill $\square$ \\[2mm]
(6.) Der Raum $PC_1(I,\R)$ der über dem Intervall $I \subseteq \R$ stückweise stetig differenzierbaren Funktionen
    \index{Raum, beschränkter Funktionen!ststetigdiff@-- stückweise stetig differenzierbarer Funktionen} ist bezüglich
    der Norm $\|x(\cdot)\|=\|x(\cdot)\|_\infty + \|\dot{x}(\cdot)\|_\infty$ nicht vollständig. \\[2mm]
    Dazu betrachte man lediglich für die Ableitung die Folge in letztem Punkt. \hfill $\square$

\subsection{Fixpunktsätze}
\begin{satz}[Fixpunktsatz von Weissinger] \label{SatzWeissinger} 
Es sei $M$ eine nichtleere abgeschlossene Teilmenge des Banachraumes $X$,
ferner seien $a_1,a_2,...$ die Glieder einer konvergenten Reihe positiver Zahlen und
$A:M \to M$ eine Selbstabbildung von $M$ mit $\|A^n x-A^n y\| \leq a_n\|x-y\|$ f"ur alle $x,y \in M$ und $n \in \N$.
Dann besitzt $A$ genau einen Fixpunkt $x \in M$, d.\,h. einen Punkt $x \in M$ mit $Ax=x$.
Dieser Fixpunkt ist Grenzwert der Iterationsfolge $x_n=Ax_{n-1}, n=1,2,...,$ bei beliebigem Startwert $x_0 \in M$. \\[2mm]
Schlie"slich gilt die Fehlerabsch"atzung
$\displaystyle\quad \|x-x_n\| \leq\sum_{k=n}^\infty a_k \cdot \|x_1-x_0\|$.
\end{satz}

{\bf Beweis} Wir betrachten die Folge $x_n=Ax_{n-1}$.
Für diese gilt
$$\|x_{n+1}-x_n\|=\|A^n x_1-A^n x_0\| \leq a_n\|x_1-x_0\|$$ und folglich für $l \geq 1$
$$\|x_{n+l}-x_n\| \leq \|x_{n+l}-x_{n+l-1}\|+...+\|x_{n+1}-x_n\| \leq \sum_{k=n}^{n+l-1} a_k \cdot \|x_1-x_0\|.$$
Da die Reihe mit Folgengliedern $a_n$ konvergiert, wird die Summe für jedes $k \in \N$ kleiner als ein beliebig vorgegebenes $\varepsilon >0$ wird,
wenn $n$ nur hinreichend groß ist.
Dies zeigt, dass die Iterationsfolge $\{x_n\}$ eine Cauchyfolge bildet.
Wegen der Vollständigkeit von $X$ konvergiert diese Folge gegen ein $x \in X$, und $x$ gehört zur Menge $M$, da $M$ abgeschlossen ist.
Es ist $x$ der gesuchte Fixpunkt.
Denn einerseits ist $\lim x_n=x$ und andererseits folgt wegen $\|x_{n+1}-Ax\|=\|Ax_n-Ax\| \leq a_1\|x_n-x\|$, dass $\lim x_n=Ax$ gilt.
Daher muss $Ax=x$ sein.
Ist $y \in M$ ein weiterer Fixpunkt von $A$, so ist
$$\|x-y\|=\|Ax-Ay\|=...=\|A^nx-A^ny\| \leq a_n\|x-y\|$$ mit einem beliebigem $n \in \N$.
Da $\{a_n\}$ notwendig eine Nullfolge sein muss, ergibt sich $x=y$.
Führen wir abschließend in obiger Ungleichung den Grenzübergang $l \to \infty$ aus,
$$\lim_{l \to \infty} \|x_{n+l}-x_n\| = \|x-x_n\| \leq \sum_{k=n}^\infty a_k \cdot \|x_1-x_0\|,$$
so ergibt sich die behauptete Fehlerabschätzung.
  \hfill $\blacksquare$

\begin{satz}[Banachscher Fixpunktsatz]
Es sei $M$ eine nichtleere abgeschlossene Teilmenge des Banachraumes $X$ und
$A:M \to M$ eine kontraktive Selbstabbildung von $M$, d.\,h. $\|A x-A y\| \leq q\|x-y\|$ f"ur alle $x,y \in M$ mit einer Zahl $q \in (0,1)$.
Dann besitzt $A$ genau einen Fixpunkt $x \in M$.
\end{satz}

{\bf Beweis} Es gilt $\|A^n x-A^n y\| \leq q^n\|x-y\|$ und der Banachsche Fixpunktsatz erweist sich als
Spezialfall des Fixpunktsatzes von Weissinger mit $a_n=q^n$. \hfill $\blacksquare$


\newpage
\subsection{Stetige lineare Abbildungen}
Es seien $X,\,Y$ normierte Räume.
Dann bezeichnet $L(X,Y)$ die Menge der stetigen linearen Abbildungen\index{Abbildung, beschränkte!stetige@--, stetige} $T:X \to Y$.
Eine Abbildung $T \in L(X,Y)$ erfüllt eine der folgenden äquivalenten Bedingungen:

\begin{enumerate}
\item[(i)] Aus $x_n \to x$ folgt $Tx_n \to Tx$.
\item[(ii)] Zu jedem $x_0 \in X$ und jedem $\varepsilon>0$ existiert ein $\delta=\delta(x_0,\varepsilon) >0$ mit
           $$\|Tx-Tx_0\| \leq \varepsilon \quad \mbox{ f"ur alle } \|x-x_0\| \leq \delta.$$
\item[(iii)] Für alle offenen $O \subseteq Y$ ist das Urbild $T^{-1}(O)=\{x \in X \,|\, Tx \in O\}$ offen in $X$.
\end{enumerate}

Die Bedingung (i) beschreibt die Folgenstetigkeit, (ii) das $(\varepsilon,\Delta)-$Kriterium und
(iii) die Stetigkeit in der Topologie.

\begin{lemma} \label{LemmaOperatornorm}
Insbesondere sind f"ur lineare Abbildungen "aquivalent:
\begin{enumerate}
\item[(a)] $T$ ist stetig, d.\,h. $T$ ist in jedem Punkt $x \in X$ stetig.
\item[(b)] $T$ ist in $x=0$ stetig.
\item[(c)] $T$ ist beschr"ankt\index{Abbildung, beschränkte}, d.\,h. es existiert eine Zahl $K\geq 0$ mit
           $\|Tx\| \leq K \|x\|$ f"ur alle $x \in X$.
\end{enumerate}
\end{lemma}

{\bf Beweis} Die Beweiskette (c)$\Rightarrow$(a)$\Rightarrow$(b) ist offensichtlich:
Ist (c) erfüllt, so ergibt sich wegen $\|Tx-Ty\|=\|T(x-y)\| \leq M \|x-y\|$ für alle $x,y\in X$ die Stetigkeit von $T$,
und die Stetigkeit in (a) impliziert unmittelbar die Stetigkeit in jedem Punkt $x \in X$. \\[2mm]
Um den Implikationskreis zu schließen, zeigen wir (b)$\Rightarrow$(c):
Angenommen, (c) würde nicht gelten.
Dann gibt es eine Folge $\{x_n\}$ in $X$ mit $\|Tx_n\|>n$.
Wir setzen $\displaystyle y_n=\frac{1}{n}\frac{x_n}{\|x_n\|}$.
Dann gelten $\displaystyle \|y_n\|=\frac{1}{n}$ und $\displaystyle \|Ty_n\|=\frac{\|Tx_n\|}{n\|x_n\|}>1$.
Obwohl $\{y_n\}$ eine Nullfolge ist, gilt $Ty_n \not\to 0$, im Widerspruch zur Stetigkeit von $T$ in $x=0$. \hfill $\blacksquare$

\begin{definition}
Für $T \in L(X,Y)$ wird die kleinste Zahl $K \geq 0$,
mit der Lemma \ref{LemmaOperatornorm} erfüllt ist, mit $\|T\|$ bezeichnet:
$$\|T\| = \inf\{K \geq 0 \,|\,  \|Tx\| \leq K \|x\| \mbox{ f"ur alle } x \in X\}.$$
\end{definition}

\begin{satz}
Die Zahl $\|T\|$ definiert eine Norm auf dem Raum $L(X,Y)$ und es gelten
$$\|T\|= \sup_{x \not=0} \frac{\|Tx\|}{\|x\|}=  \sup_{\|x\|=1} \|Tx\| = \sup_{\|x\| \leq 1} \|Tx\|.$$
\end{satz}

\begin{definition}[Operatornorm]
Für $T \in L(X,Y)$ wird die Zahl $\displaystyle \|T\|=\sup_{\|x\| \leq 1} \|Tx\|$ als die Operatornorm von $T$ bezeichnet.
\end{definition}

\newpage
\subsection{Differenzierbarkeit in normierten Räumen, implizite Funktionen}
Es seien $X,Y$ normierte R"aume, $U \subseteq X$ offen und $F:U \to Y$ eine Abbildung.

\begin{definition}[Fr\'echet-differenzierbare Abbildung] \label{Frechet}
Existiert in $x_0 \in U$ der gleichmäßige Grenzwert
\begin{equation} \label{FrechetAbl}
\lim_{\lambda \to 0} \frac{F(x_0+\lambda x)-F(x_0)}{\lambda}=Tx
\end{equation}
mit einer stetigen linearen Abbildung $T: X \to Y$,
dann hei"st die Abbildung $F$ in $x_0$ Fr\'echet-differenzierbar\index{Abbildung, beschränkte!Fr\'echet@--, Fr\'echet-differenzierbare}
und $F'(x_0)=T$ die Fr\'echet-Ableitung von $F$ in $x_0$. 
\end{definition}

Der Grenzwert (\ref{FrechetAbl}) konvergiert gleichmäßg,
falls es zu jedem $\varepsilon >0$ ein $\lambda_0 >0$  gibt mit
$$\bigg\|\frac{F(x_0+\lambda x)-F(x_0)}{\lambda}-Tx\bigg\| \leq \varepsilon \quad\mbox{ f"ur alle }
  0<\lambda<\lambda_0, \; \|x\|\leq 1.$$
  
\begin{lemma}
Eine Abbildung $F$ ist in $x_0$ genau dann Fr\'echet-differenzierbar,
falls ein $T \in L(X,Y)$ und eine Abbildung $r:X \to Y$ mit
$$F(x_0+x)=F(x_0)+Tx+r(x), \qquad \lim_{\|x\| \to 0} \frac{\|r(x)\|}{\|x\|} =0.$$
\end{lemma}

\begin{definition}[Stetig differenzierbare Abbildung]
\begin{enumerate}
\item[(a)] Die Abbildung $F$ ist in $x_0$ stetig differenzierbar\index{Abbildung, beschränkte!stetig differenzierbare@--, stetig differenzierbare},
           wenn f"ur alle Punkte einer offenen Umgebung $U_\delta(x_0) \subseteq U$ des Punktes $x_0$ die Ableitung $F'(x)$ existiert und die
           Abbildung $x \to F'(x)$ bez"uglich der Operatornorm des Raumes $L(X,Y)$ in $x_0$ stetig ist.
\item[(b)] Die Abbildung $F$ ist auf der offenen Menge $U$ stetig differenzierbar,
           wenn f"ur alle $x \in U$ die Ableitung $F'(x)$ existiert und die Abbildung $x \to F'(x)$ bez"uglich der Operatornorm des Raumes $L(X,Y)$
           auf $U$ stetig ist.
\end{enumerate}
\end{definition}

In Folgenden bezeichnet ${\rm Im\,}T=\{y \in Y \,|\, y=Tx, \; x \in X\}$ die Bildmenge der Abbildung $T:X \to Y$.

\begin{definition}[Reguläre Abbildung]
Die Abbildung $F$ nennen wir im Punkt $x_0$ regul"ar\index{Abbildung, beschränkte!reguläre@--, reguläre},
wenn sie in diesem Punkt Fr\'echet-differenzierbar ist und ${\rm Im\,}F'(x_0)=Y$ gilt.
\end{definition}

\begin{satz}[Satz über implizite Funktionen] \label{SatzImpliziteFunktionen}
Es seien $X$, $Y$ und $Z$ Banachräume, $G$ eine Umgebung des Punktes $(x_0,y_0) \in X \times Y$ und $F:G \to Z$
eine stetig differenzierbare Abbildung.
Ferner sei $F(x_0,y_0)=0$ und es sei die partielle Ableitung $F_y(x_0,y_0):Y \to Z$ ein linearer Homöomorphismus. \\
Dann gibt es auf einer Umgebung $U(x_0)$ von $x_0$ eine stetig differenzierbare Abbildung $x \to y(x)$ nach $Y$ derart,
dass $y(x_0)=y_0$ erfüllt ist, sowie für jedes $x \in U(x_0)$ die Gleichung $F\big(x,y(x)\big)=0$ und die Beziehung
$y'(x)=-\big[F_y\big(x,y(x)\big)\big]^{-1} \circ F_x\big(x,y(x)\big)$ gelten.
\end{satz}

    \newpage
\lhead[\thepage \hspace*{1mm} Differential- und Integralgleichungen]{}
\rhead[]{Differential- und Integralgleichungen \hspace*{1mm} \thepage}
    \section{Differentialgleichungen, Volterrasche Integralgleichungen} \label{AnhangDGL}
Die Lehrbücher zu gewöhnlichen Differentialgleichungen behandeln vorrangig den Fall stetiger rechter Seiten. 
Da in Steuerungsproblemen die rechte Seite der Dynamik für einen Steuerungsprozess typischerweise nicht stetig,
sondern stückweise stetig in der Zeitvariablen ist, und somit die Rahmenbedingungen von der Standardtheorie abweichen,
stellen wir die benötigten Resultate ausführlich dar. \\
Die Ideen zu den grundlegenden Sätzen für Differentialgleichung mit stückweise stetigen rechten Seiten und für Volterrasche Integralgleichungen
sind eng verknüpft mit der bekannten Theorie zu gewöhnlichen Differentialgleichungen.
Deswegen halten wir uns eng an das Lehrbuch zu gew"ohnlichen Differentialgleichungen von Heuser \cite{HeuserGD} und
an der Monographie \cite{Ioffe} von Ioffe \& Tichomirov.

\subsection{Lineare Gleichungen}
Es sei $x(\cdot):\R \to \R^n$ eine Vektorfunktion und $A(\cdot,\cdot):\R \times \R \to \R^{n \times n}$ eine matrixwertige Abbildung.
Im Weiteren betrachten wir lineare Integralgleichungssysteme der Form
\begin{eqnarray}
\label{LinDGL1} x(t) &=& \zeta(t)+ \int_\tau^t \big[A(t,s)x(s)+a(s)\big] \, ds, \\
\label{LinDGL2} x(t) &=& \zeta(t)+ \int_\tau^t \big[A(s)x(s)+a(s)\big] \, ds
\end{eqnarray}
und das lineare Differentialgleichungssystem
\begin{equation}\label{LinDGL3}
\dot{x}(t)=A(t)x(t)+a(t).
\end{equation}
Offensichtlich sind (\ref{LinDGL3}) ein Spezialfall von (\ref{LinDGL2}) und
(\ref{LinDGL2}) ein Spezialfall von (\ref{LinDGL1}).
Demzufolge basiert die eindeutige Lösbarkeit der Gleichungen (\ref{LinDGL2}) und (\ref{LinDGL3}) auf den Betrachtungen zur Gleichung (\ref{LinDGL1}).
In die Gleichung (\ref{LinDGL1}) fließt in die Abbildung $A(t,s)$ die Integrationsvariable $s$ und die obere Integrationsgrenze $t$ ein.
Die Konsequenz daraus ist der zweidimensionale Zeitbereich
\begin{equation} \label{DGLZeitbereich1}
\Theta=\{(t,s) \in \R \times \R \,|\, t_0 \leq s, t \leq t_1\}.
\end{equation}
Aufgrund ihrer Einflüsse nennen wir $t$ die äußere und $s$ die innere Zeitvariable.
In den Gleichungen (\ref{LinDGL2}), (\ref{LinDGL3}) tritt die äußere Zeitvariable nicht auf und die innerere Variable bezeichnet schlicht die Zeit. 
Im Fall der Steuerung von Integralgleichungen besteht weiterhin der Umstand stückweiser Stetigkeit,
und zwar bezüglich der inneren Zeitvariablen.
Wir werden diese Eigenschaft für die Abbildung $A(t,s)$ genauer beschreiben.
Es seien die Zeitpunkte $t_0=s_0<s_1<s_2<...<s_N<s_{N+1}=t_1$ gegeben.
Dann setzen wir
\begin{equation} \label{DGLZeitbereich2}
\Theta_i=\{(t,s) \in \R \times \R \,|\, s_{i-1} \leq s \leq s_i,\; t_0 \leq t \leq t_1\}, \quad i=1,...,N+1,
\end{equation}

\begin{figure}[h]
	\centering
	\includegraphics[height=4cm]{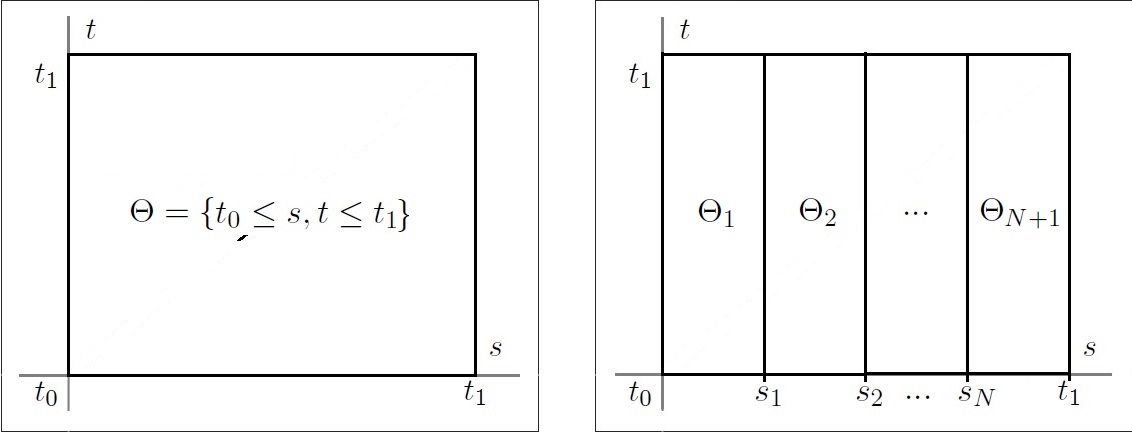}
	\caption[Zeitbereich]{Zeitbereich $\Theta$ und die Trapezaufteilung.}
\end{figure}

und die disjunkte Zerlegung $\Theta=Z_1 \cup Z_2 \cup ... \cup Z_N \cup Z_{N+1}$ des Zeitbereichs mit den Mengen
\begin{equation} \label{DGLZeitbereich3}
Z_i=\{(t,s) \in \R \times \R \,|\, s_{i-1} \leq s < s_i,\; t_0 \leq t \leq t_1\}, \quad i=1,...,N, \quad Z_{N+1}=\Theta_{N+1}.
\end{equation}

Im Weiteren bezeichnen $\Theta,\, \Theta_i,\, Z_i$ die in (\ref{DGLZeitbereich1})--(\ref{DGLZeitbereich3}) eingeführten Mengen. 
\begin{definition}[Stückweise stetige Abbildung] \label{DefinitionSTSTAbbildung}
Die Abbildung $A(t,s):\Theta \to \R^{n \times n}$ heißt über dem Zeitbereich $\Theta \subseteq \R^2$
stetig in $t$ und stückweise stetig in $s$, \index{Abbildung, beschränkte!ststetig@--, stückweise stetige}
wenn es endlich viele Stellen $t_0=s_0<s_1<s_2<...<s_N<s_{N+1}=t_1$ derart gibt,
dass für jedes $i=1,...,N+1$ über der kompakten Mengen $\Theta_i$ eine stetige Abbildung $A_i(t,s):\Theta_i \to \R^{n \times n}$ existiert mit
$A(t,s)=A_i(t,s)$ für alle $(t,s) \in Z_i$.
\end{definition}

\begin{satz} \label{LemmaDGL1}
Es sei die Abbildung $A(t,s):\Theta \to \R^{n \times n}$ stetig in $t$ und stückweise stetig in $s$.
Weiter sei $a(\cdot) \in PC([t_0,t_1],\R^n)$.
Dann existiert zu jedem $\zeta(\cdot) \in C([t_0,t_1],\R^n)$ und jedem $\tau \in [t_0,t_1]$
eine eindeutig bestimmte Lösung $x(\cdot) \in C([t_0,t_1],\R^n)$ der Gleichung
$$\zeta(t)=x(t) - \int_\tau^t \big[A(t,s)x(s)+a(s)\big] \, ds, \quad t \in [t_0,t_1].$$
\end{satz}

{\bf Beweis} Wir bemerken zunächst,
dass zu der stetigen und stückweise stetigen Abbildung $A(t,s)$ für jedes $s \in [t_0,t_1]$
der Wert $c(s)=\max\limits_{t \in [t_0,t_1]} \|A(t,s)\|$ existiert.
Weiterhin folgt aus der Stetigkeit der Abbildungen $A_i(t,s)$ über den Mengen $\Theta_i$ in Definition \ref{DefinitionSTSTAbbildung},
dass die Funktion $s \to c(s)$ stückweise stetig über $[t_0,t_1]$ ist. \\[2mm]
Im Folgenden werden wir zeigen,
dass die Fixpunktgleichung $x(\cdot) = T\big(x(\cdot)\big)$, wobei der Operator $T$ durch
$$x(\cdot) \to T\big(x(\cdot)\big), \quad
  T\big(x(\cdot)\big)(t) = \zeta(t) + \int_\tau^t \big[A(t,s)x(s)+a(s)\big] \, ds, \quad t \in [t_0,t_1],$$
gegeben wird, stets eine eindeutige L"osung besitzt.
Der Operator $T$ bildet den Raum $C([t_0,t_1],\R^n)$ in sich ab.
Mit der Funktion $c(\cdot) \in PC([t_0,t_1], \R)$ setzen wir zur abk"urzenden Schreibweise
$$C(t) = \int_\tau^t c(s) \, ds, \qquad c_0=\int_{t_0}^{t_1} c(s) \, ds.$$
Bei mehrfacher Anwendung des Operators $T$ ergeben sich f"ur $x_1(\cdot),x_2(\cdot) \in C([t_0,t_1],\R^n)$ die Beziehungen
\begin{eqnarray*}
\lefteqn{\big\| \big[T\big(x_1(\cdot) - x_2(\cdot)\big)\big](t) \big\|
         \leq \int_\tau^t c(s) \, ds \cdot \| x_1(\cdot) - x_2(\cdot) \|_\infty,} \\
\lefteqn{\big\| \big[T^2\big(x_1(\cdot) - x_2(\cdot)\big)\big](t) \big\|
         \leq \int_\tau^t c(s) \big\| \big[T\big(x_1(\cdot) - x_2(\cdot)\big)\big](s) \big\| \, ds} \\
&& \hspace*{10mm} \leq \int_\tau^t c(s) C(s) \, ds \cdot \big\| x_1(\cdot) - x_2(\cdot) \big\|_\infty
   = \frac{1}{2} C^2(t) \cdot \| x_1(\cdot) - x_2(\cdot) \|_\infty, \\
&& \hspace*{30mm}\vdots \\
\lefteqn{\big\| \big[T^m \big(x_1(\cdot) - x_2(\cdot)\big)\big](t) \big\|
         \leq \int_\tau^t c(s) \big\| \big[T^{m-1}\big(x_1(\cdot) - x_2(\cdot)\big)\big](s) \big\| \, ds} \\
&& \hspace*{10mm} \leq \int_\tau^t c(s) \frac{C^{m-1}(s)}{(m-1)!} \, ds \cdot \big\| x_1(\cdot) - x_2(\cdot) \big\|_\infty
       =\frac{C^m(t)}{m!} \cdot \| x_1(\cdot) - x_2(\cdot) \|_\infty.
\end{eqnarray*}
In der Topologie des Raumes $C([t_0,t_1],\R^n)$ gilt daher
$$\big\|T^m \big(x_1(\cdot) - x_2(\cdot)\big) \big\|_\infty \leq \frac{c_0^m}{m!} \cdot \| x_1(\cdot) - x_2(\cdot) \|_\infty.$$
Die Zahlen $a_m = c_0^m/m!$ liefern eine Folge, deren Reihe konvergiert.
Nach dem Fixpunktsatz von Weissinger (Satz \ref{SatzWeissinger}) existiert daher genau ein
$x(\cdot)$ mit $x(\cdot) = T x(\cdot)$. \hfill $\blacksquare$

\begin{lemma} \label{LemmaDGL2}
Unter den Annahmen in Satz \ref{LemmaDGL1} besitzt die Gleichung
$$\zeta(t)=x(t) - \int_\tau^t \big[A(t,s)x(s)+a(s)\big] \, ds, \quad t \in [t_0,t_1],$$
zu jedem $\zeta(\cdot) \in PC([t_0,t_1],\R^n)$ und jedem $\tau \in [t_0,t_1]$
eine eindeutig bestimmte Lösung $x(\cdot) \in PC([t_0,t_1],\R^n)$.

\end{lemma}

{\bf Beweis} In Abschnitt \ref{AbschnittNormRaum} haben wir bemerkt,
dass der Raum $PC([t_0,t_1],\R^n)$ bezüglich der Supremumsnorm nicht vollständig ist.
Um auf die gleiche Weise wie bei Satz \ref{LemmaDGL1} vorgehen zu können,
betten wir die Situation in den Raum $L_\infty([t_0,t_1],\R^n)$ der messbar und beschränkten Funktionen ein.
Da dies die einzige Stelle in dieser Ausarbeitung ist,
an der wir auf diesen Rahmen zurückgreifen,
verweisen wir hier bezüglich einer Einführung der $L_p-$Räume auf die Lehrbücher \cite{HeuserFA} von Heuser oder \cite{Werner} von Werner. \\
Der Raum $L_\infty([t_0,t_1],\R^n)$ ist versehen mit der wesentlichen Supremumsnorm $\|\cdot\|_{L_\infty}$ vollständig.
Ersetzt man im Beweis von Satz \ref{LemmaDGL1} die Norm $\|\cdot\|_\infty$ durch $\|\cdot\|_{L_\infty}$,
so ergibt sich wortwörtlich die eindeutige Lösung $x(\cdot) \in L_\infty([t_0,t_1],\R^n)$ der Gleichung
$$\zeta(t)= x(t) - \int_\tau^t \big[A(t,s)x(s)+a(s)\big] \, ds \quad \mbox{ für fast alle } t \in [t_0,t_1].$$
Darin sind die Abbildungen $\displaystyle t \to \int_\tau^t \big[A(t,s)x(s)+a(s)\big] \, ds$ stetig und $t \to \zeta(t)$ stückweise stetig.
Demzufolge muss $t \to x(t)$ über $[t_0,t_1]$ stückweise stetig sein. 
\hfill $\blacksquare$

\begin{lemma} \label{LemmaDGL3}
Es seien $A(\cdot) \in PC([t_0,t_1],\R^{n \times n})$ und $a(\cdot) \in PC([t_0,t_1],\R^n)$.
Dann existiert zu jedem $\zeta(\cdot) \in C([t_0,t_1],\R^n)$ und jedem $\tau \in [t_0,t_1]$
eine eindeutig bestimmte Lösung $x(\cdot) \in C([t_0,t_1],\R^n)$ der Gleichung
$$\zeta(t)=x(t) - \int_\tau^t \big[A(s)x(s)+a(s)\big] \, ds \quad \mbox{ für alle } t \in [t_0,t_1].$$
\end{lemma}

{\bf Beweis} Dies folgt unmittelbar aus Satz \ref{LemmaDGL1}.
\hfill $\blacksquare$

\begin{lemma} \label{LemmaDGL4}
Es seien $A(\cdot):[0,\infty) \to \R^{n \times n}$ und $a(\cdot):[0,\infty) \to \R^n$ über jedem endlichen
Intervall $[0,T]$ stückweise stetig, und es sei $t \to \|A(t)\|$ und $t \to |a(t)|$ "uber $[0,\infty)$ integrierbar.
Dann existiert zu jedem $\zeta(\cdot) \in C_{\lim}([0,\infty),\R^n)$ und jedem $\tau \in [0,\infty]$
eine eindeutig bestimmte Lösung $x(\cdot) \in C_{\lim}([0,\infty),\R^n)$ der Gleichung
$$\zeta(t)=x(t) - \int_\tau^t \big[A(s)x(s)+a(s)\big] \, ds \quad \mbox{ für alle } t \in [0,\infty].$$
\end{lemma}

{\bf Beweis} Da für jedes $x(\cdot) \in C_{\lim}([0,\infty),\R^n)$ der Grenzwert im Unendlichen existiert,
können wir $C_{\lim}([0,\infty),\R^n)$ mit dem Banachraum der stetigen Funktionen über dem Abschluss $[0,\infty]$ identifizieren
(vgl. Abschnitt \ref{AbschnittNormRaum}).
Wir ersetzen im Beweis von Satz \ref{LemmaDGL1} das Intervall $[t_0,t_1]$ durch $[0,\infty]$.
Ferner seien der Operator $T$ wie dort definiert und
$$c(t) = \|A(t)\|, \qquad C(t) = \int_\tau^t c(s) \, ds, \qquad c_0 = \int_0^\infty c(s) \, ds.$$
Bei mehrfacher Anwendung des Operators $T$ ergibt sich f"ur $x_1(\cdot),x_2(\cdot) \in C_{\lim}([0,\infty),\R^n)$:
$$\big\|T^m \big(x_1(\cdot) - x_2(\cdot)\big) \big\|_\infty \leq \frac{c_0^m}{m!} \cdot \| x_1(\cdot) - x_2(\cdot) \|_\infty.$$
Die Zahlen $a_m = c_0^m / m!$ liefern eine Folge, deren Reihe konvergiert.
Nach dem Fixpunktsatz von Weissinger existiert daher genau ein $x(\cdot)$ mit $x(\cdot) = T x(\cdot)$.
\hfill $\blacksquare$

\begin{lemma} \label{LemmaDGL5}
Es seien die Voraussetzungen in Lemma \ref{LemmaDGL3} erfüllt.
Dann existiert zu jedem $\zeta \in \R^n$ und jedem $\tau \in [t_0,t_1]$
eine eindeutige L"osung $x(\cdot) \in PC_1([t_0,t_1],\R^n)$ der Gleichung
$$\dot{x}(t)=A(t)x(t)+a(t), \qquad x(\tau)=\zeta.$$
\end{lemma}

{\bf Beweis} Wir merken an, dass in Lemma \ref{LemmaDGL3} der Integrand $t \to A(t)x(t)+a(t)$ stückweise stetig ist.
Daher folgt $\dot{x}(\cdot) \in PC([t_0,t_1],\R^n)$ unmittelbar aus Lemma \ref{LemmaDGL3} mit $\zeta(t) \equiv \zeta$. \hfill $\blacksquare$

\begin{lemma} \label{LemmaDGL6}
Es seien die Voraussetzungen in Lemma \ref{LemmaDGL4} erfüllt.
Dann existiert zu jedem $\zeta \in \R^n$ und jedem $\tau \in [0,\infty]$
eine eindeutige L"osung $x(\cdot) \in PC_1([0,\infty),\R^n)$ der Gleichung
$$\dot{x}(t)=A(t)x(t)+a(t), \qquad x(\tau)=\zeta \quad (\mbox{wobei } x(\infty):=\lim_{t \to \infty} x(t) \mbox{ gilt}).$$
\end{lemma}

{\bf Beweis} Die Behauptung folgt aus Lemma \ref{LemmaDGL4},
wenn wir $\zeta(t) \equiv \zeta$ setzen. \hfill $\blacksquare$


\subsection{Existenz und Eindeutigkeit, Abhängigkeit von Anfangsdaten}
Wir widmen unser Interesse nun den nichtlinearen Differentialgleichungen der Form
\begin{equation}\label{DGL1}
\dot{x}(t) = \varphi\big(t,x(t)\big)
\end{equation}
und den nichtlinearen Volterraschen Integralgleichungen der Gestalt
\begin{equation}\label{DGL2}
x(t)=x(t_0)+ \int_{t_0}^t \varphi\big(t,s,x(s)\big) \, ds.
\end{equation}
Da sich die Differentialgleichung (\ref{DGL1}) in der äquivalenten Form
\begin{equation}\label{DGL3}
x(t)=x(t_0)+ \int_{t_0}^t \varphi\big(s,x(s)\big) \, ds
\end{equation}
als Spezialfall von (\ref{DGL2}) erweist, konzentrieren wir uns auf die Integralgleichung.
Dementsprechend sind für die Gleichung (\ref{DGL3}) die nachstehenden Resultate lediglich ohne der äußeren Zeitvariablen $t$
zu formulieren und zu übernehmen. \\[2mm]
Wir treffen für die weiteren Untersuchungen folgende Festlegungen:
\begin{enumerate}
\item[(1)] Für $\gamma >0$ und $x(\cdot) \in C([a,b],\R^n)$ sei $\mathscr{R}$ der kompakte Umgebungsstreifen 
           $$\{(t,s,x) \in \R \times \R \times \R^n \,|\, a \leq s, t\leq b, \, \|x-x(t)\| \leq \gamma \}.$$
           Ferner bezeichne $U_\delta(x) \subset \R^n$ die offene Kugel $\{z \in \R^n \,|\, \|z-x\| < \delta\}$.
\item[(2)] Die Abbildung $\varphi(t,s,x)$ sei auf der Menge $\mathscr{R}$ in der Variablen $x$ gleichmäßig stetig und
           gleichmäßig stetig differenzierbar, d.\,h. es gibt eine positive Konstante $L$ mit
           $$\|\varphi(t,s,x_1)-\varphi(t,s,x_2)\| \leq L \|x_1-x_2\|, \quad \|\varphi_x(t,s,x_1)-\varphi_x(t,s,x_2)\| \leq L \|x_1-x_2\|$$
           für alle $(t,s,x_1),\, (t,s,x_2) \in \mathscr{R}$.
\item[(3)] Die Abbildungen $\varphi(t,s,x)$ und $\varphi_x(t,s,x)$ seien auf der Menge $\mathscr{R}$ in der Variablen $t$ stetig
           und in der Variablen $s$ stückweise stetig (vgl. Definition \ref{DefinitionSTSTAbbildung}).
\item[(4)] Die Abbildungen $\varphi(t,s,x)$ und $\varphi_x(t,s,x)$ seien über $\mathscr{R}$ beschränkt und es gelten
           $\|\varphi(t,s,x)\| \leq M,\, \|\varphi_x(t,s,x)\| \leq M$ über $\mathscr{R}$ mit der positiven Konstanten $M$.
\end{enumerate}

\begin{satz}[Lokaler Fixpunkt] \label{SatzFixpunktlokal}
Zu jedem inneren Punkt $(t_0,t_0,\xi)$ der Menge $\mathscr{R}$ gibt es Zahlen $\varepsilon >0$ und $\delta >0$ derart,
dass über dem Intervall $I=[t_0-\varepsilon,t_0+\varepsilon]$ die Abbildung
$$A: C(I,\R^n) \to C(I,\R^n), \quad
  [Ax(\cdot)](t)= \xi(t) + \int_{t_0}^t \varphi\big(t,s,x(s)\big) \, ds, \; t \in I,$$
für jedes $\xi(\cdot) \in C(I,\R^n)$ mit $\|\xi(\cdot)-\xi\|_\infty < \delta$ genau einen Fixpunkt $x_\xi(\cdot) \in C(I,\R^n)$ besitzt.
Unter diesen Voraussetzungen gilt ferner für die Fixpunkte $x_{\xi_1}(\cdot), \,x_{\xi_2}(\cdot) \in C(I,\R^n)$ der Abbildung $A$
zu den Funktionen $\xi_1(\cdot),\,\xi_2(\cdot) \in C(I,\R^n)$ die Ungleichung
$$\|x_{\xi_1}(\cdot)-x_{\xi_2}(\cdot)\|_\infty \leq  \frac{1}{1-\varepsilon L}\|\xi_1(\cdot)-\xi_2(\cdot)\|_\infty.$$
\end{satz}

{\bf Beweis} Wir wählen die positiven Zahlen $\varepsilon$ und $\delta$ so, dass 
\begin{equation} \label{DGLoffen}
\{(t,s,x) \in \R \times \R \times \R^n \,|\, a < t_0-\varepsilon \leq s, t\leq t_0+\varepsilon <b, \, \|x-\xi\| < 2\delta \}\subset \mathscr{R}
\end{equation}
gelten und weiterhin die Beziehungen $\varepsilon L < 1$ und $\varepsilon M <\delta$ erfüllt sind. \\
Wir betrachten im Raum $C(I,\R^n)$ die abgeschlossene und konvexe Teilmenge
$$K=\{ x(\cdot) \in C(I,\R^n) \,|\, \|x(t)-\xi\|\leq 2\delta \mbox{ für alle } t \in I = [t_0-\varepsilon,t_0+\varepsilon]\}.$$
Es sei $\xi(\cdot) \in C(I,\R^n)$ mit $\|\xi(\cdot)-\xi\|_\infty < \delta$.
Für $x(\cdot) \in K$ und $t \in I$ gilt dann
\begin{eqnarray*}
        \|[Ax(\cdot)](t) -\xi\|
&\leq&  \|[Ax(\cdot)](t) -\xi(t)\| + \|\xi(t)-\xi\| < \bigg\| \int_{t_0}^t \varphi\big(t,s,x(s)\big) \, ds \bigg\| + \delta \\
&\leq& \int_{t_0}^t M \, ds + \delta \leq \varepsilon M + \delta < 2\delta,
\end{eqnarray*}
d.\,h., die Abbildung $A$ bildet die Menge $K$ in sich ab.
Ferner gilt für $x_1(\cdot),\, x_2(\cdot) \in K$
\begin{eqnarray*}
    \|Ax_1(\cdot) - Ax_2(\cdot)\|_\infty
&=& \max_{t \in I} \bigg\| \int_{t_0}^t \big[\varphi\big(t,s,x_1(s)\big)-\varphi\big(t,s,x_2(s)\big)\big] \, ds \bigg\| \\
&\leq& \max_{t \in I} \int_{t_0}^t L \|x_1(s)-x_2(s)\| \, ds
       \leq \varepsilon L \cdot \big\| x_1(\cdot)- x_2(\cdot)\big\|_\infty,
\end{eqnarray*}
d.\,h., die Abbildung $A$ ist wegen $\varepsilon L < 1$ kontrahierend.
Nach dem Banachschen Fixpunktsatz besitzt die Abbildung $A$ einen eindeutigen Fixpunkt $x(\cdot) \in K$. \\ 
Es seien die Funktion $\xi_1(\cdot),\,\xi_2(\cdot)  \in C(I,\R^n)$ mit
$\|\xi_1(\cdot)-\xi\|_\infty \leq \delta$, $\|\xi_2(\cdot)-\xi\|_\infty \leq \delta$ und 
es seien $x_{\xi_1}(\cdot), \,x_{\xi_2}(\cdot) \in C(I,\R^n)$ die zugehörigen Fixpunkte der Abbildung $A$.
Dann gilt
\begin{eqnarray*}
    \|x_{\xi_1}(\cdot)-x_{\xi_2}(\cdot)\|_\infty
&=& \max_{t \in I}
    \bigg\| \xi_1(t)-\xi_2(t)+\int_{t_0}^t \big[\varphi\big(t,s,x_{\xi_1}(s)\big)-\varphi\big(t,s,x_{\xi_2}(s)\big)\big] \, ds \bigg\| \\
&\leq&  \|\xi_1(\cdot)-\xi_2(\cdot)\|_\infty + \varepsilon L \cdot \big\|x_{\xi_1}(\cdot)-x_{\xi_2}(\cdot)\big\|_\infty.
\end{eqnarray*}
Daraus ergibt sich abschlie"send
$$\big\|x_{\xi_1}(\cdot)-x_{\xi_2}(\cdot)\big\|_\infty \leq  \frac{1}{1-\varepsilon L}\|\xi_1(\cdot)-\xi_2(\cdot)\|_\infty.$$
Der Satz \ref{SatzFixpunktlokal} ist damit vollständig nachgewiesen. \hfill $\blacksquare$

\begin{satz}[Lokaler Existenz-, Eindeutigkeits- und Stetigkeitssatz] \label{SatzEElokal}
Zu jedem inneren Punkt $(t_0,t_0,\xi)$ der Menge $\mathscr{R}$ gibt es Zahlen $\varepsilon >0$ und $\delta >0$,
so dass zu jedem $z \in U_\delta(\xi)$ genau eine stetige L"osung $x_z(\cdot)$ der Gleichung (\ref{DGL2}) zur Anfangsbedingung $x(t_0)=z$
auf dem Intervall $I=[t_0-\varepsilon,t_0+\varepsilon]$ existiert. \\
Ist ferner $\{z_k\}$ eine Folge aus $U_\delta(\xi)$, die gegen $z \in U_\delta(\xi)$ konvergiert,
so gilt
$$\max_{t \in I} \big\|x_{z_k}(t)-x_z(t)\big\|_\infty \leq  \frac{1}{1-\varepsilon L}\|z_k-z\|$$
und $\{x_{z_k}(\cdot)\}$ konvergiert auf $I$ gleichm"a"sig gegen $x_z(\cdot)$.
\end{satz}

{\bf Beweis} Nach Satz \ref{SatzFixpunktlokal} können wir $\varepsilon > 0$ und $\delta >0$ so wählen,
dass die Abbildung 
$$[Ax(\cdot)](t)= z + \int_{t_0}^t \varphi\big(t,s,x(s)\big) \, ds, \quad t \in I=[t_0-\varepsilon,t_0+\varepsilon],$$
für jedes $z \in U_\delta(\xi)$ genau einen Fixpunkt $x_z(\cdot) \in C(I,\R^n)$ besitzt, d.\,h.
$$Ax_z(\cdot) = x_z(\cdot) \quad\Leftrightarrow\quad 
  x_z(t)=z + \int_{t_0}^t \varphi\big(t,s,x_z(s)\big) \, ds \mbox{ für alle } t \in I=[t_0-\varepsilon,t_0+\varepsilon].$$
Die Existenz und Eindeutigkeit der Lösung ist damit gezeigt.\\
Weiter sei $\delta >0$ so gewählt, dass (\ref{DGLoffen}) erfüllt ist.
Es seien nun $z_1,z_2 \in U_\delta(\xi)$.
Dann existieren die Lösungen $x_{z_1}(\cdot)$, $x_{z_2}(\cdot)$ der Gleichung (\ref{DGL2})
über $I$ und nehmen Werte in der Menge $U_{2\delta}(\xi)$ an.
Nach Satz \ref{SatzFixpunktlokal} ist die Ungleichung
$$\big\|x_{z_1}(\cdot)-x_{z_2}(\cdot)\big\|_\infty \leq  \frac{1}{1-\varepsilon L}\|z_1-z_2\|$$
erfüllt,
welche die gleichmäßige Konvergenz von $x_{z_k}(\cdot)$ gegen $x_z(\cdot)$ für $z_k \to z$ bedeutet.
Satz \ref{SatzEElokal} ist somit vollständig bewiesen. \hfill $\blacksquare$ \\[2mm]
Ist $x(\cdot)$ über $[t_0,t_1]$ eine Lösung der Gleichung (\ref{DGL3}),
so lässt sich $x(\cdot)$ formal durch 
$$x(t)=x(t_1)+ \int_{t_1}^t \varphi\big(s,x(s)\big) \, ds$$
zu einer Lösung über $[t_0,t_1+\varepsilon]$ fortsetzen.
Im Fall von Gleichung (\ref{DGL2}) führt
$$x(t)=x(t_1)+ \int_{t_1}^t \varphi\big(t,s,x(s)\big) \, ds$$
im Allgemeinen nicht zu der gewünschten Fortsetzung.
Stattdessen erfolgt eine Fortsetzung der Gleichung (\ref{DGL2}) an der Stelle $t=t_1$ nicht mit dem Anfangswert $x(t_1)$,
sondern mit einer Funktion $\xi(t)$, für die $\xi(t_1)=x(t_1)$ gilt.
Dieser Sachverhalt ist der Gegenstand der nachstehen Untersuchung:

\begin{lemma}[Lokale Fortsetzbarkeit] \label{LemmaDGLFortsetzung}
Es sei $x(\cdot)$ über $[t_0,t_1] \subset (a,b)$ eine Lösung von (\ref{DGL2}) und es sei $\big(t_1,t_1,x(t_1)\big)$
ein innerer Punkt der Menge $\mathscr{R}$.
Dann existiert eine Zahl $\varepsilon >0$ derart,
dass $x(\cdot)$ eindeutig zu einer stetigen Lösung der Gleichung (\ref{DGL2}) über dem Intervall $[t_0,t_1+\varepsilon]$ fortgesetzt werden kann.
Die Fortsetzung ergibt sich für $t \in [t_1,t_1+\varepsilon]$ als die eindeutige Lösung der Gleichung
\begin{equation}\label{DGLFortsetzung}
x(t) = \displaystyle \xi(t)+ \int_{t_1}^t \varphi\big(t,s,x(s)\big) \, ds \quad\mbox{mit}\quad
\xi(t)=x(t_0)  + \int_{t_0}^{t_1} \varphi\big(t,s,x(s)\big) \, ds.
\end{equation}
\end{lemma}

{\bf Beweis} Als Lösung der Gleichung (\ref{DGL2}) gilt für $x(\cdot)$ in $t=t_1$ die Beziehung
$$x(t_1)= x(t_0) + \int_{t_0}^{t_1} \varphi\big(t_1,s,x(s)\big) \, ds.$$
Damit bringen wir die Abbildung $t \to \xi(t)$ in die Form
\begin{equation} \label{DGLUmformung}
\xi(t) = x(t_1)+ \int_{t_0}^{t_1} \big[ \varphi\big(t,s,x(s)\big) - \varphi\big(t_1,s,x(s)\big) \big] \, ds.
\end{equation}
Nun sieht man unmittelbar,
dass $t \to \xi(t)$ für hinreichend kleine $\varepsilon >0$ auf $[t_1,t_1+\varepsilon]$ stetig ist und $\xi(t_1)=x(t_1)$ gilt.
Daher kann zu jedem $\delta >0$ ein $\varepsilon >0$ gewählt werden, dass $\|\xi(t)-x(t_1)\| < \delta$ über $[t_1,t_1+\varepsilon]$ erfüllt ist.
Nach Satz \ref{SatzFixpunktlokal} besitzt die Abbildung
$$[Ay(\cdot)](t)= \xi(t) + \int_{t_1}^t \varphi\big(t,s,y(s)\big) \, ds, \qquad t \in [t_1,t_1+\varepsilon],$$
für hinreichend kleine $\varepsilon > 0$ und $\delta >0$ genau einen Fixpunkt $y(\cdot) \in C([t_1,t_1+\varepsilon],\R^n)$.
Es gilt $y(t_1)=x(t_1)$.
Setzen wir zudem $y(t)=x(t)$ für $t \in [t_0,t_1]$,
so ist $y(\cdot) \in C([t_0,t_1+\varepsilon],\R^n)$ eine Fortsetzung der Lösung $x(\cdot)$ von (\ref{DGL2}).
Denn offenbar erfüllt $y(\cdot)$ für jedes $t \in [t_0,t_1]$ die Gleichung (\ref{DGL2}).
Für $t \in (t_1,t_1+\varepsilon]$ folgt
\begin{eqnarray*}
y(t) &=& \xi(t) + \int_{t_1}^t \varphi\big(t,s,y(s)\big) \, ds
         = x(t_0)  + \int_{t_0}^{t_1} \varphi\big(t,s,x(s)\big) \, ds + \int_{t_1}^t \varphi\big(t,s,y(s)\big) \, ds \\
     &=& x(t_0)  + \int_{t_0}^{t_1} \varphi\big(t,s,y(s)\big) \, ds + \int_{t_1}^t \varphi\big(t,s,y(s)\big) \, ds
         = x(t_0 ) + \int_{t_0}^t \varphi\big(t,s,y(s)\big) \, ds.
\end{eqnarray*}
Damit löst $y(\cdot)$ die Gleichung $\displaystyle y(t) = x(t_0 ) + \int_{t_0}^t \varphi\big(t,s,y(s)\big) \, ds$ für alle $t \in [t_0,t_1+\varepsilon]$. \\
Ist umgekehrt $\zeta(\cdot)$ eine fortgesetzte Lösung der Gleichung (\ref{DGL2}) von $x(\cdot)$,
d.\,h. $\zeta(t)=x(t)$ für alle $t \in [t_0,t_1]$,
so gilt für $t > t_1$:
$$\zeta(t) = x(t_0 ) + \int_{t_0}^t \varphi\big(t,s,\zeta(s)\big) \, ds
           = \underbrace{x(t_0 ) + \int_{t_0}^{t_1} \varphi\big(t,s,x(s)\big) \, ds}_{=\xi(t)} + \int_{t_1}^t \varphi\big(t,s,\zeta(s)\big) \, ds.$$
Damit wird $x(\cdot)$ nur durch die eindeutige Lösung von (\ref{DGLFortsetzung}) fortgesetzt. \hfill $\blacksquare$ 

\begin{satz}[Globaler Existenz-, Eindeutigkeits- und Stetigkeitssatz] \label{SatzEEglobal}
Über $[t_0,t_1]$ sei $x(\cdot)$ eine stetige Lösung der Gleichung (\ref{DGL2}) zum Anfangswert $x(t_0)=x_0$,
für die die Menge $\{(t,s,x) \in \R \times \R \times \R^n \,|\, t_0 \leq s, t\leq t_1, \, x=x(s) \}$ dem Inneren der Menge $\mathscr{R}$ angehört. \\
Dann existiert ein $\delta>0$ derart,
dass für jedes $z \in U_\delta(x_0)$ die Gleichung
\begin{equation} \label{AnhangDGLMax}
x(t) = z+ \int_{t_0}^t \varphi\big(t,s,x(s)\big) \, ds
\end{equation}
eine eindeutige stetige Lösung $x_z(\cdot)$ über $[t_0,t_1]$ besitzt und
die Funktionen $x_z(\cdot)$ gleichmäßig gegen $x(\cdot)$ über $[t_0,t_1]$ für $z \to x_0$ konvergieren.
\end{satz}

{\bf Beweis} Nach Satz \ref{SatzEElokal} gibt es ein $\delta >0$ und ein $\tau >t_0$ derart,
dass für jedes $z \in U_\delta(x_0)$ die Gleichung (\ref{AnhangDGLMax})
eine Lösung über $[t_0,\tau]$ besitzt und die Funktionen $x_z(\cdot)$ gleichmäßig gegen $x(\cdot)$ über $[t_0,\tau]$ für $z \to x_0$
konvergieren.
Im Weiteren sei $\overline{\tau}$ durch das Supremum über alle $\tau > t_0$ mit dieser Eigenschaft definiert.
D.\,h. es gibt ein $\delta >0$ derart,
dass für jedes $z \in U_\delta(x_0)$ eine Lösung $x_z(\cdot)$ von (\ref{AnhangDGLMax}) über $[t_0,\tau]$ existiert und
die Funktionen $x_z(\cdot)$ gleichmäßig gegen $x(\cdot)$ über $[t_0,\tau]$ für $z \to x_0$ konvergieren.
Zum Beweis des Satzes ist es ausreichend, die Beziehung $\overline{\tau}>t_1$ nachzuweisen. \\[1mm]
Angenommen, es ist $\overline{\tau} \leq t_1$.
Dann ist $\big(\overline{\tau},\overline{\tau},x(\overline{\tau}) \big)$ ein innerer Punkt der Menge $\mathscr{R}$.
Nach Satz \ref{SatzEElokal} lässt sich zu jedem $\varrho >0$ ein $\delta >0$ derart angeben,
dass $\|x_z(\cdot)-x(\cdot)\|_\infty \leq \varrho$ für alle $z \in U_\delta(x_0)$ gilt.
Insbesondere konvergiert $x_z(\cdot)$ gegen $x(\cdot)$ gleichmäßig für $z \to x_0$
und es ist $x_z(\overline{\tau}) \in U_\varrho\big(x(\overline{\tau})\big)$. \\

\begin{figure}[h]
	\centering
	\includegraphics[height=5cm]{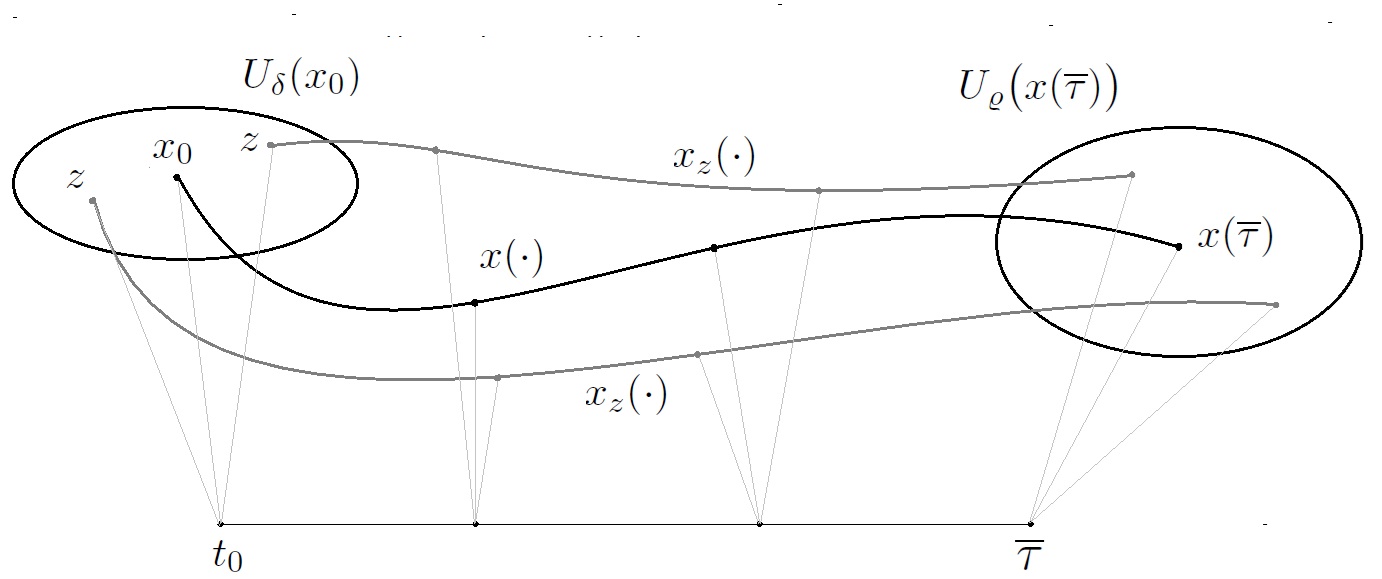}
	\caption[DGLglobal1]{Für jedes $z \in U_\delta(x_0)$ gehört $x_z(\overline{\tau})$ zu $U_\varrho\big(x(\overline{\tau})\big)$.}
\end{figure}

Sei zunächst $x_z(\cdot)$ die stetige Lösung von (\ref{AnhangDGLMax}) über $[t_0,\overline{\tau}]$ zu einem festen $z \in U_\delta(x_0)$.
Dann gibt es nach Lemma \ref{LemmaDGLFortsetzung} eine Zahl $\varepsilon > 0$ und eine stetige Fortsetzung $y_z(\cdot)$ von $x_z(\cdot)$,
welche die Gleichung (\ref{AnhangDGLMax}) über $[t_0,\overline{\tau}+\varepsilon]$ löst.
Die Fortsetzung $y_z(\cdot)$ ist für $t \in [\overline{\tau},\overline{\tau}+\varepsilon]$ die Lösung der Gleichung
$$y_z(t) = \xi_z(t)+ \int_{\overline{\tau}}^t \varphi\big(t,s,y_z(s)\big)\big) \, ds \quad\mbox{mit}\quad
\xi_z(t)= z  + \int_{t_0}^{\overline{\tau}} \varphi\big(t,s,x_z(s)\big) \, ds.$$
Es ist nun zu zeigen, dass zu jedem $\varrho>0$ eine Kugel $U_{\delta'}(x_0) \subseteq U_\delta(x_0)$ und ein $\varepsilon > 0$ so angeben lassen,
dass für alle $z \in U_{\delta'}(x_0)$ sich die Lösungen $x_z(\cdot)$ durch $y_z(\cdot)$ auf $t \in [\overline{\tau},\overline{\tau}+\varepsilon]$
fortsetzen lassen und die Ungleichung
$$\max_{t \in [\overline{\tau},\overline{\tau}+\varepsilon]} \|y_z(t)-y(t) \| \leq \varrho$$
erfüllt ist.
Denn dann konvergiert $x_z(\cdot)$ zusammen mit der zugehörigen Fortsetzung $y_z(\cdot)$ gleichmäßig über $[t_0,\overline{\tau}+\varepsilon]$ gegen
die Lösung $x(\cdot)$ und deren Fortsetzung $y(\cdot)$.

\newpage
\begin{figure}[h]
	\centering
	\includegraphics[height=5cm]{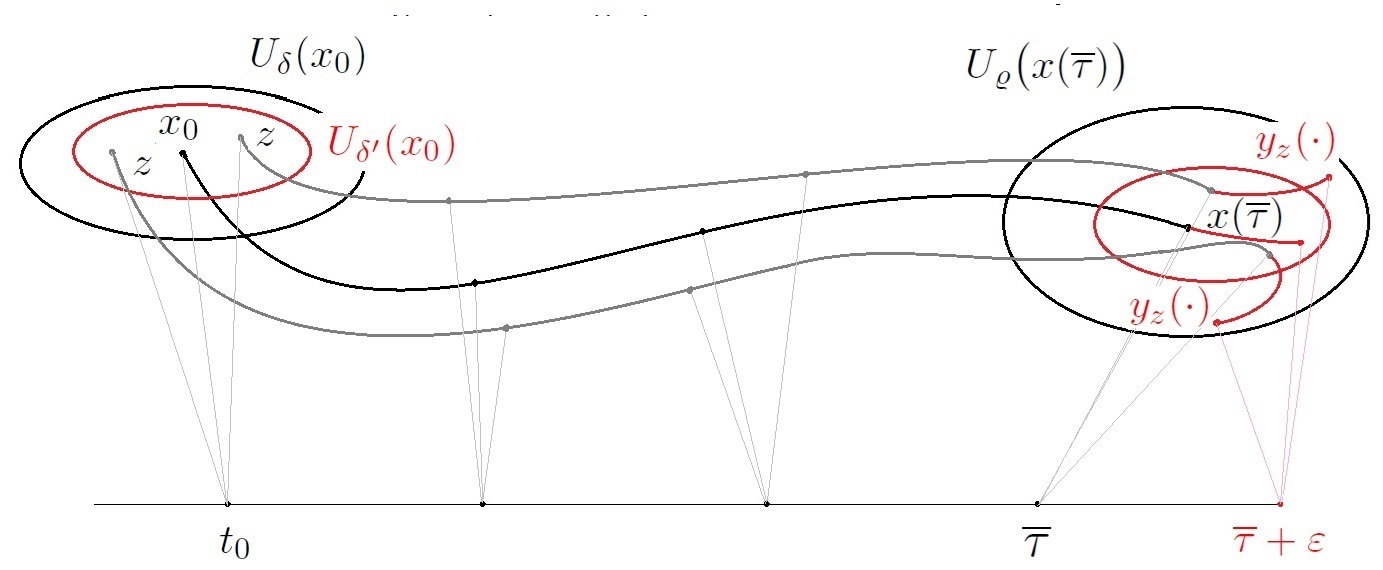}
	\caption[DGLglobal2]{Die Fortsetzungen $y_z(\cdot)$ von $x_z(\cdot)$.}
\end{figure}

Nach Satz \ref{SatzFixpunktlokal} existieren $\varepsilon >0$ und $\varrho >0$ derart,
dass die Abbildung $A$ für jedes stetige $\xi(\cdot)$ mit $\|\xi(t)-x(\overline{\tau})\| < \varrho$ für alle $t \in [\overline{\tau},\overline{\tau}+\varepsilon]$
genau einen Fixpunkt $x_\xi(\cdot)$ über $[\overline{\tau},\overline{\tau}+\varepsilon]$ besitzt.
Weiterhin wählen wir $\varepsilon$ mit $\varepsilon L <1$ und den Radius $\delta'$ der Kugel $U_{\delta'}(x_0)$ mit
$$\max_{t \in [t_0,\overline{\tau}]} \|x_z(t)-x(t) \| < \frac{1-\varepsilon L}{1 + 2L ( \overline{\tau}-t_0 ) } \varrho \quad
  \mbox{ für alle } z \in U_{\delta'}(x_0).$$
Die positiven Zahlen $\varepsilon, \, \varrho,\, \delta'$ seien im Weiteren so vorgegeben.
Ferner bezeichnen $\xi(\cdot)$ bzw. $\xi_z(\cdot)$ die Funktion,
mit deren Hilfe die Fortsetzungen für $t > \overline{\tau}$ gebildet werden:
$$\xi(t)= x(t_0)  + \int_{t_0}^{\overline{\tau}} \varphi\big(t,s,x(s)\big) \, ds,
  \quad \xi_z(t)= z  + \int_{t_0}^{\overline{\tau}} \varphi\big(t,s,x_z(s)\big) \, ds.$$
Für $z \in U_{\delta'}(x_0)$ ergeben sich mit den Darstellungen von $\xi(\cdot)$ und $\xi_z(\cdot)$ in der Form (\ref{DGLUmformung}):
\begin{eqnarray*}
\|\xi_z(t)-\xi(t)\|
&=& \bigg\|x_z(\overline{\tau}) - x(\overline{\tau})  
       + \int_{t_0}^{\overline{\tau}} \big[ \varphi\big(t,s,x_z(s)\big) - \varphi\big(t,s,x(s)\big) \big]\, ds \\
&&  + \int_{t_0}^{\overline{\tau}} \big[ \varphi\big(\overline{\tau},s,x_z(s)\big) - \varphi\big(\overline{\tau},s,x(s)\big) \big]\, ds \bigg\|\\
&\leq& \|x_z(\overline{\tau}) - x(\overline{\tau}) \| + 
       2L(\overline{\tau}-t_0) \cdot \max_{t \in [t_0,\overline{\tau}]} \|x_z(t)-x(t)\|.
\end{eqnarray*}
Mit Hilfe von Satz \ref{SatzFixpunktlokal} ergibt sich für die Fortsetzungen $y_z(\cdot), \, y(\cdot)$ nun der Abstand
\begin{eqnarray*}
\max_{t \in [\overline{\tau},\overline{\tau}+\varepsilon]} \|y_z(t)-y(t)\|
&\leq& \frac{1}{1-\varepsilon L} \cdot \max_{t \in [\overline{\tau},\overline{\tau}+\varepsilon]} \|\xi_z(t)-\xi(t)\| \\
&\leq& \frac{1 + 2L(\overline{\tau}-t_0)}{1-\varepsilon L} \cdot \max_{t \in [t_0,\overline{\tau}]} \|x_z(t)-x(t)\|
       < \varrho.
\end{eqnarray*}
Dies zeigt, dass sich zu jedem $\varrho >0$ positive Zahlen $\varepsilon$ und $\delta'$ so angeben lassen,
dass über $[t_0,\overline{\tau}+\varepsilon]$ die Ungleichung $\|y_z(t)-y(t)\| < \varrho$ für alle $z \in U_{\delta'}(x_0)$
erfüllt ist;
also die fortgesetzten Lösungen $x_z(\cdot)$ gleichmäßig gegen $x(\cdot)$ über $[t_0,\overline{\tau}+\varepsilon]$ für $z \to x_0$ konvergieren. 
Daher muss $\overline{\tau} > t_1$ gelten. 
Der Satz \ref{SatzEEglobal} ist damit bewiesen. \hfill $\blacksquare$

\begin{lemma} \label{LemmaDGLDifferenzierbarkeit}
Es bezeichnet weiterhin $\mathscr{R}=\{(t,s,x) \,|\, a \leq s, t\leq b, \, \|x-x_0\| \leq \gamma \}$.
Ferner seien $a \leq t_0<t_1 \leq b$.
Dann ist die Abbildung
$$G:C([t_0,t_1],\R^n) \to C([t_0,t_1],\R^n), \quad
  \big[G\big(x(\cdot)\big)\big](t)= x(t_0) - \int_{t_0}^t \varphi\big(t,s,x(s)\big) \, ds, \; t \in [t_0,t_1],$$
über der offenen Menge
$U=\{ x(\cdot) \in C([t_0,t_1],\R^n) \,|\, \|x(\cdot)-x_0\|_\infty < \gamma\}$
stetig differenzierbar und es gilt für jedes $x(\cdot) \in U$ und jedes $y(\cdot) \in C([t_0,t_1],\R^n)$:
$$\big[G'\big(x(\cdot)\big)y(\cdot)\big](t)= \int_{t_0}^t \varphi_x\big(t,s,x(s)\big) y(s) \, ds, \quad t \in [t_0,t_1].$$
\end{lemma}

{\bf Beweis} Nach Voraussetzung ist die Abbildung $G$ für jedes $x(\cdot) \in U$ auf einer Umgebung wohldefiniert.
Weiterhin ist der lineare Operator $G'\big(x(\cdot)\big)$ wegen
$$\big\|G'\big(x(\cdot)\big)\big\| =\sup_{\|y(\cdot)\|_\infty=1} \big\|G'\big(x(\cdot)\big)y(\cdot)\big\|
  \leq \int_{t_0}^{t_1} \big\|\varphi_x\big(t,s,x(s)\big)\big\| \, ds \cdot \|y(\cdot)\|_\infty \leq M(t_1-t_0)$$
beschränkt und somit stetig. 
Da nach unseren Annahmen die Abbildung $\varphi$ in $x$ gleichmäßig stetig differenzierbar ist,
gibt es zu $x(\cdot) \in U$ und $\varepsilon >0$ eine Zahl $\lambda_0>0$ mit 
$$\big\|\varphi_x\big(t,s,x(t)+\lambda y\big)-\varphi_x\big(t,s,x(t)\big)\big\| \leq L \cdot \lambda \|y\| \leq \frac{\varepsilon}{t_1-t_0}$$
f"ur alle $t_0 \leq s, t \leq t_1$, f"ur alle $\|y\| < 1$ und alle $0 < \lambda < \lambda_0$.
Weiter nutzen wir 
$$g(1)-g(0)=\int_0^1 g'(\tau)\, d\tau \mbox{ für }
  g(\tau)=\frac{\varphi(t,s,x+\tau\lambda y) - \varphi(t,s,x)}{\lambda}, \; g'(\tau)=\varphi_x(t,s,x+ \tau \lambda y)y.$$
Damit erhalten wir f"ur alle $\|y(\cdot) \|_\infty < 1$ und alle $0 < \lambda < \lambda_0$:
\begin{eqnarray*}
\lefteqn{\bigg\|
     \frac{G\big(x(\cdot)+\lambda y(\cdot)\big) - G\big(x(\cdot)\big)}{\lambda} - G'\big(x(\cdot)\big) y(\cdot)\bigg\|_\infty } \\
&\leq& \int_{t_0}^{t_1} \bigg[ \int_0^1
       \big\|\varphi_x\big(t,s,x(t)+\tau \lambda y(t)\big)-\varphi_x\big(t,s,x(t)\big)\big\|\, d\tau \bigg]  \, dt \leq \varepsilon,
\end{eqnarray*}
d.\,h., die Abbildung $G$ ist über der Menge $U$ Fr\'echet-differenzierbar. 
Wir zeigen noch die Stetigkeit der Abbildung $x(\cdot) \to G'\big(x(\cdot)\big)$ über der Menge $U$.
Für $x_1(\cdot),\,x_2(\cdot) \in U$ gilt
\begin{eqnarray*}
\lefteqn{    \big\|G'\big(x_1(\cdot)\big) - G'\big(x_2(\cdot)\big) \big\|
= \sup_{\|y(\cdot)\|_\infty=1} \big\| \big[G'\big(x_1(\cdot)\big) - G'\big(x_2(\cdot)\big)\big] y(\cdot) \big\|_\infty} \\
&\leq& \int_{t_0}^{t_1} \big\|  \varphi_x\big(t,s,x_1(t)\big) - \varphi_x\big(t,s,x_2(t)\big) \big\| \, ds
       \leq L (t_1-t_0) \|x_1(\cdot)-x_2(\cdot) \|_\infty.
\end{eqnarray*}
Damit ist $G$ über der Menge $U$ stetig differenzierbar. \hfill $\blacksquare$

\begin{satz}[Differenzierbarkeitssatz] \label{SatzDGLDifferenzierbarkeit}
Es seien $x_0(\cdot)$ und $x_z(\cdot)$ die Lösungen von (\ref{AnhangDGLMax}) zu den Anfangswerten $x_0$ bzw. $z \in U_\delta(x_0)$ in Satz \ref{SatzEEglobal}.
Weiter bezeichne $\Phi$ die Abbildung, die jedem Anfangswert $z \in U_\delta(x_0)$ die Lösung $x_z(\cdot)$ von (\ref{AnhangDGLMax}) zuordnet. \\
Dann ist $\Phi:U_\delta(x_0) \to C([t_0,t_1],\R^n)$ mit $\Phi(z)=x_z(\cdot)$ in $x_0$ stetig differenzierbar,
und für jedes $y \in \R^n$ genügt die Ableitung
$$\lim_{\lambda \to 0} \frac{\Phi(x_0+\lambda y)-\Phi(x_0)}{\lambda} = \Phi'(x_0)y=\xi_y(\cdot) \in C([t_0,t_1],\R^n)$$
der linearen Integralgleichung
$$\xi_y(t)= y + \int_{t_0}^t \varphi_x\big(t,s,x_0(s)\big) \xi_y(s) \, ds, \quad t \in [t_0,t_1].$$ 
\end{satz}

{\bf Beweis} Wir betrachten die Abbildung $F:\R^n \times C([t_0,t_1],\R^n) \to C([t_0,t_1],\R^n)$, die durch
$$\big[F\big(z,x(\cdot)\big)\big](t)= x(t)- z- \int_{t_0}^t \varphi\big(t,s,x(s)\big) \, ds = \big[G\big(x(\cdot)\big)\big](t) - z$$
definiert ist.
Diese ist nach Lemma \ref{LemmaDGLDifferenzierbarkeit} in einer Umgebung des Punktes $x_0(\cdot)$ nach $x(\cdot)$ stetig differenzierbar und besitzt die Ableitung
$$\big[F_{x(\cdot)}\big(z,x(\cdot)\big)y(\cdot)\big](t)= y(t)- \int_{t_0}^t \varphi_x\big(t,s,x(s)\big) y(s) \, ds.$$
Es gilt $F\big(x_0,x_0(\cdot)\big) = 0$,
da $x_0(\cdot)$ Lösung (\ref{AnhangDGLMax}) zum Anfangswert $x_0$ ist.
Ferner stellt nach Lemma \ref{LemmaDGL1} der Operator $F_{x(\cdot)}\big(x_0,x_0(\cdot)\big)$ eine lineare homöomorphe Abbildung des Raumes $C([t_0,t_1],\R^n)$ in sich dar.
Gemä"s Satz \ref{SatzImpliziteFunktionen} über implizite Funktionen wird in einer hinreichend kleinen Umgebung $U_\delta(x_0)$ des Punktes $x_0$ eine
stetig differenzierbare Abbildung $z \to x_z(\cdot)$ in den Raum $C([t_0,t_1],\R^n)$ definiert,
für die $F\big(z,x_z(\cdot)\big)=0$ gilt.
Diese Abbildung ist in $x_0$ stetig differenzierbar. 
Die Bedingung $F\big(z,x_z(\cdot)\big)=0$ bedeutet
$$x_z(t)= z + \int_{t_0}^t \varphi\big(t,s,x_z(s)\big) \, ds,$$
d.\,h., dass diese Abbildung jedem Anfangswert $z \in U_\delta(x_0)$ die Lösung $x_z(\cdot)$ von (\ref{AnhangDGLMax}) zuordnet.
Demnach wird dadurch die Abbildung $\Phi$ beschrieben und $\Phi$ ist in $x_0$ stetig differenzierbar.
Ihre Ableitung ordnet nach Satz \ref{SatzImpliziteFunktionen} jedem $y \in \R^n$ die Vektorfunktion
$$\Phi'(x_0)y=\xi_y(\cdot) = - F^{-1}_{x(\cdot)}\big(x_0,x_0(\cdot)\big) \circ \big[ F_z\big(x_0,x_0(\cdot)\big)y \big] = -F^{-1}_{x(\cdot)}\big(x_0,x_0(\cdot)\big) y(\cdot)$$
zu, wobei $F_z\big(x_0,x_0(\cdot)\big)y$ die Funktion $y(\cdot)$ mit $y(t) \equiv y$ liefert.
Weiter lässt sich die Beziehung $\xi_y(\cdot) = F^{-1}_{x(\cdot)}\big(x_0,x_0(\cdot)\big) y(\cdot)$
in die Form
$$F_{x(\cdot)}\big(x_0,x_0(\cdot)\big)\xi_y(\cdot) =  y(\cdot) \quad\Leftrightarrow\quad
  \xi_y(t) - \int_{t_0}^t \varphi_x\big(t,s,x_0(s)\big) \xi_y(s) \, ds = y$$
bringen. Der Satz ist damit bewiesen. \hfill $\blacksquare$
    \newpage
\lhead[\thepage \hspace*{1mm} Konvexe Funktionen]{}
\rhead[]{Konvexe Funktionen \hspace*{1mm} \thepage}
    \section{Konvexe Mengen und Funktionen}
\begin{definition}[Konvexe Menge]
Die Menge $M \subseteq \R^n$ heißt konvex,
wenn für alle $x_1,x_2 \in M$ und $\lambda \in [0,1]$ der Punkt $x=\lambda x_1 + (1- \lambda) x_2$
ebenfalls zu $M$ gehört.
\end{definition} 

\begin{satz}[Trennungssatz] \label{Trennungssatz}
Es sei $M \subset \R$ eine abgeschlossene und konvexe Menge.
Ferner sei $x_0$ ein Randpunkt der Menge $M$.
Dann exisitiert ein nichttrivialer Vektor $a \in \R^n\setminus \{0\}$ mit
$\langle a, x \rangle \geq \langle a, x_0 \rangle$
für alle $x \in M$.
\end{satz}

\begin{definition}[Epigraph]
Es sei $f: X \subseteq \R^m \to \R$ eine Funktion.
Dann heißt die Menge ${\rm epi\,}f = \big\{ (\alpha,x) \in \R \times X \,\big|\, \alpha \geq f(x) \big\}$
der Epigraph der Funktion $f$.
\end{definition}

\begin{definition}[Konvexe Funktion]
Die Funktion $f:X \to \R$ ist konvex, wenn der Epigraph ${\rm epi\,}f$ eine konvexe Menge ist.
\end{definition} 

\begin{lemma}[Jensensche Ungleichung]
Eine Funktion $f:X \to \R$ ist genau dann konvex,
wenn die Ungleichung $f\big(\lambda x_1 + (1-\lambda) x_2 \big) \leq \lambda f(x_1) + (1-\lambda) f(x_2)$
f"ur alle $x_1,x_2 \in X$ und alle $\lambda \in [0,1]$ gilt.
\end{lemma}

\begin{satz}[Stetigkeit konvexer Funktionen] \label{SatzStetigkeitKonvexerFunktionen}
Es sei $f$ eine konvexe Funktion über der offenen Menge $X \subseteq \R^n$.
Die Funktion $f$ ist genau dann in dem Punkt $x_0 \in X$ stetig, wenn sie auf einer Umgebung dieses Punktes nach oben beschr"ankt ist.
\end{satz}

{\em Beweis:} Ist $f$ in $x_0$ stetig, so ist sie nach Definition der Stetigkeit auf einer Umgebung von $x_0$ beschränkt.
Umgekehrt nehmen wir zur Vereinfachung $x_0=0$ und $f(0)=0$ an.
Es sei $U_\varrho(0) \subseteq X$ eine offene Kugel mit Radius $\varrho >0$ um $x_0=0$,
auf der $f$ nach oben durch $0<c < \infty$ beschränkt ist:
$f(x) < c$ f"ur alle $x \in U_\varrho(0)$. \\
Wir w"ahlen $0 < \varepsilon \leq c$ hinreichend klein und betrachten
$$U_\varepsilon = \frac{\varepsilon}{c} \cdot U_\varrho(0) \subseteq X \quad\mbox{mit}\quad
  \frac{\varepsilon}{c} \cdot U_\varrho(0)=\Big\{ y \in \R^n \,\Big|\, \frac{c}{\varepsilon} \cdot y \in U_\varrho(0)\Big\}.$$
Wir zeigen nun $|f(y)| < \varepsilon$ für alle $y \in U_\varepsilon$.
Da $\varepsilon$ beliebig klein sein kann, bedeutet dies die Stetigkeit von $f$ in $x_0=0$. \\[2mm]
Es sei $y \in U_\varepsilon$.
Dann ist $\displaystyle x=\frac{c}{\varepsilon} \cdot y \in U_\varrho(0)$ und wegen der Konvexit"at von $f$ folgt mit der Jensenschen Ungleichung einerseits
$$f(y) = f \bigg( \frac{\varepsilon}{c} \cdot x +   \bigg(1-\frac{\varepsilon}{c} \bigg) \cdot 0 \bigg)
  \leq \frac{\varepsilon}{c} \cdot f(x) +   \bigg(1-\frac{\varepsilon}{c} \bigg) \cdot f(0)
  < \frac{\varepsilon}{c} \cdot c = \varepsilon.$$
Andererseits ist $\displaystyle x = -\frac{c}{\varepsilon} \cdot y \in U_\varrho(0)$,
da $U_\varrho(0) \subset \R^n$ als eine Kugel mit Radius $\varrho >0$ um den Nullpunkt gewählt wurde.
Es ergibt sich die Beziehung
$$0 = f(0) =f\bigg( \frac{1}{1+\frac{\varepsilon}{c}} \cdot y + \frac{\frac{\varepsilon}{c}}{1+\frac{\varepsilon}{c}} \bigg(-\frac{c}{\varepsilon}\cdot y \bigg)\bigg)
    \leq \frac{1}{1+\frac{\varepsilon}{c}} \, f(y) + \frac{\frac{\varepsilon}{c}}{1+\frac{\varepsilon}{c}}\, f\left(-\frac{c}{\varepsilon}\,y\right),$$
also    
$$f(y) \geq -\frac{\varepsilon}{c} \, f\bigg(-\frac{c}{\varepsilon}\,y\bigg) =  -\frac{\varepsilon}{c} \, f(x)
       > -\frac{\varepsilon}{c} \cdot c = - \varepsilon$$
und damit die Stetigkeit von $f$ in $x_0=0$.  \hfill $\blacksquare$
\cleardoublepage
\end{appendix}
\cleardoublepage

\addcontentsline{toc}{part}{Literatur- und Sachverzeichnis}

\addcontentsline{toc}{section}{Literatur}
\lhead[\thepage \hspace*{1mm}  Literatur]{}
\rhead[]{Literatur \hspace*{1mm} \thepage}



\end{document}